\pdfoutput=1
\documentclass{amsart}

\usepackage{lmodern}
\usepackage{amsfonts}
\usepackage{amsxtra}
\usepackage{amssymb}
\usepackage{mathdots}
\usepackage{array}
\usepackage[reversemp,
    paperwidth=210mm,
    paperheight=297mm,
    top={26mm},
    headheight={5.5pt},
    headsep={5.6mm},
    text={31pc,245mm},
    marginparsep=5mm,
    marginparwidth=12mm,
    bindingoffset=6mm,
    footskip=10mm]{geometry}
\usepackage{xcolor}
\definecolor{cite}{rgb}{0.30,0.60,1.00}
\definecolor{url}{rgb}{0.00,0.00,0.80}
\definecolor{link}{rgb}{0.40,0.10,0.20}
\usepackage[colorlinks,linkcolor=link,urlcolor=url,citecolor=cite,pagebackref,breaklinks]{hyperref}
\usepackage{bbm}
\usepackage{mathtools}
\usepackage{mathrsfs}
\usepackage{appendix}
\usepackage[all]{xy}
\usepackage[lite,abbrev,msc-links]{amsrefs}
\usepackage{graphicx}
\usepackage{multirow}
\usepackage{pstricks}
\usepackage{pst-pdf}
\usepackage{enumitem}
\usepackage{marvosym}

\numberwithin{equation}{section}

\theoremstyle{plain}
\newtheorem{proposition}{Proposition}[section]

\newtheorem{corollary}[proposition]{Corollary}
\newtheorem{lem}[proposition]{Lemma}
\newtheorem{theorem}[proposition]{Theorem}

\newtheorem{hypothesis}[proposition]{Hypothesis}

\theoremstyle{definition}
\newtheorem{definition}[proposition]{Definition}
\newtheorem{construction}[proposition]{Construction}
\newtheorem{notation}[proposition]{Notation}

\newtheorem{assumption}[proposition]{Assumption}

\theoremstyle{remark}
\newtheorem{remark}[proposition]{Remark}
\newtheorem{example}[proposition]{Example}

\renewcommand{\b}[1]{\mathbf{#1}}
\renewcommand{\c}[1]{\mathcal{#1}}
\renewcommand{\d}[1]{\mathbb{#1}}
\newcommand{\f}[1]{\mathfrak{#1}}
\renewcommand{\r}[1]{\mathrm{#1}}
\newcommand{\s}[1]{\mathscr{#1}}
\renewcommand{\sf}[1]{\mathsf{#1}}
\renewcommand{\(}{\left(}
\renewcommand{\)}{\right)}
\newcommand{\res}{\mathbin{|}}
\newcommand{\ol}[1]{\overline{#1}{}}

\newcommand{\ul}{\underline}
\renewcommand{\leq}{\leqslant}
\renewcommand{\geq}{\geqslant}
\renewcommand{\clubsuit}{\et}

\newcommand{\bA}{\b A}

\newcommand{\bC}{\b C}
\newcommand{\bD}{\b D}

\newcommand{\bG}{\b G}
\newcommand{\bH}{\b H}

\newcommand{\bM}{\b M}

\newcommand{\bw}{\b w}

\newcommand{\cA}{\c A}

\newcommand{\cC}{\c C}
\newcommand{\cD}{\c D}
\newcommand{\cE}{\c E}
\newcommand{\cF}{\c F}
\newcommand{\cG}{\c G}
\newcommand{\cH}{\c H}

\newcommand{\cM}{\c M}

\newcommand{\cP}{\c P}

\newcommand{\cT}{\c T}
\newcommand{\cU}{\c U}
\newcommand{\cV}{\c V}
\newcommand{\cW}{\c W}
\newcommand{\cX}{\c X}
\newcommand{\cY}{\c Y}
\newcommand{\cZ}{\c Z}

\newcommand{\dA}{\d A}
\newcommand{\dB}{\d B}
\newcommand{\dC}{\d C}
\newcommand{\dD}{\d D}
\newcommand{\dE}{\d E}
\newcommand{\dF}{\d F}

\newcommand{\dI}{\d I}
\newcommand{\dJ}{\d J}
\newcommand{\dK}{\d K}
\newcommand{\dL}{\d L}
\newcommand{\dM}{\d M}
\newcommand{\dN}{\d N}
\newcommand{\dO}{\d O}

\newcommand{\dQ}{\d Q}
\newcommand{\dR}{\d R}
\newcommand{\dS}{\d S}
\newcommand{\dT}{\d T}

\newcommand{\dW}{\d W}

\newcommand{\dZ}{\d Z}

\newcommand{\fD}{\f D}

\newcommand{\fO}{\f O}

\newcommand{\fS}{\f S}
\newcommand{\fT}{\f T}

\newcommand{\fY}{\f Y}
\newcommand{\fZ}{\f Z}

\newcommand{\fg}{\f g}

\newcommand{\fm}{\f m}

\newcommand{\fp}{\f p}

\newcommand{\rB}{\r B}
\newcommand{\rC}{\r C}
\newcommand{\rD}{\r D}
\newcommand{\rE}{\r E}
\newcommand{\rF}{\r F}
\newcommand{\rG}{\r G}
\newcommand{\rH}{\r H}
\newcommand{\rI}{\r I}

\newcommand{\rL}{\r L}
\newcommand{\rM}{\r M}
\newcommand{\rN}{\r N}

\newcommand{\rP}{\r P}

\newcommand{\rR}{\r R}

\newcommand{\rT}{\r T}
\newcommand{\rU}{\r U}
\newcommand{\rV}{\r V}

\newcommand{\rZ}{\r Z}

\newcommand{\rc}{\r c}
\newcommand{\rd}{\,\r d}

\newcommand{\rh}{\r h}

\newcommand{\rs}{\r s}
\newcommand{\rt}{\r t}

\newcommand{\sA}{\s A}
\newcommand{\sB}{\s B}
\newcommand{\sC}{\s C}
\newcommand{\sD}{\s D}
\newcommand{\sE}{\s E}
\newcommand{\sF}{\s F}
\newcommand{\sG}{\s G}

\newcommand{\sJ}{\s J}

\newcommand{\sL}{\s L}

\newcommand{\sN}{\s N}
\newcommand{\sO}{\s O}

\newcommand{\sR}{\s R}
\newcommand{\sS}{\s S}
\newcommand{\sT}{\s T}

\newcommand{\sW}{\s W}
\newcommand{\sX}{\s X}

\newcommand{\sfA}{\sf A}

\newcommand{\sfC}{\sf C}

\newcommand{\sfF}{\sf F}
\newcommand{\sfG}{\sf G}
\newcommand{\sfH}{\sf H}

\newcommand{\sfP}{\sf P}
\newcommand{\sfQ}{\sf Q}

\newcommand{\sfT}{\sf T}
\newcommand{\sfU}{\sf U}
\newcommand{\sfV}{\sf V}
\newcommand{\sfW}{\sf W}
\newcommand{\sfX}{\sf X}
\newcommand{\sfY}{\sf Y}
\newcommand{\sfZ}{\sf Z}

\newcommand{\sfu}{\sf u}

\newcommand{\tP}{\mathtt{P}}

\newcommand{\tR}{\mathtt{R}}
\newcommand{\tS}{\mathtt{S}}
\newcommand{\tT}{\mathtt{T}}
\newcommand{\tU}{\mathtt{U}}
\newcommand{\tV}{\mathtt{V}}

\newcommand{\tc}{\mathtt{c}}

\newcommand{\tf}{\mathtt{f}}

\newcommand{\ts}{\mathtt{s}}

\newcommand{\tw}{\mathtt{w}}

\newcommand{\bomega}{\boldsymbol{\omega}}
\newcommand{\bsigma}{\boldsymbol{\Sigma}}
\newcommand{\bphi}{\boldsymbol{\phi}}
\newcommand{\bPhi}{\boldsymbol{\Phi}}

\newcommand{\bpi}{\boldsymbol{\pi}}
\newcommand{\bbh}{\boldsymbol{h}}
\newcommand{\bbi}{\boldsymbol{i}}
\newcommand{\bbq}{\boldsymbol{q}}

\newcommand{\bbA}{\boldsymbol{A}}

\newcommand{\pres}[2]{\prescript{#1}{}{{\hskip-0.02em\relax}#2}}

\newcommand{\ab}{\r{ab}}

\newcommand{\ad}{\r{ad}}

\newcommand{\ac}{\r{ac}}
\newcommand{\an}{\r{an}}

\newcommand{\can}{\r{can}}
\newcommand{\CF}{\mathbf{1}}

\newcommand{\cont}{\r{cont}}
\newcommand{\conv}{\r{conv}}
\newcommand{\cris}{\r{cris}}

\newcommand{\der}{\r{der}}

\newcommand{\dr}{\r{dR}}
\newcommand{\dtm}{\r{det}\:}

\newcommand{\et}{{\acute{\r{e}}\r{t}}}
\newcommand{\qet}{{\r{q}\acute{\r{e}}\r{t}}}

\newcommand{\fin}{\r{fin}}
\newcommand{\free}{\r{fr}}
\newcommand{\gr}{\r{gr}}

\newcommand{\Herm}{\r{Herm}}

\newcommand{\hol}{\r{hol}}
\newcommand{\id}{\r{id}}
\newcommand{\inc}{\r{inc}}
\newcommand{\inert}{\r{int}}

\newcommand{\ram}{\r{ram}}
\newcommand{\rec}{\r{rec}}
\newcommand{\reg}{\r{reg}}

\newcommand{\rig}{\r{rig}}
\newcommand{\RE}{\r{Re}\,}

\newcommand{\sel}{\r{Sel}}
\newcommand{\st}{\r{st}}
\newcommand{\std}{\r{std}}

\newcommand{\SF}{\r{SF}}

\newcommand{\spl}{\r{spl}}
\newcommand{\sph}{\r{sph}}

\newcommand{\tor}{\r{tor}}
\newcommand{\triv}{\r{triv}}

\newcommand{\unr}{\r{unr}}

\newcommand{\pr}{\mathrm{pr}}

\DeclareMathOperator{\AJ}{AJ}

\DeclareMathOperator{\Aut}{Aut}

\DeclareMathOperator{\BC}{BC}
\DeclareMathOperator{\CH}{CH}
\DeclareMathOperator{\Cor}{Cor}

\DeclareMathOperator{\Diff}{Diff}

\DeclareMathOperator{\End}{End}

\DeclareMathOperator{\Frac}{Frac}

\DeclareMathOperator{\Gal}{Gal}

\DeclareMathOperator{\GL}{GL}

\DeclareMathOperator{\Hom}{Hom}
\DeclareMathOperator{\IM}{Im}

\DeclareMathOperator{\Ind}{Ind}

\DeclareMathOperator{\Ker}{Ker}

\DeclareMathOperator{\Lie}{Lie}
\DeclareMathOperator{\Map}{Map}
\DeclareMathOperator{\Mat}{Mat}

\DeclareMathOperator{\modulo}{mod}
\DeclareMathOperator{\Nm}{Nm}

\DeclareMathOperator{\Res}{Res}

\DeclareMathOperator{\Sch}{Sch}

\DeclareMathOperator{\Spec}{Spec}

\DeclareMathOperator{\supp}{supp}
\DeclareMathOperator{\Sym}{Sym}

\DeclareMathOperator{\tr}{tr}
\DeclareMathOperator{\Tr}{Tr}

\DeclareMathOperator{\vol}{vol}

\begin{document}

\title{A $p$-adic arithmetic inner product formula}

\author{Daniel Disegni}
\address{Department of Mathematics, Ben-Gurion University of the Negev, Be'er Sheva 84105, Israel}
\address{Aix-Marseille University, CNRS, I2M - Institut de Math\'ematiques de Marseille, campus de Luminy, 13288 Marseille, France}
\email{daniel.disegni@univ-amu.fr}

\author{Yifeng Liu}
\address{Institute for Advanced Study in Mathematics, Zhejiang University, Hangzhou 310058, China}
\email{liuyf0719@zju.edu.cn}

\date{\today}
\subjclass[2010]{11G18, 11G40, 11G50, 11R34}

\begin{abstract}
  Fix a prime number $p$ and let $E/F$ be a CM extension of number fields in which $p$ splits relatively. Let $\pi$ be an automorphic representation of a quasi-split unitary group of even rank with respect to $E/F$ such that $\pi$ is ordinary above $p$ with respect to the Siegel parabolic subgroup. We construct the cyclotomic $p$-adic $L$-function of $\pi$, and a certain generating series of Selmer classes of special cycles on Shimura varieties. We show, under some conditions, that if the vanishing order of the $p$-adic $L$-function is $1$, then our generating series is modular and yields explicit nonzero classes (called Selmer theta lifts) in the Selmer group of the Galois representation of $E$ associated with $\pi$; in particular, the rank of this Selmer group is at least $1$. In fact, we prove a precise formula relating the $p$-adic heights of Selmer theta lifts to the derivative of the $p$-adic $L$-function. In parallel to Perrin-Riou's $p$-adic analogue of the Gross--Zagier formula, our formula is the $p$-adic analogue of the arithmetic inner product formula recently established by Chao~Li and the second author.
\end{abstract}

\maketitle

\tableofcontents

\section{Introduction}
\label{ss:introduction}

In 1986, Gross and Zagier published a groundbreaking formula relating the heights of Heegner points on modular curves to derivatives of $L$-functions, known as the Gross--Zagier formula \cite{GZ86}. For a cuspidal eigenform $f$ of weight $2$, an imaginary quadratic field $K$ and an unramified Dirichlet character $\xi$ of $K$, the formula shows, under the so-called Heegner condition (which implies that the Rankin--Selberg $L$-function $L(s,f,\xi)$ vanishes at the center $1$), that up to some explicit constant, $L'(1,f,\xi)$ equals the N\'{e}ron--Tate height of $H_\xi(f)$ -- the $f$-isotypic component of the $K$-Heegner point weighted by $\xi$ on a modular curve. Shortly after, Perrin-Riou found an analogous result in the $p$-adic world \cite{PR87}. Namely, she constructed a $p$-adic analogue of the (complex) $L$-function as a $p$-adic measure $\sL_p(f,\xi)$ in the Iwasawa algebra that interpolates $L(1,f\otimes\chi,\xi)$ where $\chi$ is a Dirichlet character ramified only at $p$, assuming that $f$ is ordinary at $p$ and $p$ splits in $K$. Then she proved that under the same Heegner condition, up to some explicit constant, the derivative of the $p$-adic $L$-function $\sL_p(f,\xi)$ at the trivial character equals the $p$-adic height of $H_\xi(f)$ -- this is known as the \emph{$p$-adic Gross--Zagier formula}.

Since the original work of Gross and Zagier, the Gross--Zagier formula and its $p$-adic avatar have been extended to various settings but all (essentially) for curves or fibrations/local systems over curves (see Remark \ref{re:history} below for a brief review of the $p$-adic results), until the very recent works by Chao~Li and one of us \cites{LL,LL2}. There, the authors proved a formula computing central $L$-derivatives for unitary groups of higher ranks in terms of the Beilinson--Bloch heights of special cycles. This originates from a program initiated by Kudla \cites{Kud02,Kud03,Kud04} and can be regarded as a Gross--Zagier formula in higher dimensions, as well as an arithmetic analogue of Rallis' inner product formula in the theory of the theta correspondence \cite{Ral82}. The current work contains a \emph{$p$-adic} avatar of the \emph{arithmetic inner product formula} in \cites{LL,LL2}; this is likewise the first generalization of the $p$-adic Gross--Zagier formula to genuinely higher dimensional varieties. A secondary aim of this article is to develop some foundational results in the theory of $p$-adic heights of algebraic cycles (in the two appendices); in particular, we prove a crystalline property of bi-extensions, which generalizes the fact that $p$-adic regulators take values in Selmer groups.

In the rest of this introduction, we explain our results in more detail. Throughout the article, we fix a prime number $p$, an algebraic closure $\ol\dQ_p$ of $\dQ_p$, and a CM extension $E/F$ of number fields such that \emph{every $p$-adic place of $F$ splits in $E$}. Denote by
\begin{itemize}
  \item $\tc\in\Gal(E/F)$ the Galois involution,

  \item $\tV_F^{(\lozenge)}$ the set of places of $F$ above a finite set $\lozenge$ of places of $\dQ$,\footnote{When $\lozenge=\{w\}$ is a singleton, we simply write $\tV_F^{(w)}$ for $\tV_F^{(\{w\})}$.}

  \item $\tV_F^\fin$ the set of non-archimedean places of $F$,

  \item $\tV_F^\spl$, $\tV_F^\inert$ and $\tV_F^\ram$ the subsets of $\tV_F^\fin$ of those that are split, inert and ramified in $E$, respectively.
\end{itemize}

For every number field $K$, we denote by $\Gamma_{K,p}$ the $p$-completion of
\[
K^\times\backslash\dA_K^{\infty,\times}\left/\(O_K\otimes\prod_{w\neq p}\dZ_w\)^\times\right.,
\]
which is naturally a finitely generated $\dZ_p$-module; and let $\sX_{K,p}$ be the rigid analytic space over $\dQ_p$ such that for every $\dQ_p$-affinoid algebra $R$, $\sX_{K,p}(R)$ is the set of continuous characters from $\Gamma_{K,p}$ to $R^\times$.

\subsection{Cyclotomic $p$-adic $L$-function}
\label{ss:cyclotomic}

Take an integer $r\geq 1$ and put $n=2r$. We equip $W_r\coloneqq E^n$ with the skew-hermitian form (with respect to $\tc$) given by the matrix $\tw_r\coloneqq\(\begin{smallmatrix}&1_r\\ -1_r &\end{smallmatrix}\)$. Put $G_r\coloneqq\rU(W_r)$, the unitary group of $W_r$, which is a quasi-split reductive group over $F$. Denote by $\dag$ the involution of $G_r$ given by the conjugation by the element $\(\begin{smallmatrix}1_r & \\ & -1_r \end{smallmatrix}\)$ inside $\Res_{E/F}\GL_n$. For $v\in\tV_F^\fin$, let $K_{r,v}\subseteq G_r(F_v)$ be the stabilizer of the lattice $O_{E_v}^n$.

\begin{definition}\label{de:relevant}
Let $\dL$ be a field embeddable into $\dC$. A \emph{relevant $\dL$-representation} of $G_r(\dA_F^\infty)$ is a representation $\pi$ with coefficients in $\dL$ satisfying that for every embedding $\iota\colon\dL\to\dC$,
\[
\pres{\iota}\pi\coloneqq\(\otimes_{v\in\tV_F^{(\infty)}}\pi^{[r]}_v\)\otimes\iota\pi
\]
is a tempered cuspidal automorphic representation of $G_r(\dA_F)$. Here, for $v\in\tV_F^{(\infty)}$, $\pi^{[r]}_v$ denotes the (unique up to isomorphism) holomorphic discrete series representation of $G_r(F_v)=G_r(\dR)$ with the Harish--Chandra parameter
\[
(\tfrac{n-1}{2},\tfrac{n-3}{2},\dots,\tfrac{1}{2};-\tfrac{1}{2},\dots,-\tfrac{n-3}{2},-\tfrac{n-1}{2}).
\]
In particular, $\pi$ is admissible and absolutely irreducible.
\end{definition}

We consider a finite extension $\dL/\dQ_p$ contained in $\ol\dQ_p$ and a relevant $\dL$-representation $\pi$ of $G_r(\dA_F^\infty)$. By Lemma \ref{le:dual}, $\hat\pi\coloneqq(\pi^\vee)^\dag$ is a relevant $\dL$-representation of $G_r(\dA_F^\infty)$ as well.

\begin{definition}
For $v\in\tV_F^{(p)}$, let $\tP_v$ be the set (of two elements) of places of $E$ above $v$. For $u\in\tP_v$, we have the representation $\pi_u$ of $\GL_n(F_v)$ as a local component of $\pi$ via the isomorphism $G_r(F_v)\simeq\GL_n(E_u)=\GL_n(F_v)$. In particular, $\pi_u^\vee\simeq\pi_{u^\tc}$. We say that $\pi_u$ is \emph{Panchishkin unramified} if
\begin{enumerate}
  \item $\pi_u$ is unramified;

  \item if we write the Satake polynomial of $\pi_u$, which makes sense by (1), as
     \begin{align*}
     \sfP_{\pi_u}(T)=T^n+\beta_{u,1}\cdot T^{n-1}+\cdots + \beta_{u,r}\cdot q_v^{\frac{r(r-1)}{2}} T^r+\dots+ \beta_{u,n}\cdot q_v^{\frac{n(n-1)}{2}}
     \in\dL[T]
     \end{align*}
     (see Definition \ref{de:satake1} for more details), then $\beta_{u,r}\in O_\dL^\times$, where $q_v$ is the residue cardinality of $F_v$.
\end{enumerate}
\end{definition}

We collect two important facts about Panchishkin unramified representations:
\begin{itemize}
  \item The representation $\pi_u$ is Panchishkin unramified if and only if $\pi_{u^\tc}$ is (Lemma \ref{le:satake2}). In particular, it makes sense to say that $\pi_v$ is Panchishkin unramified.

  \item If $\pi_u$ is Panchishkin unramified, then there is a unique polynomial $\sfQ_{\pi_u}(T)\in\dL[T]$ that divides $\sfP_{\pi_u}(T)$ and has the form
      \[
      \sfQ_{\pi_u}(T)=T^r+\gamma_{u,1}\cdot T^{r-1}+ \gamma_{u,2}\cdot q_vT^{r-2}+\cdots + \gamma_{u,r}\cdot q_v^{\frac{r(r-1)}{2}}
      \]
      with $\gamma_{u,r}\in O_\dL^\times$ (Proposition \ref{pr:satake}). In particular, we have an unramified principal series $\ul{\pi_u}$ of $\GL_r(F_v)$ defined over $\dL$ whose Satake polynomial is $\sfQ_{\pi_u}(T)$.
\end{itemize}

\begin{remark}
In fact, $\pi_v$ is Panchishkin unramified if and only if $\pi_v$ is unramified and $\pi$ is ordinary at $v$ with respect to the standard Siegel parabolic subgroup of $G_r$ in the sense of Hida \cite{Hid98}.
\end{remark}

\begin{theorem}\label{th:pL}
Under the above setup, suppose that $\pi_v$ is Panchishkin unramified for every $v\in\tV_F^{(p)}$. For every finite set $\lozenge$ of places of $\dQ$ containing $\{\infty,p\}$ such that $\pi_v$ is unramified for every $v\in\tV_F^\fin\setminus\tV_F^{(\lozenge)}$, there is a unique bounded analytic function $\sL_p^\lozenge(\pi)$ on the rigid analytic space $\sX_{F,p}\otimes_{\dQ_p}\dL$ such that for every finite (continuous) character $\chi\colon\Gamma_{F,p}\to\ol\dQ_p^\times$ and every embedding $\iota\colon\ol\dQ_p\to\dC$, we have
\begin{align*}
\iota\sL_p^\lozenge(\pi)(\chi)&=
\frac{1}{\rP^\iota_\pi}\cdot\frac{Z_r^{[F:\dQ]}}{b_{2r}^\lozenge(\CF)} \cdot \prod_{v\in\tV_F^{(p)}}\prod_{u\in\tP_v}
\gamma(\tfrac{1+r}{2},\iota(\ul{\pi_u}\otimes\chi_v),\psi_{F,v})^{-1} \\
&\quad\times L(\tfrac{1}{2},\BC(\iota\pi^\lozenge)\otimes(\iota\chi^\lozenge\circ\Nm_{E/F})),
\end{align*}
where
\begin{itemize}
  \item $\rP^\iota_\pi\in\dC^\times$ is a certain period for $\pi$ with respect to $\iota$ for every embedding $\iota\colon\dL\to\dC$, satisfying $\rP^\iota_\pi=\rP^\iota_{\hat\pi}$;

  \item $Z_r\coloneqq(-1)^r2^{-r^2-r}\bpi^{r^2}\frac{\Gamma(1)\cdots\Gamma(r)}{\Gamma(r+1)\cdots\Gamma(2r)}$ is the value of a certain explicit archimedean local doubling zeta integral;

  \item $b_{2r}^\lozenge(\CF)=\prod_{i=1}^{2r}L^\lozenge(i,\eta_{E/F}^i)$ is defined in \S\ref{ss:run}(F4);

  \item $\gamma(s,\iota(\ul{\pi_u}\otimes\chi_v),\psi_{F,v})$ is the gamma factor \cite{Jac79} in which $\psi_F\coloneqq\psi_\dQ\circ\Tr_{F/\dQ}$ with $\psi_\dQ\colon\dA_\dQ\to\dC^\times$ the standard automorphic additive character;\footnote{By definition, for every $v\in\tV_F^{(p)}$ and $u\in\tP_v$, $\iota\ul{\pi_u}\otimes|\;|_{F_v}^{\frac{r}{2}}$ is tempered, which implies that $\gamma(\tfrac{1+r}{2},\iota(\ul{\pi_u}\otimes\chi_v),\psi_{F,v})\in\dC^\times$.}

  \item $L(s,\BC(\iota\pi^\lozenge)\otimes(\iota\chi^\lozenge\circ\Nm_{E/F}))$ is the complex $L$-function of the (complex) representation $\BC(\iota\pi^\lozenge)\otimes(\iota\chi^\lozenge\circ\Nm_{E/F})$ of $\GL_n(\dA_E^\lozenge)$, hence is an Euler product away from $\lozenge$.
\end{itemize}
\end{theorem}

\begin{remark}\label{re:padic}
We have the following remarks concerning Theorem \ref{th:pL}.
\begin{enumerate}
  \item A bounded analytic function on the rigid analytic space $\sX_{F,p}\otimes_{\dQ_p}\dL$ is equivalent to an element in $\dZ_p[[\Gamma_{F,p}]]\otimes_{\dZ_p}\dL$, that is, an $\dL$-valued $p$-adic measure on $\Gamma_{F,p}$. In particular, the uniqueness of $\sL_p^\lozenge(\pi)$ is clear.

  \item The collection of periods $(\rP^\iota_\pi)_\iota$ is only well-defined up to a common factor in $\dL^\times$ (see Notation \ref{no:period}). In particular, the $p$-adic $L$-function $\sL_p^\lozenge(\pi)$ is only well-defined up to a factor in $\dL^\times$.

  \item The vanishing order of $\sL_p^\lozenge(\pi)$ at $\b1$ does not depend on $\lozenge$. From the interpolation formula, we have $\sL^\lozenge_p(\pi)=\sL^\lozenge_p(\hat\pi)$.

  \item Our $p$-adic $L$-function is defined over the $p$-adic field of definition of the representation and interpolates complex $L$-values along \emph{all} isomorphisms $\ol\dQ_p\simeq\dC$; this is a rationality property stronger than the one under a fixed isomorphism $\ol\dQ_p\simeq\dC$ as in the setup of many previous works in this field.

  \item Among other technical assumptions, at least when $\pi$ is ordinary at $p$ in the usual sense (that is, for every $u\in\tP$, $\beta_{u,m}\in O_\dL^\times$ for every $1\leq m \leq n$ in the Satake polynomial of $\pi_u$), our $p$-adic $L$-function has already been constructed in \cite{EHLS} up to some constant (and with a weaker rationality property). In fact, in \cite{EHLS} the authors construct more generally a multi-variable $p$-adic $L$-function in which $\pi$ is allowed to vary in an ordinary Hida family as well. In this article, we will give a (relatively) self-contained construction of our $p$-adic $L$-function independent of \cite{EHLS} since first, the process of the construction itself is an ingredient for the $p$-adic height formula; and second, our construction is technically much simpler to follow.
\end{enumerate}

\end{remark}

\subsection{Modularity of generating functions in Selmer groups}

In this subsection, we construct a Selmer group analogue of Kudla's generating functions and state a theorem on its modularity. We now suppose that $F\neq\dQ$. Fix an embedding $E\hookrightarrow\dC$ and regard $E$ as a subfield of $\dC$. For the simplicity of the introduction, we fix an embedding $\ol\dQ_p\hookrightarrow\dC$ and will not pay attention to the rationality of the constructions below, while the full details with full generality can be found in \S\ref{ss:height} and \S\ref{ss:theta}.

Let $V$ be a hermitian space over $E$ of rank $n=2r$ that has signature $(n-1,1)$ along the induced inclusion $F\subseteq\dR$ and signature $(n,0)$ at other archimedean places of $F$. Put $H\coloneqq\rU(V)$. We then have a system of Shimura varieties $\{X_L\}_L$ indexed by neat open compact subgroups $L$ of $H(\dA_F^\infty)$, which are smooth projective schemes over $E$ of dimension $n-1$.
Take a neat open compact subgroup $L\subseteq H(\dA_F^\infty)$. Let $\rV_{\pi,L}$ be the $\theta(\pi)$-isotypic subspace of $\rH^{2r-1}(X_L\otimes_E\ol{E},\ol\dQ_p(r))$ (which could be zero), where $\theta(\pi)$ denotes the (product of) local theta lifting of $\pi$. We have a canonical map
\[
\wp_\pi\colon\rH^{2r}(X_L,\ol\dQ_p(r))\to\rH^1(E,\rV_{\pi,L})
\]
from Lemma \ref{le:theta2}.

For every Schwartz function $\phi\in\sS(V^r\otimes_F\dA_F^\infty)^L$ and every $g\in G_r(\dA_F^\infty)$, we have Kudla's generating function
\[
Z_{\phi,L}(g)\coloneqq\sum_{T\in\Herm_r(F)^+}\sum_{\substack{x\in L\backslash V^r\otimes_F\dA_F^\infty\\ T(x)=T}}
(\omega_r(g)\phi)(x) Z(x)_L\cdot q^T
\]
as a formal power series indexed by totally semi-positive definite hermitian matrices $T$ over $E/F$ of rank $r$, with coefficients that are special cycles $Z(x)_L\in\CH^r(X_L)$ indexed by elements $x\in L\backslash V^r\otimes_F\dA_F^\infty$ with moment matrix $T$. Denote by $Z^\pi_{\phi,L}(g)$ its image under the composition of the absolute cycle class map $\CH^r(X_L)\to\rH^{2r}(X_L,\ol\dQ_p(r))$ and the canonical map $\wp_\pi$ mentioned above. We say that $\pi$ satisfies the \emph{Modularity Hypothesis} if:

\emph{There exists a (unique) holomorphic automorphic form $\cZ^\pi_{\phi,L}$ on $G_r(\dA_F)$ valued in the Bloch--Kato Selmer group $\rH^1_f(E,\rV_{\pi,L})$ \cite{BK90} such that the $q$-expansion of $g\cdot\cZ^\pi_{\phi,L}$ coincides with $Z^\pi_{\phi,L}(g)$ for every $g\in G_r(\dA_F^\infty)$.}

Our first result concerns the Modularity Hypothesis under certain assumptions.

\begin{assumption}\label{st:main}
Suppose that $F\neq\dQ$, that $\tV_F^\spl$ contains all 2-adic (and $p$-adic) places, and that every prime in $\tV_F^\ram$ is unramified over $\dQ$. Suppose that the relevant $\dL$-representation $\pi$ of $G_r(\dA_F^\infty)$ (with $\dL/\dQ_p$ a finite extension contained in $\ol\dQ_p$) satisfies:
\begin{enumerate}
  \item For every $v\in\tV_F^\ram$, $\pi_v$ is spherical with respect to $K_{r,v}$, that is, $\pi_v^{K_{r,v}}\neq\{0\}$.

  \item For every $v\in\tV_F^\inert$, $\pi_v$ is either unramified or almost unramified (see \cite{LL}*{Remark~1.4(3)}) with respect to $K_{r,v}$; moreover, if $\pi_v$ is almost unramified, then $v$ is unramified over $\dQ$.

  \item We have $\tR_\pi\cup\tS_\pi\subseteq\tV_F^\heartsuit$ (see below), where
     \begin{itemize}
         \item $\tR_\pi\subseteq\tV_F^\spl$ denotes the (finite) subset for which $\pi_v$ is ramified,

         \item $\tS_\pi\subseteq\tV_F^\inert$ denotes the (finite) subset for which $\pi_v$ is almost unramified.
     \end{itemize}

  \item For every $v\in\tV_F^{(p)}$, $\pi_v$ is Panchishkin unramified.
\end{enumerate}
Here, we recall from \cite{LL2} (and refer to \cite{LL2}*{Remark~1.2} for its technical nature) that $\tV_F^\heartsuit$ is the subset of $\tV_F^\spl\cup\tV_F^\inert$ consisting of $v$ satisfying that for every $v'\in\tV_F^{(p_v)}\cap\tV_F^\ram$, the subfield of $\ol{F_v}$ generated by $F_v$ and the Galois closure of $E_{v'}$ is unramified over $F_v$. In particular, $\tV_F^\heartsuit$ contains $\tV_F^{(p)}$.
\end{assumption}

\begin{theorem}[Theorem \ref{th:modularity}]\label{th:main2}
Suppose that we are in the situation of Assumption \ref{st:main} and $n<p$. If the vanishing order of $\sL_p^\lozenge(\pi)$ at $\b1$ is one, then $\pi$ satisfies the Modularity Hypothesis.
\end{theorem}

\subsection{A $p$-adic arithmetic inner product formula}

In this subsection, we construct a Selmer group analogue of the (arithmetic) theta lift, and state a corresponding inner product formula for it, which we call the \emph{$p$-adic arithmetic inner product formula}. The details can be found in \S\ref{ss:theta} and \S\ref{ss:aipf}. We keep the setup from the previous subsection.

Suppose that both $\pi$ and $\hat\pi$ satisfy the Modularity Hypothesis. For every $\varphi\in\hat\pi$, we define $\Theta_\phi^\sel(\varphi)_L$ to be the convolution of $\varphi^\dag$ and $\cZ^\pi_{\phi,L}$, which is an element of $\rH^1_f(E,\rV_{\pi,L})$. The element $\Theta_\phi^\sel(\varphi)_L$ is the Selmer group analogue of the arithmetic theta lift constructed in \cites{Liu11,LL}, which we call a \emph{Selmer theta lift}.

The Poincar\'{e} duality for $X_L$ induces a pairing $\rV_{\pi,L}\times\rV_{\hat\pi,L}\to\ol\dQ_p(1)$. By Nekov\'{a}\v{r}'s theory \cite{Nek93}, we have a $p$-adic height pairing
\[
\langle\;,\;\rangle_E\colon\rH^1_f(E,\rV_{\pi,L})\times\rH^1_f(E,\rV_{\hat\pi,L})\to\Gamma_{E,p}\otimes_{\dZ_p}\ol\dQ_p
\]
using certain canonical Hodge splitting at $p$-adic places. For every $\varphi_1\in\hat\pi$, every $\varphi_2\in\pi$ and every pair $\phi_1,\phi_2\in\sS(V^r\otimes_F\dA_F^\infty)^L$, the height
\[
\vol^\natural(L)\cdot\langle\Theta_{\phi_1}^\sel(\varphi_1)_L,\Theta_{\phi_2}^\sel(\varphi_2)_L\rangle_E
\in\Gamma_{E,p}\otimes_{\dZ_p}\ol\dQ_p
\]
is independent of $L$, where $\vol^\natural(L)$ denotes a certain canonical volume of $L$ introduced in \cite{LL}*{Definition~3.8}. We will denote the above canonical value as $\langle\Theta_{\phi_1}^\sel(\varphi_1),\Theta_{\phi_2}^\sel(\varphi_2)\rangle_{\pi,E}^\natural$.

\begin{theorem}[$p$-adic arithmetic inner product formula, Theorem \ref{th:aipf}]\label{th:main1}
Suppose that we are in the situation of Assumption \ref{st:main} and $n<p$.
\begin{enumerate}
  \item If the vanishing order of $\sL_p^\lozenge(\pi)$ at $\CF$ is one (so that both $\pi$ and $\hat\pi$ satisfy the Modularity Hypothesis by Theorem \ref{th:main2} and Remark \ref{re:padic}(3)), then for every choice of elements
     \begin{itemize}
        \item $\varphi_1=\otimes_v\varphi_{1,v}\in\hat\pi$ and $\varphi_2=\otimes_v\varphi_{2,v}\in\pi$ such that for every $v\not\in\tV_F^{(\lozenge)}$, $\varphi_{1,v}$ and $\varphi_{2,v}$ are fixed by $K_{r,v}$ such that $\langle\varphi_{1,v},\varphi_{2,v}\rangle_{\pi_v}=1$,

        \item $\phi_1=\otimes_v\phi_{1,v},\phi_2=\otimes_v\phi_{2,v}\in\sS(V^r\otimes_F\dA_F^\infty)$ with $\phi_1^\lozenge=\phi_2^\lozenge$ being the characteristic function of $(\Lambda^\lozenge)^r$ in which $\Lambda^\lozenge$ is a self-dual lattice of $V\otimes_F\dA_F^\lozenge$,
     \end{itemize}
     the identity
     \begin{align*}
     \Nm_{E/F}\langle\Theta_{\phi_1}^\sel(\varphi_1),\Theta_{\phi_2}^\sel(\varphi_2)\rangle_{\pi,E}^\natural
     &=\partial\sL_p^\lozenge(\pi)(\b1)\cdot
     \prod_{v\in\tV_F^{(p)}}\prod_{u\in\tP_v}\gamma(\tfrac{1+r}{2},\ul{\pi_u},\psi_{F,v}) \\
     &\quad\times \prod_{v\in\tV_F^{(\lozenge\setminus\{\infty\})}}
     Z(\varphi_{1,v}^\dag\otimes\varphi_{2,v},f^{\r{SW}}_{\phi_{1,v}\otimes\phi_{2,v}})
     \end{align*}
     holds in $\Gamma_{F,p}\otimes_{\dZ_p}\dC$, where the term $Z(\varphi_{1,v}^\dag\otimes\varphi_{2,v},f^{\r{SW}}_{\phi_{1,v}\otimes\phi_{2,v}})$ is the local doubling zeta integral with respect to the Siegel--Weil section $f^{\r{SW}}_{\phi_{1,v}\otimes\phi_{2,v}}$ associated with $\phi_{1,v}\otimes\phi_{2,v}$.

  \item If the vanishing order of $\sL_p^\lozenge(\pi)$ at $\CF$ is not one, then assuming that both $\pi$ and $\hat\pi$ satisfy the Modularity Hypothesis, we have
      \[
      \Nm_{E/F}\langle\Theta_{\phi_1}^\sel(\varphi_1),\Theta_{\phi_2}^\sel(\varphi_2)\rangle_{\pi,E}^\natural=0
      \]
      for every $\varphi_1\in\hat\pi$, $\varphi_2\in\pi$, and $\phi_1,\phi_2\in\sS(V^r\otimes_F\dA_F^\infty)$. (See Theorem \ref{th:aipf}(2) for a version of this part that does not rely on the Modularity Hypothesis.)
\end{enumerate}
\end{theorem}

The above theorem is only nontrivial when $r[F:\dQ]+|\tS_\pi|$ is odd (Remark \ref{re:aipf}(3)).

The theorem has applications to the \emph{$p$-adic} Beilinson--Bloch--Kato conjecture. Associated with $\pi$, we have a semisimple continuous representation $\rho_\pi$ of $\Gal(\ol{E}/E)$ of dimension $n$ with coefficients in $\ol\dQ_p$, satisfying $\rho_\pi^\tc\simeq\rho_\pi^\vee(1-n)$ (Lemma \ref{le:galois}). Then in the interpolation property of $\sL_p^\lozenge(\pi)$ in Theorem \ref{th:pL}, we have
\[
L(\tfrac{1}{2},\BC(\iota\pi^\lozenge)\otimes(\iota\chi^\lozenge\circ\Nm_{E/F}))
=L^\lozenge(0,\iota(\rho_\pi(r)\otimes\chi\res_{\Gal(\ol{E}/E)})),
\]
where on the right-hand side we view $\chi$ as a $\ol\dQ_p$-valued character of $\Gal(\ol{E}/F)$ via the global class field theory. The following corollary provides evidence toward the $p$-adic Beilinson--Bloch--Kato conjecture for (genuinely) higher-dimensional motives, whose deduction is provided after Remark \ref{re:aipf}.

\begin{corollary}\label{th:main}
Suppose that we are in the situation of Assumption \ref{st:main} and $n<p$. If the vanishing order of $\sL_p^\lozenge(\pi)$ at $\b1$ is one, then
\[
\dim_{\ol\dQ_p}\rH^1_f(E,\rho_\pi(r))\geq 1.
\]
\end{corollary}

\begin{remark}\label{re:history}
When $n=2$, this result is a variant of the main application of the $p$-adic Gross--Zagier formula of \cite{PR87}, as generalized to totally real fields by one of us \cite{Dis17} following the development of \cite{GZ86} in \cite{YZZ}. In different directions, Perrin-Riou's results had been generalized to the case of higher-weight modular forms by Nekov\'{a}\v{r} \cite{Nek95} and further to the case with twists by higher-weight Hecke characters by Shnidman \cite{Shn16}, to the supersingular case by Kobayashi \cite{Kob13}, and to the case where $p$ is not necessarily relative split by one of us \cite{Dis22}.\footnote{In fact, in \cite{Kob13}, a formula in the nonsplit case is deduced from the split case by making use of some special features of the setup under consideration.} A common generalization of \cites{Nek95,Shn16,Dis17,Dis22} was developed in \cite{Dis23}.
\end{remark}

\begin{remark}
Strictly speaking, Theorem \ref{th:main1} (together with Corollary \ref{th:main} and Corollary \ref{co:main} below) relies on a hypothesis on the characterization of the tempered part of the cohomology of certain unitary Shimura varieties (see Hypothesis \ref{hy:galois} and Remark \ref{re:hypothesis}), which is expected to be verified in a sequel of the work \cite{KSZ}.
\end{remark}

\subsection{Application to symmetric power of elliptic curves}

The above results have applications to the motives of symmetric power of elliptic curves. We consider a \emph{modular} elliptic curve $A$ over $F$ without complex multiplication that has \emph{ordinary good reduction} at every $p$-adic place of $F$. Denote by $\tV_F^A\subseteq\tV_F^\fin$ the subset consisting of places over which $A$ has bad reduction.

By the very recent breakthrough on the automorphy of symmetric powers of Hilbert modular forms \cite{NT}, there exists a unique cuspidal automorphic representation $\Pi(\Sym^{n-1}A)$ of $\GL_n(\dA_F)$ satisfying
\begin{itemize}
  \item for every $v\in\tV_F^{(\infty)}$, the base change of $\Pi(\Sym^{n-1}A)_v$ to $\GL_n(\dC)$ is the principal series representation of characters $(\arg^{1-n},\arg^{3-n},\dots,\arg^{n-3},\arg^{n-1})$, where $\arg\colon\dC^\times\to\dC^\times$ is the character given by $\arg(z)\coloneqq z/\sqrt{z\ol{z}}$;

  \item for every $v\in\tV_F^\fin\setminus\tV_F^A$, $\Pi(\Sym^{n-1}A)_v$ is unramified with the Satake polynomial
      \[
      \prod_{j=0}^{n-1}\(T-\alpha_{v,1}^j\alpha_{v,2}^{n-1-j}\)\in\dQ[T],
      \]
      where $\alpha_{v,1}$ and $\alpha_{v,2}$ are the two roots of the polynomial $T^2-a_v(A)T+q_v$ (with $q_v$ the residue cardinality of $F_v$).
\end{itemize}
Let $\Pi(\Sym^{n-1}A_E)$ be the (solvable) base change of $\Pi(\Sym^{n-1}A)$ to $E$, which is a cuspidal automorphic representation of $\GL_n(\dA_E)$. The representation $\Pi(\Sym^{n-1}A_E)$ satisfies $\Pi(\Sym^{n-1}A_E)^\vee\simeq\Pi(\Sym^{n-1}A_E)\simeq\Pi(\Sym^{n-1}A_E)^\tc$, hence is a relevant representation in the sense of \cite{LTXZZ}*{Definition~1.1.3}. By \cite{LTXZZ}*{Remark~1.1.4} and the endoscopic classification for quasi-split unitary groups \cite{Mok15}, there exists a cuspidal automorphic representation $\pi(\Sym^{n-1}A_E)$ of $G_r(\dA_F)$ satisfying
\begin{itemize}
  \item for every $v\in\tV_F^{(\infty)}$, $\pi(\Sym^{n-1}A_E)_v$ is isomorphic to $\pi^{[r]}_v$;

  \item for every $v\in\tV_F^\fin\setminus\tV_F^A$, $\pi(\Sym^{n-1}A_E)_v$ is spherical with respect to $K_{r,v}$ and its base change to $\GL_n(E_v)$ is isomorphic to $\Pi(\Sym^{n-1}A_E)_v$.
\end{itemize}
In particular, there exists a relevant $\dQ$-representation $\pi$ in the sense of Definition \ref{de:relevant} such that $\pi\otimes_\dQ\dC\simeq\pi(\Sym^{n-1}A_E)^\infty$. Moreover, for every $v\in\tV_F^{(p)}$, $\pi_v\otimes_\dQ\dQ_p$ is Panchishkin unramified. Applying Theorem \ref{th:pL} to $\pi$ (or rather $\pi\otimes_\dQ\dQ_p$), we obtain a bounded analytic function $\sL_p^\lozenge(\pi)$ on $\sX_{F,p}$ for every finite set $\lozenge$ of places of $\dQ$ containing $\{\infty,p\}$ and every prime number underlying $\tV_F^\ram\cup\tV_F^A$. For every $v\in\tV_F^{(p)}$ and $u\in\tP_v$, the unramified representation $\ul{\pi_u}$ of $\GL_r(F_v)$ is the one with the Satake polynomial
\[
\prod_{j=r}^{n-1}\(T-\alpha_{v,1}^j\alpha_{v,2}^{n-1-j}\)\in\dQ_p[T],
\]
where we have ordered $\alpha_{v,1},\alpha_{v,2}\in\dQ_p^\times$ in the way that $\alpha_{v,i}\in q_v^{i-1}\dZ_p^\times$. The following is an immediate consequence of Corollary \ref{th:main} in which $\tS_\pi=\emptyset$.

\begin{corollary}\label{co:main}
Under the above setup, we further assume that
\begin{itemize}
  \item $n<p$,

  \item $[F:\dQ]>1$,

  \item $r[F:\dQ]$ is odd,

  \item every prime in $\tV_F^\ram$ is unramified over $\dQ$,

  \item $\tV_F^A\cup\tV_F^{(2)}$ is contained in $\tV_F^\spl$.
\end{itemize}
Then $\sL_p^\lozenge(\pi)(\CF)=0$. Moreover, if $\partial\sL_p^\lozenge(\pi)(\b{1})\neq 0$, then
\[
\dim_{\dQ_p}\rH^1_f(E,\Sym^{n-1}\rH^1_{\et}(A_{\ol{E}},\dQ_p)(r))\geq 1.
\]
\end{corollary}

\subsection{Structure and strategy}

We explain the structure of the article and the strategy for the proofs. Before that, we point out that throughout the article, we have restricted ourselves to only use $p$-adic measures valued in finite products of finite extensions of $\dQ_p$ to reduce the technical burden such as infinite dimensional $p$-adic Banach spaces.

In Section \ref{ss:2}, we make preparation for proving the rationality property of our $p$-adic $L$-function. In \S\ref{ss:run}, we collect two sets of more specialized notation that will be used throughout the main part of the article. In \S\ref{ss:hermitian}, we introduce the notion of Siegel hermitian varieties which are over $\dQ_p$ and are the main stage to characterize the rationality of automorphic forms on the unitary group $G_r$. In \S\ref{ss:pel}, we review the construction of an auxiliary Shimura variety over $\dQ$ that is of PEL type in the sense of Kottwitz, which is needed to prove the rationality of certain Eisenstein series used in the doubling method. The main reason we pass to this auxiliary one is that the theory of algebraic $q$-expansions is only available for such Shimura varieties. However, if the reader is satisfied with fixing an isomorphism $\ol\dQ_p\simeq\dC$ from the beginning and does not care about the field of definition of the $p$-adic $L$-function, then there is no need to use those parts of \S\ref{ss:hermitian} that are related to Shimura varieties and the entire \S\ref{ss:pel}.

In Section \ref{ss:3}, we construct the $p$-adic $L$-function. The main strategy is to use the doubling method for an ``analytic'' family of sections in the degenerate principal series of the doubling unitary group $G_{2r}$, similar to \cite{EHLS}. However, it is worth pointing out that our computation makes no use of Weil representations (or their twisted versions). In particular, we do not need any explicit Schwartz functions on hermitian spaces. In fact, we do not even need an explicit formula for the sections in the degenerate principal series at $p$-adic places -- what we need is just their Fourier transforms, which have very simple forms. The main reason we can simplify the computation is a formula obtained in the previous work \cite{LL} for computing the local doubling zeta integral (see Lemma \ref{le:zeta}). Using this formula, the gamma factor in Theorem \ref{th:pL} appears naturally and immediately. In \S\ref{ss:degenerate}, we review the doubling degenerate principal series and collect some facts on their Siegel--Fourier coefficients. In \S\ref{ss:eisenstein}, we review the doubling Eisenstein series and prove a certain rationality property of their pullbacks to the diagonal block. In \S\ref{ss:relevant}, we make all the representational-theoretical preparations; in particular, we study Panchishkin unramified representations. In \S\ref{ss:zeta}, we prove several formulae for local doubling zeta integrals. In \S\ref{ss:construction}, we complete the construction of the $p$-adic $L$-function by defining it as an inner product of a specific element in $\hat\pi\boxtimes\pi$ and the pullback of the family of doubling Eisenstein series with respect to a careful choice of sections in degenerate principal series. In \S\ref{ss:measure}, we collect some basic facts about $p$-adic measures that will be used later.

In Section \ref{ss:4}, we construct the so-called Selmer theta lifts, which are Selmer group analogues of the classical theta lifts, and study their $p$-adic heights. In \S\ref{ss:setup}, we introduce further notation for the whole section and study the rationality of local theta liftings. In \S\ref{ss:height}, we construct a canonical projection from the absolute $p$-adic cohomology to the Galois cohomology $\rH^1(E,\rV_{\pi,L})$ or $\rH^1(E,\rV_{\hat\pi,L})$, and define $p$-adic height pairings on the latter. In \S\ref{ss:theta}, we state a theorem (Theorem \ref{th:modularity}) on the modularity of Kudla's generating functions valued in the Selmer groups of the above Galois cohomology, and then construct Selmer theta lifts which belong to the Selmer groups. In \S\ref{ss:aipf}, we state the precise version of our $p$-adic arithmetic inner formula (Theorem \ref{th:aipf}), which is slightly stronger than Theorem \ref{th:main1} by taking rationality into account. In \S\ref{ss:strategy}, we present our strategy of reducing Theorem \ref{th:modularity} to the problem of computing $p$-adic height itself, and explain that it makes sense and suffices to consider the $p$-adic height pairing $\langle Z_{T_1}(\phi_1),Z_{T_2}(\phi_2)\rangle_E$ between (weighted) special cycles for a certain pool of Schwartz functions, together with a formula decomposing the (global) $p$-adic height pairing into local ones. In \S\ref{ss:height1}, we compute local $p$-adic height pairings between special cycles at (nonarchimedean) places of $E$ not above $p$, based on a result of Scholl that relates local $p$-adic heights to Beilinson's local indices and the formulae for the latter from previous works \cites{LL,LL2}. In \S\ref{ss:height2}, we study local $p$-adic height pairings between special cycles at $p$-adic places of $E$. With a crucial ingredient (Theorem \ref{th:crystalline}) on the crystalline property of the corresponding bi-extensions, we show that the local $p$-adic heights approach $0$ $p$-adically when one repeatedly applies a certain operator $\rU_p$ to the Schwartz functions. In \S\ref{ss:proof}, we finish the proof of Theorem \ref{th:main} assuming a nonvanishing result (Proposition \ref{pr:modularity3}), by using the previous formulae on local $p$-adic heights together with certain limit processes. In \S\ref{ss:proof1}, we prove Proposition \ref{pr:modularity3} by a variant of the $p$-adic doubling formula from \S\ref{ss:construction}; then we complete the proof of Theorem \ref{th:aipf}.

The article has two appendices. In Appendix \ref{ss:a}, we develop further the theory of $p$-adic heights on \emph{general varieties}, after Nekov\'{a}\v{r}. For local $p$-adic heights above $p$, we state a key theorem (Theorem \ref{th:crystalline}) on the crystalline property for certain bi-extensions, whose proof occupies the entire Appendix \ref{ss:b}.

\subsection{Notation and conventions}
\label{ss:notation}

\begin{itemize}
  \item We denote $\dN\coloneqq\{0,1,2,\dots\}$.

  \item We denote by $\ol\dZ_p$ the ring of integers of $\ol\dQ_p$.

  \item We write $\bpi=3.1415926\dots$, to be distinguished from the representation $\pi$. We also write $\bbi$ for the imaginary unit in $\dC$, to be distinguished from the commonly used index $i$.

  \item When we have a function $f$ on a product set $A_1\times\cdots\times A_s$, we will write $f(a_1,\dots,a_s)$ instead of $f((a_1,\dots,a_s))$ for its value at an element $(a_1,\dots,a_s)\in A_1\times\cdots\times A_s$.

  \item For a set $S$, we denote by $\CF_S$ the characteristic function of $S$.

  \item All rings are commutative and unital; and ring homomorphisms preserve units. However, we use the word \emph{algebra} in the general sense, which is not necessarily commutative or unital.

  \item If a base ring is not specified in the tensor operation $\otimes$, then it is $\dZ$.

  \item For an abelian group $A$ and a ring $R$, we put $A_R\coloneqq A\otimes R$ as an $R$-module.

  \item For an abelian group $A$, we denote by $A^\free$ its free quotient.

  \item For a ring $R$, we denote by $\Sch'_{/R}$ the category of locally Noetherian schemes over $R$.

  \item We denote by $\bG$ the multiplicative group scheme, that is, $\Spec\dZ[X,X^{-1}]$.

  \item For an integer $m\geq 0$, we denote by $0_m$ and $1_m$ the null and identity matrices of rank $m$, respectively, and by $\tw_m$ the matrix $\(\begin{smallmatrix}&1_m\\ -1_m &\end{smallmatrix}\)$.

  \item Let $\psi_\dQ\colon \dQ\backslash\dA_\dQ\to\dC^\times$ be the standard automorphic additive character that sends $w^{-1}$ at a prime $w$ to $\exp(-2\bpi\bbi/w)$, and put $\psi_K\coloneqq\psi_\dQ\circ\Tr_{K/\dQ}$ for every number field $K$.

  \item For a subring $R\subseteq\dC$ and a positive integer $\Delta$, we denote by $R\langle\Delta\rangle\subseteq\dC$ the subring generated by $\Delta^l$-th roots of unity for all $l\geq 0$.

  \item For a locally compact totally disconnected space $X$ and a ring $R$, we denote by $\sS(X,R)$ the $R$-module of $R$-valued locally constant compactly supported functions on $X$. We omit $R$ from the notation when $R=\dC$.
\end{itemize}

\subsubsection*{Acknowledgements}

D.~D. would like to thank Ellen~Eischen and Zheng~Liu for correspondence on $p$-adic $L$-functions. Y.~L. would like to thank Yichao~Tian for some general discussion related to Theorem \ref{th:crystalline}. We would like to thank Marc-Hubert~Nicole and Congling~Qiu for useful comments. Finally, we are grateful to the anonymous referee for many helpful comments and suggestions. The research of D.~D. is partially supported by
ISF grant 1963/20 and BSF grant 2018250.

\section{Siegel hermitian varieties}
\label{ss:2}

Recall that we have fixed the CM extension $E/F$ of number fields with the Galois involution $\tc$, such that every $p$-adic place of $F$ splits in $E$.

\subsection{Running notation}
\label{ss:run}

We introduce two sets of more specialized notation that will be used throughout the main part of the article.

\begin{enumerate}[label=(F\arabic*)]
  \item We denote by
      \begin{itemize}
        \item $\tV_F$ and $\tV_F^\fin$ the set of all places and non-archimedean places of $F$, respectively;

        \item $\tV_F^\spl$, $\tV_F^\inert$ and $\tV_F^\ram$ the subsets of $\tV_F^\fin$ of those that are split, inert and ramified in $E$, respectively;

        \item $\tV_F^{(\lozenge)}$ the subset of $\tV_F$ of places above a finite set $\lozenge$ of places of $\dQ$.
      \end{itemize}
      Moreover,
      \begin{itemize}
        \item for every $v\in\tV_F$, we put $E_v\coloneqq E\otimes_FF_v$;

        \item for every finite set $\lozenge$ of places of $\dQ$, we put $F_\lozenge\coloneqq\prod_{v\in\tV_F^{(\lozenge)}}F_v$;

        \item for every $v\in\tV_F^\fin$, we denote by $p_v$ the underlying rational prime of $v$ and by $\fp_v$ the maximal ideal of $O_{F_v}$, put $q_v\coloneqq|O_{F_v}/\fp_v|$ which is a power of $p_v$, and let $d_v\geq 0$ be the integer such that $\fp_v^{d_v}$ generates the different ideal of $F_v/\dQ_{p_v}$.
      \end{itemize}

  \item For every $v\in\tV_F^{(p)}$, let $\tP_v$ be the set of places of $E$ above $v$. Put $\tP\coloneqq\bigcup_{v\in\tV_F^{(p)}}\tP_v$. We fix a subset $\tP_{\r{CM}}$ of $\tP$ satisfying that $\tP_{\r{CM}}\cap\tP_v$ is a singleton for every $v\in\tV_F^{(p)}$.

  \item Let $m\geq 0$ be an integer.
      \begin{itemize}
        \item We denote by $\Herm_m$ the subscheme of $\Res_{O_E/O_F}\Mat_{m,m}$ of $m$-by-$m$ matrices $b$ satisfying $\pres{\rt}{b}^\tc=b$. Put $\Herm_m^\circ\coloneqq\Herm_m\cap\Res_{O_E/O_F}\GL_m$.

        \item For every (ordered) partition $m=m_1+\cdots+m_s$ with $m_i$ a positive integer, we denote by
            \[
            \partial_{m_1,\dots,m_s}\colon\Herm_m\to\Herm_{m_1}\times\cdots\times\Herm_{m_s}
            \]
            the morphism that extracts the diagonal blocks with corresponding ranks.

        \item We denote by $\Herm_m(F)^+$ (resp.\ $\Herm^\circ_m(F)^+$) the subset of $\Herm_m(F)$ of elements that are totally semi-positive definite (resp.\ totally positive definite).
      \end{itemize}

  \item Let $\eta_{E/F}\colon F^\times\backslash\dA_F^\times\to\dC^\times$ be the quadratic character associated with $E/F$. For every finite character $\chi\colon F^\times\backslash\dA_F^\times\to\dC^\times$ and every integer $m\geq 1$, we put
      \begin{itemize}
        \item for every $v\in\tV_F$,
          \[
          b_{m,v}(\chi)\coloneqq\prod_{i=1}^m L(i,\chi_v\eta_{E/F,v}^{m-i});
          \]

        \item for a finite set $\lozenge$ of places of $\dQ$,
            \begin{align*}
            b_{m,\lozenge}(\chi)\coloneqq\prod_{v\in\tV_F^{(\lozenge)}}b_{m,v}(\chi),\quad
            b_m^\lozenge(\chi)\coloneqq\prod_{v\in\tV_F\setminus\tV_F^{(\lozenge)}}b_{m,v}(\chi),
            \end{align*}
            in which the latter product is absolutely convergent when $m$ is even or $\chi\neq\CF$.
      \end{itemize}
\end{enumerate}

Let $m\geq 1$ be an integer. We equip $W_m=E^{2m}$ and $\bar{W}_m=E^{2m}$ with the skew-hermitian forms (that are $E$-linear in the first variable) given by the matrices $\tw_m$ and $-\tw_m$, respectively.

\begin{enumerate}[label=(G\arabic*)]
  \item Let $G_m$ be the unitary group of both $W_m$ and $\bar{W}_m$. We write elements of $W_m$ and $\bar{W}_m$ in the row form, on which $G_m$ acts from the right. Denote by $\dag$ the involution of $G_m$ given by the conjugation by the element $\(\begin{smallmatrix}1_m & \\ & -1_m \end{smallmatrix}\)$ inside $\Res_{E/F}\GL_{2m}$.

  \item We denote by $\{e_1,\dots,e_{2m}\}$ and $\{\bar{e}_1,\dots,\bar{e}_{2m}\}$ the natural bases of $W_m$ and $\bar{W}_m$, respectively.

  \item Let $P_m\subseteq G_m$ be the parabolic subgroup stabilizing the subspace generated by $\{e_{m+1},\dots,e_{2m}\}$, and $N_m\subseteq P_m$ its unipotent radical.

  \item We have
     \begin{itemize}
       \item a homomorphism $m\colon\Res_{E/F}\GL_m\to P_m$ sending $a$ to
          \[
          m(a)\coloneqq
          \begin{pmatrix}
              a &  \\
               & \pres{\rt}{a}^{\tc,-1} \\
          \end{pmatrix}
          ,
          \]
          which identifies $\Res_{E/F}\GL_m$ as a Levi factor of $P_m$, denoted by $M_m$.

       \item a homomorphism $n\colon\Herm_{m,F}\to N_m$ sending $b$ to
          \[
          n(b)\coloneqq
          \begin{pmatrix}
              1_m & b \\
               & 1_m \\
          \end{pmatrix}
          ,
          \]
          which is an isomorphism.
     \end{itemize}

  \item We define a maximal compact subgroup $K_m=\prod_{v\in\tV_F}K_{m,v}$ of $G_m(\dA_F)$ in the following way:
    \begin{itemize}
      \item for $v\in\tV_F^\fin$, $K_{m,v}$ is the stabilizer of the lattice $O_{E_v}^{2m}$;

      \item for $v\in\tV_F^{(\infty)}$, $K_{m,v}$ is the subgroup of the form
         \[
         [k_1,k_2]\coloneqq\frac{1}{2}
         \begin{pmatrix}
           k_1+k_2   & -\bbi k_1+\bbi k_2 \\
           \bbi k_1-\bbi k_2   & k_1+k_2 \\
         \end{pmatrix}
         ,
         \]
         in which $k_i\in\GL_m(\dC)$ satisfies $k_i\pres{\rt}{\ol{k_i}}=1_m$ for $i=1,2$.\footnote{Here, we choose a complex embedding of $E$ above $v$ to identify $G_m(F_v)$ as a subgroup of $\GL_{2m}(\dC)$. However, neither $K_{m,v}$ nor the character $\kappa_{m,v}$ in (G6) depends on such a choice.}
    \end{itemize}
    Moreover,
    \begin{itemize}
      \item for every place $w$ of $\dQ$, put $K_{m,w}\coloneqq\prod_{v\in\tV_F^{(w)}}K_{m,v}$;

      \item for a set $\lozenge$ of places of $\dQ$, put $K_m^\lozenge\coloneqq\prod_{w\not\in\lozenge}K_{m,w}$.
    \end{itemize}

  \item For every $v\in\tV_F^{(\infty)}$, we have a character $\kappa_{m,v}\colon K_{m,v}\to\dC^\times$ that sends $[k_1,k_2]$ to $\dtm k_1/\dtm k_2$.

  \item For every $v\in\tV_F$, we define a Haar measure $\r{d}g_v$ on $G_m(F_v)$ as follows:
    \begin{itemize}
      \item for $v\in\tV_F^\fin$, $\r{d}g_v$ is the Haar measure under which $K_{m,v}$ has volume $1$;

      \item for $v\in\tV_F^{(\infty)}$, $\r{d}g_v$ is the product of the Haar measure on $K_{m,v}$ under which $K_{m,v}$ has volume $1$ and the standard hyperbolic measure on $G_m(F_v)/K_{m,v}$ (see, for example, \cite{EL}*{Section~2.1}).
    \end{itemize}
    Put $\r{d}g=\prod_{v}\r{d}g_v$, which is a Haar measure on $G_m(\dA_F)$.

  \item Let $m_1,\dots,m_s$ be finitely many positive integers. Put
       \[
       G_{m_1,\dots,m_s}\coloneqq G_{m_1}\times\cdots\times G_{m_s}.
       \]
       We denote by $\cA_{m_1,\dots,m_s}$ the space of both $\cZ(\fg_{m_1,\dots,m_s,\infty})$-finite and $K_{m_1,\infty}\times\cdots\times K_{m_s,\infty}$-finite automorphic forms (in the sense of \cite{BJ79}*{\S4.2}) on $G_{m_1,\dots,m_s}(\dA_F)$, where $\cZ(\fg_{m_1,\dots,m_s,\infty})$ denotes the center of the complexified universal enveloping algebra of the Lie algebra $\fg_{m_1,\dots,m_s,\infty}$ of $G_{m_1,\dots,m_s}\otimes_\dQ\dR$. For every integer $w\geq 0$ (as weight), we denote by
      \begin{itemize}
        \item $\cA_{m_1,\dots,m_s}^{[w]}$ the maximal subspace of $\cA_{m_1,\dots,m_s}$ on which for every $v\in\tV_F^{(\infty)}$ and every $1\leq j\leq s$, $K_{m_j,v}$ acts by the character $\kappa_{m_j,v}^w$,

        \item $\cA_{m_1,\dots,m_s,\hol}^{[w]}\subseteq\cA_{m_1,\dots,m_s}^{[w]}$ the subspace of holomorphic ones.
      \end{itemize}

  \item For every vector space $\cH$ on which $G_{m_1,\dots,m_s}(\dA_F^\infty)$ acts, we put $\cH(K)\coloneqq\cH^K$ for every open compact subgroup $K\subseteq G_{m_1,\dots,m_s}(\dA_F^\infty)$.
\end{enumerate}

\subsection{Siegel hermitian varieties and line bundles of automorphy}
\label{ss:hermitian}

We first recall the construction of a CM moduli problem following \cite{LTXZZ}*{Section~3.5}. Let $T$ be the subtorus of $\Res_{E/\dQ}\bG$ that is the inverse image of $\bG_\dQ$ under the norm map $\Nm_{E/F}\colon\Res_{E/\dQ}\bG\to\Res_{F/\dQ}\bG$.

For every nonzero element $\delta\in E^{\tc=-1}$, we denote by $W^\delta$ the $E$-vector space $E$ (itself) together with a pairing $\langle\;,\;\rangle^\delta\colon E\times E\to\dQ$ given by $\langle x,y\rangle^\delta=\Tr_{E/\dQ}(\delta x y^\tc)$. For every $\dQ$-ring $R$, we have
\[
T(R)=\{t\in(E\otimes_\dQ R)^\times\res\langle tx,ty\rangle^\delta=c(t)\langle x,y\rangle^\delta\text{ for some $c(t)\in R^\times$}\}.
\]

For every neat open compact subgroup $K_T$ of $T(\dA^\infty)$, we define a moduli problem $\bsigma^\delta(K_T)$ on $\Sch'_{/\dQ_p}$ as follows: for every $S\in\Sch'_{/\dQ_p}$, $\bsigma^\delta(K_T)(S)$ is the set of equivalence classes of quadruples $(A_0,i_0,\lambda_0,\eta_0)$ in which
\begin{itemize}
  \item $A_0$ is an abelian scheme over $S$ of relative dimension $[F:\dQ]$,

  \item $i_0\colon E\to\End_S(A_0)\otimes\dQ$ is an $E$-action such that for every $x\in E$,
      \[
      \tr(i_0(x)\res\Lie_S(A_0))=\sum_{u\in\tP_{\r{CM}}}\Tr_{E_u/\dQ_p}(x)
      \]
      holds, where $\tP_{\r{CM}}$ is the fixed subset of $\tP$ (\S\ref{ss:run}(F2)),

  \item $\lambda_0\colon A_0\to A_0^\vee$ is a quasi-polarization under which the Rosati involution coincides with the complex conjugation on $E$ under $i_0$,

  \item $\eta_0\colon W^\delta\otimes_\dQ\dA^\infty\to\rH^{\et}_1(A_0,\dA^\infty)$ is a $K_T$-level structure (see, for example, \cite{LTXZZ}*{Definition~3.5.4}).\footnote{In this article, we have been vague in writing level structures: Strictly speaking, one should choose a geometric point $s$ on every connected component of $S$ and the level structure is a $\pi_1(S,s)$-invariant orbit (with respect to the level subgroup) of isometries concerning the fiber at $s$.}
\end{itemize}
It is known that $\bsigma^\delta(K_T)$ is a nonempty scheme finite \'{e}tale over $\dQ_p$, which admits a natural action by the finite group $T(\dA^\infty)/T(\dQ)K_T$ such that each orbit is a finite Galois \'{e}tale scheme over $\Spec\dQ_p$ with the Galois group $T(\dA^\infty)/T(\dQ)K_T$. We fix such an orbit $\bsigma^\delta_0(K_T)$.

For every neat open compact subgroup $K\subseteq G_m(\dA^\infty_F)$, we consider the moduli problem $\bsigma_m^\delta(K,K_T)$ on $\Sch'_{/\dQ_p}$ as follows: for every $S\in\Sch'_{/\dQ_p}$, $\bsigma_m^\delta(K,K_T)(S)$ is the set of equivalencep classes of octuples $(A_0,i_0,\lambda_0,\eta_0;A,i,\lambda,\eta)$ in which
\begin{itemize}
  \item $(A_0,i_0,\lambda_0,\eta_0)$ is an element of $\bsigma^\delta_0(K_T)(S)$,

  \item $A$ is an abelian scheme over $S$ of relative dimension $2m[F:\dQ]$,

  \item $i\colon E\to\End_S(A)\otimes\dQ$ is an $E$-action such that for every $x\in E$, $\tr(i(x)\res\Lie_S(A))=m\Tr_{E/\dQ}(x)$,

  \item $\lambda\colon A\to A^\vee$ is a quasi-polarization under which the Rosati involution coincides with the complex conjugation on $E$ under $i$,

  \item $\eta\colon W_m^\delta\otimes_E\dA_E^\infty\to\Hom_{\dA_E^\infty}(\rH^{\et}_1(A_0,\dA^\infty),\rH^{\et}_1(A,\dA^\infty))$ is a $K$-level structure, where $W_m^\delta$ denotes the space $E^{2m}$ equipped with the \emph{hermitian} form $\delta^{-1}\cdot\tw_m$ (see, for example, \cite{LTXZZ}*{Definition~4.2.2}).
\end{itemize}
It is known that $\bsigma_m^\delta(K,K_T)$ is a scheme finite type over $\bsigma^\delta_0(K_T)$, which admits a natural lift of the action of $T(\dA^\infty)/T(\dQ)K_T$. We denote by $\bsigma_m^\delta(K,K_T)^\flat$ the quotient of $\bsigma_m^\delta(K,K_T)$ by $T(\dA^\infty)/T(\dQ)K_T$, as a presheaf on $\Sch'_{/\dQ_p}$.

Now we discuss the relation between $\bsigma_m^\delta(K,K_T)^\flat$ and usual Shimura varieties. For every CM type $\bPhi$, we have the Deligne homomorphism
\begin{align*}
\rh_m^{\bPhi}\colon\Res_{\dC/\dR}\bG &\to(\Res_{F/\dQ}G_m)\otimes_\dQ\dR \\
z &\mapsto
\([1_m,(\ol{z}/z)1_m],\cdots,[1_m,(\ol{z}/z)1_m]\)\in K_{m,\infty},
\end{align*}
in which for every archimedean place $v$ of $F$, the notation $[1_m,(\ol{z}/z)1_m]$ is understood via the unique complex embedding of $E$ in $\bPhi$ inducing $v$. Then we obtain a projective system of Shimura varieties $\{\bsigma_m^{\bPhi}(K)\}_K$ associated with the Shimura data $(\Res_{F/\dQ}G_m,\rh_m^{\bPhi})$ indexed by neat open compact subgroups $K\subseteq G_m(\dA_F^\infty)$, which are smooth quasi-projective complex schemes of dimension $m^2[F:\dQ]$, with the complex analytification
\[
\bsigma_m^{\bPhi}(K)^\an= G_m(F)\backslash G_m(\dA_F)/K_{m,\infty}K.
\]

For every embedding $\iota\colon\dQ_p\to\dC$, we denote by $\bPhi_\iota$ the set of complex embeddings $i\colon E\to\dC$ such that the $p$-adic place induced by the embedding $i\colon E\hookrightarrow i(E).\iota(\dQ_p)$ belongs to $\tP_{\r{CM}}$ (\S\ref{ss:run}(F2)). Then $\bPhi_\iota$ is a CM type of $E$.

\begin{lem}
The presheaf $\bsigma_m^\delta(K,K_T)^\flat$ is a scheme over $\dQ_p$ independent of the choices of $K_T$, $\delta$, and the orbit $\bsigma^\delta_0(K_T)$.\footnote{But $\bsigma_m^\delta(K,K_T)^\flat$ depends on the fixed subset $\tP_{\r{CM}}$.} Moreover, for every embedding $\iota\colon\dQ_p\to\dC$, we have a canonical isomorphism
\[
\bsigma_m^\delta(K,K_T)^\flat\otimes_{\dQ_p,\iota}\dC\xrightarrow{\sim}\bsigma_m^{\bPhi_\iota}(K).
\]
\end{lem}

\begin{proof}
By definition, the reflex field $E_{\bPhi_\iota}\subseteq\dC$ of $\bPhi_\iota$ is contained in $\iota(\dQ_p)$. Then there is a canonical isomorphism
\[
\(X_K\otimes_{E_{\bPhi_\iota}}Y_{K_T}\)\otimes_{E_{\bPhi_\iota},\iota^{-1}}\dQ_p\simeq\bsigma_m^\delta(K,K_T)
\]
of schemes over $\dQ_p$, where $X_K$ and $Y_{K_T}$ are the usual Shimura varieties for $G_m$ and $T$ of level $K$ and $K_T$, respectively, over their common reflex field $E_{\bPhi_\iota}$. Under such isomorphism, $T(\dA^\infty)/T(\dQ)K_T$ acts on the left side via the second factor $Y_{K_T}$ whose quotient is nothing but $\Spec E_{\bPhi_\iota}$. Thus, we obtain a canonical isomorphism $X_K\otimes_{E_{\bPhi_\iota},\iota^{-1}}\dQ_p\simeq\bsigma_m^\delta(K,K_T)^\flat$. The lemma follows.
\end{proof}

\begin{definition}
We define the Siegel hermitian variety (of genus $m$ and level $K$) over $\dQ_p$, denoted as $\bsigma_m(K)$, to be $\bsigma_m^\delta(K,K_T)^\flat$, which makes sense by the lemma above.\footnote{By construction, $\bsigma_m(K)$ also depends on the choice of the subset $\tP_{\r{CM}}$ of $\tP$ (\S\ref{ss:run}(F2)).}
\end{definition}

Now we define the \emph{line bundle of automorphy} on $\bsigma_m(K)$. Denote by $\bbA$ (the $A$ part of) the universal object over $\bsigma_m^\delta(K,K_T)$. Then $\Lie(\bbA)$ is a projective $\sO\otimes_\dQ E$-module of rank $m$, where $\sO=\sO_{\bsigma_m^\delta(K,K_T)}$ is the structure sheaf. Put
\[
\bomega_m^\delta\coloneqq\r{det}_\sO\(\r{det}_{\sO\otimes_\dQ E}\Lie(\bbA)^\vee\),
\]
which is a line bundle on $\bsigma_m^\delta(K,K_T)$. Since $T(\dA^\infty)/T(\dQ)K_T$ acts trivially on $\bomega_m^\delta$, $\bomega_m^\delta$ descends to a line bundle $\bomega_m$ on $\bsigma_m(K)$. It is easy too see that $\bomega_m$ does not depend on the choices of $K_T$, $\delta$, and the orbit $\bsigma^\delta_0(K_T)$.

Now suppose that we are given a partition $m=m_1+\cdots+m_s$ of $m$ into positive integers. We have a natural isometry
\begin{align}\label{eq:partition1}
W_{m_1}\oplus\cdots\oplus W_{m_s}\simeq W_m
\end{align}
such that if we write $\{e^j_1,\dots,e^j_{2m_j}\}$ as the standard bases for $W_{m_j}$ for $1\leq j\leq s$, then the standard basis of $W_m$ is identified with
\[
\{e^1_1,\dots,e^1_{m_1},\dots,e^s_1,\dots,e^s_{m_s},e^1_{m_1+1},\dots,e^1_{2m_1},\dots,e^s_{m_s+1},\dots,e^s_{2m_s}\}.
\]
In particular, we may regard $G_{m_1,\dots,m_s}= G_{m_1}\times\cdots\times G_{m_s}$ as a subgroup of $G_m$. We obtain a map
\begin{align}\label{eq:descent4}
\rho_{m_1,\dots,m_s}\colon\cA^{[w]}_{m,\hol}\to\cA^{[w]}_{m_1,\dots,m_s,\hol}
\end{align}
(see \S\ref{ss:run}(G8)) given by the restriction to the subgroup $G_{m_1,\dots,m_s}(\dA_F)$.

For neat open compact subgroups $K_j\subseteq G_{m_j}(\dA_F^\infty)$ for $1\leq j\leq s$, we put
\begin{align*}
\bsigma_{m_1,\dots,m_s}^\delta(K_1\times\cdots\times K_s,K_T)&\coloneqq
\bsigma_{m_1}^\delta(K_1,K_T)\times_{\bsigma^\delta_0(K_T)}\cdots\times_{\bsigma^\delta_0(K_T)}\bsigma_{m_s}^\delta(K_s,K_T), \\
\bomega_{m_1,\dots,m_s}^\delta&\coloneqq\bomega_{m_1}^\delta\boxtimes\cdots\boxtimes\bomega_{m_s}^\delta;
\end{align*}
and
\begin{align*}
\bsigma_{m_1,\dots,m_s}(K_1\times\cdots\times K_s)&\coloneqq\bsigma_{m_1}(K_1)\times_{\dQ_p}\cdots\times_{\dQ_p}\bsigma_{m_s}(K_s), \\
\bomega_{m_1,\dots,m_s} &\coloneqq\bomega_{m_1}\boxtimes\cdots\boxtimes\bomega_{m_s}.
\end{align*}
We have the natural quotient map
\[
\xi_{m_1,\dots,m_s}\colon\bsigma_{m_1,\dots,m_s}^\delta(K_1\times\cdots\times K_s,K_T)\to\bsigma_{m_1,\dots,m_s}(K_1\times\cdots\times K_s)
\]
under which $\xi_{m_1,\dots,m_s}^*\bomega_{m_1,\dots,m_s}\simeq\bomega_{m_1,\dots,m_s}^\delta$.

For a neat open compact subgroup $K\subseteq G_m(\dA_F^\infty)$ containing $K_1\times\cdots\times K_s$, there is a natural morphism
\[
\sigma_{m_1,\dots,m_s}^\delta\colon\bsigma_{m_1,\dots,m_s}^\delta(K_1\times\cdots\times K_s,K_T)\to\bsigma_m^\delta(K,K_T)
\]
sending $((A_0,i_0,\lambda_0,\eta_0;A_j,i_j,\lambda_j,\eta_j))_{1\leq j\leq s}$ to
\[
(A_0,i_0,\lambda_0,\eta_0;A_1\times\cdots\times A_s,(i_1,\dots,i_s),\lambda_1\times\cdots\times\lambda_s,(\eta_1,\dots,\eta_s)).
\]
It is clear that $\sigma_{m_1,\dots,m_s}^\delta$ descends to a morphism
\[
\sigma_{m_1,\dots,m_s}\colon\bsigma_{m_1,\dots,m_s}(K_1\times\cdots\times K_s)\to\bsigma_m(K)
\]
rendering the following diagram
\begin{align}\label{eq:descent3}
\xymatrix{
\bsigma_{m_1,\dots,m_s}^\delta(K_1\times\cdots\times K_s,K_T) \ar[rr]^-{\sigma_{m_1,\dots,m_s}^\delta}\ar[d]_-{\xi_{m_1,\dots,m_s}}
&& \bsigma_m^\delta(K,K_T) \ar[d]^-{\xi_m} \\
\bsigma_{m_1,\dots,m_s}(K_1\times\cdots\times K_s) \ar[rr]^-{\sigma_{m_1,\dots,m_s}}
&& \bsigma_m(K)
}
\end{align}
in $\Sch'_{/\dQ_p}$ commutative. It is independent of the choices of $K_T$, $\delta$, and the orbit $\bsigma^\delta_0(K_T)$. For the line bundles of automorphy, we have $(\sigma_{m_1,\dots,m_s}^\delta)^*\bomega_m^\delta\simeq\bomega_{m_1,\dots,m_s}^\delta$, and hence $\sigma_{m_1,\dots,m_s}^*\bomega_m\simeq\bomega_{m_1,\dots,m_s}$.

For every integer $w\geq 0$, put
\begin{align*}
\cH_{m_1,\dots,m_s}^w(K_1\times\cdots\times K_s)&\coloneqq\rH^0(\bsigma_{m_1,\dots,m_s}(K_1\times\cdots\times K_s),\bomega_{m_1,\dots,m_s}^{\otimes w}), \\
\cH_{m_1,\dots,m_s}^w&\coloneqq\varinjlim_{K_1,\dots,K_s}\cH_{m_1,\dots,m_s}^w(K_1\times\cdots\times K_s).
\end{align*}
For every embedding $\iota\colon\dQ_p\to\dC$, we have an injective map
\begin{align}\label{eq:holomorphic1}
\bbh_{m_1,\dots,m_s}^\iota\colon\cA^{[w]}_{m_1,\dots,m_s,\hol}\to\cH_{m_1,\dots,m_s}^w\otimes_{\dQ_p,\iota}\dC,
\end{align}
which fits into the following commutative diagram
\[
\xymatrix{
\cA^{[w]}_{m,\hol} \ar[rr]^-{\rho_{m_1,\dots,m_s}}_-{\eqref{eq:descent4}}\ar[d]_-{\bbh_m^\iota} && \cA^{[w]}_{m_1,\dots,m_s,\hol} \ar[d]^-{\bbh_{m_1,\dots,m_s}^\iota} \\
\cH_m^w\otimes_{\dQ_p,\iota}\dC \ar[rr]^-{\sigma_{m_1,\dots,m_s}^*} && \cH_{m_1,\dots,m_s}^w\otimes_{\dQ_p,\iota}\dC
}
\]
of complex vector spaces.

\begin{definition}\label{de:holomorphic}
Let the notation be as above.
\begin{enumerate}
  \item We define $\cH_{m_1,\dots,m_s}^{[w]}$ to be the maximal subspace of $\cH_{m_1,\dots,m_s}^w$ such that for every embedding $\iota\colon\dQ_p\to\dC$, $\cH_{m_1,\dots,m_s}^{[w]}\otimes_{\dQ_p,\iota}\dC$ is contained in the image of $\cA^{[w]}_{m_1,\dots,m_s,\hol}$ under $\bbh_{m_1,\dots,m_s}^\iota$.

  \item For every $\varphi\in\cH_{m_1,\dots,m_s}^{[w]}$ and every embedding $\iota\colon\dQ_p\to\dC$, we denote by $\varphi^\iota$ the unique element in $\cA^{[w]}_{m_1,\dots,m_s,\hol}$ such that $\bbh_{m_1,\dots,m_s}^\iota(\varphi^\iota)=\iota\varphi$.
\end{enumerate}
\end{definition}

\begin{remark}\label{re:holomorphic}
We have the following remarks concerning $\cH_{m_1,\dots,m_s}^{[w]}$.
\begin{enumerate}
  \item The inclusion $\cH_{m_1,\dots,m_s}^{[w]}\subseteq\cH_{m_1,\dots,m_s}^w$ is proper in general since in the definition of $\cH_{m_1,\dots,m_s}^w$, we do not impose any growth condition along the boundary.

  \item It is clear that the subspace $\cH_{m_1,\dots,m_s}^{[w]}$ is closed under the action of $G_{m_1,\dots,m_s}(\dA_F^\infty)$. Moreover, in its definition, it suffices to check for \emph{some} embedding $\iota$.

  \item The natural map $\cH_{m_1}^{[w]}\otimes_{\dQ_p}\cdots\otimes_{\dQ_p}\cH_{m_s}^{[w]}\to\cH_{m_1,\dots,m_s}^{[w]}$ given by exterior product is an isomorphism. Indeed, it suffices to check it at every finite level, which is then an isomorphism of \emph{finite-dimensional} $\dQ_p$-vector spaces.
\end{enumerate}
\end{remark}

To end this subsection, we review the notion of analytic $q$-expansion (or Siegel--Fourier expansion).

\begin{definition}\label{de:fourier}
For every ring $R$, we denote by $\SF_{m_1,\dots,m_s}(R)$ the $R$-module of formal power series
\[
\sum_{(T_1,\dots,T_s)\in\Herm_{m_1}(F)^+\times\cdots\times\Herm_{m_s}(F)^+}a_{T_1,\dots,T_s}q^{T_1,\dots,T_s},\quad a_{T_1,\dots,T_s}\in R
\]
in which $a_{T_1,\dots,T_s}$ vanishes unless the entries of $T_1,\dots,T_s$ are in a fractional ideal of $E$. We have a restriction map
\[
\varrho_{m_1,\dots,m_s}\colon\SF_m(R)\to\SF_{m_1,\dots,m_s}(R)
\]
sending
\[
\sum_{T\in\Herm_m(F)^+}a_T q^T
\]
to
\[
\sum_{(T_1,\dots,T_s)\in\Herm_{m_1}(F)^+\times\cdots\times\Herm_{m_s}(F)^+}
\(\sum_{\substack{T\in\Herm_m(F)^+ \\ \partial_{m_1,\dots,m_s}T=(T_1,\dots,T_s)}}a_T\)q^{T_1,\dots,T_s},
\]
where $\partial_{m_1,\dots,m_s}$ is the map from \S\ref{ss:run}(F3). It is an easy exercise to show that the interior summation is always a finite sum.
\end{definition}

For every integer $w\geq 0$, we have a map
\begin{align}\label{eq:expansion1}
\bbq_{m_1,\dots,m_s}^\an\colon\cA_{m_1,\dots,m_s,\hol}^{[w]}&\to\SF_{m_1,\dots,m_s}(\dC) \\
\varphi&\mapsto\sum_{(T_1,\dots,T_s)\in\Herm_{m_1}(F)^+\times\cdots\times\Herm_{m_s}(F)^+}a_{T_1,\dots,T_s}(\varphi) q^{T_1,\dots,T_s} \notag
\end{align}
in which $a_{T_1,\dots,T_s}(\varphi)$ equals
\begin{align*}
\int_{\Herm_{m_1}(F)\backslash\Herm_{m_1}(\dA_F)}\cdots\int_{\Herm_{m_s}(F)\backslash\Herm_{m_s}(\dA_F)}
&\varphi(n(b_1),\dots,n(b_s))\psi_F(\tr T_1 b_1)^{-1} \\
&\times\psi_F(\tr T_s b_s)^{-1}\rd b_1\cdots\rd b_s
\end{align*}
with $\rd b_1,\dots,\rd b_s$ being the Tamagawa measures.

We have the following commutative diagram
\[
\xymatrix{
\cA^{[w]}_{m,\hol} \ar[rr]^-{\rho_{m_1,\dots,m_s}}_-{\eqref{eq:descent4}}\ar[d]_-{\bbq_m^\an} && \cA^{[w]}_{m_1,\dots,m_s,\hol} \ar[d]^-{\bbq_{m_1,\dots,m_s}^\an} \\
\SF_m(\dC) \ar[rr]^-{\varrho_{m_1,\dots,m_s}} && \SF_{m_1,\dots,m_s}(\dC)
}
\]
under restrictions.

We also need an equivariant version of the above constructions for use in \S\ref{ss:4}.

\begin{definition}\label{de:sf1}
For every ring $R$, we denote by $\c{SF}_{m_1,\dots,m_s}(R)$ the $R[G_{m_1,\dots,m_s}(\dA_F^\infty)]$-module
\[
\Map\(G_{m_1,\dots,m_s}(\dA_F^\infty),\SF_{m_1,\dots,m_s}(R)\)
\]
in which $G_{m_1,\dots,m_s}(\dA_F^\infty)$ acts via the right translation. We have an injective $G_{m_1,\dots,m_s}(\dA_F^\infty)$-equivariant complex linear map
\[
\bbq_{m_1,\dots,m_s}^\infty\colon\cA_{m_1,\dots,m_s,\hol}^{[w]}\to\c{SF}_{m_1,\dots,m_s}(C)
\]
such that $\bbq_{m_1,\dots,m_s}^\infty(\varphi)$ sends $g$ to $\bbq_{m_1,\dots,m_s}^\an(g\cdot\varphi)$.
\end{definition}

\subsection{Relation with PEL type moduli spaces}
\label{ss:pel}

In order to show the rationality of some Eisenstein series later, we need the theory of algebraic $q$-expansions. However, since such theory was only developed for PEL type Shimura varieties (in the sense of Kottwitz), we need to study its relation with our Siegel hermitian varieties.

Let $\widetilde{W}_m$ be the space $E^{2m}$ equipped with the pairing
\[
\Tr_{E/\dQ}\circ\langle\;,\;\rangle_{W_{2m}}\colon E^{2m}\times E^{2m}\to\dQ.
\]
Let $\widetilde{G}_m$ be the similitude group of $\widetilde{W}_m$, which is a reductive group over $\dQ$. Let $\widetilde{P}_m\subseteq\widetilde{G}_m$ be the parabolic subgroup stabilizing the subspace generated by $\{e_{m+1},\dots,e_{2m}\}$,

Consider a partition $m=m_1+\cdots+m_s$ of $m$ into positive integers. We denote by $\widetilde{G}_{m_1,\dots,m_s}$ the subgroup of $\widetilde{G}_{m_1}\times\cdots\times\widetilde{G}_{m_s}$ of common similitudes; in other words, it fits into a Cartesian diagram
\[
\xymatrix{
\widetilde{G}_{m_1,\dots,m_s} \ar[r]\ar[d] & \widetilde{G}_{m_1}\times\cdots\times\widetilde{G}_{m_s} \ar[d] \\
\bG_\dQ \ar[r]^{\r{diagonal}} & \bG_\dQ^s
}
\]
in which the vertical arrows are similitude maps. In particular, we may regard $\widetilde{G}_{m_1,\dots,m_s}$ as a subgroup of $\widetilde{G}_m$. Put $\widetilde{P}_{m_1,\dots,m_s}\coloneqq\widetilde{G}_{m_1,\dots,m_s}\cap\widetilde{P}_m$.

For every neat open compact subgroup $\widetilde{K}_{m_1,\dots,m_s}\subseteq\widetilde{G}_{m_1,\dots,m_s}(\dA^\infty)$, we consider the PEL type moduli problem $\widetilde\bsigma_{m_1,\dots,m_s}(\widetilde{K}_{m_1,\dots,m_s})$ on $\Sch'_{/\dQ}$ as follows: for every $S\in\Sch'_{/\dQ}$, $\widetilde\bsigma_{m_1,\dots,m_s}(\widetilde{K}_{m_1,\dots,m_s})(S)$ is the set of equivalence classes of $s$-tuples of quadruples $\((A_1,i_1,\lambda_1,\widetilde\eta_1),\dots,(A_s,i_s,\lambda_s,\widetilde\eta_s)\)$ in which
\begin{itemize}
  \item for $1\leq j\leq s$, $A_j$ is an abelian scheme over $S$ of relative dimension $2m_j[F:\dQ]$,

  \item for $1\leq j\leq s$, $i_j\colon E\to\End_S(A_j)\otimes\dQ$ is an $E$-action such that for every $x\in E$, $\tr(i_j(x)\res\Lie_S(A_j))=m_j\Tr_{E/\dQ}(x)$,

  \item for $1\leq j\leq s$, $\lambda_j\colon A_j\to A_j^\vee$ is a quasi-polarization under which the Rosati involution coincides with the complex conjugation on $E$ under $i_j$,

  \item $\left\{\widetilde\eta_j\colon\widetilde{W}_m\otimes_\dQ\dA^\infty\to\rH^{\et}_1(A_j,\dA^\infty)\right\}_{1\leq j\leq s}$ is a $\widetilde{K}_{m_1,\dots,m_s}$-orbit of skew-hermitian similitudes with similitude factors independent of $j$.
\end{itemize}
Then $\widetilde\bsigma_{m_1,\dots,m_s}(\widetilde{K}_{m_1,\dots,m_s})$ is a scheme of finite type over $\dQ$. Now for a neat open compact subgroup $\widetilde{K}\subseteq\widetilde{G}_m(\dA^\infty)$ containing $\widetilde{K}_{m_1,\dots,m_s}$, we have an obvious morphism
\[
\widetilde\sigma_{m_1,\dots,m_s}\colon\widetilde\bsigma_{m_1,\dots,m_s}(\widetilde{K}_{m_1,\dots,m_s})
\to\widetilde\bsigma_m(\widetilde{K})
\]
over $\dQ$ by ``taking the product of all factors''. For neat open compact subgroups $\widetilde{K}_j\subseteq\widetilde{G}_{m_j}(\dA^\infty)$ containing the image of $\widetilde{K}_{m_1,\dots,m_s}$ under the natural projection map $\widetilde{G}_{m_1,\dots,m_s}\to\widetilde{G}_{m_j}$, we have another obvious map
\[
\tau_{m_1,\dots,m_s}\colon\widetilde\bsigma_{m_1,\dots,m_s}(\widetilde{K}_{m_1,\dots,m_s})\to
\widetilde\bsigma_{m_1}(\widetilde{K}_1)\times_\dQ\cdots\times_\dQ\widetilde\bsigma_{m_s}(\widetilde{K}_s)
\]
over $\dQ$. On $\widetilde\bsigma_m(\widetilde{K})$, we have the line bundle of automorphy $\widetilde\bomega_m$ similar to $\bomega_m^\delta$, which satisfies
\[
\widetilde\sigma_{m_1,\dots,m_s}^*\widetilde\bomega_m\simeq
\tau_{m_1,\dots,m_s}^*\(\widetilde\bomega_{m_1}\boxtimes\cdots\boxtimes\widetilde\bomega_{m_s}\).
\]
Put $\widetilde\bomega_{m_1,\dots,m_s}\coloneqq\widetilde\sigma_{m_1,\dots,m_s}^*\widetilde\bomega_m$ for future use.

\begin{remark}\label{no:pel}
For every $1\leq j\leq s$, we have an isometry $W_{m_j}^\delta\otimes_E W^\delta\xrightarrow\sim\widetilde{W}_{m_j}$. These isometries induce a homomorphism
\begin{align*}
\zeta_{m_1,\dots,m_s}\colon\Res_{F/\dQ}G_{m_1,\dots,m_s}\times T\to\widetilde{G}_{m_1,\dots,m_s}
\end{align*}
sending $(g_1,\dots,g_s,t)$ to $(g_1t,\dots,g_st)$, which is independent of the choice of $\delta$. Using this map, we regard $\Res_{F/\dQ}G_{m_1,\dots,m_s}$ as a subgroup of $\widetilde{G}_{m_1,\dots,m_s}$ in what follows.
\end{remark}

For neat open compact subgroups $K_j\subseteq G_{m_j}(\dA_F^\infty)$ for $1\leq j\leq s$ and $K_T\subseteq T(\dA^\infty)$ such that  $K_1\times\cdots\times K_s\times K_T$ is contained in $\widetilde{K}_{m_1,\dots,m_s}$, we have a natural morphism
\[
\zeta_{m_1,\dots,m_s}\colon\bsigma_{m_1,\dots,m_s}^\delta(K_1\times\cdots\times K_s,K_T)\to
\widetilde\bsigma_{m_1,\dots,m_s}(\widetilde{K}_{m_1,\dots,m_s})\otimes_\dQ\dQ_p
\]
sending $((A_0,i_0,\lambda_0,\eta_0;A_j,i_j,\lambda_j,\eta_j))_{1\leq j\leq s}$ to $((A_j,i_j,\lambda_j,\widetilde\eta_j))_{1\leq j\leq s}$, where $\widetilde\eta_j$ sends $w\otimes v$ to $\eta_j(w)(\eta_0(v))$. The morphism $\zeta_{m_1,\dots,m_s}$ is finite \'{e}tale.

In summary, for every neat open compact subgroup $K\subseteq G_m(\dA_F^\infty)$ containing $K_1\times\cdots\times K_s$ and such that $\zeta_m(K\times K_T)$ is contained in $\widetilde{K}$, we have a diagram
\begin{align}\label{eq:holomorphic5}
\resizebox{\hsize}{!}{
\xymatrix{
\widetilde\bsigma_{m_1}(\widetilde{K}_1)_{\dQ_p}\times_{\dQ_p}\cdots
\times_{\dQ_p}\widetilde\bsigma_{m_s}(\widetilde{K}_s)_{\dQ_p}
& \widetilde\bsigma_{m_1,\dots,m_s}(\widetilde{K}_{m_1,\dots,m_s})_{\dQ_p} \ar[l]_-{\tau_{m_1,\dots,m_s}}\ar[rr]^-{\widetilde\sigma_{m_1,\dots,m_s}}
&& \widetilde\bsigma_m(\widetilde{K})_{\dQ_p} \\
\bsigma_{m_1}^\delta(K_1,K_T)\times_{\bsigma^\delta_0(K_T)}\cdots\times_{\bsigma^\delta_0(K_T)}\bsigma_{m_s}^\delta(K_s,K_T)
\ar[u]^-{\zeta_{m_1}\times\cdots\times\zeta_{m_s}} \ar@{=}[r]^-{\text{def}} \ar[d]_-{\xi_{m_1}\times\cdots\times\xi_{m_s}}
& \bsigma_{m_1,\dots,m_s}^\delta(K_1\times\cdots\times K_s,K_T) \ar[u]^-{\zeta_{m_1,\dots,m_s}}\ar[rr]^-{\sigma_{m_1,\dots,m_s}^\delta}\ar[d]_-{\xi_{m_1,\dots,m_s}}
&& \bsigma_m^\delta(K,K_T) \ar[u]_-{\zeta_m} \ar[d]^-{\xi_m} \\
\bsigma_{m_1}(K_1)\times_{\dQ_p}\cdots\times_{\dQ_p}\bsigma_{m_s}(K_s)  \ar@{=}[r]^-{\text{def}}
& \bsigma_{m_1,\dots,m_s}(K_1\times\cdots\times K_s) \ar[rr]^-{\sigma_{m_1,\dots,m_s}}
&& \bsigma_m(K)
}
}
\end{align}
in $\Sch'_{/\dQ_p}$ expanding \eqref{eq:descent3} as the lower-right square, in which various line bundles of automorphy are compatible under pullbacks.

Similar to $\cA_{m_1,\dots,m_s,\hol}^{[w]}$ (\S\ref{ss:run}(G8)), we define a space $\widetilde\cA_{m_1,\dots,m_s,\hol}^{[w]}$ of certain automorphic forms on $\widetilde{G}_{m_1,\dots,m_s}(\dA)$ with the additional requirement that $(t 1_{m_1},\dots,t 1_{m_s})$ acts trivially for every $t\in T(\dR)$. We have a map
\begin{align}\label{eq:descent10}
\widetilde\rho_{m_1,\dots,m_s}\colon\widetilde\cA^{[w]}_{m,\hol}\to\widetilde\cA^{[w]}_{m_1,\dots,m_s,\hol}
\end{align}
given by the restriction to the subgroup $\widetilde{G}_{m_1,\dots,m_s}(\dA)$.

For every integer $w\geq 0$, put
\begin{align*}
\widetilde\cH_{m_1,\dots,m_s}^w(\widetilde{K}_{m_1,\dots,m_s})
&\coloneqq\rH^0(\widetilde\bsigma_{m_1,\dots,m_s}(\widetilde{K}_{m_1,\dots,m_s}),\widetilde\bomega_{m_1,\dots,m_s}^{\otimes w}), \\
\widetilde\cH_{m_1,\dots,m_s}^w&
\coloneqq\varinjlim_{\widetilde{K}_{m_1,\dots,m_s}}\widetilde\cH_{m_1,\dots,m_s}^w(\widetilde{K}_{m_1,\dots,m_s}).
\end{align*}

\begin{definition}\label{de:pel}
Similar to \eqref{eq:holomorphic1}, we have an injective map
\begin{align*}
\widetilde{\bbh}_{m_1,\dots,m_s}\colon\widetilde\cA_{m_1,\dots,m_s,\hol}^{[w]}
\to\widetilde\cH_{m_1,\dots,m_s}^w\otimes_\dQ\dC
\end{align*}
for $w\geq 0$. We define $\widetilde\cH_{m_1,\dots,m_s}^{[w]}$ to be the subspace of $\widetilde\cH_{m_1,\dots,m_s}^w$ such that the image of $\widetilde{\bbh}_{m_1,\dots,m_s}$ coincides with $\widetilde\cH_{m_1,\dots,m_s}^{[w]}\otimes_\dQ\dC$. Thus, we obtain an isomorphism
\begin{align}\label{eq:holomorphic3}
\widetilde{\bbh}_{m_1,\dots,m_s}\colon\widetilde\cA_{m_1,\dots,m_s,\hol}^{[w]}
\xrightarrow\sim\widetilde\cH_{m_1,\dots,m_s}^{[w]}\otimes_\dQ\dC.
\end{align}
\end{definition}

Now we review the algebraic theory of $q$-expansion for $\widetilde\bsigma_{m_1,\dots,m_s}$ from \cite{Lan12}. Take an open compact subgroup $\widetilde{K}_{m_1,\dots,m_s}\subseteq\widetilde{G}_{m_1,\dots,m_s}(\dA^\infty)$. We choose a smooth projective toroidal compactification $\widetilde\bsigma_{m_1,\dots,m_s}(\widetilde{K}_{m_1,\dots,m_s})^\tor$ of $\widetilde\bsigma_{m_1,\dots,m_s}(\widetilde{K}_{m_1,\dots,m_s})$ over $\dQ$, and let $\widetilde\bomega_{m_1,\dots,m_s}^\tor$ be the canonical extension of $\widetilde\bomega_{m_1,\dots,m_s}$ to $\widetilde\bsigma_{m_1,\dots,m_s}(\widetilde{K}_{m_1,\dots,m_s})^\tor$. Then by \cite{Lan12}*{Definition~5.3.4}, for every $w\geq0$, we have the \emph{algebraic $q$-expansion map}
\begin{align*}
\rH^0(\widetilde\bsigma_{m_1,\dots,m_s}(\widetilde{K}_{m_1,\dots,m_s})^\tor,(\widetilde\bomega_{m_1,\dots,m_s}^\tor)^{\otimes w})\otimes_\dQ\dC\to
\SF_{m_1,\dots,m_s}(\dC)
\end{align*}
(Definition \ref{de:fourier}) at the cusp ``at infinity''. We remark that the map $\bbq_{m_1,\dots,m_s}$ is not necessarily injective, since we only expand the section on the connected component of $\widetilde\bsigma_{m_1,\dots,m_s}(\widetilde{K}_{m_1,\dots,m_s})^\tor\otimes_\dQ\dC$ that contains the cusp ``at infinity''. By \cite{Lan12}*{Remark~5.2.14}, the natural map
\[
\widetilde\cH_{m_1,\dots,m_s}^{[w]}(\widetilde{K}_{m_1,\dots,m_s})
\to\rH^0(\widetilde\bsigma_{m_1,\dots,m_s}(\widetilde{K}_{m_1,\dots,m_s}),(\widetilde\bomega_{m_1,\dots,m_s})^{\otimes w})
\]
(here we adopt a similar notation as in \S\ref{ss:run}(G9)) factors through a map
\begin{align*}
\widetilde\cH_{m_1,\dots,m_s}^{[w]}(\widetilde{K}_{m_1,\dots,m_s})\to
\rH^0(\widetilde\bsigma_{m_1,\dots,m_s}(\widetilde{K}_{m_1,\dots,m_s})^\tor,(\widetilde\bomega_{m_1,\dots,m_s}^\tor)^{\otimes w}),
\end{align*}
hence we obtain a map
\begin{align}\label{eq:expansion}
\bbq_{m_1,\dots,m_s}\colon\widetilde\cH_{m_1,\dots,m_s}^{[w]}(\widetilde{K}_{m_1,\dots,m_s})\otimes_\dQ\dC
\to\SF_{m_1,\dots,m_s}(\dC),
\end{align}
which is independent of the choice of the toroidal compactification. Thus, by passing to the colimit, we obtain a map
\begin{align}\label{eq:holomorphic2}
\bbq_{m_1,\dots,m_s}\colon\widetilde\cH_{m_1,\dots,m_s}^{[w]}\otimes_\dQ\dC\to
\SF_{m_1,\dots,m_s}(\dC),
\end{align}
which fits into the following commutative diagram
\begin{align}\label{eq:holomorphic4}
\xymatrix{
\widetilde\cA^{[w]}_{m,\hol} \ar[rr]^-{\widetilde\rho_{m_1,\dots,m_s}}_-{\eqref{eq:descent10}} \ar[d]_-{\widetilde\bbh_m}
&& \widetilde\cA^{[w]}_{m_1,\dots,m_s,\hol} \ar[d]^-{\widetilde\bbh_{m_1,\dots,m_s}} \\
\widetilde\cH_{m}^{[w]}\otimes_\dQ\dC \ar[rr]^-{\widetilde\sigma_{m_1,\dots,m_s}^*}_-{} \ar[d]_-{\bbq_m}
&& \widetilde\cH_{m_1,\dots,m_s}^{[w]}\otimes_\dQ\dC \ar[d]^-{\bbq_{m_1,\dots,m_s}} \\
\SF_m(\dC) \ar[rr]^-{\varrho_{m_1,\dots,m_s}}_-{\text{Def.~\ref{de:fourier}}} && \SF_{m_1,\dots,m_s}(\dC)
}
\end{align}
of complex vector spaces.

\begin{definition}\label{de:lattice}
Denote by $\fD_E\subseteq O_E$ the different ideal of $E/\dQ$. The (projective) $O_E$-lattice $\cW_m\coloneqq(O_E)^m\oplus (\fD_E^{-1})^m$ of $W_m$ defines an integral model $\cG_m$ (resp.\ $\widetilde\cG_m$) of $G_m$ (resp.\ $\widetilde{G}_m$) over $O_F$ (resp.\ $\dZ$).\footnote{For $v\in\tV_F^\fin$, $\cG_m(O_{F_v})=K_{m,v}$ if and only if $d_v=0$ and $v\not\in\tV_F^\ram$.} Similarly, we have $\cG_{m_1,\dots,m_s}$ and $\widetilde\cG_{m_1,\dots,m_s}$ and their parabolic subgroups $\cP_{m_1,\dots,m_s}$ and $\widetilde\cP_{m_1,\dots,m_s}$, respectively.
\end{definition}

\begin{notation}\label{no:open}
For future use, we introduce some standard open compact subgroups. Take two positive integers $\Delta$ and $\Delta'$ that are coprime to each other. We put
\begin{align*}
\widetilde{K}_{m_1,\dots,m_s}(\Delta,\Delta')&\coloneqq\widetilde\cG_{m_1,\dots,m_s}(\widehat\dZ)
\times_{\widetilde\cG_{m_1,\dots,m_s}(\dZ/\Delta\Delta')}\widetilde\cP_{m_1,\dots,m_s}(\dZ/\Delta),\\
K_{m_1,\dots,m_s}(\Delta,\Delta')&\coloneqq G_{m_1,\dots,m_s}(\dA_F^\infty)\cap\widetilde{K}_{m_1,\dots,m_s}(\Delta,\Delta')
\end{align*}
in view of Remark \ref{no:pel}.
\end{notation}

\begin{lem}\label{le:rational1}
When $\widetilde{K}_{m_1,\dots,m_s}=\widetilde{K}_{m_1,\dots,m_s}(\Delta,\Delta')$, the map \eqref{eq:expansion} is equivariant under $\Aut(\dC/\dQ\langle\Delta'\rangle)$, where we recall that $\dQ\langle\Delta'\rangle\subseteq\dC$ is the subfield generated by $\Delta'^l$-th roots of unity for all $l\geq 1$.
\end{lem}

\begin{proof}
This follows from the fact that the cusp ``at infinity'' is defined over the subfield $\dQ\langle\Delta'\rangle$ at this level structure. See \cite{Lan12} for more details.
\end{proof}

\begin{remark}\label{re:adjoint}
Denote by $\widetilde{G}_{m_1,\dots,m_s}^\der$ the derived subgroup of $\widetilde{G}_{m_1,\dots,m_s}$ and consider the maximal abelian quotient $\widetilde{G}_{m_1,\dots,m_s}^\ab\coloneqq\widetilde{G}_{m_1,\dots,m_s}/\widetilde{G}_{m_1,\dots,m_s}^\der$. Since $\widetilde{G}_{m_1,\dots,m_s}^\der$ is simply connected, for every open compact subgroup $\widetilde{K}_{m_1,\dots,m_s}\subseteq\widetilde{G}_{m_1,\dots,m_s}(\dA^\infty)$, the natural map
\[
\bsigma_{m_1,\dots,m_s}(\widetilde{K}_{m_1,\dots,m_s})(\dC)\to \widetilde{G}_{m_1,\dots,m_s}^\ab(\dQ)\backslash \widetilde{G}_{m_1,\dots,m_s}^\ab(\dA^\infty)/\widetilde{K}_{m_1,\dots,m_s}^\ab
\]
has connected fibers, where $\widetilde{K}_{m_1,\dots,m_s}^\ab$ denotes the image of $\widetilde{K}_{m_1,\dots,m_s}$ in $\widetilde{G}_{m_1,\dots,m_s}^\ab(\dA^\infty)$. It is clear that $\widetilde{K}_{m_1,\dots,m_s}(\Delta,\Delta')^\ab$ depends only on $\Delta'$, which we denote by $\widetilde{K}_{m_1,\dots,m_s}^\ab(\Delta')$.
\end{remark}

\section{Cyclotomic $p$-adic $L$-function}
\label{ss:3}

In this section, we construct the $p$-adic $L$-function. We fix an even positive integer $n=2r$.

\subsection{Doubling space and degenerate principal series}
\label{ss:degenerate}

We have the doubling skew-hermitian space $W_r^\Box\coloneqq W_r\oplus\bar{W}_r$. Let $G_r^\Box$ be the unitary group of $W_r^\Box$, which admits a canonical embedding $\imath\colon G_r\times G_r\hookrightarrow G_r^\Box$. We now take a basis $\{e^\Box_1,\dots,e^\Box_{4r}\}$ of $W_r^\Box$ by the formula
\[
e^\Box_i=e_i,\quad e^\Box_{r+i}=-\bar{e}_i,\quad e^\Box_{2r+i}=e_{r+i}, \quad e^\Box_{3r+i}=\bar{e}_{r+i}
\]
for $1\leq i\leq r$, under which we may identify $W_r^\Box$ with $W_{2r}$ and $G_r^\Box$ with $G_{2r}$. Put $\tw_r^\Box\coloneqq\tw_{2r}$, $P_r^\Box\coloneqq P_{2r}$ and $N_r^\Box\coloneqq N_{2r}$. We denote by
\[
\delta^\Box_r\colon P_r^\Box\to\bG_F
\]
the composition of the Levi quotient map $P_r^\Box=P_{2r}\to M_{2r}$, the isomorphism $m^{-1}\colon M_{2r}\to\Res_{E/F}\GL_{2r}$, the determinant $\Res_{E/F}\GL_{2r}\to\Res_{E/F}\bG$ and the norm $\Nm_{E/F}\colon\Res_{E/F}\bG\to\bG_F$. Put
\begin{align}\label{eq:wr}
\bw_r\coloneqq
\begin{pmatrix}
     &  & 1_r &  \\
     & 1_r &  &  \\
    -1_r & 1_r &  &  \\
     &  & 1_r & 1_r \\
\end{pmatrix}
\in G_r^\Box(F).
\end{align}
Then $P_r^\Box\cdot\bw_r\cdot\imath(G_r\times G_r)$ is Zariski open in $G_r^\Box$.

In what follows, we will regard $G_r\times G_r$ as a subgroup of $G_{2r}=G_r^\Box$ via the isometry \eqref{eq:partition1}, which is precisely the embedding
\begin{align}\label{eq:embedding}
\(
\begin{pmatrix}
  a_1 & b_1 \\
  c_1 & d_1
\end{pmatrix},
\begin{pmatrix}
  a_2 & b_2 \\
  c_2 & d_2
\end{pmatrix}
\)
\mapsto
\begin{pmatrix}
    a_1 &  & b_1 &  \\
     & a_2 &  & b_2 \\
    c_1 &  & d_1 &  \\
     & c_2 & & d_2 \\
\end{pmatrix}.
\end{align}

\begin{remark}\label{re:dagger}
The embedding $\imath\colon G_r\times G_r\hookrightarrow G^\Box_r=G_{2r}$ coincides with the embedding \eqref{eq:embedding} twisted by the involution $\id\times\dag$ on $G_r\times G_r$.
\end{remark}

Let $\chi\colon \Gamma_{F,p}\to\dC^\times$ be a finite character, regarded as an automorphic character of $\dA_F^\times$. For every place $v$ of $F$, we have the degenerate principal series of $G_r^\Box(F_v)$, which is defined as the normalized induced representation
\[
\rI^\Box_{r,v}(\chi_v)\coloneqq\Ind_{P_r^\Box(F_v)}^{G_r^\Box(F_v)}(\chi_v\circ\delta^\Box_{r,v})
\]
of $G_r^\Box(F_v)$ with complex coefficients. For every $f\in\rI^\Box_{r,v}(\chi_v)$ and every $T^\Box\in\Herm_{2r}^\circ(F_v)$, we can regularize the following integral
\begin{align}\label{eq:whittaker}
W_{T^\Box}(f)\coloneqq\int_{\Herm_{2r}(F_v)}f(\tw_r^\Box n(b))\psi_{F,v}(\tr T^\Box b)^{-1}\rd b,
\end{align}
where $\r{d}b$ is the self-dual measure on $\Herm_{2r}(F_v)$ with respect to $\psi_{F,v}$. Indeed, one has a family of integrals $W_{T^\Box}(f_s)$ for $s\in\dC$, where $f_s\in\rI^\Box_{r,v}(\chi_v|\;|_{F_v}^s)$ is the standard section induced by $f$; it is absolutely convergent when $\RE s$ is large enough and has an analytic continuation to $\dC$. Then $W_{T^\Box}(f)$ is defined as the value at $0$ of this analytic continuation. See \cite{Wal88}*{Theorem~8.1} and \cite{Kar79}*{Corollary~3.6.1} for more details.

In order to show the rationality of our $p$-adic $L$-function, we need to extend the degenerate principal series to $\widetilde{G}_{2r}$ along the natural inclusion $\Res_{F/\dQ}G^\Box_r=\Res_{F/\dQ}G_{2r}\hookrightarrow\widetilde{G}_{2r}$. We have a map
\[
s\colon\bG_\dQ\to\widetilde{G}_{2r}
\]
sending $c$ to $\(\begin{smallmatrix} c1_{2r}& \\  & 1_{2r} \end{smallmatrix}\)$. Then the natural map $\Res_{F/\dQ}P_{2r}\times s(\bG_\dQ)\to\widetilde{P}_{2r}$ is an isomorphism.

Take a place $w$ of $\dQ$. Put
\begin{align*}
\psi_{F,w}&\coloneqq\prod_{v\in\tV_F^{(w)}}\psi_{F,v},\\
\chi_w&\coloneqq\prod_{v\in\tV_F^{(w)}}\chi_v,\\
\rI_{r,w}^\Box(\chi_w)&\coloneqq\bigotimes_{v\in\tV_F^{(w)}}\rI^\Box_{r,v}(\chi_v),
\end{align*}
and
\[
\delta^\Box_{r,w}\coloneqq\prod_{v\in\tV_F^{(w)}}\delta^\Box_{r,v}\colon
\prod_{v\in\tV_F^{(w)}}P^\Box_r(F_v)=(\Res_{F/\dQ}P^\Box_r)(\dQ_w)\to(F_w)^\times.
\]
The map $\delta^\Box_{r,w}$ extends uniquely to a map $\widetilde\delta^\Box_{r,w}$ along the inclusion $(\Res_{F/\dQ}P^\Box_r)(\dQ_w)=(\Res_{F/\dQ}P_{2r})(\dQ_w)\subseteq\widetilde{P}_{2r}(\dQ_w)$ that sends $s(c)$ to $c^{2r}$ for $c\in\dQ_w^\times$. Then we have a canonical isomorphism
\[
\rI_{r,w}^\Box(\chi_w)\simeq
\Ind_{\widetilde{P}_{2r}(\dQ_w)}^{\widetilde{G}_{2r}(\dQ_w)}(\chi_w\circ\widetilde\delta^\Box_{r,w})
\]
so that $\rI_{r,w}^\Box(\chi_w)$ becomes a representation of $\widetilde{G}_{2r}(\dQ_w)$. For every $T^\Box\in\Herm_{2r}^\circ(F_w)$, we define the functional $W_{T^\Box}(-)$ on $\rI_{r,w}^\Box(\chi_w)$ to be the product of the corresponding ones over $v\in\tV_F^{(w)}$.

\begin{lem}\label{le:section6}
For every $v\in\tV_F^{(\infty)}$, denote by $f_v^{[r]}\in\rI^\Box_{r,v}(\chi_v)=\rI^\Box_{r,v}(\b{1})$ the unique section whose restriction to $K_{2r,v}$ is the character $\kappa_{2r,v}^r$. Put $f_\infty^{[r]}\coloneqq\otimes_{v\in\tV_F^{(\infty)}}f_v^{[r]}$. Then there exists $W_{2r}\in\dQ_{>0}$ such that
\[
W_{T^\Box}(f_\infty^{[r]})=W_{2r}\cdot b_{2r}^\infty(\CF)\cdot \exp(-2\bpi\Tr_{F/\dQ}\tr T^\Box)
\]
for every $T^\Box\in\Herm_{2r}^\circ(F)^+$.
\end{lem}

\begin{proof}
For two elements $x,y\in\dC^\times$, we write $x\sim y$ if their quotient is rational.

By \cite{Liu11}*{Proposition~4.5(2)}, we have
\[
W_{T^\Box}(f_\infty^{[r]})=\(\frac{(2\bpi)^{r(2r+1)}}{\Gamma(1)\Gamma(2)\cdots\Gamma(2r)}\)^{[F:\dQ]}
\exp(-2\bpi\Tr_{F/\dQ}\tr T^\Box)
\]
for every $T^\Box\in\Herm_{2r}^\circ(F)^+$. The positivity of $W_{2r}$ then follows. Thus, it remains to show that $b_{2r}^\infty(\CF)\sim\bpi^{r(2r+1)[F:\dQ]}$.

Write $L(s,\eta_{E/F}^i)$ for the complete $L$-function for the self-dual character $\eta_{E/F}^i$. Then by the functional equation, we have
\[
\prod_{i=1}^{2r}L(i,\eta_{E/F}^i)\sim\prod_{i=1}^{2r}L(1-i,\eta_{E/F}^i).
\]
By a well-known result of Siegel, $\prod_{i=1}^{2r}L^\infty(1-i,\eta_{E/F}^i)$ is rational. It follows that
\begin{align*}
b_{2r}^\infty(\CF)&\sim\frac{\prod_{i=1}^{2r}L_\infty(1-i,\eta_{E/F}^i)}{\prod_{i=1}^{2r}L_\infty(i,\eta_{E/F}^i)}
=\(\frac{\prod_{i=1}^{2r}L_\dR(1-i,\r{sgn}^i)}{\prod_{i=1}^{2r}L_\dR(i,\r{sgn}^i)}\)^{[F:\dQ]} \\
&\sim\(\frac{\bpi^{r^2}}{\bpi^{-r(r+1)}}\)^{[F:\dQ]}=\bpi^{r(2r+1)[F:\dQ]}.
\end{align*}
The lemma follows.
\end{proof}

From now to the end of this subsection, we assume $w\neq\infty$.

\begin{lem}\label{le:section1}
We have
\begin{enumerate}
  \item For $v\in\tV_F^{(w)}$ and $b\in\Herm_{2r}(F_v)$, the relation
      \[
      W_{T^\Box}(n(b)f)=\psi_{F,v}(\tr T^\Box b)\cdot W_{T^\Box}(f)
      \]
      holds for every $f\in\rI^\Box_{r,v}(\chi_v)$ and every $T^\Box\in\Herm_{2r}^\circ(F_v)$.

  \item For $v\in\tV_F^{(w)}$ and $a\in\GL_{2r}(E_v)$, the relation
      \[
      W_{T^\Box}(m(a)f)=\chi_v(\Nm_{E_v/F_v}\det a)^{-1}|\dtm a|_{E_v}^r\cdot W_{\pres{\rt}{a}^\tc T^\Box a}(f)
      \]
      holds for every $f\in\rI^\Box_{r,v}(\chi_v)$ and every $T^\Box\in\Herm_{2r}^\circ(F_v)$.

  \item For $c\in\dQ_w^\times$, the relation
      \[
      W_{T^\Box}(s(c)f)=\chi_w(c)^{-2r}|c|_{F_v}^{2r^2}\cdot W_{cT^\Box}(f)
      \]
      holds for every $f\in\rI^\Box_{r,w}(\chi_w)$ and every $T^\Box\in\Herm_{2r}^\circ(F_w)$.
\end{enumerate}
\end{lem}

\begin{proof}
This is well-known. For readers' convenience, we give a (formal) proof.

For (1), we have
\begin{align*}
&\quad W_{T^\Box}(n(b)f) \\
&=\int_{\Herm_{2r}(F_v)}f(\tw_r^\Box n(b')n(b))\psi_{F,v}(\tr T^\Box b')^{-1}\rd b' \\
&=\psi_{F,v}(\tr T^\Box b)\int_{\Herm_{2r}(F_v)}f(\tw_r^\Box n(b'+b))\psi_{F,v}(\tr T^\Box(b'+b))^{-1}\rd b' \\
&=\psi_{F,v}(\tr T^\Box b)\cdot W_{T^\Box}(f).
\end{align*}

For (2), we have
\begin{align*}
&\quad W_{T^\Box}(m(a)f) \\
&=\int_{\Herm_{2r}(F_v)}f(\tw_r^\Box n(b)m(a))\psi_{F,v}(\tr T^\Box b)^{-1}\rd b \\
&=\int_{\Herm_{2r}(F_v)}f(\tw_r^\Box m(a) n(a^{-1}b \pres{\rt}{a}^{\tc,-1}))\psi_{F,v}(\tr T^\Box b)^{-1}\rd b \\
&=\int_{\Herm_{2r}(F_v)}f(m(\pres{\rt}{a}^{\tc,-1})\tw_r^\Box n(a^{-1}b \pres{\rt}{a}^{\tc,-1}))
\psi_{F,v}(\tr (\pres{\rt}{a}^\tc T^\Box a)(a^{-1}b \pres{\rt}{a}^{\tc,-1}))^{-1}\rd b \\
&=\chi_v(\Nm_{E_v/F_v}\det a)^{-1}|\dtm a|_{E_v}^r
\int_{\Herm_{2r}(F_v)}f(\tw_r^\Box n(b))\psi_{F,v}(\tr (\pres{\rt}{a}^\tc T^\Box a)b)^{-1}\rd b \\
&=\chi_v(\Nm_{E_v/F_v}\det a)^{-1}|\dtm a|_{E_v}^r\cdot W_{\pres{\rt}{a}^\tc T^\Box a}(f).
\end{align*}

The proof for (3) is similar to (2) and we omit it. The lemma is proved.
\end{proof}

\begin{notation}\label{no:section1}
Let $v\in\tV_F^\fin$ be a finite place.
\begin{enumerate}
  \item We denote by $\rI^\Box_{r,v}(\chi_v)^\circ$ the subspace of $\rI^\Box_{r,v}(\chi_v)$ consisting of sections that are supported on the big Bruhat cell $P_r^\Box(F_v)\cdot\tw_r^\Box\cdot N^\Box_r(F_v)$.

  \item When $v\in\tV_F^\fin\setminus\tV_F^{(p)}$, we denote by $f_{\chi_v}^\sph\in\rI^\Box_{r,v}(\chi_v)$ the unique section that takes value $1$ on $K_{2r,v}$.
\end{enumerate}
\end{notation}

It is clear that $\rI^\Box_{r,v}(\chi_v)^\circ$ is stable under the action of $P_r^\Box(F_v)$. For $f\in\rI^\Box_{r,v}(\chi_v)^\circ$ and $T^\Box\in\Herm_{2r}(F_v)$, we put
\[
W_{T^\Box}(f)\coloneqq\int_{\Herm_{2r}(F_v)}f(\tw_r^\Box n(b))\psi_{F,v}(\tr T^\Box b)^{-1}\rd b,
\]
which is in fact a finite sum and coincides with \eqref{eq:whittaker} for $T^\Box\in\Herm_{2r}^\circ(F_v)$. It is clear that the assignment $T^\Box\mapsto W_{T^\Box}(f)$ is a Schwartz function on $\Herm_{2r}(F_v)$. Conversely, using the Fourier inversion formula, we know that for every $\tf\in\sS(\Herm_{2r}(F_v))$, there exists a unique section $\tf^{\chi_v}\in\rI^\Box_{r,v}(\chi_v)^\circ$ such that $W_{T^\Box}(\tf^{\chi_v})=\tf(T^\Box)$ holds for every $T^\Box\in\Herm_{2r}(F_v)$. In other words, we obtain a bijection
\begin{align}\label{eq:section}
-^{\chi_v}\colon \sS(\Herm_{2r}(F_v))\xrightarrow{\sim}\rI^\Box_{r,v}(\chi_v)^\circ.
\end{align}
Put $\rI_{r,w}^\Box(\chi_w)^\circ\coloneqq\bigotimes_{v\in\tV_F^{(w)}}\rI^\Box_{r,v}(\chi_v)^\circ$ and
we obtain an isomorphism
\[
-^{\chi_w}\colon \sS(\Herm_{2r}(F_w))\xrightarrow{\sim}\rI^\Box_{r,w}(\chi_w)^\circ
\]
by taking product over $v\in\tV_F^{(w)}$.

\begin{lem}\label{le:section4}
Suppose that (the rational prime) $w\neq p$.
\begin{enumerate}
  \item For every $v\in\tV_F^{(w)}\setminus\tV_F^\ram$ and every $g\in G_{2r}(F_v)$, there exists a finitely generated ring $\dO_g$ contained in $\dZ_{(p)}\langle w\rangle$ such that for every $T^\Box\in\Herm_{2r}^\circ(F_v)$, there exists a unique element $\pres{g}\sfW_{T^\Box,v}^\sph\in\dO_g[X,X^{-1}]$ such that
      \[
      \pres{g}\sfW_{T^\Box,v}^\sph(\chi_v(\varpi_v))=b_{2r,v}(\chi)\cdot W_{T^\Box}(g\cdot f_{\chi_v}^\sph)
      \]
      holds for every finite character $\chi\colon \Gamma_{F,p}\to\dC^\times$, where $\varpi_v$ is an arbitrary uniformizer of $F_v$. Moreover, $\sfW_{T^\Box,v}^\sph\coloneqq\pres{1_{4r}}\sfW_{T^\Box,v}^\sph\in\dZ[X]$.

  \item For every $f\in\rI^\Box_{r,w}(\chi_w)$ and every $T^\Box\in\Herm_{2r}^\circ(F_w)$, we have
      \[
      W_{T^\Box}(\sigma f)=\sigma W_{T^\Box}(f)
      \]
      for $\sigma\in\Aut(\dC/\dQ\langle w\rangle)$.

  \item For every $f\in\rI^\Box_{r,w}(\chi_w)^\circ$ that is fixed by $\widetilde\cP_{2r}(\dZ_w)$ and every $T^\Box\in\Herm_{2r}^\circ(F_w)$, we have
      \[
      W_{T^\Box}(\sigma f)=\sigma W_{T^\Box}(f)
      \]
      for $\sigma\in\Aut(\dC/\dQ)$.
\end{enumerate}
\end{lem}

\begin{proof}
For (1), by Lemma \ref{le:section1}(1,2) and the Iwasawa decomposition $G_{2r}(F_v)=P_{2r}(F_v)K_{2r,v}$, it suffices to consider the case where $g=1_{4r}$. Then the statement follows from \cite{LZ}*{Theorem~3.5.1}, together with the discussion in \cite{LZ}*{Sections~3.2~\&~3.3}.\footnote{Though \cite{LZ} only treats the case where $v$ is inert in $E$, the same argument works in the case where $v$ splits in $E$ as well.}

Part (2) follows from the proof of \cite{Kar79}*{Corollary~3.6.1} and the fact that $\psi_{F,w}$ takes values in $\dQ\langle w\rangle$.

For (3), put $\fO_b\coloneqq\{c\pres{\rt}{a}^\tc b a\res a\in\GL_{2r}(O_{E_w}),c\in\dZ_w^\times\}$ for every $b\in\Herm_{2r}(F_w)$, which is an open compact subset of $\Herm_{2r}(F_w)$. It follows easily that
\[
c_{T^\Box,\fO_b}\coloneqq\int_{\fO_b}\psi_{F,v}(\tr T^\Box b')^{-1} \rd b'\in\dQ.
\]
Since $\chi_w$ is unramified, the assignment $b'\mapsto f(\tw_r^\Box n(b'))$ is constant on each subset $\fO_b$, which we denote as $f_{\fO_b}$. Then $(\sigma f)_{\fO_b}=\sigma f_{\fO_b}$. It follows that
\begin{align*}
W_{T^\Box}(\sigma f)&=\sum_{\fO}c_{T^\Box,\fO}\cdot(\sigma f)_{\fO}=\sum_{\fO}c_{T^\Box,\fO}\cdot\sigma f_{\fO} \\
&=\sigma\sum_{\fO}c_{T^\Box,\fO}\cdot f_{\fO}=\sigma W_{T^\Box}(f)
\end{align*}
in which the sum is taken over a finite set of disjoint open compact subset of $\Herm_{2r}(F_w)$ of the form $\fO_b$. Thus, (3) follows.
\end{proof}

\begin{lem}\label{le:section2}
The representation $\rI^\Box_{r,w}(\chi_w)$ is semisimple and of finite length as a representation of $\widetilde{G}_{2r}(\dQ_w)$. When $w\neq p$, every irreducible summand of $\rI^\Box_{r,w}(\chi_w)$ contains a nonzero element $f$ in $\rI^\Box_{r,w}(\chi_w)^\circ$ that is fixed by $\widetilde\cP_{2r}(\dZ_w)$.
\end{lem}

\begin{proof}
The first statement follows since it is the parabolic induction of a unitary character.

Now we show the second statement. For every $v\in\tV_F^{(v)}$, by \cite{KS97}*{Theorem~1.2 \& Theorem~1.3}, $\rI^\Box_{r,v}(\chi_v)$ is an irreducible representation of $G_{2r}(F_v)$ unless $\chi_v^2=\b{1}$. Moreover, when $\chi_v^2=\b{1}$, each direct summand of $\rI^\Box_{r,v}(\chi_v)$ is of the form $\rI(V_v)$ for some (nondegenerate) hermitian space $V_v$ over $E_v$ of rank $2r$. Here, $\rI(V_v)$ is the image of the Siegel--Weil section map $\sS(V_v^{2r})\to\rI^\Box_{r,v}(\chi_v)$ under the Weil representation with respect to (the standard additive character $\psi_{F,v}$ and) the splitting character $\chi_v\circ\Nm_{E_v/F_v}$ (again see \cite{KS97}). Put $\tV\coloneqq\{v\in\tV_F^{(w)}\res\chi_v^2=\b{1}\}$.

Now let $\rI$ be an irreducible summand of $\rI^\Box_{r,w}(\chi_w)$ as a representation of $\widetilde{G}_{2r}(\dQ_w)$. One can find a collection of hermitian spaces $V_v$ over $E_v$ of rank $2r$ for $v\in\tV$ such that $\rI$ contains
\[
\(\bigotimes_{v\in\tV}\rI(V_v)\)\otimes\(\bigotimes_{v\in\tV_F^{(v)}\setminus\tV}\rI^\Box_{r,v}(\chi_v)\).
\]
For every $v\in\tV_F^{(w)}$, we define a subset $\fT_v$ of $\Herm_{2r}^\circ(F_v)$ as follows. If $v\in\tV$, then we define $\fT_v$ to be the intersection of $\Herm_{2r}^\circ(F_v)$ and the image of the moment map $V_v^{2r}\to\Herm_{2r}(F_v)$ (see \S\ref{ss:setup}(H1) if one needs recall). If $v\not\in\tV$, then we define $\fT_v$ to be $\Herm_{2r}^\circ(F_v)$. Take any open compact subset $\fT$ of $\Herm_{2r}(F_w)=\prod_{v\in\tV_F^{(w)}}\Herm_{2r}(F_v)$ that is contained in $\prod_{v\in\tV_F^{(w)}}\fT_v\cap\Herm_{2r}(O_{F_v})$ satisfying that $c\pres{\rt}{a}^\tc \fT a=\fT$ for every $a\in\GL_{2r}(O_{E,w})$ and every $c\in\dZ_w^\times$. Then $(\CF_\fT)^{\chi_w}\in\rI^\Box_{r,w}(\chi_w)^\circ$ is a nonzero element of $\rI$. Moreover, by Lemma \ref{le:section1}, it is fixed by $\widetilde\cP_{2r}(\dZ_w)$.

The lemma is proved.
\end{proof}

In the rest of this subsection, we construct some explicit sections in $\rI^\Box_{r,p}(\chi_p)^\circ$.

\begin{notation}\label{no:section2}
For every place $v\in\tV_F^{(p)}$, we
\begin{itemize}
  \item fix a uniformizer $\varpi_v$ of $F_v$,

  \item for every element $e=(e_u)_u\in\dZ^{\tP_v}$, put $|e|\coloneqq\sum_{u\in\tP_v}e_u$ and denote by $\varpi_v^e$ the element in $E_v=\prod_{u\in\tP_v}E_u$ whose component in $E_u$ is $\varpi_v^{e_u}$,

  \item for $u\in\tP_v$, denote by $1_u\in\dZ^{\tP_v}$ the element that takes values $1$ at $u$ and $0$ at $u^\tc$,

  \item for every $u\in\tP_v$, introduce an element
      \[
      \rU_u\coloneqq
      \sum_{b\in\Herm_r(O_{F_v}/\varpi_v)}
      \left[
      \begin{pmatrix}
      1_r & \varpi_v^{-d_v} b^\sharp \\
      & 1_r
      \end{pmatrix}
      \begin{pmatrix}
      \varpi_v^{1_u}\cdot 1_r & \\
      & \varpi_v^{-1_{u^\tc}}\cdot 1_r
      \end{pmatrix}
      \right]
      \in\dZ[G_r(F_v)],
      \]
      where $b^\sharp\in\Herm_r(O_{F_v})$ denotes the Teichm\"{u}ller lift of $b$,

  \item for every $e=(e_u)_u\in\dN^{\tP_v}$, define
      \[
      \rU_v^e\coloneqq\prod_{u\in\tP_v}\rU_u^{e_u}\in\dZ[G_r(F_v)],
      \]
      where we note that the subalgebra of $\dZ[G_r(F_v)]$ generated by $\rU_u$ for $u\in\tP_v$ is commutative.
\end{itemize}
\end{notation}

\begin{construction}\label{co:section}
For $v\in\tV_F^{(p)}$ and every element $e\in\dZ^{\tP_v}$, put
\[
\fT_v^{[e]}\coloneqq\left\{\left.
T^\Box=\begin{pmatrix} T^\Box_{11}& T^\Box_{12} \\ T^\Box_{21} & T^\Box_{22} \end{pmatrix}\right|
T^\Box_{11},T^\Box_{22}\in\Herm_r(O_{F_v}),T^\Box_{12}\in\varpi_v^{-e}\cdot\GL_r(O_{E_v})
\right\}
\]
as a subset of $\Herm_{2r}(F_v)$. Define a function $\tf_{\chi_v}^{[e]}\in\sS(\Herm_{2r}(F_v))$ by the formula
\[
\tf_{\chi_v}^{[e]}(T^\Box)\coloneqq\chi_v(\Nm_{E_v/F_v}\det T^\Box_{12})\cdot \CF_{\fT_v^{[e]}}(T^\Box).
\]
In particular, we obtain a section $(\tf_{\chi_v}^{[e]})^{\chi_v}\in\rI^\Box_{r,v}(\chi_v)^\circ$ by \eqref{eq:section}.
\end{construction}

In what follows, we will identity $\dZ^\tP$ and $\dN^\tP$ with $\prod_{v\in\tV_F^{(p)}}\dZ^{\tP_v}$ and $\prod_{v\in\tV_F^{(p)}}\dN^{\tP_v}$, respectively. For $e\in\dZ^\tP$, we put
\begin{align*}
\|e\|&\coloneqq\max_{v\in\tV_F^{(p)}}|e_v|,\\
\fT_p^{[e]}&\coloneqq\prod_{v\in\tV_F^{(p)}}\fT_v^{[e_v]},\\
\tf_{\chi_p}^{[e]}&\coloneqq\bigotimes_{v\in\tV_F^{(p)}}\tf_{\chi_v}^{[e_v]}.
\end{align*}
For $e\in\dN^\tP$, we put
\[
\rU_p^e\coloneqq\bigotimes_{v\in\tV_F^{(p)}}\rU_v^{e_v}\in
\bigotimes_{v\in\tV_F^{(p)}}\dZ[G_r(F_v)]=\dZ[G_r(F\otimes\dZ_p)].
\]
For two elements $e_1,e_2\in\dN^\tP$, we have the element $\rU_p^{e_1}\times\rU_p^{e_2}$ as the image of $\rU_p^{e_1}\otimes\rU_p^{e_2}$ under the natural map $\dZ[G_r(F\otimes\dZ_p)]\otimes\dZ[G_r(F\otimes\dZ_p)]\to\dZ[G_{2r}(F\otimes\dZ_p)]$ induced by the embedding \eqref{eq:embedding}.

\begin{example}\label{ex:operator}
Suppose that $F=\dQ$ and write $\tP=\{u,u^\tc\}$. If we take $\varpi_p=p$ and identify $G_{2r}(\dQ_p)$ with $\GL_{4r}(\dQ_p)$ via $u$, then
\begin{align*}
\rU_p^{1_u}\times\rU_p^0&=
\sum_{b\in\Herm_r(\dF_p)}
\left[
\begin{pmatrix}
 1_r &  & b^\sharp &  \\
     & 1_r &  &  \\
     &  & 1_r &  \\
     &  & & 1_r \\
\end{pmatrix}
\begin{pmatrix}
 p1_r &  &  &  \\
     & 1_r &  &  \\
     &  & 1_r &  \\
     &  & & 1_r \\
\end{pmatrix}
\right] \\
\rU_p^{1_{u^\tc}}\times\rU_p^0&=
\sum_{b\in\Herm_r(\dF_p)}
\left[
\begin{pmatrix}
 1_r &  & b^\sharp &  \\
     & 1_r &  &  \\
     &  & 1_r &  \\
     &  & & 1_r \\
\end{pmatrix}
\begin{pmatrix}
 1_r &  &  &  \\
     & 1_r &  &  \\
     &  & p^{-1}1_r &  \\
     &  & & 1_r \\
\end{pmatrix}
\right] \\
\rU_p^0\times\rU_p^{1_u}&=
\sum_{b\in\Herm_r(\dF_p)}
\left[
\begin{pmatrix}
 1_r &  &  &  \\
     & 1_r &  & b^\sharp \\
     &  & 1_r &  \\
     &  & & 1_r \\
\end{pmatrix}
\begin{pmatrix}
 1_r &  &  &  \\
     & p1_r &  &  \\
     &  & 1_r &  \\
     &  & & 1_r \\
\end{pmatrix}
\right] \\
\rU_p^0\times\rU_p^{1_{u^\tc}}&=
\sum_{b\in\Herm_r(\dF_p)}
\left[
\begin{pmatrix}
 1_r &  & &  \\
     & 1_r &  & b^\sharp \\
     &  & 1_r &  \\
     &  & & 1_r \\
\end{pmatrix}
\begin{pmatrix}
 1_r &  &  &  \\
     & 1_r &  &  \\
     &  &1_r &  \\
     &  & & p^{-1}1_r \\
\end{pmatrix}
\right]
\end{align*}
and the general ones $\rU_p^{e_1}\times\rU_p^{e_2}$ are the composition of the above four.
\end{example}

\begin{lem}\label{le:section3}
For every element $e\in\dZ^\tP$, the section $(\tf_{\chi_p}^{[e]})^{\chi_p}\in\rI^\Box_{r,p}(\chi_p)^\circ$ is invariant under the subgroup $\widetilde\cP_{r,r}(\dZ_p)$ (Definition \ref{de:lattice}) of $\widetilde{G}_{2r}(\dQ_p)$.
\end{lem}

\begin{proof}
This follows immediately from the construction of $\tf_{\chi_p}^{[e]}$.
\end{proof}

\begin{lem}\label{le:section5}
For every element $e\in\dZ^\tP$ and every $e_1,e_2\in\dN^\tP$, we have
\[
(\rU_p^{e_1}\times\rU_p^{e_2})(\tf_{\chi_p}^{[e]})^{\chi_p}=
(\tf_{\chi_p}^{[e+e_1^\tc+e_2]})^{\chi_p},
\]
where $e_1^\tc\coloneqq e_1\circ\tc$.
\end{lem}

\begin{proof}
By induction, we may assume either $e_1=0$ or $e_2=0$. We consider the case where $e_2=0$ and leave the other similar case to the reader. Again by induction, we may assume $e_1=1_u$ for some $u\in\tP$, with $v\in\tV_F^{(p)}$ its underlying place.

For two square matrices $a$ and $b$, we write $[a,b]$ for the block diagonal matrix. As an element in $\dZ[\widetilde{G}_{2r}(\dQ_p)]$, we have
\[
\rU_p^{e_1}\times\rU_p^{e_2}=\sum_{b\in\Herm_r(O_{F_{v}}/\varpi_{v})}
\left[n([\varpi_{v}^{-d_{v}}\cdot b^\sharp,0_r])\cdot m([\varpi_v^{1_u}\cdot 1_r,1_r])\right]
\]
in which all components away from $v$ are $1_{4r}$. By Lemma \ref{le:section1}, we have
\begin{align}\label{eq:section1}
&\quad W_{T^\Box}((\rU_p^{e_1}\times\rU_p^{e_2})(\tf_{\chi_p}^{[e]})^{\chi_p}) \\
&=\sum_{b\in\Herm_r(O_{F_{v}}/\varpi_{v})}
W_{T^\Box}(n([\varpi_{v}^{-d_{v}}\cdot b^\sharp,0_r])\cdot m([\varpi_v^{1_u}\cdot 1_r,1_r])\cdot(\tf_{\chi_p}^{[e]})^{\chi_p}) \notag\\
&=\(\sum_{b\in\Herm_r(O_{F_{v}}/\varpi_{v})}\psi_{F,v}(\varpi_v^{-d_v}\tr T^\Box_{11,v} b^\sharp)\) \notag\\
&\quad\times W_{T^\Box}(m([\varpi_v^{1_u}\cdot 1_r,1_r])\cdot(\tf_{\chi_p}^{[e]})^{\chi_p}) \notag\\
&=\(\sum_{b\in\Herm_r(O_{F_{v}}/\varpi_{v})}\psi_{F,v}(\varpi_v^{-d_v}\tr T^\Box_{11,v} b^\sharp)\) \notag\\
&\quad\times \chi_v(\varpi_v^r)^{-1}q_v^{-r^2}\cdot\tf_{\chi_p}^{[e]}([\varpi_v^{1_{u^\tc}}\cdot 1_r,1_r] T^\Box [\varpi_v^{1_u}\cdot 1_r,1_r]) \notag\\
&=\(\sum_{b\in\Herm_r(O_{F_{v}}/\varpi_{v})}\psi_{F,v}(\varpi_v^{-d_v}\tr T^\Box_{11,v} b^\sharp)\) \notag\\
&\quad\times \chi_v(\varpi_v^r)^{-1}q_v^{-r^2}\cdot\tf_{\chi_p}^{[e]}
\(
\begin{pmatrix}
\varpi_v\cdot T^\Box_{11} & \varpi_v^{1_{u^\tc}}\cdot T^\Box_{12} \\
\varpi_v^{1_u}\cdot T^\Box_{21} & T^\Box_{22}
\end{pmatrix}
\). \notag
\end{align}
Since
\begin{align*}
&\quad\sum_{b\in\Herm_r(O_{F_{v}}/\varpi_{v})}\psi_{F,v}(\varpi_v^{-d_v}\tr T^\Box_{11,v} b^\sharp) \\
&=
\begin{dcases}
q_v^{r^2} & \text{if $T^\Box_{11,v}\in\Herm_{2r}(O_{F_v})$,} \\
0 & \text{if $T^\Box_{11,v}\in\varpi_v^{-1}\Herm_{2r}(O_{F_v})\setminus\Herm_{2r}(O_{F_v})$,}
\end{dcases}
\end{align*}
we have
\begin{align*}
\eqref{eq:section1}
&=\chi_v(\varpi_v^r)^{-1}\chi_p(\Nm_{E_v/F_v}\det \varpi_v^{1_{u^\tc}}T^\Box_{12})\cdot\CF_{\fT_p^{[e+e_1^\tc]}}(T^\Box) \\
&=\chi_p(\Nm_{E_p/F_p}\det T^\Box_{12})\cdot\CF_{\fT_p^{[e+e_1^\tc]}}(T^\Box) \\
&=\tf_{\chi_p}^{[e+e_1^\tc]}(T^\Box).
\end{align*}
The lemma follows.
\end{proof}

\subsection{Siegel hermitian Eisenstein series}
\label{ss:eisenstein}

Let $\chi\colon \Gamma_{F,p}\to\dC^\times$ be a finite character, regarded as an automorphic character of $\dA_F^\times$. We define $\rI^\Box_r(\chi)$ to be the restricted tensor product of $\rI^\Box_{r,v}(\chi_v)$ over all places $v$ of $F$, which is a smooth representation of $\widetilde{G}_{2r}(\dA)$. For $f_\chi\in\rI^\Box_r(\chi)$, we have the Siegel hermitian Eisenstein series\footnote{We remind the reader that the sums in the following expressions are not absolutely convergent in general; they are rather defined by analytic continuation. For example, one has a family of Eisenstein series $E(g,f_{\chi,s})$ for $s\in\dC$, where $f_{\chi,s}\in\rI^\Box_r(\chi|\;|_{\dA_F}^s)$ is the standard section induced by $f_\chi$; it is absolutely convergent when $\RE s$ is large enough and has a meromorphic continuation to $\dC$. Then $E(g,f_\chi)$ is defined as the value at $0$ of this continuation, known to be analytic.}
\begin{align*}
E(g,f_\chi)&\coloneqq\sum_{\gamma\in P_{2r}(F)\backslash G_{2r}(F)}f_\chi(\gamma g),\quad g\in G_{2r}(\dA_F),\\
\widetilde{E}(g,f_\chi)&\coloneqq\sum_{\gamma\in\widetilde{P}_{2r}(\dQ)\backslash\widetilde{G}_{2r}(\dQ)}f_\chi(\gamma g),\quad g\in \widetilde{G}_{2r}(\dA).
\end{align*}

For a finite set $\lozenge$ of places of $\dQ$ containing $\{\infty,p\}$, an element $e\in\dZ^\tP$ and a section $f\in\rI^\Box_r(\chi)^{\infty p}$, we put
\begin{align}\label{eq:eisenstein2}
\widetilde{E}^{[e]}_{\lozenge}(-,\chi,f)\coloneqq b_{2r}^\lozenge(\CF)^{-1} \cdot b_{2r}^\lozenge(\chi)\cdot
\widetilde{E}(-,f_\infty^{[r]}\otimes (\tf_{\chi_p}^{[e]})^{\chi_p}\otimes f),
\end{align}
where $b_{2r}^\lozenge(\CF)$ is defined in \S\ref{ss:run}(F4); $f_\infty^{[r]}$ is introduced in Lemma \ref{le:section6}; and $(\tf_{\chi_p}^{[e]})^{\chi_p}$ is introduced in Construction \ref{co:section}. It is clear that $\widetilde{E}^{[e]}_{\lozenge}(-,\chi,f)$ belongs to $\widetilde\cA^{[r]}_{2r,\hol}$. Put
\begin{align}\label{eq:eisenstein1}
W_{2r}^\lozenge\coloneqq W_{2r}\cdot b_{2r,\lozenge\setminus\{\infty\}}(\CF)\in\dQ^\times,
\end{align}
where $W_{2r}$ is the constant in Lemma \ref{le:section6}.

\begin{lem}\label{le:rational}
Suppose that $\|e\|>0$. Then for $f=\otimes_{w\nmid\infty p} f_w$ that is a pure tensor,
\begin{align*}
&\quad\bbq_{2r}\widetilde\bbh_{2r}\(\widetilde{E}^{[e]}_{\lozenge}(-,\chi,f)\) \\
&=W_{2r}^\lozenge\sum_{T^\Box\in\Herm_{2r}^\circ(F)^+}
\(\chi_p(\Nm_{E_p/F_p}\det T^\Box_{12})\CF_{\fT_p^{[e]}}(T^\Box)\cdot\prod_{w\nmid\infty p}W_{T^\Box}^\lozenge(f_w)\) q^{T^\Box}
\end{align*}
in which the product is finite. Here, $\widetilde\bbh_{2r}$ is the map \eqref{eq:holomorphic3}; $\bbq_{2r}$ is the map \eqref{eq:holomorphic2}; and
\[
W_{T^\Box}^\lozenge(f_w)\coloneqq
\begin{dcases}
W_{T^\Box}(f_w) & \text{if $w\in\lozenge$,} \\
b_{2r,w}(\chi)\cdot W_{T^\Box}(f_w) & \text{if $w\not\in\lozenge$.}
\end{dcases}
\]
\end{lem}

\begin{proof}
First, one notices that when $\|e\|>0$, we have $\tf_{\chi_p}^{[e]}(T^\Box)=0$ for every element $T^\Box\in\Herm_{2r}(F)\setminus\Herm_{2r}^\circ(F)$. By the discussion in \cite{Liu12}*{Section~2B} and Lemma \ref{le:section6}, the analytic $q$-expansion \eqref{eq:expansion1} of $\widetilde{E}(-,f_\infty^{[r]}\otimes(\tf_{\chi_p}^{[e]})^{\chi_p}\otimes f)$ equals
\[
W_{2r}\cdot b_{2r}^\infty(\CF)\sum_{T^\Box\in\Herm_{2r}^\circ(F)^+}
\(\chi_p(\Nm_{E_p/F_p}\det T^\Box_{12})\CF_{\fT_p^{[e]}}(T^\Box)\cdot\prod_{w\nmid\infty p}W_{T^\Box}(f_w)\) q^{T^\Box}
\]
in which we recall that $b_{2r}^\infty(\CF)$ is absolutely convergent as in \S\ref{ss:run}(F4). It follows that the analytic $q$-expansion of $\widetilde{E}^{[e]}_{\lozenge}(-,\chi,f)$ equals
\[
W_{2r}^\lozenge\sum_{T^\Box\in\Herm_{2r}^\circ(F)^+}
\(\chi_p(\Nm_{E_p/F_p}\det T^\Box_{12})\CF_{\fT_p^{[e]}}(T^\Box)\cdot\prod_{w\nmid\infty p}W_{T^\Box}^\lozenge(f_w)\) q^{T^\Box}
\]
in which the product is actually finite by \cite{Tan99}*{Proposition~3.2}. The lemma follows by the coincidence of the analytic and the algebraic $q$-expansions \cite{Lan12}*{Theorem~5.3.5}.
\end{proof}

Put
\begin{align}\label{eq:eisenstein}
\widetilde{D}^{[e]}_{\lozenge}(-,\chi,f)\coloneqq
\widetilde\rho_{r,r}\(\widetilde{E}^{[e]}_{\lozenge}(-,\chi,f)\)
\in\widetilde\cA_{r,r,\hol}^{[r]}
\end{align}
(see \eqref{eq:descent10} for the map $\widetilde\rho_{r,r}$).\footnote{The letter $D$ stands for pullback along the \emph{diagonal} block.} The following proposition concerns the rationality of $\widetilde{D}^{[e]}_{\lozenge}(-,\chi,f)$, which is the main result of this subsection.

\begin{proposition}\label{pr:rational}
Suppose that $\|e\|>0$ and let $f\in\rI^\Box_r(\chi)^{\infty p}$ be a section. For every $\sigma\in\Aut(\dC/\dQ)$, we have
\[
\widetilde\bbh_{r,r}\(\widetilde{D}^{[e]}_{\lozenge}(-,\sigma\chi,\sigma f)\)
=\sigma\widetilde\bbh_{r,r}\(\widetilde{D}^{[e]}_{\lozenge}(-,\chi,f)\),
\]
where $\widetilde\bbh_{r,r}$ is the map \eqref{eq:holomorphic3}.
\end{proposition}

Note that for $f\in\rI^\Box_r(\chi)^{\infty p}$, $\sigma f\in\rI^\Box_r(\sigma\chi)^{\infty p}$. Thus, the statement of the proposition makes sense.

\begin{proof}
Take an integer $d\geq 1$ such that $(\tf_{\chi_p}^{[e]})^{\chi_p}$ is fixed by the kernel of the map $\widetilde\cG_{2r}(\dZ_p)\to\widetilde\cG_{2r}(\dZ/p^d)$. The proof consists of two steps.

\textbf{Step 1.} We first show that for every $f\in\rI^\Box_r(\chi)^{\infty p}$ and every $\sigma\in\Aut(\dC/\dQ\langle p\rangle)$, the relation
\begin{align}\label{eq:rational0}
\widetilde\bbh_{2r}\(\widetilde{E}^{[e]}_{\lozenge}(-,\sigma\chi,\sigma f)\)
=\sigma\widetilde\bbh_{2r}\(\widetilde{E}^{[e]}_{\lozenge}(-,\chi,f)\)
\end{align}
holds.

Take an irreducible summand $\rI$ of $\rI^\Box_r(\chi)^{\infty p}$ (as a representation of $\widetilde{G}_{2r}(\dA^{\infty p})$). Choose a positive integer $\Delta=\Delta_\rI>1$ that is coprime to $p$ such that
\begin{itemize}
  \item[(1)] for every rational prime $w$ not dividing $p\Delta$, $\rI_w$ has nonzero invariants under $\widetilde\cG_{2r}(\dZ_w)$;

  \item[(2)] one can write $\Delta=\Delta_1\cdot\Delta_2$ with $(\Delta_1,\Delta_2)=1$ such that for $i=1,2$, $\prod_{w\mid\Delta_i}\widetilde{P}_{2r}(\dQ_w)$ maps surjectively to $\widetilde{G}_{2r}^\ab(\dQ)\backslash \widetilde{G}_{2r}^\ab(\dA^\infty)/\widetilde{K}_{2r}^\ab(p^d)$ (Remark \ref{re:adjoint}).
\end{itemize}
For every $\sigma\in\Aut(\dC/\dQ\langle p\rangle)$, since the map
\[
f\mapsto\widetilde\bbh_{2r}\(\widetilde{E}^{[e]}_{\lozenge}(-,\sigma\chi,\sigma f)\)-
\sigma\widetilde\bbh_{2r}\(\widetilde{E}^{[e]}_{\lozenge}(-,\chi,f)\)
\]
is $\widetilde{G}_{2r}(\dA^{\infty p})$-equivariant, it suffices to show that there exists a nonzero element $f=f_\sigma\in\rI$ such that
\begin{align}\label{eq:rational1}
\widetilde\bbh_{2r}\(\widetilde{E}^{[e]}_{\lozenge}(-,\sigma\chi,\sigma f)\)
-\sigma\widetilde\bbh_{2r}\(\widetilde{E}^{[e]}_{\lozenge}(-,\chi,f)\)=0.
\end{align}
In practice below, we will first check \eqref{eq:rational1} for $\sigma\in\Aut(\dC/\dQ\langle\Delta_1\rangle)$ and then for $\sigma\in\Aut(\dC/\dQ\langle\Delta_2\rangle)$.

Choose a nonzero element $f=\otimes_{w\nmid\infty p} f_w\in\rI$ such that $f_w$ satisfies the condition in Lemma \ref{le:section2} (that is, it belongs to $\rI^\Box_{r,w}(\chi_w)^\circ$ and is fixed by $\widetilde\cP_{2r}(\dZ_w)$) for $w\mid\Delta$ and that $f_w$ is the unique section that is fixed by $\widetilde\cG_{2r}(\dZ_w)$ and satisfies $f_w(1_{4r})=1$ for $w\nmid \Delta$. Replacing $\Delta$ by a power of $\Delta$, we may assume that $f_w$ is invariant under $\widetilde\cG_{2r}(\dZ_w)\times_{\widetilde\cG_{2r}(\dZ_w/\Delta)}\widetilde\cP_{2r}(\dZ_w/\Delta)$ for every $w\nmid\infty p$. In particular, we have
\[
\widetilde{E}^{[e]}_{\lozenge}(-,\chi,f)\in\widetilde\cA_{2r,\hol}^{[r]}(\widetilde{K}_{2r}(\Delta,p^d))
\]
(Notation \ref{no:open}). For such $f$, we first show that \eqref{eq:rational1} holds for $\sigma\in\Aut(\dC/\dQ\langle p\Delta_1\rangle)$. By property (2) for $\Delta$ and Remark \ref{re:adjoint}, the $q$-expansions of $\bbh_{2r}\(\widetilde{E}^{[e]}_{\lozenge}(-,\chi,g\cdot f)\)$ for all $g\in\prod_{w\mid\Delta_1}\widetilde{P}_{2r}(\dQ_w)$ determines $\widetilde{E}^{[e]}_{\lozenge}(-,\chi,f)$. For every $g\in\prod_{w\mid\Delta_1}\widetilde{P}_{2r}(\dQ_w)$, there exists an integer $d_g\geq 1$ such that $\widetilde{E}^{[e]}_{\lozenge}(-,\chi,g\cdot f)$ belongs to $\widetilde\cA_{2r,\hol}^{[r]}(\widetilde{K}_{2r}(\Delta_2,p^d\Delta_1^{d_g}))$. Then by Lemma \ref{le:rational1}, to check \eqref{eq:rational1} for $\sigma\in\Aut(\dC/\dQ\langle p\Delta_1\rangle)$, it suffices to check that
\[
\bbq_{2r}\widetilde\bbh_{2r}\(\widetilde{E}^{[e]}_{\lozenge}(-,\sigma\chi,g\cdot\sigma f)\)
-\sigma\bbq_{2r}\widetilde\bbh_{2r}\(\widetilde{E}^{[e]}_{\lozenge}(-,\chi,g\cdot f)\)=0
\]
for every $g\in\prod_{w\mid\Delta_1}\widetilde{P}_{2r}(\dQ_w)$. However, this follows from Lemma \ref{le:rational} and Lemma \ref{le:section4}(2,3). Since the roles of $\Delta_1$ and $\Delta_2$ are symmetric, \eqref{eq:rational1} also holds for $\sigma\in\Aut(\dC/\dQ\langle p\Delta_2\rangle)$. Together, \eqref{eq:rational1} holds for $\sigma\in\Aut(\dC/\dQ\langle p\rangle)$. Thus, \eqref{eq:rational0} holds.

\textbf{Step 2.} By Step 1 and the upper square of the functorial diagram \eqref{eq:holomorphic4}, for the proposition, it suffices to show that for every $f\in\rI^\Box_r(\chi)^{\infty p}$, there exists a positive integer $\Delta$ that is coprime to $p$ such that
\begin{align}\label{eq:rational2}
\widetilde\bbh_{r,r}\(\widetilde{D}^{[e]}_{\lozenge}(-,\sigma\chi,\sigma f)\)-
\sigma\widetilde\bbh_{r,r}\(\widetilde{D}^{[e]}_{\lozenge}(-,\chi,f)\)=0
\end{align}
holds for every $\sigma\in\Aut(\dC/\dQ\langle\Delta\rangle)$.

By linearity, we may assume that $f=\otimes_{w\nmid\infty p} f_w$ is a pure tensor. Let $\Delta$ be a positive integer that is coprime to $p$ such that
\begin{itemize}
  \item[(3)] for every prime $w$ not dividing $p\Delta$, $f_w$ is the unique section that is fixed by $\widetilde\cG_{2r}(\dZ_w)$ and satisfies $f_w(1_{4r})=1$;

  \item[(4)] for every prime $w$ dividing $\Delta$, $f_w$ is fixed by the kernel of the map $\widetilde\cG_{2r}(\dZ_w)\to\widetilde\cG_{2r}(\dZ_w/\Delta)$.
\end{itemize}
Combining with Lemma \ref{le:section3}, we see that $\widetilde{D}^{[e]}_{\lozenge}(-,\chi,f)$ belongs to $\widetilde\cA_{r,r,\hol}^{[r]}(\widetilde{K}_{r,r}(p^d,\Delta))$. Thus, by Lemma \ref{le:rational1} (with $\Delta=p$ and $\Delta'=\Delta$), for \eqref{eq:rational2}, it suffices to show that
\begin{align}\label{eq:rational3}
\bbq_{r,r}\widetilde\bbh_{r,r}\(\widetilde{D}^{[e]}_{\lozenge}(-,\sigma\chi,\sigma f)\)-
\sigma\bbq_{r,r}\widetilde\bbh_{r,r}\(\widetilde{D}^{[e]}_{\lozenge}(-,\chi,f)\)=0
\end{align}
holds for every $\sigma\in\Aut(\dC/\dQ\langle\Delta\rangle)$ (this time, we only need to consider the $q$-expansion on one connected component since we argue for all $f$). By Lemma \ref{le:rational} and Lemma \ref{le:section4}(2), \eqref{eq:rational3} holds for $\sigma\in\Aut(\dC/\dQ\langle\Delta\rangle)$.

The proposition is proved.
\end{proof}

\subsection{Relevant representations}
\label{ss:relevant}

\begin{lem}\label{le:dual}
Let $\dL/\dQ_p$ be a finite extension and let $\pi$ be a relevant $\dL$-representation of $G_r(\dA_F^\infty)$ (Definition \ref{de:relevant}).
\begin{enumerate}
  \item The representation $\hat\pi\coloneqq(\pi^\vee)^\dag$ is also a relevant $\dL$-representation of $G_r(\dA_F^\infty)$, where $\dag$ is the involution introduced at the beginning of \S\ref{ss:cyclotomic}.

  \item The $\dL$-vector space $\Hom_{G_r(\dA_F^\infty)}(\pi,\cH_r^{[r]}\otimes_{\dQ_p}\dL)$ has dimension $1$.
\end{enumerate}
See Definition \ref{de:holomorphic} for the notation $\cH_r^{[r]}$.
\end{lem}

\begin{proof}
Part (1) follows from the fact that for every $v\in\tV_F^{(\infty)}$, $((\pi^{[r]}_v)^\vee)^\dag$ is isomorphic to $\pi^{[r]}_v$.

For (2), it suffices to show that for every embedding $\iota\colon\dL\to\dC$, the complex vector space $\Hom_{G_r(\dA_F^\infty)}(\iota\pi,\cA^{[r]}_{r,\hol})$ has dimension $1$. However, this follows from Arthur's multiplicity one property proved in \cite{Mok15}.
\end{proof}

Now we fix a relevant $\dL$-representation $\pi$ of $G_r(\dA_F^\infty)$ for some finite extension $\dL/\dQ_p$ contained in $\ol\dQ_p$ such that $\pi_v$ is \emph{unramified} for every $v\in\tV_F^{(p)}$. After Lemma \ref{le:dual}, we let $\cV_\pi$ and $\cV_{\hat\pi}$ be the unique subspaces of $\cH_r^{[r]}\otimes_{\dQ_p}\dL$ that are isomorphic to $\pi$ and $\hat\pi$, respectively.

\begin{notation}\label{no:period}
We fix a $G_r(\dA_F^\infty)$-invariant pairing $\langle\;,\;\rangle_\pi\colon\cV_{\hat\pi}^\dag\times\cV_\pi\to\dL$. For every embedding $\iota\colon\dL\to\dC$, since $\pi$ is absolutely irreducible, there is a unique element $\rP^\iota_\pi\in\dC^\times$ such that
\[
\int_{G_r(F)\backslash G_r(\dA_F)}\varphi_1^\iota(g^\dag)\varphi_2^\iota(g)\rd g
=\rP^\iota_\pi\cdot\langle\varphi_1^\dag,\varphi_2\rangle_\pi
\]
for every $\varphi_1\in\cV_{\hat\pi}$ and $\varphi_2\in\cV_\pi$. See Definition \ref{de:holomorphic} for the notation $\varphi_i^\iota$.
\end{notation}

\begin{remark}\label{re:period}
Since $\hat{\hat\pi}=\pi$, the pairing $\langle\;,\;\rangle_\pi$ is equivalent to a similar pairing $\langle\;,\;\rangle_{\hat\pi}\colon\cV_\pi^\dag\times\cV_{\hat\pi}\to\dL$ for $\hat\pi$, for which we have $\rP^\iota_{\hat\pi}=\rP^\iota_\pi$ for every embedding $\iota\colon\dL\to\dC$.
\end{remark}

\begin{lem}\label{le:projection}
There is a unique $\dL$-linear map
\begin{align*}
\pr_\pi\colon\cH_r^{[r]}\otimes_{\dQ_p}\dL\to\cV_\pi
\end{align*}
satisfying that for every $\cZ\in\cH_r^{[r]}\otimes_{\dQ_p}\dL$, every $\varphi\in\cV_{\hat\pi}$ and every embedding $\iota\colon\dL\to\dC$,
\begin{align*}
\int_{G_r(F)\backslash G_r(\dA_F)}\varphi^\iota(g^\dag)\cZ^\iota(g)\rd g
=\rP^\iota_\pi\cdot\iota\langle \varphi^\dag, \pr_\pi(\cZ)\rangle_\pi
\end{align*}
holds.
\end{lem}

\begin{proof}
Take an open compact subgroup $K$ of $G_r(\dA_F^\infty)$. The $\dL$-vector space $\cH_r^{[r]}(K)\otimes_{\dQ_p}\dL$ is a semisimple module over $\dL[K\backslash G_r(\dA_F^\infty)/K]$, in which $\cV_\pi(K)$ is the unique summand that is isomorphic to $\pi^{K}$. We denote by $\cH_r^{[r]}(K)^\pi\subseteq\cH_r^{[r]}(K)\otimes_{\dQ_p}\dL$ the direct sum of the remaining summands. Then we have a direct sum decomposition
\[
\cH_r^{[r]}(K)\otimes_{\dQ_p}\dL=\cV_\pi(K)\oplus\cH_r^{[r]}(K)^\pi
\]
of $\dL[K\backslash G_r(\dA_F^\infty)/K]$-modules. Denote by $\pr_\pi^K\colon\cH_r^{[r]}(K)\otimes_{\dQ_p}\dL\to\cV_\pi(K)\subseteq\cV_\pi$ the corresponding projection map. It is clear that the maps $\pr_\pi^K$ are compatible with each other for different $K$, hence defining a map $\pr_\pi\colon\cH_r^{[r]}\otimes_{\dQ_p}\dL\to\cV_\pi$ which satisfies the property of the lemma by construction. The lemma is proved as the uniqueness is obvious.
\end{proof}

In the rest of this subsection, we take an element $v\in\tV_F^{(p)}$. For every $u\in\tP_v$, we have the representation $\pi_u$ of $\GL_n(F_v)$ as a local component of $\pi$ via the isomorphism $G_r(F_v)\simeq\GL_n(E_u)=\GL_n(F_v)$ (recall that $n=2r$). In particular, $\pi_u^\vee\simeq\pi_{u^\tc}\simeq(\pi^\vee)_u$. Note that we will also speak of $\pi_v$, a representation of $G_r(F_v)$ without any identification with $\GL_n(F_v)$.

\begin{definition}\label{de:satake1}
We let $\{\alpha_{u,1},\dots,\alpha_{u,n}\}\subseteq\ol\dQ_p^\times$ (as a multi-subset) be the Satake parameter of $\pi_u$.
\begin{enumerate}
  \item Define the \emph{Satake polynomial} of $\pi_u$ to be
      \[
      \sfP_{\pi_u}(T)\coloneqq\prod_{j=1}^{n}\(T-\alpha_{u,j}\sqrt{q_v}^{n-1}\).
      \]

  \item For every integer $1\leq m\leq n$, put
      \[
      \sfA(\pi_u,m)\coloneqq\left\{\left.\(\prod_{j\in J}\alpha_{u,j}\)\sqrt{q_v}^{m(n-m)}\right| J\subseteq\{1,\dots,n\},|J|=m\right\}
      \]
      as a multi-subset of $\ol\dQ_p$.
\end{enumerate}
Note that to define the Satake parameter, one needs to choose a square root of $q_v$ in $\ol\dQ_p$. However, both $\sfP_{\pi_u}(T)$ and $\sfA(\pi_u,m)$ are independent of such choice.
\end{definition}

\begin{lem}\label{le:satake}
We have
\begin{enumerate}
  \item There exist $\beta_{u,1},\dots,\beta_{u,n}\in O_\dL$ such that
     \[
     \sfP_{\pi_u}(T)=T^n+\sum_{s=1}^n\beta_{u,s}\cdot q_v^{\frac{s(s-1)}{2}} T^{n-s}.
     \]

  \item For every integer $1\leq m\leq n$, $\sfA(\pi_u,m)$ is contained in $\ol\dZ_p$ and contains at most one element (with multiplicity one) in $\ol\dZ_p^\times$. Moreover, $\sfA(\pi_u,m)\cap\ol\dZ_p^\times\neq\emptyset$ if and only if $\beta_{u,m}\in O_\dL^\times$.

  \item We have that $\sfA(\pi_u,m)\cap\ol\dZ_p^\times\neq\emptyset$ is equivalent to that $\sfA(\pi_{u^\tc},n-m)\cap\ol\dZ_p^\times\neq\emptyset$.
\end{enumerate}
\end{lem}

\begin{proof}
Part (1) follows from Definition \ref{de:relevant} and \cite{Hid98}*{Theorem~8.1(3)}, that is, the Newton polygon is above the Hodge polygon.

For (2), we may order the multi-set $\{\alpha_{u,1},\dots,\alpha_{u,n}\}$ in the way that $\alpha_{u,j+1}/\alpha_{u,j}\in\ol\dZ_p$ for $1\leq j<n$. Then it follows by (1) and induction that for every $1\leq m\leq n$, $\prod_{j=1}^m(\alpha_{u,j}\sqrt{q_v}^{n-1})\in q_v^{\frac{m(m-1)}{2}}\ol\dZ_p$. Thus, (2) follows.

Part (3) follows from the fact that $\prod_{j=1}^n\alpha_{u,j}$ is a root of unity and the fact that $\{\alpha_{u^\tc,1},\dots,\alpha_{u^\tc,n}\}=\{\alpha_{u,1}^{-1},\dots,\alpha_{u,n}^{-1}\}$.
\end{proof}

Put
\begin{align}\label{eq:iwahori}
I_v\coloneqq\cG_r(O_{F_v})\times_{\cG_r(O_{F_v}/\varpi_v)}\cP_r(O_{F_v}/\varpi_v)
\end{align}
which is an open compact subgroup of $G_r(F_v)$. For every $u\in\tP_v$, define two Hecke operators
\[
\rT_u^\pm\coloneqq I_v
\begin{pmatrix}
\varpi_v^{\pm 1_u}\cdot 1_r &  \\
 & \varpi_v^{\mp 1_{u^\tc}} \cdot 1_r \\
\end{pmatrix}
I_v
\]
(in which the volume of $I_v$ is regarded as $1$). In particular, $\rT_u^+=\rU_u\cdot I_v$ (Notation \ref{no:section2}).

\begin{lem}\label{le:satake1}
For every $u\in\tP_v$, the multisets of generalized eigenvalues of the actions of $\rT_u^+$ and $\rT_u^-$ on $\pi_v^{I_v}$ are $\sfA(\pi_u,r)$ and $\sfA(\pi_{u^\tc},r)$, respectively.
\end{lem}

The proof of this lemma will be given at the end of this subsection.

\begin{definition}\label{de:satake}
We say that the (unramified) representation $\pi_v$ of $G_r(F_v)$ is \emph{Panchishkin} if $\beta_{u,r}\in O_\dL^\times$ for every $u\in\tP_v$ under the notation in Lemma \ref{le:satake}.
\end{definition}

\begin{lem}\label{le:satake2}
The following statements are equivalent:
\begin{enumerate}
  \item $\pi_v$ is Panchishkin unramified;

  \item $\hat\pi_v$ is Panchishkin unramified;

  \item $\sfA(\pi_u,r)$ contains a unique element in $O_\dL^\times$ for some $u\in\tP_v$.
\end{enumerate}
\end{lem}

\begin{proof}
This is an immediate consequence of Lemma \ref{le:satake}. The fact that the unique element in $\sfA(\pi_u,r)\cap\ol\dZ_p^\times$ belongs to $\dL$ follows from the Galois action and the uniqueness.
\end{proof}

\begin{lem}\label{le:satake3}
Suppose that $\pi_v$ is Panchishkin unramified.
\begin{enumerate}
  \item The one-dimensional subspace of $\pi_v^{I_v}$ that is the eigenspace of the operator $\rT_u^+$ for the eigenvalue that is the unique element in $\sfA(\pi_u,r)\cap O_\dL^\times$ is independent of $u\in\tP_v$.

  \item The one-dimensional subspace of $\pi_v^{I_v}$ that is the eigenspace of the operator $\rT_u^-$ for the eigenvalue that is the unique element in $\sfA(\pi_{u^\tc},r)\cap O_\dL^\times$ is independent of $u\in\tP_v$.
\end{enumerate}
\end{lem}

The proof of this lemma will be given at the end of this subsection.

\begin{notation}\label{no:satake}
Suppose that $\pi_v$ is Panchishkin unramified.
\begin{enumerate}
  \item For every $u\in\tP_v$, we denote by $\alpha(\pi_u)\in O_\dL^\times$ the unique element in Lemma \ref{le:satake3}(3).

  \item In view of Lemma \ref{le:satake1} and Lemma \ref{le:satake3}, we denote by $\pi_v^+$ and $\pi_v^-$ the one-dimensional subspaces of $\pi_v^{I_v}$ that are the eigenspaces of the operators $\rT_u^+$ and $\rT_u^-$ for the eigenvalues $\alpha(\pi_u)$ and $\alpha(\pi_{u^\tc})$ for every $u\in\tP_v$, respectively.
\end{enumerate}
\end{notation}

\begin{proposition}\label{pr:satake}
Suppose that $\pi_v$ is Panchishkin unramified.
\begin{enumerate}
  \item For every $u\in\tP_v$, there is a unique polynomial $\sfQ_{\pi_u}(T)\in\dL[T]$ that divides $\sfP_{\pi_u}(T)$ and has the form
      \[
      \sfQ_{\pi_u}(T)=T^r+\gamma_{u,1}\cdot T^{r-1}+ \gamma_{u,2}\cdot q_v T^{r-2}+\cdots + \gamma_{u,r}\cdot q_v^{\frac{r(r-1)}{2}}
      \]
      with $\gamma_{u,r}\in O_\dL^\times$.

  \item There is a unique $\dL$-valued locally constant function $\xi_{\pi_v}$ on $G_r(F_v)$ such that
     \begin{enumerate}
       \item there exist $\varphi_v^\vee\in(\pi_v^\vee)^-$ and $\varphi_v\in\pi_v^-$ such that $\xi_{\pi_v}=\langle\pi_v^\vee(g)\varphi_v^\vee,\varphi_v\rangle_{\pi_v}$ for $g\in G_r(F_v)$;

       \item $\xi_{\pi_v}(\tw_r)=1$.
     \end{enumerate}
     In particular, $\xi_{\pi_v}$ is bi-$I_v$-invariant.

  \item For $u\in\tP_v$, denote by $\ul{\pi_u}$ the unramified principal series of $\GL_r(F_v)$ with $\sfQ_{\pi_u}(T)$ as its Satake polynomial, which is defined over $\dL$. Then there exist $\GL_r(O_{F_v})$-invariant vectors $\phi_u\in\ul{\pi_u}$ and $\phi_u^\vee\in(\ul{\pi_u})^\vee$ for every $u\in\tP_v$ such that
      \[
      \xi_{\pi_v}(m(a)\tw_r)=\prod_{u\in\tP_v}
      \langle\ul{\pi_u}(a_{u^\tc})\phi_u,\phi_u^\vee\rangle_{(\ul{\pi_u})^\vee}
      \]
      holds for every $a=(a_u)_u\in\GL_r(E_v)=\prod_{u\in\tP_v}\GL_r(E_u)$.
\end{enumerate}
\end{proposition}

In this rest of this subsection, we prove Lemma \ref{le:satake1}, Lemma \ref{le:satake3} and Proposition \ref{pr:satake}. To ease the notation, we will suppress the subscript $v$ hence $F=F_v$, $\tP=\tP_v$, and $\pi=\pi_v$ temporarily. It is easy to see that for these three statements, we may replace $\dL$ by a finite extension (inside $\ol\dQ_p$). Thus, without loss of generality, we may assume that $\dL$ contains both $\sqrt{q}$ and $\alpha_{u,j}$ for $u\in\tP$ and $1\leq j\leq n$. We need some preparation on Jacquet modules.

For every subset $J\subseteq\{1,\dots,n\}$, put $\ol{J}\coloneqq\{1,\dots,n\}\setminus J$. For every subset $J\subseteq\{1,\dots,n\}$ of cardinality $r$, every $u\in\tP$ and every sign $\epsilon\in\{+,-\}$, we denote by $\rI(\alpha_{u,j}\sqrt{q}^{\epsilon r}\res j\in J)$ the unramified principal series of $\GL_r(F)$ with the Satake parameter $\{\alpha_{u,j}\sqrt{q}^{\epsilon r}\res j\in J\}$, with coefficients in $\dL$.

Put $\ol{P}_r\coloneqq\tw_r^{-1}P_r\tw_r$ and let $\ol{N}_r$ be its unipotent radical. We identify both Levi quotients $P_r/N_r$ and $\ol{P}_r/\ol{N}_r$ with $\Res_{E/F}\GL_r$ via the map $m$ in \S\ref{ss:run}(G4). We define the Jacquet modules
\begin{align*}
\pi_{N_r}&\coloneqq\pi/\{\varphi-\pi(n)\varphi\res n\in N_r(F),\varphi\in\pi\}, \\
\pi_{\ol{N}_r}&\coloneqq\pi/\{\varphi-\pi(n)\varphi\res n\in\ol{N}_r(F),\varphi\in\pi\},
\end{align*}
which are admissible representations of $\GL_r(E)$ of finite length. Fix an order $\{u_1,u_2\}$ of $\tP$. Recall that $\{\alpha_{u_1,1},\dots,\alpha_{u_1,n}\}=\{\alpha_{u_2,1}^{-1},\dots,\alpha_{u_2,n}^{-1}\}$. Without loss of generality, we may assume $\alpha_{u_1,j}\alpha_{u_2,j}=1$ for $1\leq j\leq n$. It is well-known that
\begin{align*}
\pi_{N_r}^{\r{ss}}&\simeq\bigoplus_{\substack{J\subseteq\{1,\dots,n\} \\ |J|=r}}
\rI(\alpha_{u_1,j}\sqrt{q}^{-r}\res j\in J)\boxtimes\rI(\alpha_{u_2,j}\sqrt{q}^{-r}\res j\in\ol{J}), \\
\pi_{\ol{N}_r}^{\r{ss}}&\simeq\bigoplus_{\substack{J\subseteq\{1,\dots,n\} \\ |J|=r}}
\rI(\alpha_{u_1,j}\sqrt{q}^r\res j\in\ol{J})\boxtimes\rI(\alpha_{u_2,j}\sqrt{q}^r\res j\in J),
\end{align*}
as representations of $\GL_r(E)=\GL_r(E_{u_1})\times\GL_r(E_{u_2})$. Since $\tw_r$ conjugates $m(a_1,a_2)\in G_{2r}(F)$ to $m(\pres{\rt}{a}_2^{-1},\pres{\rt}{a}_1^{-1})$, the isomorphism $\tw_r\colon\pi\xrightarrow\sim\pi$ descends to an isomorphism $\pi_{N_r}^{\r{ss}}\to\pi_{\ol{N}_r}^{\r{ss}}$ that sends
\[
\rI(\alpha_{u_1,j}\sqrt{q}^{-r}\res j\in J)\boxtimes\rI(\alpha_{u_2,j}\sqrt{q}^{-r}\res j\in\ol{J})
\]
to
\[
\rI(\alpha_{u_1,j}\sqrt{q}^r\res j\in\ol{J})\boxtimes\rI(\alpha_{u_2,j}\sqrt{q}^r\res j\in J).
\]

\begin{proof}[Proofs of Lemma \ref{le:satake1} and Lemma \ref{le:satake3}]
The element
\[
\begin{pmatrix}
 & \varpi^{-1_{u^c}}\cdot 1_r \\
-\varpi^{1_u}\cdot 1_r & \\
\end{pmatrix}\in G_r(F)
\]
normalizes $I$ and induces an operator on $\pi^I$ that switches $\rT_u^+$ and $\varpi^{1_u-1_{u^\tc}}\cdot\rT_u^-$. In particular, if the multiset of generalized eigenvalues of $\rT_u^+$ on $\pi^I$ is $\sfA(\pi_u,r)$, then the multiset for $\rT_u^-$ is
\[
\left\{\left.\(\prod_{j\in J}\alpha_{u,j}^{-1}\)\sqrt{q}^{r^2}\right| J\subseteq\{1,\dots,n\},|J|=r\right\},
\]
which is nothing but $\sfA(\pi_{u^\tc},r)$. Thus, it suffices to study $\rT_u^+$ in both lemmas.

The quotient map $\pi\to\pi_{N_r}$ induces an isomorphism
\[
\pi^I\xrightarrow{\sim}\pi_{N_r}^{\GL_r(O_E)}
\]
under which the operator $\rT_u^+$ (which is nothing but the operator $\rU_u$ in Notation \ref{no:section2}) corresponds to the operator $q^{r^2}\cdot\(\begin{smallmatrix}\varpi^{1_u}\cdot 1_r & \\ & 1_r\end{smallmatrix}\)$.

Since the operator $\(\begin{smallmatrix}\varpi^{1_u}\cdot 1_r & \\ & 1_r\end{smallmatrix}\)$ acts on
\[
\rI(\alpha_{u,j}\sqrt{q}^{-r}\res j\in J)\boxtimes\rI(\alpha_{u^\tc,j}\sqrt{q}^{-r}\res j\in\ol{J})
\]
by the scalar $\prod_{j\in J}\alpha_{u,j}\sqrt{q}^{-r^2}$, the multiset of (generalized) eigenvalues of $\rT_u^+$ on $\pi^I$ is $\sfA(\pi_u,r)$. Lemma \ref{le:satake1} is proved.

Now we consider Lemma \ref{le:satake3}, for which it suffices to show (1). For $i=1,2$, let $J_i$ be the unique subset of $\{1,\dots,n\}$ of cardinality $r$ such that $\prod_{j\in J_i}\alpha_{u_i,j}\sqrt{q}^{-r^2}\in\ol\dZ_p^\times$. Then $J_1\cup J_2=\{1,\dots,n\}$. Thus, for both $i=1,2$, the one-dimensional subspace of $\pi^I$ that is the eigenspace of the operator $\rT_{u_i}^+$ for the eigenvalue that is the unique element in $\sfA(\pi_{u_i},r)\cap O_\dL^\times$ is the $\GL_r(O_E)$-invariant subspace of
\[
\rI(\alpha_{u_1,j}\sqrt{q}^{-r}\res j\in J_1)\boxtimes\rI(\alpha_{u_2,j}\sqrt{q}^{-r}\res j\in J_2).
\]
Lemma \ref{le:satake3} is proved.
\end{proof}

\begin{proof}[Proof of Proposition \ref{pr:satake}]
Without loss of generality, by Lemma \ref{le:satake2}, we may assume that the unique subset $J$ of $\{1,\dots,n\}$ with $|J|=r$ such that $\sqrt{q}^{r^2}\prod_{j\in J}\alpha_{u_1,j}\in O_\dL^\times$ is $\{1,\dots,r\}$. It follows that the unique subset $J$ of $\{1,\dots,n\}$ with $|J|=r$ such that $\sqrt{q}^{r^2}\prod_{j\in J}\alpha_{u_2,j}\in O_\dL^\times$ is $\{r+1,\dots,n\}$.

For (1), note that every factor of $\sfP_{\pi_u}(T)$ in $\dL[T]$ that is monic of degree $r$ has the form
\[
\prod_{j\in J}\(T-\alpha_{u,j}\sqrt{q}^{n-1}\)
\]
for some $J\subseteq\{1,\dots,n\}$ with $|J|=r$. In particular, the corresponding term $\gamma_{u,r}$ equals $\sqrt{q}^{r^2}\prod_{j\in J}\alpha_{u,j}$. Thus, we must have
\begin{align*}
\sfQ_{\pi_{u_1}}(T)&=\prod_{j=1}^r\(T-\alpha_{u_1,j}\sqrt{q}^{n-1}\),\\
\sfQ_{\pi_{u_2}}(T)&=\prod_{j=r+1}^n\(T-\alpha_{u_2,j}\sqrt{q}^{n-1}\).
\end{align*}

For (2) and (3), it suffices to show the following claim: For nonzero vectors $\varphi^\vee\in(\pi^\vee)^-$ and $\varphi\in\pi^-$, there exist nonzero $\GL_r(O_F)$-invariant vectors $\phi_1\in\ul{\pi_{u_1}}$, $\phi_1^\vee\in(\ul{\pi_{u_1}})^\vee$, $\phi_2\in\ul{\pi_{u_2}}$, $\phi_2^\vee\in(\ul{\pi_{u_2}})^\vee$ such that
\[
\langle\pi^\vee(m(a_1,a_2)\tw_r)\varphi^\vee,\varphi\rangle_{\pi}=\prod_{i=1}^2
\langle\ul{\pi_{u_i}}(a_{3-i})\phi_i,\phi_i^\vee\rangle_{(\ul{\pi_{u_i}})^\vee}
\]
holds for every $(a_1,a_2)\in\GL_r(E)=\GL_r(E_{u_1})\times\GL_r(E_{u_2})$.

Again by Lemma \ref{le:satake2}, the two factors
\begin{align*}
&\rI(\alpha_{u_1,j}\sqrt{q}^{-r}\res 1\leq j\leq r)\boxtimes\rI(\alpha_{u_2,j}\sqrt{q}^{-r}\res r+1\leq j\leq n), \\
&\rI(\alpha_{u_1,j}\sqrt{q}^{-r}\res r+1\leq r\leq n)\boxtimes\rI(\alpha_{u_2,j}\sqrt{q}^{-r}\res 1\leq j\leq r)
\end{align*}
are direct summands of $\pi_{N_r}$. We see from the proof of Lemma \ref{le:satake3} that under the projection $\pi\to\pi_{N_r}$, the one-dimensional subspaces $\pi^+,\pi^-\subseteq\pi^I$ map to
\begin{align*}
&\rI(\alpha_{u_1,j}\sqrt{q}^{-r}\res 1\leq j\leq r)^{\GL_r(O_F)}\boxtimes\rI(\alpha_{u_2,j}\sqrt{q}^{-r}\res r+1\leq j\leq n)^{\GL_r(O_F)},\\
&\rI(\alpha_{u_1,j}\sqrt{q}^{-r}\res r+1\leq r\leq n)^{\GL_r(O_F)}\boxtimes\rI(\alpha_{u_2,j}\sqrt{q}^{-r}\res 1\leq j\leq r)^{\GL_r(O_F)},
\end{align*}
respectively. However, we observe that
\begin{align*}
\rI(\alpha_{u_1,j}\sqrt{q}^{-r}\res r+1\leq r\leq n)&\simeq(\ul{\pi_{u_2}})^\vee,\\
\rI(\alpha_{u_2,j}\sqrt{q}^{-r}\res 1\leq j\leq r)&\simeq(\ul{\pi_{u_1}})^\vee.
\end{align*}
The claim follows.

The proposition is proved.
\end{proof}

\subsection{Local doubling zeta integral}
\label{ss:zeta}

Let $\pi$ be as in \S\ref{ss:relevant}. Let $\chi\colon \Gamma_{F,p}\to\dC^\times$ be a finite character, regarded as an automorphic character of $\dA_F^\times$.

Take a finite place $v\in\tV_F^\fin$. For every $\varphi_v^\vee\in\pi_v^\vee$, $\varphi_v\in\pi_v$ and $f\in\rI^\Box_{r,v}(\chi_v)$, we have the local doubling zeta integral
\begin{align*}
Z^\iota(\varphi_v^\vee\otimes\varphi_v,f)\coloneqq
\int_{G_r(F_v)}\iota\langle\pi_v^\vee(g)\varphi_v^\vee,\varphi_v\rangle_{\pi_v}\cdot f(\bw_r(g,1_{2r}))\rd g
\end{align*}
for every embedding $\iota\colon\dL\to\dC$. Here, $\bw_r$ is defined in \eqref{eq:wr} and $(g,1_{2r})$ is an element of $G_{2r}(F_v)$ via the embedding \eqref{eq:embedding}. Since $\iota\pi_v$ is tempered, the above integral is absolutely convergent by \cite{Yam14}*{Lemma~7.2}.

\begin{lem}\label{le:zeta}
Define a map $\varsigma\colon(\Res_{E/F}\GL_r)\times\Herm_F\times\Herm_F\to G_r$ by the formula
\[
\varsigma(a,u_1,u_2)\coloneqq
\begin{pmatrix}
1_r & u_2 \\
0 & 1_r \\
\end{pmatrix}
\begin{pmatrix}
-a & 0 \\
0 & -\pres{\rt}{a}^{\tc,-1} \\
\end{pmatrix}
\tw_r
\begin{pmatrix}
1_r & u_1 \\
0 & 1_r \\
\end{pmatrix},
\]
whose image is contained in the big Bruhat cell $P_r\tw_r N_r$. Then $Z^\iota(\varphi_v^\vee\otimes\varphi_v,f)$ equals
\begin{align*}
\int_{P_r(F_v)\tw_r N_r(F_v)}
&\iota\langle\pi_v^\vee(\varsigma(a,u_1,u_2))\varphi_v^\vee,\varphi_v\rangle_{\pi_v}
\cdot\chi(\Nm_{E_v/F_v}\det a)|\dtm a|_{E_v}^r \\
&\times f\(\tw_r^\Box\cdot
n\begin{pmatrix}
u & \pres{\rt}{a}^\tc \\
a & v
\end{pmatrix}\)\rd\varsigma(a,u_1,u_2),
\end{align*}
where the integral is absolutely convergent. Here, we recall that $\tw_r^\Box=\tw_{2r}=\(\begin{smallmatrix}&1_{2r}\\ -1_{2r} &\end{smallmatrix}\)$ from \S\ref{ss:notation}.
\end{lem}

\begin{proof}
This formula is deduced in the proof of \cite{LL}*{Proposition~3.13}.
\end{proof}

\begin{definition}\label{de:typical}
For a pair $\varphi_v^\vee\in\pi_v^\vee$ and $\varphi_v\in\pi_v$, we say that an element $\tf\in\sS(\Herm_{2r}(F_v))$ is \emph{$(\varphi_v^\vee,\varphi_v)$-typical} if its Fourier transform $\widehat{\tf}\in\sS(\Herm_{2r}(F_v))$ with respect to $\psi_{F,v}$ (recall from \S\ref{ss:notation}) satisfies
\begin{enumerate}
  \item $\widehat{\tf}$ takes values in $\dQ$;

  \item $\widehat{\tf}$ is supported on the subset
       \[
       \left\{
       \left.
       \begin{pmatrix}
       u_1 & \pres{\rt}{a}^\tc \\
       a & u_2
       \end{pmatrix}
       \right|
       a\in\GL_r(O_{E_v}),
       u_1,u_2\in\Herm_r(O_{F_v})
       \right\}
       \subseteq\Herm_{2r}(F_v);
       \]

  \item $\widehat{\tf}$ satisfies
       \[
       \int_{G_r(F_v)}\langle\pi_v^\vee(\varsigma(a,u_1,u_2))\varphi_v^\vee,\varphi_v\rangle_{\pi_v}\cdot
       \widehat{\tf}
       \(
       \begin{pmatrix}
       u_1 & \pres{\rt}{a}^\tc \\
       a & u_2
       \end{pmatrix}
       \)\cdot\rd\varsigma(a,u_1,u_2)=1,
       \]
       where the integration is in fact a finite sum by (2) and $\varsigma$ is the map in Lemma \ref{le:zeta}.
\end{enumerate}
\end{definition}

\begin{remark}\label{re:typical}
It is clear that $(\varphi_v^\vee,\varphi_v)$-typical element exists if $\langle\pi_v^\vee(\tw_r)\varphi_v^\vee,\varphi_v\rangle_{\pi_v}\in\dQ^\times$.
\end{remark}

\begin{lem}\label{le:typical}
Consider elements $\varphi_v^\vee\in\pi_v^\vee$, $\varphi_v\in\pi_v$, and a $(\varphi_v^\vee,\varphi_v)$-typical element $\tf\in\sS(\Herm_{2r}(F_v))$. If $\chi_v$ is unramified, then
\[
Z^\iota(\varphi_v^\vee\otimes\varphi_v,\tf^{\chi_v})=1
\]
(see \eqref{eq:section} for $\tf^{\chi_v}$) holds for every $\iota\colon\dL\to\dC$.
\end{lem}

\begin{proof}
This is immediate from Lemma \ref{le:zeta} and Definition \ref{de:typical}.
\end{proof}

This following lemma will not be used until Section \ref{ss:4}.

\begin{lem}\label{le:zeta1}
For every $\varphi_v^\vee\in\pi_v^\vee$, $\varphi_v\in\pi_v$ and a $\dQ$-valued section $f\in\rI^\Box_{r,v}(\CF)$, there exists a unique element
\[
Z(\varphi_v^\vee\otimes\varphi_v,f)\in\dL
\]
such that for every embedding $\iota\colon\dL\to\dC$, $\iota Z(\varphi_v^\vee\otimes\varphi_v,f)$ coincides with $Z^\iota(\varphi_v^\vee\otimes\varphi_v,f)$.
\end{lem}

\begin{proof}
We may regard $\rI^\Box_{r,v}(\CF)$ as a representation with coefficients in $\dQ$. Let $\Omega$ be the set of all embeddings $\iota\colon\dL\to\dC$. The assignment
\[
(\varphi_v^\vee\otimes\varphi_v,f)\mapsto\left\{Z^\iota(\varphi_v^\vee\otimes\varphi_v,f)\right\}_{\iota\in\Omega}
\]
defines an element
\[
\fZ\in\Hom_{G_r(F_v)\times G_r(F_v)}\((\pi_v^\vee\boxtimes\pi_v)\otimes\rI^\Box_{r,v}(\CF),\dC^{\Omega}\).
\]
We need to show that $\fZ$ takes values in $\dL$, which is tautologically a subring of $\dC^\Omega$. By \cite{LL}*{Proposition~3.6(1)}, it suffices to find one pair of elements $(\varphi_v^\vee\otimes\varphi_v,f)$ such that $\fZ(\varphi_v^\vee\otimes\varphi_v,f)\in\dL^\times$. Indeed, choose $\varphi_v^\vee\in\pi_v^\vee$, $\varphi_v\in\pi_v$ such that $\langle\pi_v^\vee(\tw_r)\varphi_v^\vee,\varphi_v\rangle_{\pi_v}=1$, and a $(\varphi_v^\vee,\varphi_v)$-typical element $\tf\in\sS(\Herm_{2r}(F_v))$ (which exists by Remark \ref{re:typical}). Then $\tf^\CF$ is $\dQ$-valued. By Lemma \ref{le:typical}, $\fZ(\varphi_v^\vee\otimes\varphi_v,\tf^\CF)=1\in\dL^\times$. The lemma is proved.
\end{proof}

\begin{lem}\label{le:zeta2}
Suppose that $v\not\in\tV_F^{(p)}$. If $\pi_v$ is unramified with respect to $K_{r,v}$ and $\varphi_v^\vee,\varphi_v$ are both $K_{r,v}$-invariant such that $\langle\varphi_v^\vee,\varphi_v\rangle_{\pi_v}=1$, then for every $\iota\colon\dL\to\dC$,
\[
Z^\iota(\varphi_v^\vee\otimes\varphi_v,f_{\chi_v}^\sph)
=\frac{L(\tfrac{1}{2},\BC(\iota\pi_v)\otimes(\chi_v\circ\Nm_{E_v/F_v}))}{b_{2r,v}(\chi)},
\]
where $f_{\chi_v}^\sph$ is defined in Notation \ref{no:section1}(2).
\end{lem}

\begin{proof}
This is a well-known calculation of Piatetski-Shapiro and Rallis. See \cite{Li92}*{Theorem~3.1} for a full account including our case.
\end{proof}

\begin{proposition}\label{pr:zeta}
Suppose that $v\in\tV_F^{(p)}$ and that $\pi_v$ is Panchishkin unramified. For every embedding $\iota\colon\dL\to\dC$, we have
\[
\int_{G_r(F_v)}\iota\xi_{\pi_v}(g)\cdot(\tf_{\chi_v}^{[0]})^{\chi_v}(\bw_r(g,1_{2r}))\rd g=
q_v^{d_vr^2}\prod_{u\in\tP_v}\gamma(\tfrac{1+r}{2},\iota\ul{\pi_u}\otimes\chi_v,\psi_{F,v})^{-1},
\]
where $\xi_{\pi_v}$ and $\ul{\pi_u}$ are introduced in Proposition \ref{pr:satake}.
\end{proposition}

Note that the left-hand side is a local doubling zeta integral.

\begin{proof}
To ease notation, we omit $v$ and $\iota$ in the proof. In particular, $\varpi^d$ generates the different ideal of $F/\dQ_p$, and $\xi_\pi$ is $\dC$-valued.

By Lemma \ref{le:zeta}, we have
\begin{align}\label{eq:zeta1}
&\quad\int_{G_r(F)}\xi_\pi(g)\cdot(\tf_\chi^{[0]})^\chi(\bw_r(g,1_{2r}))\rd g \\
&=\int_{G_r(F)}\xi_\pi(\varsigma(a,u_1,u_2))\cdot\chi(\Nm_{E/F}\det a)\cdot
|\Nm_{E/F}\det a|_{F_v}^r \notag\\
&\qquad\qquad\times\widehat{\tf^{[0]}_\chi}(a,u_1,u_2)\cdot\rd\varsigma(a,u_1,u_2), \notag
\end{align}
where
\[
\widehat{\tf^{[0]}_\chi}(a,u_1,u_2)\coloneqq
\int_{\Herm_{2r}(F)}\tf^{[0]}_\chi(T^\Box)\psi_F
\(\tr
\begin{pmatrix}
u_1 & \pres{\rt}a^c \\
a & u_2 \\
\end{pmatrix}
\begin{pmatrix}
T^\Box_{11} & T^\Box_{12} \\
T^\Box_{21} & T^\Box_{22} \\
\end{pmatrix}
\)
\rd T^\Box
\]
in which $\rd T^\Box$ is the self-dual measure with respect to $\psi_F$.

It follows easily that
\begin{align}\label{eq:zeta2}
&\quad\widehat{\tf^{[0]}_\chi}(a,u_1,u_2)=\\
&
\begin{dcases}
q^{-dr^2}\int_{\GL_r(O_E)}\chi(\Nm_{E/F}\det T)\psi_F\(\Tr_{E/F}\tr aT\)\rd T & \text{if $u_1,u_2\in\varpi^{-d}\Herm_r(O_F)$,} \\
0 & \text{otherwise,}
\end{dcases} \notag
\end{align}
in which $\rd T$ is the self-dual measure on $\Mat_{r,r}(E)$ with respect to $\psi_F$.

Since $\xi_\pi$ is bi-$I$-invariant, \eqref{eq:zeta2} implies that
\begin{align*}
\eqref{eq:zeta1}=q^{dr^2}\int_{\GL_r(E)}
&\xi_\pi(m(a)\tw_r)\cdot\chi(\Nm_{E/F}\det a)\cdot|\Nm_{E/F}\det a|_F^r \\
&\times\(\int_{\GL_r(O_E)}\chi(\Nm_{E/F}\det T)\psi_F\(\Tr_{E/F}\tr aT\)\rd T\)\rd a,
\end{align*}
which, by Proposition \ref{pr:satake}(3), equals
\begin{align*}
&\quad q^{dr^2}\prod_{u\in\tP}
\(\int_{\GL_r(F)}\langle\ul{\pi_u}(a)\phi_u,\phi_u^\vee\rangle_{(\ul{\pi_u})^\vee}\cdot\chi(\det a)\cdot|\dtm a|_F^r
\cdot X(a)\rd a\)\\
&=q^{dr^2}
\prod_{u\in\tP}\(\int_{\GL_r(F)}\langle(\ul{\pi_u}\otimes\chi)(a)\phi_u,\phi_u^\vee\rangle_{(\ul{\pi_u}\otimes\chi)^\vee}
\cdot|\dtm a|_F^r\cdot X(a)\rd a\)
\end{align*}
with
\begin{align*}
X(a)\coloneqq\int_{\GL_r(O_F)}\chi(\det T)\psi_F\(\tr aT\)\rd T.
\end{align*}
Applying \cite{Jac79}*{Proposition~1.2(3)} with $\Phi=(\chi\circ\det)\cdot\CF_{\GL_r(O_F)}$, we have
\begin{align*}
&\quad\int_{\GL_r(F)}\langle(\ul{\pi_u}\otimes\chi)(a)\phi_u,\phi_u^\vee\rangle_{(\ul{\pi_u}\otimes\chi)^\vee}\cdot|\dtm a|_F^r\cdot X(a)\rd a \\
&=\gamma(\tfrac{1-r}{2},(\ul{\pi_u}\otimes\chi)^\vee,\psi_F) \\
&\quad\int_{\GL_r(F)}\langle\phi_u,(\ul{\pi_u}\otimes\chi)^\vee(a)\phi_u^\vee\rangle_{(\ul{\pi_u}\otimes\chi)^\vee}\cdot
|\dtm a|_F^{\frac{1-r}{2}}\cdot\chi(\det a)\cdot\CF_{\GL_r(O_F)}(a)\rd a \\
&=\gamma(\tfrac{1-r}{2},(\ul{\pi_u}\otimes\chi)^\vee,\psi_F)\cdot\langle\phi_u^\vee,\phi_u\rangle_{\ul{\pi_u}} \\
&=\gamma(\tfrac{1+r}{2},\ul{\pi_u}\otimes\chi,\psi_F)^{-1}\cdot\langle\phi_u^\vee,\phi_u\rangle_{\ul{\pi_u}}.
\end{align*}
Together, we have
\begin{align*}
\eqref{eq:zeta1}&=q^{dr^2}
\prod_{u\in\tP}\gamma(\tfrac{1+r}{2},\ul{\pi_u}\otimes\chi,\psi_F)^{-1}\cdot\langle\phi_u^\vee,\phi_u\rangle_{\ul{\pi_u}} \\
&=q^{dr^2}\xi_\pi(\tw_r)\prod_{u\in\tP}\gamma(\tfrac{1+r}{2},\ul{\pi_u}\otimes\chi,\psi_F)^{-1} \\
&=q^{dr^2}\prod_{u\in\tP}\gamma(\tfrac{1+r}{2},\ul{\pi_u}\otimes\chi,\psi_F)^{-1}.
\end{align*}
The proposition is proved.
\end{proof}

\subsection{Construction of the $p$-adic $L$-function}
\label{ss:construction}

Let $\pi$ be a relevant $\dL$-representation of $G_r(\dA_F^\infty)$ for some finite extension $\dL/\dQ_p$ contained in $\ol\dQ_p$ such that $\pi_v$ is \emph{Panchishkin unramified} for every $v\in\tV_F^{(p)}$.

Choose a finite set $\lozenge$ of places of $\dQ$ containing $\{\infty,p\}$ such that $\pi_v$ is unramified (hence $v\not\in\tV_F^\ram$) for every $v\in\tV_F\setminus\tV_F^{(\lozenge)}$.

We choose decomposable elements $\varphi_1=\otimes_v\varphi_{1,v}\in\cV_{\hat\pi}$ and $\varphi_2=\otimes_v\varphi_{2,v}\in\cV_\pi$ satisfying
\begin{enumerate}[label=(T\arabic*)]
  \item $\langle\pi_v^\vee(\tw_r)\varphi^\dag_{1,v},\varphi_{2,v}\rangle_{\pi_v}=1$ for $v\in\tV_F^{(\lozenge\setminus\{\infty,p\})}$,

  \item $\varphi^\dag_{1,v}\in(\pi_v^\vee)^-$, $\varphi_{2,v}\in\pi_v^-$ and $\langle\pi_v^\vee(\tw_r)\varphi^\dag_{1,v},\varphi_{2,v}\rangle_{\pi_v}=q_v^{-d_vr^2}$ for $v\in\tV_F^{(p)}$,

  \item $\varphi^\dag_{1,v}\in(\pi_v^\vee)^{K_{r,v}}$, $\varphi_{2,v}\in\pi_v^{K_{r,v}}$ and $\langle\varphi^\dag_{1,v},\varphi_{2,v}\rangle_{\pi_v}=1$ for $v\in\tV_F\setminus\tV_F^{(\lozenge)}$.
\end{enumerate}
Note that (T2) is possible by Proposition \ref{pr:satake}(2). We also choose a $(\varphi^\dag_{1,v},\varphi_{2,v})$-typical element $\tf_v\in\sS(\Herm_{2r}(F_v))$ (Definition \ref{de:typical}, which exists by (T1) and Remark \ref{re:typical}) for $v\in\tV_F^{(\lozenge\setminus\{\infty,p\})}$.

For every finite character $\chi\colon\Gamma_{F,p}\to\dC^\times$, put
\[
f_{\chi^{\infty p}}\coloneqq\bigotimes_{v\in\tV_F^\fin\setminus\tV_F^{(p)}}f_{\chi_v}\in\rI_r^\Box(\chi)^{\infty p},
\]
where $f_{\chi_v}\in\rI_{r,v}^\Box(\chi_v)$ is the section $\tf_v^{\chi_v}$ (resp.\ $f_{\chi_v}^\sph$) for $v\in\tV_F^{(\lozenge\setminus\{\infty,p\})}$ (resp.\ $v\in\tV_F\setminus\tV_F^{(\lozenge)}$).

Consider an open compact subset $\Omega\subseteq\Gamma_{F,p}$. By the linear independence of characters, one can write
\[
\CF_\Omega=\sum_{i}c_i\cdot\chi_i
\]
as a finite sum in a unique way with $c_i\in\dC$ and finite characters $\chi_i\colon\Gamma_{F,p}\to\dC^\times$. For an element $e\in\dZ^\tP$, we put
\begin{align}\label{eq:omega}
\widetilde{D}^{[e]}_{\lozenge}(-,\Omega)\coloneqq
\sum_{i}c_i\widetilde{D}^{[e]}_{\lozenge}(-,\chi_i,f_{\chi_i^{\infty p}}),
\end{align}
where $\widetilde{D}^{[e]}_{\lozenge}(-,\chi_i,f_{\chi_i^{\infty p}})$ is defined in \eqref{eq:eisenstein}.

For every $w\not\in\lozenge$, we choose a nonnegative power $\Delta_w$ of $w$ such that
\[
G_{2r}(F_w)\bigcap\(\widetilde\cG_{2r}(\dZ_w)\times_{\widetilde\cG_{2r}(\dZ_w/\Delta_w)}\widetilde\cP_{2r}(\dZ_w/\Delta_w)\)
\]
is contained in $K_{2r,w}$ (and we may take $\Delta_w=1$ when $w$ is unramified in $E$). For every $w\in\lozenge\setminus\{\infty,p\}$, we may choose a nonnegative power $\Delta_w$ of $w$ such that $\otimes_{v\in\tV_F^{(w)}}f_{\chi_v}$ is fixed by the kernel of the map $\widetilde\cG_{2r}(\dZ_w)\to\widetilde\cG_{2r}(\dZ_w/\Delta_w)$ for every finite character $\chi\colon\Gamma_{F,p}\to\dC^\times$. Indeed, by Definition \ref{de:typical}(2), the restriction of $f_{\chi_v}$ to $K_{2r,v}$ is independent of $\chi_v$, which implies the existence of $\Delta_w$. Finally, put
\[
\Delta\coloneqq\prod_{w\not\in\lozenge}\Delta_w,\quad
\Delta'\coloneqq\prod_{w\in\lozenge\setminus\{\infty,p\}}\Delta_w.
\]

\begin{lem}\label{le:rational3}
For every open compact subset $\Omega\subseteq\Gamma_{F,p}$, if $\|e\|>0$, then
\[
\widetilde\bbh_{r,r}\(\widetilde{D}^{[e]}_{\lozenge}(-,\Omega)\)
\in\varinjlim_{d\in\dN}\widetilde\cH_{r,r}^{[r]}(\widetilde{K}_{r,r}(p^d\Delta,\Delta'))
\]
(Notation \ref{no:open}).
\end{lem}

\begin{proof}
By construction and Lemma \ref{le:section3}, it is clear that
\[
\widetilde\bbh_{r,r}\(\widetilde{D}^{[e]}_{\lozenge}(-,\Omega)\)
\in\varinjlim_{d\in\dN}\widetilde\cH_{r,r}^{[r]}(\widetilde{K}_{r,r}(p^d\Delta,\Delta'))\otimes_\dQ\dC.
\]
It remains to show the rationality when $\|e\|>0$.

Take an arbitrary element $\sigma\in\Aut(\dC/\dQ)$. We have $\sigma f_{\chi_v}=f_{\sigma\chi_v}$ for every $v\in\tV_F^\fin\setminus\tV_F^{(p)}$ and every finite character $\chi\colon\Gamma_{F,p}\to\dC^\times$ by construction. By Proposition \ref{pr:rational}, we have
\[
\sigma\widetilde\bbh_{r,r}\(\widetilde{D}^{[e]}_{\lozenge}(-,\Omega)\)=
\widetilde\bbh_{r,r}\(\sum_{i}\sigma(c_i)\widetilde{D}^{[e]}_{\lozenge}
(-,\sigma\chi_i,f_{\sigma\chi_i^{\infty p}})\).
\]
On the other hand, we have $\CF_\Omega=\sigma\CF_\Omega=\sum_i\sigma(c_i)\cdot\sigma\chi_i$, which implies that
\[
\widetilde\bbh_{r,r}\(\sum_{i}\sigma(c_i)\widetilde{D}^{[e]}_{\lozenge}
(-,\sigma\chi_i,f_{\sigma\chi_i^{\infty p}})\)
=\widetilde\bbh_{r,r}\(\widetilde{D}^{[e]}_{\lozenge}(-,\Omega)\).
\]
The lemma is proved.
\end{proof}

\begin{lem}\label{le:rational4}
For every open compact subset $\Omega\subseteq\Gamma_{F,p}$, if $\|e\|>0$, then there is a unique element
\[
\cD^{[e]}_{\lozenge}(-,\Omega)\in\varinjlim_{d\in\dN}\cH_{r,r}^{[r]}(K_{r,r}(p^d\Delta,\Delta'))
\]
(Notation \ref{no:open}) such that
\[
\xi_{r,r}^*\cD^{[e]}_{\lozenge}(-,\Omega)=\zeta_{r,r}^*\widetilde\bbh_{r,r}\(\widetilde{D}^{[e]}_{\lozenge}(-,\Omega)\)
\]
in terms of the diagram \eqref{eq:holomorphic5}.
\end{lem}

\begin{proof}
Since the center of $\widetilde{G}_{2r}(\dA^\infty)$ (as a subgroup of $\widetilde{G}_{r,r}(\dA^\infty)$) acts trivially on $\widetilde{D}^{[e]}_{\lozenge}(-,\Omega)$, the element $\zeta_{r,r}^*\widetilde{D}^{[e]}_{\lozenge}(-,\Omega)$ descends to the desired element $\cD^{[e]}_{\lozenge}(-,\Omega)$.
\end{proof}

\begin{notation}\label{no:doubling}
By Remark \ref{re:holomorphic}(2), we have a map
\[
\pr_{\pi,\hat\pi}\coloneqq\pr_\pi\otimes\pr_{\hat\pi}\colon
\cH_{r,r}^{[r]}=\cH_r^{[r]}\otimes_{\dQ_p}\cH_r^{[r]}\to\cV_\pi\otimes_\dL\cV_{\hat\pi}
\]
that is the tensor product of $\pr_\pi$ and $\pr_{\hat\pi}$ from Lemma \ref{le:projection}. In what follows, for every $\Psi\in\cH_{r,r}^{[r]}$, $\varphi_1\in\cV_{\hat\pi}$ and $\varphi_2\in\cV_\pi$, we put
\[
\langle\varphi_1\otimes\varphi_2,\Psi\rangle_{\pi,\hat\pi}\coloneqq
\left\langle\varphi_2^\dag,\left\langle \varphi_1^\dag, \pr_{\pi,\hat\pi}\Psi\right\rangle_\pi\right\rangle_{\hat\pi}.
\]
\end{notation}

\begin{definition}\label{de:measure}
We define an $\dL$-valued distribution $\sL_p^\lozenge(\pi)$ on $\Gamma_{F,p}$ to be the following assignment
\[
\Omega\subseteq\Gamma_{F,p}\mapsto
\(\prod_{u\in\tP}\alpha(\pi_u)\)^{-1}
\left\langle\varphi_1\otimes\varphi_2, \cD^{[1]}_{\lozenge}(-,\Omega)\right\rangle_{\pi,\hat\pi},
\]
which is additive from the construction. Here, $1$ is regarded as a constant tuple in $\dN^\tP$.
\end{definition}

\begin{theorem}\label{th:measure}
The distribution $\sL_p^\lozenge(\pi)$ on $\Gamma_{F,p}$ in Definition \ref{de:measure} is a $p$-adic measure. Moreover, for every finite (continuous) character $\chi\colon\Gamma_{F,p}\to\ol\dQ_p^\times$ and every embedding $\iota\colon\ol\dQ_p\to\dC$, we have
\begin{align*}
\iota\sL_p^\lozenge(\pi)(\chi)&=
\frac{1}{\rP^\iota_\pi}\cdot\frac{Z_r^{[F:\dQ]}}{b_{2r}^\lozenge(\CF)} \cdot \prod_{v\in\tV_F^{(p)}}\prod_{u\in\tP_v}
\gamma(\tfrac{1+r}{2},\iota(\ul{\pi_u}\otimes\chi_v),\psi_{F,v})^{-1} \\
&\quad\times L(\tfrac{1}{2},\BC(\iota\pi^\lozenge)\otimes(\iota\chi^\lozenge\circ\Nm_{E/F})),
\end{align*}
where
\[
Z_r\coloneqq(-1)^r2^{-2r^2}\cdot 2^{r^2-r}\bpi^{r^2}\frac{\Gamma(1)\cdots\Gamma(r)}{\Gamma(r+1)\cdots\Gamma(2r)}
\]
and $\ul{\pi_u}$ is introduced in Proposition \ref{pr:satake}. In particular, in terms of the data chosen from this subsection, $\sL_p^\lozenge(\pi)$ depends on $\lozenge$ only, justifying its notation.
\end{theorem}

\begin{proof}
For the first statement, it amounts to showing that the map
\[
\Omega\mapsto\int_\Omega\rd\sL_p^\lozenge(\pi)\coloneqq
\(\prod_{u\in\tP}\alpha(\pi_u)\)^{-1}
\left\langle\varphi_1\otimes\varphi_2, \cD^{[1]}_{\lozenge}(-,\Omega)\right\rangle_{\pi,\hat\pi}\in\dL
\]
is uniformly bounded.

Now we show the uniform boundedness. By Lemma \ref{le:open} below, for every $\Omega$, there is an integer $e=e_\Omega\geq 1$, regarded as a constant tuple in $\dN^\tP$, such that
\[
D^{[e]}_{\lozenge}(-,\Omega)\in\cH_{r,r}^{[r]}(K_{r,r}(p\Delta,\Delta')).
\]
By (T2) and Lemma \ref{le:section5}, we have
\begin{align*}
\int_\Omega\rd\sL_p^\lozenge(\pi)&=\(\prod_{u\in\tP}\alpha(\pi_u)\)^{-e}
\left\langle\varphi_1\otimes(\rT_p^-)^{e-1}\varphi_2,\cD^{[1]}_{\lozenge}(-,\Omega)\right\rangle_{\pi,\hat\pi} \\
&=\(\prod_{u\in\tP}\alpha(\pi_u)\)^{-e}
\left\langle \varphi_1\otimes\varphi_2,(1\times\rU_p^{e-1})\cD^{[1]}_{\lozenge}(-,\Omega)\right\rangle_{\pi,\hat\pi} \\
&=\(\prod_{u\in\tP}\alpha(\pi_u)\)^{-e}
\left\langle \varphi_1\otimes\varphi_2,\cD^{[e]}_{\lozenge}(-,\Omega)\right\rangle_{\pi,\hat\pi},
\end{align*}
where $\rT_p^-\coloneqq\prod_{u\in\tP}\rT_u^-$. Since $\alpha(\pi_u)\in O_\dL^\times$ for every $u\in\tP$ and $\cH_{r,r}^{[r]}(K(p\Delta,\Delta'))$ is a finite-dimensional $\dQ_p$-vector space, it suffices to show that there exists an integer $M\geq 0$ such that
\[
p^M\bbq_{r,r}^\an(g\cdot D^{[e]}_{\lozenge}(-,\Omega))\in\SF_{r,r}(\ol\dZ_{(p)})
\]
holds for every $g\in G_{r,r}(\dA^\lozenge)$, every $e\geq 1$ and every $\Omega$. By \eqref{eq:holomorphic4}, it suffices to study $\bbq_{2r}^\an(g\cdot E^{[e]}_{\lozenge}(-,\Omega))$, where
\[
E^{[e]}_{\lozenge}(-,\Omega)\coloneqq\sum_{i}c_iE^{[e]}_{\lozenge}(-,\chi_i,f_{\chi_i^{\infty p}})
\]
is similarly defined as \eqref{eq:omega}. By Lemma \ref{le:rational} and Lemma \ref{le:section4}(1), we know that for every $g\in G_{r,r}(\dA^\lozenge)$, there exists a finitely generated ring $\dO_g$ contained in $\ol\dZ_{(p)}$ and a collection of $\dO_g$-valued functions $\{\pres{g}\sfW_{T^\Box}^\lozenge\res T^\Box\in \Herm_{2r}^\circ(F)\}$ on $\dA_F^{\lozenge,\times}/\(O_F\otimes\prod_{v\not\in\lozenge}\dZ_w\)^\times$ with finite support such that
\begin{align*}
&\bbq_{2r}^\an(g\cdot E^{[e]}_{\lozenge}(-,\Omega)) =W_{2r}^\lozenge\sum_{T^\Box\in\Herm_{2r}^\circ(F)^+} \\
&
\(\prod_{v\in\tV_F^{(\lozenge\setminus\{\infty,p\})}}\widehat{\tf_v}(T^\Box)
\CF_{\fT_p^{[e]}}(T^\Box)\sum_{x}
\CF_\Omega((\Nm_{E_p/F_p}\det T^\Box_{12},x))\cdot\pres{g}\sfW_{T^\Box}^\lozenge(x)\) q^{T^\Box}
\end{align*}
in which $x$ is taken over the set $\dA_F^{\lozenge,\times}/\(O_F\otimes\prod_{v\not\in\lozenge}\dZ_w\)^\times$ and $(\Nm_{E_p/F_p}\det T^\Box_{12},x)$ is regarded as an element in $\Gamma_{F,p}$. We may choose an integer $M$ such that
\[
p^M W_{2r}^\lozenge\cdot\prod_{v\in\tV_F^{(\lozenge\setminus\{\infty,p\})}}\widehat{\tf_v}(T^\Box)\in\dZ_{(p)}
\]
for every $T^\Box\in\Herm_{2r}^\circ(F)$. It follows that $p^M\bbq_{2r}^\an(g\cdot E^{[e]}_{\lozenge}(-,\Omega))\in\SF_{2r}(\ol\dZ_{(p)})$ holds for every $\Omega$, every $e\geq 1$, and every $g\in G_{r,r}(\dA^\lozenge)$. Thus, $\rd\sL_p^\lozenge(\pi)$ is a $p$-adic measure.

Next, we show the second statement, that is, the interpolation property. By construction, Remark \ref{re:period} and Lemma \ref{le:section5}, for every finite character $\chi\colon\Gamma_{F,p}\to\ol\dQ_p^\times$ and embedding $\iota\colon\ol\dQ_p\to\dC$,
\begin{align*}
&\quad\iota\sL_p^\lozenge(\pi)(\chi) \\
&=\(\iota\prod_{u\in\tP}\alpha(\pi_u)\)^{-1}
\frac{1}{(\rP^\iota_\pi)^2}\iint\limits_{\(G_r(F)\backslash G_r(\dA_F)\)^2}\varphi_1^\iota(g_1^\dag)\varphi_2^\iota(g_2^\dag)E^{[1]}_{\lozenge}
((g_1,g_2),\iota\chi,f_{\iota\chi^{\infty p}})
\rd g_1\rd g_2 \notag\\
&=\frac{1}{(\rP^\iota_\pi)^2}\iint\limits_{\(G_r(F)\backslash G_r(\dA_F)\)^2}
\varphi_1^\iota(g_1^\dag)\varphi_2^\iota(g_2^\dag)E^{[0]}_{\lozenge}
((g_1,g_2),\iota\chi,f_{\iota\chi^{\infty p}})
\rd g_1\rd g_2 \notag\\
&=\frac{1}{(\rP^\iota_\pi)^2}
\iint\limits_{\(G_r(F)\backslash G_r(\dA_F)\)^2}(\varphi_1^\dag)^\iota(g_1)\varphi_2^\iota(g_2)E^{[0]}_{\lozenge}((g_1,g_2^\dag),\iota\chi,
f_{\iota\chi^{\infty p}})\rd g_1\rd g_2. \notag
\end{align*}
By \eqref{eq:eisenstein2},
\begin{align*}
&\quad\iota\sL_p^\lozenge(\pi)(\chi) \\
&=\frac{1}{(\rP^\iota_\pi)^2}\cdot b_{2r}^\lozenge(\CF)^{-1} \cdot b_{2r}^\lozenge(\iota\chi) \\
&\quad\iint\limits_{\(G_r(F)\backslash  G_r(\dA_F)\)^2}(\varphi_1^\dag)^\iota(g_1)\varphi_2^\iota(g_2)E((g_1,g_2^\dag),f_\infty^{[r]}\otimes (\tf_{\iota\chi_p}^{[0]})^{\iota\chi_p}\otimes
f_{\iota\chi^{\infty p}})\rd g_1\rd g_2 \\
&=\frac{1}{(\rP^\iota_\pi)^2}\cdot b_{2r}^\lozenge(\CF)^{-1} \cdot b_{2r}^\lozenge(\iota\chi)  \\ &\quad\iint\limits_{\(G_r(F)\backslash G_r(\dA_F)\)^2}(\varphi_1^\dag)^\iota(g_1)\varphi_2^\iota(g_2)E(\imath(g_1,g_2),f_\infty^{[r]}\otimes (\tf_{\iota\chi_p}^{[0]})^{\iota\chi_p}\otimes
f_{\iota\chi^{\infty p}})\rd g_1\rd g_2,
\end{align*}
where we have used $(g_1,g_2^\dag)=\imath(g_1,g_2)$ as in Remark \ref{re:dagger}. By the well-known doubling integral expansion (see \cite{Ral82} or \cite{Liu11}*{Section~2B} in the case of unitary groups) and Lemma \ref{le:zeta2}, we have
\begin{align*}
&\quad\iint\limits_{\(G_r(F)\backslash G_r(\dA_F)\)^2}(\varphi_1^\dag)^\iota(g_1)\varphi_2^\iota(g_2)E(\imath(g_1,g_2),f_\infty^{[r]}\otimes (\tf_{\iota\chi_p}^{[0]})^{\iota\chi_p}
f_{\iota\chi^{\infty p}})\rd g_1\rd g_2 \\
&=\frac{L(\tfrac{1}{2},\BC(\iota\pi^\lozenge)\otimes(\iota\chi^\lozenge\circ\Nm_{E/F}))}
{b_{2r}^\lozenge(\iota\chi)}\cdot Z((\varphi_1^\dag)^\iota_\infty\otimes(\varphi_2^\iota)_\infty,f_\infty^{[r]})\\
&\quad\times
\prod_{v\in\tV_F^{(p)}}Z^\iota(\varphi_{1,v}^\dag\otimes\varphi_{2,v},(\tf_{\iota\chi_v}^{[0]})^{\iota\chi_v})
\cdot
\prod_{v\in\tV_F^{(\lozenge\setminus\{\infty,p\})}}Z^\iota(\varphi_{1,v}^\dag\otimes\varphi_{2,v},f_{\iota\chi_v}).
\end{align*}
There are three cases:
\begin{itemize}
  \item By \cite{EL}*{Theorem~1.3 \& Proposition~3.3.2} (with $n=k=2r$, $a=b=r$, $\tau_1=\cdots=\tau_r=r$, $\nu_1=\cdots=\nu_r=-r$, and $\chi_{\r{ac}}^r=1$), we have (see the proof of \cite{LL}*{Proposition~3.7} for more details)
      \[
      Z((\varphi_1^\dag)^\iota_\infty\otimes(\varphi_2^\iota)_\infty,f_\infty^{[r]})=\rP^\iota_\pi\cdot Z_r^{[F:\dQ]}.
      \]

  \item By (T2) and Proposition \ref{pr:zeta}, for $v\in\tV_F^{(p)}$, we have
      \[
      Z^\iota(\varphi_{1,v}^\dag\otimes\varphi_{2,v},(\tf_{\iota\chi_v}^{[0]})^{\iota\chi_v})=
      \prod_{u\in\tP_v}\gamma(\tfrac{1+r}{2},\iota(\ul{\pi_u}\otimes\chi_v),\psi_{F,v})^{-1}.
      \]

  \item By Lemma \ref{le:typical}, for $v\in\tV_F^{(\lozenge\setminus\{\infty,p\})}$, we have $Z^\iota(\varphi_{1,v}^\dag\otimes\varphi_{2,v},f_{\iota\chi_v})=1$.
\end{itemize}
Putting together, we have
\begin{align*}
\iota\sL_p^\lozenge(\pi)(\chi)
&=\frac{1}{\rP^\iota_\pi}\cdot\frac{Z_r^{[F:\dQ]}}{b_{2r}^\lozenge(\CF)} \cdot \prod_{v\in\tV_F^{(p)}}\prod_{u\in\tP_v}
\gamma(\tfrac{1+r}{2},\iota(\ul{\pi_u}\otimes\chi_v),\psi_{F,v})^{-1} \\
&\quad\times L(\tfrac{1}{2},\BC(\iota\pi^\lozenge)\otimes(\iota\chi^\lozenge\circ\Nm_{E/F})).
\end{align*}
The theorem is proved.
\end{proof}

\begin{lem}\label{le:open}
For every finite character $\chi\colon\Gamma_{F,p}\to\dC^\times$, there exists $e_\chi\in\dN$ such that for every $e\in\dN^\tP$ satisfying $e_u\geq e_\chi$ for every $u\in\tP$, the section $(\tf_{\chi_p}^{[e]})^{\chi_p}$ is invariant under $\prod_{v\in\tV_F^{(p)}}\cG_{r,r}(O_{F_v})\times_{\cG_{r,r}(O_{F_v}/\varpi_v)}\cP_{r,r}(O_{F_v}/\varpi_v)$ (Definition \ref{de:lattice}).
\end{lem}

\begin{proof}
It is well-known that for every $v\in\tV_F^{(p)}$ and $u\in\tP_u$, we have $I_v^d\rU_u I_v^d=I_v^{d-1}\rU_u I_v^d$ for every integer $d\geq 2$, where $I_v^d\coloneqq\cG_r(O_{F_v})\times_{\cG_r(O_{F_v}/\varpi_v)}\cP_r(O_{F_v}/\varpi_v^d)$. Since $(\tf_{\chi_p}^{[0]})^{\chi_p}$ is fixed by $\prod_{v\in\tV_F^{(p)}}\cP_{r,r}(O_{F_v})$, it follows that there exists a pair $(e_{\chi1},e_{\chi2})\in\dN\times\dN$ such that for every $(e_1,e_2)\in\dN^\tP\times\dN^\tP$ satisfying $e_{1,u}\geq e_{\chi1}$ and $e_{2,u}\geq e_{\chi2}$ for every $u\in\tP$, $(\rU_p^{e_1}\times\rU_p^{e_2})(\tf_{\chi_p}^{[0]})^{\chi_p}$ is invariant under $\prod_{v\in\tV_F^{(p)}}\cG_{r,r}(O_{F_v})\times_{\cG_{r,r}(O_{F_v}/\varpi_v)}\cP_{r,r}(O_{F_v}/\varpi_v)$. By Lemma \ref{le:section5}, we have
\[
(\rU_p^{e_1}\times\rU_p^{e_2})(\tf_{\chi_p}^{[0]})^{\chi_p}=(\tf_{\chi_p}^{[e_1^\tc+e_2]})^{\chi_p}.
\]
Thus, the lemma follows by taking $e_\chi=e_{\chi1}+e_{\chi2}$.
\end{proof}

To end this subsection, we discuss the parity of the vanishing order of $\sL^\lozenge_p(\pi)$ at $\b1$. For every $v\in\tV_F^\fin$, the root number $\epsilon(\BC(\pi_v)\otimes\iota\chi_v\circ\Nm_{E_v/F_v})$ does not depend on the finite character $\chi\colon\Gamma_{F,p}\to\ol\dQ_p^\times$ and the embedding $\iota\colon\ol\dQ_p\to\dC$, which we denote by $\epsilon(\pi_v)$. Put $\epsilon(\pi)\coloneqq\prod_{v\in\tV_F^\fin}\epsilon(\pi_v)$, which is indeed a finite product.

\begin{proposition}\label{pr:complete}
The vanishing order of $\sL^\lozenge_p(\pi)$ at $\b1$ has the same parity as
\[
r[F:\dQ]+\frac{1-\epsilon(\pi)}{2}.
\]
\end{proposition}

\begin{proof}
Denote by $\Gamma_{F,p}^\lozenge$ the subgroup of $\Gamma_{F,p}$ generated by uniformizers above $\lozenge\setminus\{\infty,p\}$. For every $v\in\tV_F^{(\lozenge\setminus\{\infty,p\})}$, there is a unique element $\sL_p(\pi_v)\in\dL(\Gamma_{F,p}^\lozenge)$ such that for every finite character $\chi\colon\Gamma_{F,p}\to\ol\dQ_p^\times$ and every embedding $\iota\colon\ol\dQ_p\to\dC$,
\[
\iota\sL_p(\pi_v)(\chi)=L(\tfrac{1}{2},\BC(\iota\pi_v)\otimes(\iota\chi_v\circ\Nm_{E_v/F_v})).
\]
In particular, $\sL_p(\pi_v)$ has neither poles nor zeros at points corresponding to finite characters.

Put
\[
\sL_p(\pi)\coloneqq\sL_p^\lozenge(\pi)\prod_{v\in\tV_F^{(\lozenge\setminus\{\infty,p\})}}b_{2r,v}(\CF)^{-1}\cdot\sL_p(\pi_v),
\]
regarded as an element of $\dZ_p[[\Gamma_{F,p}]]\otimes_{\dZ_p[\Gamma_{F,p}^\lozenge]}\dL(\Gamma_{F,p}^\lozenge)$. Then $\sL_p(\pi)$ is the unique element such that for every finite character $\chi\colon\Gamma_{F,p}\to\ol\dQ_p^\times$ and every embedding $\iota\colon\ol\dQ_p\to\dC$,
\begin{align*}
\iota\sL_p(\pi)(\chi)&=
\frac{1}{\rP^\iota_\pi}\cdot\frac{Z_r^{[F:\dQ]}}{b_{2r}^{\infty p}(\CF)} \cdot \prod_{v\in\tV_F^{(p)}}\prod_{u\in\tP_v}
\gamma(\tfrac{1+r}{2},\iota(\ul{\pi_u}\otimes\chi_v),\psi_{F,v})^{-1} \\
&\quad\times L(\tfrac{1}{2},\BC(\iota\pi^p)\otimes(\iota\chi^p\circ\Nm_{E/F}))
\end{align*}
holds. As $\BC(\iota\hat\pi_v)\simeq\BC(\iota\pi^\vee_v)\simeq\BC(\iota\pi_v)\circ\tc$ for every $v\in\tV_F$, we have the functional equation
\begin{align}\label{eq:complete}
\iota\sL_p(\pi)(\chi)=\epsilon(\BC(\pres{\iota}\pi)\otimes\iota\chi\circ\Nm_{E/F})\cdot\iota\sL_p(\pi)(\chi^{-1})
\end{align}
(Definition \ref{de:relevant}) for the root number
\[
\epsilon(\BC(\pres{\iota}\pi)\otimes\iota\chi\circ\Nm_{E/F})=
\prod_{v\in\tV_F}\epsilon(\BC(\pres{\iota}\pi_v)\otimes\iota\chi_v\circ\Nm_{E_v/F_v})\in\{\pm 1\}.
\]
It is clear that for $v\in\tV_F^{(\infty)}$, $\epsilon(\BC(\pres{\iota}\pi_v)\otimes\iota\chi_v\circ\Nm_{E_v/F_v})=(-1)^r$; and for $v\in\tV_F^\fin$, $\epsilon(\BC(\pres{\iota}\pi_v)\otimes\iota\chi_v\circ\Nm_{E_v/F_v})=\epsilon(\pi_v)$ by definition.

To summarize, if we denote by $\vee$ the involution on $\dZ_p[[\Gamma_{F,p}]]\otimes_{\dZ_p[\Gamma_{F,p}^\lozenge]}\dL(\Gamma_{F,p}^\lozenge)$ induced by the inverse homomorphism of $\Gamma_{F,p}$, then \eqref{eq:complete} implies the functional equation
\[
\sL_p(\pi)=(-1)^{r[F:\dQ]}\epsilon(\pi)\cdot\sL_p(\pi)^\vee.
\]
It follows that the vanishing order of $\sL_p(\pi)$ at $\b1$ has the same parity as $r[F:\dQ]+\frac{1-\epsilon(\pi)}{2}$. The proposition is then proved since the vanishing order of $\sL^\lozenge_p(\pi)$ at $\b1$ is same as that of $\sL_p(\pi)$.
\end{proof}

\begin{remark}\label{re:complete}
We expect that the $p$-adic $L$-function $\sL_p(\pi)$ constructed in the proof of Proposition \ref{pr:complete} is again a $p$-adic measure, that is, an element of $\dZ_p[[\Gamma_{F,p}]]\otimes_{\dZ_p}\dL$.
\end{remark}

\subsection{Remarks on $p$-adic measures}
\label{ss:measure}

In this subsection, we review some facts about derivatives of $p$-adic measures and make some remarks that will be used in the next section. For $d\geq 1$, we denote by $U_d$ the image of $1+O_F\otimes p^d\dZ_p$ in $\Gamma_{F,p}$, which is an open subgroup of finite index.

Let $\mu$ be an $\dL$-valued $p$-adic measure on $\Gamma_{F,p}$ (for a finite extension $\dL/\dQ_p$). For every continuous character $\chi\colon\Gamma_{F,p}\to R^\times$ for an $\dL$-affinoid algebra $R$, we put
\[
\mu(\chi)\coloneqq\int_{\Gamma_{F,p}}\chi\rd\mu\coloneqq\lim_{d\to\infty}\sum_{x\in\Gamma_d}\chi(x)\vol(xU_d,\mu)
\]
where $\Gamma_1\subseteq\Gamma_2\subseteq\cdots$ is an arbitrary increasing chain of sets of representatives of $\Gamma_{F,p}/U_d$ for $d=1,2,\dots$. Then $\mu(\chi)$ does not depend on $(\Gamma_d)_d$; and hence $\mu$ is nothing but a bounded rigid analytic function on $\sX_{F,p}\otimes_{\dQ_p}\dL$, or equivalently, an element in $\dZ_p[[\Gamma_{F,p}]]\otimes_{\dZ_p}\dL$.

We consider its derivative $\partial\mu(\CF)$ at $\CF$, which is an element in $\Gamma_{F,p}\otimes_{\dZ_p}\dL$ -- the cotangent space of $\sX_{F,p}\otimes_{\dQ_p}\dL$ at $\CF$. Since in this article, we only consider derivatives at $\CF$, we will simply write $\partial\mu$ for $\partial\mu(\CF)$ in the rest of the writing.

More precisely, $\partial\mu$ is the linear functional in $\Hom_{\dZ_p}(\Hom_{\dZ_p}(\Gamma_{F,p},\dZ_p),\dL)$ that sends $\lambda\in\Hom_{\dZ_p}(\Gamma_{F,p},\dZ_p)$ to
\begin{align*}
\partial_\lambda\mu
&\coloneqq\lim_{c\to\infty}\frac{1}{p^c}\(\mu(\exp(p^c\lambda))-\mu(\CF)\) \\
&=
\lim_{c\to\infty}\frac{1}{p^c}\lim_{d\to\infty}\sum_{x\in\Gamma_d}\(\exp(p^c\lambda(x))-1\)\vol(xU_d,\mu).
\end{align*}
Since
\[
\frac{1}{p^c}\(\exp(p^c\lambda(x))-1\)=\lambda(x)+\frac{p^c\lambda(x)^2}{2!}+\frac{p^{2c}\lambda(x)^3}{3!}+\cdots,
\]
and $\vol(xU_d,\mu)$ is bounded independent of $x$ and $d$, we have
\begin{align}\label{eq:derivative}
\partial_\lambda\mu=\lim_{d\to\infty}\sum_{x\in\Gamma_d}\lambda(x)\vol(xU_d,\mu).
\end{align}

\begin{definition}\label{de:integral}
We say that an $\dL$-valued $p$-adic measure $\mu$ on $\Gamma_{F,p}$ is \emph{integral} if $\vol(\Omega,\mu)\in O_\dL$ for every open compact subset $\Omega\subseteq\Gamma_{F,p}$, that is, $\mu$ belongs to $O_\dL[[\Gamma_{F,p}]]$.
\end{definition}

\begin{lem}\label{le:integral}
Let $\mu$ be an integral $\dL$-valued $p$-adic measure on $\Gamma_{F,p}$. Then for every $\lambda\in\Hom_{\dZ_p}(\Gamma_{F,p},\dZ_p)$ and every $d\geq 1$, we have
\[
\partial_\lambda\mu-\sum_{x\in\Gamma_d}\lambda(x)\vol(xU_d,\mu)\in p^d O_\dL.
\]
In particular, $\partial\mu\in\Gamma_{F,p}^\free\otimes_{\dZ_p}O_\dL$.
\end{lem}

\begin{proof}
Since $U_d\subseteq p^d\Gamma_{F,p}$, we have $\lambda(x)-\lambda(x')\in p^d O_\dL$ if $x=x'$ in $\Gamma_{F,p}/U_d$. Then the lemma follows from \eqref{eq:derivative} since $\mu$ is integral.
\end{proof}

The discussion of this subsection can be easily generalized to $p$-adic measures valued in a finite product of finite extensions of $\dL$.

\section{Selmer theta lifts and their $p$-adic heights}
\label{ss:4}

In this section, we introduce Selmer theta lifts and study their $p$-adic heights. We fix an embedding $E\hookrightarrow\dC$ and regard $E$ as a subfield of $\dC$ and regard $\ol{E}$ as the algebraic closure of $E$ in $\dC$. Fix an even positive integer $n=2r$.

\subsection{Hermitian spaces and Weil representations}
\label{ss:setup}

Let $\pi$ be a relevant $\dL$-representation of $G_r(\dA_F^\infty)$ for some finite extension $\dL/\dQ_p$ contained in $\ol\dQ_p$.

Choose a finite set $\lozenge$ of places of $\dQ$ containing $\{\infty,p\}$ such that $\pi_v$ is unramified (hence $v\not\in\tV_F^\ram$) for every $v\in\tV_F\setminus\tV_F^{(\lozenge)}$.

Let $V,(\;,\;)_V$ be a hermitian space (that is nondegenerate and $E$-linear in the second variable) over $E$ of rank $n$ that is split at every $v\in\tV_F\setminus\tV_F^{(\lozenge)}$, has signature $(n-1,1)$ along the induced inclusion $F\subseteq\dR$ and signature $(n,0)$ at other archimedean places of $F$. We introduce the following sets of notation.

\begin{enumerate}[label=(H\arabic*)]
  \item For every $F$-ring $R$ and every integer $m\geq 0$, we denote by
      \[
      T(x)\coloneqq\((x_i,x_j)_V\)_{i,j}\in\Herm_m(R)
      \]
      the moment matrix of an element $x=(x_1,\dots,x_m)\in V^m\otimes_FR$.

  \item For every $v\in\tV_F$, we put $\epsilon_v\coloneqq\eta_{E/F}((-1)^r\dtm V_v)\in\{\pm 1\}$. In particular, $\epsilon_v=1$ for $v\not\in\tV_F^{(\lozenge)}$.

  \item Let $v\in\tV_F^\fin$ be an element and $m\geq 0$ an integer.
      \begin{itemize}
        \item For $T\in\Herm_m(F_v)$, we put $(V^m_v)_T\coloneqq\{x\in V^m_v\res T(x)=T\}$, and
            \[
            (V^m_v)_\reg\coloneqq\bigcup_{T\in\Herm_m^\circ(F_v)}(V^m_v)_T,
            \]
            where we recall $\Herm_m^\circ$ from \S\ref{ss:run}(F3).

        \item For every $\dZ[p_v^{-1}]\langle p_v\rangle$-ring $\dM$, we have a Fourier transform map
            \[
            \widehat{\phantom{a}}\colon \sS(V_v^m,\dM)\to \sS(V_v^m,\dM)
            \]
            sending $\phi$ to $\widehat\phi$ defined by the formula
            \[
            \widehat\phi(x)\coloneqq\int_{V_v^m}\phi(y)\psi_{F,v}\(\Tr_{E_v/F_v}\sum_{i=1}^m(x_i,y_i)_V\)\rd y,
            \]
            which is in fact a finite sum, where $\r{d}y$ is the self-dual Haar measure on $V_v^m$ with respect to $\psi_{F,v}$. In what follows, we will always use this self-dual Haar measure on $V_v^m$.
      \end{itemize}

  \item Put $H\coloneqq\rU(V)$, which is a reductive group over $F$.

  \item For $v\in\tV_F^\fin\setminus\{v\in\tV_F^\ram\res\text{either $\epsilon_v=-1$ or $v\mid 2$}\}$, a \emph{good lattice} of $V_v$ is an $O_{E_v}$-lattice $\Lambda_v$ of $V_v$ that is a subgroup of $\Lambda_v^\vee$ of index $q_v^{1-\epsilon_v}$, where
      \[
      \Lambda_v^\vee\coloneqq\{x\in V_v\res\Tr_{E_v/F_v}(x,y)_V\in \fp_v^{-d_v}\text{ for every }y\in\Lambda_v\}.
      \]
      We say that
      \begin{itemize}
        \item an open compact subgroup $L^\lozenge$ of $H(\dA_F^{\lozenge})$ is \emph{good} if it is the product of the stabilizers of good lattices at $v\not\in\tV_F\setminus\tV_F^{(\lozenge)}$;

        \item a Schwartz function $\phi^\lozenge\in\sS(V_v^m\otimes_F\dA_F^\lozenge)$ is \emph{good} if the it is the product of $\CF_{\Lambda_v^m}$ in which $\Lambda_v$ is a good lattice of $V_v$.
      \end{itemize}

  \item Denote by $\dT^\lozenge$ the (abstract) spherical Hecke algebra of rank-$n$ unitary groups over $E$ away from $\lozenge$, and $\dS^\lozenge$ its subring consisting of Hecke operators supported on places split in $E$. In particular, for every good open compact subgroup $L^\lozenge$ of $H(\dA_F^{\lozenge})$, we have canonical isomorphisms
      \begin{align*}
      \dT^\lozenge&=\dZ[L^\lozenge\backslash H(\dA_F^{\lozenge})/L^\lozenge],\\
      \dS^\lozenge&=\varinjlim_{\substack{\tT\subseteq\tV_F^\spl\setminus\tV_F^{(\lozenge)}\\|\tT|<\infty}}
      \dZ[L^\lozenge_\tT\backslash H(F_\tT)/L^\lozenge_\tT]\otimes\CF_{L^{\lozenge\tT}}
      \end{align*}
      of commutative rings.

  \item For every integer $m\geq 1$, every $v\in\tV_F^\fin$ and every $\dZ[p_v^{-1}]\langle p_v\rangle$-ring $\dM$, we have the Weil representation $\omega_{m,v}$ of $G_m(F_v)\times H(F_v)$ on $\sS(V_v^m,\dM)$ given by the following formulae:
      \begin{itemize}
         \item for $a\in\GL_m(E_v)$ and $\phi\in\sS(V_v^m,\dM)$, we have
           \[
           \omega_{m,v}(m(a))\phi(x)=|\dtm a|_{E_v}^r\cdot\phi(x a);
           \]

         \item for $b\in\Herm_m(F_v)$ and $\phi\in\sS(V_v^m,\dM)$, we have
           \[
           \omega_{m,v}(n(b))\phi(x)=\psi_{F,v}(\tr b T(x))\phi(x);
           \]

        \item for $\phi\in\sS(V_v^m,\dM)$, we have
           \[
           \omega_{m,v}\(\tw_m\)\phi(x)=\gamma_{V_v,\psi_{F,v}}^m\cdot\widehat\phi(x),
           \]
           where $\gamma_{V_v,\psi_{F,v}}\in\{\pm 1\}$ is the Weil constant of $V_v$ with respect to $\psi_{F,v}$;

        \item for $h\in H(F_v)$ and $\phi\in \sS(V_v^m,\dM)$, we have
           \[
           \omega_{m,v}(h)\phi(x)=\phi(h^{-1}x).
           \]
      \end{itemize}

  \item When $m=n=2r$, we have the \emph{Siegel--Weil section} map
      \[
      f^{\r{SW}}_-\colon \sS(V_v^{2r})\to\rI^\Box_{r,v}(\CF)
      \]
      for $v\in\tV_F^\fin$ sending $\Phi$ to $f_\Phi^{\r{SW}}$ defined by the formula
      \[
      f_\Phi^{\r{SW}}(g)=\(\omega_{2r,v}(g)\Phi\)(0), \quad g\in G_{2r}(F_v)=G_r^\Box(F_v).
      \]

  \item For every $v\in\tV_F^\fin$, there is a unique $\dQ$-valued Haar measure $\rd h_v$ on $H(F_v)$, called the \emph{Siegel--Weil measure}, satisfying that for every $T^\Box\in\Herm_{2r}^\circ(F_v)$ and every $\Phi\in\sS(V_v^{2r})$,\footnote{We recall our convention from \S\ref{ss:notation} that $\sS(V_v^{2r})$ means $\sS(V_v^{2r},\dC)$.}
      \[
      I_{T^\Box}(\Phi)\coloneqq\int_{H(F_v)}\Phi(h_v^{-1}x)\rd h_v=b_{2r,v}(\CF)\cdot W_{T^\Box}(f^{\r{SW}}_\Phi),
      \]
      where $x$ is an arbitrary element in $(V_v^{2r})_{T^\Box}$. When $v$ is unramified over $\dQ$ and $H\otimes_FF_v$ is unramified, the measure $\rd h_v$ gives volume $1$ to every hyperspecial maximal subgroup of $H(F_v)$. For $\Phi\in\sS(V_v^{2r},R)$ with $R$ a general $\dQ$-ring, $I_{T^\Box}(\Phi)$ is well-defined and belongs to $R$.

  \item Let $\iota\colon\dL\to\dC$ be an embedding. For every $v\in\tV_F^\fin$, put
      \[
      \theta(\iota\pi_v)\coloneqq\Hom_{G_r(F_v)}(\sS(V_{\pi_v}^r),\iota\pi_v)
      \]
      as a complex representation of $H(F_v)$. Then put
      \[
      \theta(\iota\pi)\coloneqq\otimes'_{v\in\tV_F^\fin}\theta(\iota\pi_v)
      \]
      as a complex representation of $H(\dA_F^\infty)$.
\end{enumerate}

\begin{lem}\label{le:theta}
For every $v\in\tV_F^\fin$, there exists a (unique up to isomorphism) hermitian space $V_{\pi_v}$ over $E_v$ of rank $n$ such that for every embedding $\iota\colon\dL\to\dC$, $\theta(\iota\pi_v)\neq 0$ if and only if $V_v\simeq V_{\pi_v}$. When $V_v\simeq V_{\pi_v}$, $\theta(\iota\pi_v)$ is a tempered irreducible admissible representation of $H(F_v)$, satisfying
\[
\Hom_{H(F_v)}\(\sS(V_{\pi_v}^r),\theta(\iota\pi_v)\)\simeq\iota\pi_v
\]
as $\dC[G_r(F_v)]$-modules.
\end{lem}

\begin{proof}
For every fixed embedding $\iota\colon\dL\to\dC$, the existence and the uniqueness of $V_{\pi_v}$ follow from the local theta dichotomy \cite{GG11}*{Theorem~3.10} (see also \cite{HKS96}*{Corollary~4.4} and \cite{Har07}*{Theorem~2.1.7}). As $V_{\pi_v}$ does not change if we twist the additive character $\psi_{F,v}$ by automorphisms of $\dC$, it is the same for all $\iota$.

Since $\iota\pi_v$ is tempered, the irreducibility and the temperedness of $\theta(\iota\pi_v)$ follow from \cite{GI16}*{Theorem~4.1(v)} and (the same argument for) \cite{GS12}*{Theorem~1.3(ii)}, respectively. The last isomorphism follows from the dual statements.
\end{proof}

\begin{definition}\label{de:coherent}
We say that $V$ (as above) is \emph{$\pi$-coherent} if $V_v\simeq V_{\pi_v}$ for every $v\in\tV_F^\fin$.
\end{definition}

\begin{remark}\label{re:coherent}
We have the following remarks.
\begin{enumerate}
  \item The following conditions are equivalent for $V$ as above: $\pi$-coherent, $\pi^\vee$-coherent, $\pi^\dag$-coherent, $\hat\pi$-coherent.

  \item There exists a hermitian space $V$ as above that is $\pi$-coherent if and only if
      \begin{align}\label{eq:incoherent}
      \prod_{v\in\tV_F^{(\lozenge\setminus\{\infty\})}}\eta_{E/F}\((-1)^r\det V_{\pi_v}\)=-(-1)^{r[F:\dQ]}.
      \end{align}
      Moreover, when \eqref{eq:incoherent} holds, $\sL_p^\lozenge(\pi)$ vanishes at $\CF$.
\end{enumerate}

\end{remark}

In the rest of this subsection, we discuss the rationality of local theta liftings. For readers who are willing to fix an embedding $\ol\dQ_p\hookrightarrow\dC$ and do not care about the rationality of the coefficients of the Selmer theta lifts below, this discussion may be ignored.

Take a place $v\in\tV_F^\fin$. We say that $\pi_v$ is \emph{symmetric} if for every element $a\in F_v^\times$, $\pi_v^{\dag_a}\simeq\pi_v$, where $\dag_a$ is the automorphism of $G_r(F_v)$ given by the conjugation of the element $\(\begin{smallmatrix}1_r & \\ & a1_r \end{smallmatrix}\)\in\GL_{2r}(E_v)$.\footnote{Note that $\pi_v^{\dag_a}\simeq\pi_v$ when $a\in\Nm_{E_v/F_v}E_v^\times$.} Denote by $\tU_\pi$ the subset of $\tV_F^\fin$ consisting of $v$ such that $\pi_v$ is \emph{not} symmetric.

It is easy to see that
\begin{align}\label{eq:weil}
\left\{(a_w)_w\in\prod_{w<\infty}\dZ_w^\times=\widehat\dZ^\times \left| a_w\in \bigcap_{v\in\tU_\pi\cap\tV_F^{(w)}}\Nm_{E_v/F_v}E_v^\times\right.\right\}
\end{align}
is an open subgroup of $\widehat\dZ^\times$. Thus, we may define $\dQ_\pi$ to be the finite abelian extension of $\dQ$ contained in $\dC$ determined by this subgroup via the global class field theory.

\begin{remark}\label{re:weil}
We have the following remarks concerning $\dQ_\pi$.
\begin{enumerate}
  \item It is clear that $\dQ_{\hat\pi}=\dQ_{\pi}$ since $\tU_{\hat\pi}=\tU_{\pi^\dag}=\tU_{\pi}$.

  \item It is clear that in \eqref{eq:weil}, we may replace $\tU_\pi$ by $\tU_\pi\cap\tV_F^\ram$.

  \item Suppose that we are in the situation of Assumption \ref{st:main}. For every $v\in\tV_F^\ram$ and every $a\in O_{F_v}^\times$, since $\dag_a$ preserves $K_{r,v}$, we have $\pi_v^{\dag_a}\simeq\pi_v$, that is, $\pi_v$ is symmetric. In other words, $\tU_\pi\cap\tV_F^\ram=\emptyset$ hence $\dQ_\pi=\dQ$. For this reason, the readers may just assume $\dQ_\pi=\dQ$ for further reading.

  \item Since every $p$-adic place of $F$ splits in $E$, $p$ is unramified in $\dQ_\pi$.
\end{enumerate}

\end{remark}

\begin{lem}\label{le:theta3}
For every $v\in\tV_F^\fin$, every embedding $\iota\colon\dL\to\dC$, and every $\sigma\in\Aut(\dC/\dQ_\pi)$, $\theta(\sigma\iota\pi_v)$ is isomorphic to $\sigma\theta(\iota\pi_v)$.
\end{lem}

\begin{proof}
We have $\sigma\theta(\iota\pi_v)=\Hom_{G_r(F_v)}(\sigma\sS(V_{\pi_v}^r),\sigma\iota\pi_v)$. The representation $\sigma\sS(V_{\pi_v}^r)$ has the same formulae of definition as $\sS(V_{\pi_v}^r)$ except that $n(b)$ sends $\phi$ to the function $x\mapsto\psi_{F,v}(a\tr b T(x))\phi(x)$ for some element $a\in\dZ_{p_v}^\times$ (resp. $a\in\dZ_{p_v}^\times\cap\Nm_{E_v/F_v}E_v^\times$) if $v\not\in\tU_\pi$ (resp.\ $v\in\tU_\pi$). It follows that
\[
\Hom_{G_r(F_v)}(\sigma\sS(V_{\pi_v}^r),\sigma\iota\pi_v)\simeq\Hom_{G_r(F_v)}(\sS(V_{\pi_v}^r),\sigma\iota\pi_v^{\dag_a})
=\theta(\sigma\iota\pi_v^{\dag_a}).
\]
By definition, we have $\pi_v^{\dag_a}\simeq\pi_v$. The lemma follows.
\end{proof}

\subsection{$p$-adic height pairing on unitary Shimura varieties}
\label{ss:height}

From this subsection, we will assume $F\neq\dQ$. Put $\dL_\pi\coloneqq\dL\otimes_\dQ\dQ_\pi$.

Back to the setup in \S\ref{ss:setup}, we have the projective system of Shimura varieties $X_L$ associated with $\Res_{F/\dQ}H$ indexed by neat open compact subgroups $L\subseteq H(\dA_F^\infty)$, which are smooth \emph{projective} schemes over $E$ of dimension $n-1$. Put $\ol{X}_L\coloneqq X_L\otimes_E\ol{E}$.

\begin{lem}\label{le:theta1}
For every $L$, there exists a unique decomposition
\[
\rH^{2r-1}_{\et}(\ol{X}_L,\dL_\pi(r))=
\rH^{2r-1}_{\et}(\ol{X}_L,\dL_\pi(r))[\theta_\pi]\oplus\rH^{2r-1}_{\et}(\ol{X}_L,\dL_\pi(r))[\widehat{\theta_\pi}]
\]
of $\dL_\pi[L\backslash H(\dA_F^\infty)/L]$-modules such that for every homomorphism $\iota\colon\dL_\pi\to\dC$ extending the inclusion $\dQ_\pi\subseteq\dC$, $\iota\rH^{2r-1}_{\et}(\ol{X}_L,\dL_\pi(r))[\theta_\pi]$ is isomorphic to a finite sum of copies of $\theta(\iota\pi)^L$ (\S\ref{ss:setup}(H10)) and $\iota\rH^{2r-1}_{\et}(\ol{X}_L,\dL_\pi(r))[\widehat{\theta_\pi}]$ does not contain $\theta(\iota\pi)^L$ as a subquotient.\footnote{We warn the readers that the statement could be wrong if we replace $\dL_\pi$ by $\dL$.}
\end{lem}

In what follow, we put $\rV_{\pi,L}\coloneqq\rH^{2r-1}_{\et}(\ol{X}_L,\dL_\pi(r))[\theta_\pi]$. It is clear that $\rV_{\pi,L}$ is nonzero only if $V$ is $\pi$-coherent (Definition \ref{de:coherent}).

\begin{proof}
For every given $\iota$, the existence of such a decomposition follows from Matsushima's formula. It follows from Lemma \ref{le:theta3} that these decompositions are the same for all $\iota$.
\end{proof}

The Hochschild--Serre spectral sequence in \cite{Jan88}*{Corollary~3.4} induces a decreasing filtration
\[
\cdots\subseteq\rF^1\rH^{2r}_{\et}(X_L,\dL_\pi(r))
\subseteq\rF^0\rH^{2r}_{\et}(X_L,\dL_\pi(r))=\rH^{2r}_{\et}(X_L,\dL_\pi(r))
\]
of $\rH^{2r}_{\et}(X_L,\dL_\pi(r))$ in the category of $\dL_\pi[L\backslash H(\dA_F^\infty)/L]$-modules so that there is a canonical inclusion\footnote{Indeed, the above Hochschild--Serre spectral sequence should degenerate from the second page, as pointed out in \cite{Jan00}*{Page~262}, which implies that $\frac{\rF^i\rH^{2r}_{\et}(X_L,\dL_\pi(r))}{\rF^{i+1}\rH^{2r}_{\et}(X_L,\dL_\pi(r))}
=\rH^i(E,\rH^{2r-i}_{\et}(\ol{X}_L,\dL_\pi(r)))$ for every $i\geq 0$.}
\[
\frac{\rF^1\rH^{2r}_{\et}(X_L,\dL_\pi(r))}{\rF^2\rH^{2r}_{\et}(X_L,\dL_\pi(r))}
\subseteq\rH^1(E,\rH^{2r-1}_{\et}(\ol{X}_L,\dL_\pi(r))).
\]

\begin{lem}\label{le:theta2}
There exists a unique map of $\dL_\pi[L\backslash H(\dA_F^\infty)/L]$-modules
\[
\wp_\pi\colon\rH^{2r}_{\et}(X_L,\dL_\pi(r))\to\rH^1(E,\rH^{2r-1}_{\et}(\ol{X}_L,\dL_\pi(r))[\theta_\pi])
\]
that vanishes on $\rF^2\rH^{2r}_{\et}(X_L,\dL_\pi(r))$ and induces the inclusion map on
\[
\rH^1(E,\rH^{2r-1}_{\et}(\ol{X}_L,\dL_\pi(r))[\theta_\pi])\cap
\frac{\rF^1\rH^{2r}_{\et}(X_L,\dL_\pi(r))}{\rF^{2}\rH^{2r}_{\et}(X_L,\dL_\pi(r))}.
\]
\end{lem}

\begin{proof}
By Lemma \ref{le:theta}, $\theta_\pi$ is tempered, so that it does not appear in $\rH^i_{\et}(\ol{X}_L,\dL_\pi(r))$ for $i\neq 2r-1$ as a subquotient. In particular, since $\frac{\rF^{i}\rH^{2r}_{\et}(X_L,\dL_\pi(r))}{\rF^{i+1}\rH^{2r}_{\et}(X_L,\dL_\pi(r))}$ is a subquotient of $\rH^i(E,\rH^{2r-i}_{\et}(\ol{X}_L,\dL_\pi(r)))$, we have
\[
\Hom_{\dL_\pi[L\backslash H(\dA_F^\infty)/L]}\(\rF^{2}\rH^{2r}_{\et}(X_L,\dL_\pi(r)),\rH^{2r-1}_{\et}(\ol{X}_L,\dL_\pi(r))[\theta_\pi]\)=0.
\]
Then by Lemma \ref{le:theta1}, we have a unique map
\[
\wp_\pi^1\colon\rF^1\rH^{2r}_{\et}(X_L,\dL_\pi(r))\to\rH^1(E,\rH^{2r-1}_{\et}(\ol{X}_L,\dL_\pi(r))[\theta_\pi])
\]
satisfying the property in the lemma. It remains to show that $\wp_\pi^1$ extends uniquely to $\rH^{2r}_{\et}(X_L,\dL_\pi(r))$. The uniqueness follows from the same reason that $\theta_\pi$ does not appear in $\rH^{2r}_{\et}(\ol{X}_L,\dL_\pi(r))$. For the existence, note that \cite{LL}*{Proposition~6.9(1)} actually implies that there exists an element $\rs\in\dS^\lozenge_\dL$ (by possibly enlarging $\lozenge$) such that $\rs^*$ annihilates $\rH^{2r}_{\et}(\ol{X}_L,\dL_\pi(r))$ and acts by the identity map on $\rH^{2r-1}_{\et}(\ol{X}_L,\dL_\pi(r))[\theta_\pi]$. In particular, $\wp_\pi\coloneqq\wp_\pi^1\circ\rs^*$ is such an extension.
\end{proof}

Denote by $\dS^\lozenge_{\pi,L}$ the subset of $\dS^\lozenge_{\dL_\pi}$ consisting of elements $\rs$ such that $\rs^*$ annihilates $\rF^2\rH^{2r}_{\et}(X_L,\dL_\pi(r))$ and the induced endomorphism of $\rH^{2r}_{\et}(X_L,\dL_\pi(r))/\rF^2\rH^{2r}_{\et}(X_L,\dL_\pi(r))$ has image in $\rH^1(E,\rV_{\pi,L})$. It is clear that $\dS^\lozenge_{\pi,L}$ is an ideal. On the other hand, we have the Hecke character
\[
\chi^\lozenge_\pi\colon\dS^\lozenge_\dL\to\dL
\]
given by its action on $\pi$.

\begin{lem}\label{le:hodge}
Suppose that $L$ is of the form $L_{\lozenge} L^\lozenge$ in which $L^\lozenge$ is good (\S\ref{ss:setup}(H5)).
\begin{enumerate}
  \item For every $\rs\in\dS^\lozenge_{\pi,L}$, we have
      \[
      \rs^*=\chi^\lozenge_{\hat\pi}(\rs)\cdot\wp_\pi\colon\rH^{2r}_{\et}(X_L,\dL_\pi(r))\to\rH^1(E,\rV_{\pi,L}).
      \]

  \item If $L_{\lozenge}$ is of the form $\prod_{v\in\tV_F^{(\lozenge\setminus\{\infty\})}}L_v$ in which for every $v\in\tV_F^{(\lozenge\setminus\{\infty\})}\setminus\tV_F^\spl$, $L_v$ is special maximal and $\theta(\iota\pi_v)^{L_v}\neq 0$ for every embedding $\iota\colon\dL\to\dC$, then the restriction of $\chi^\lozenge_{\hat\pi}$ to $\dS^\lozenge_{\pi,L}$ is surjective.
\end{enumerate}
\end{lem}

\begin{proof}
By \cite{Liu11}*{Corollary~A.6(2)}, for every embedding $\iota\colon\dL\to\dC$ and every $v\in\tV_F^\spl\setminus\tV_F^{(\lozenge)}$, we have $\theta(\iota\pi_v)\simeq\hat\pi_v\otimes_{\dL,\iota}\dC$. This already implies (1).

For (2), For every embedding $\iota\colon\dL\to\dC$ and every $i\in\dZ$, we have
\begin{align}\label{eq:matsushima}
&\quad\rH^i_{\et}(\ol{X}_L,\dL(r))\otimes_{\dL,\iota}\dC \\
&\simeq\bigoplus_{\pi'}(\pi')^L\otimes_\dC\Hom_{\dC[L\backslash H(\dA_F^\infty)/L]}\((\pi')^L,\rH^i_{\et}(\ol{X}_L,\dL(r))\otimes_{\dL,\iota}\dC\) \notag
\end{align}
in which $\pi'$ runs over all irreducible admissible (complex) representations of $H(\dA_F^\infty)$, by Matsushima's formula. By \cite{LL}*{Proposition~6.9(1)}, we may find $\rs\in\dS^\lozenge_\dL$ that annihilates $\rH^i_{\et}(\ol{X}_L,\dL(r))$ for every $i\neq 2r-1$ and such that $\chi^\lozenge_{\hat\pi}(\rs)=1$. It remains to show that if $\pi'$ contributes nontrivially to $\rH^{2r-1}_{\et}(\ol{X}_L,\dL(r))$ in \eqref{eq:matsushima} satisfying that $\pi'_v\simeq\theta(\iota\pi_v)$ for every $v\in\tV_F^\spl\setminus\tV_F^{(\lozenge)}$, then the same must hold for every $v\in\tV_F^\fin$. Indeed, by the strong multiplicity one property \cite{Ram}*{Theorem~A} and the local-global compatibility of base change \cite{KMSW}, we already have the isomorphism for $v\in\tV_F^\spl$ and that $\BC(\pi'_v)\simeq\BC(\theta(\iota\pi_v))$ for $v\in\tV_F^\fin\setminus\tV_F^\spl$. Now take an element $v\in\tV_F^\fin\setminus\tV_F^\spl$. Since both $\pi'_v$ and $\theta(\iota\pi_v)$ have nontrivial $L_v$-invariants and $L_v$ is special maximal, they are constituents of the same principal series representation $\rho$ of $H(F_v)$. Since $\rho^{L_v}$ is one-dimensional, they must be the same constituent. Thus, (2) follows.
\end{proof}

\begin{lem}\label{le:hodge4}
For $v\in\tV_F^{(p)}$, if $\pi_v$ is unramified, then $\rV_{\pi,L}$ is crystalline at every place $u$ of $E$ above $v$.
\end{lem}

\begin{proof}
If $\pi_v$ is unramified, then its local theta lift is also an unramified representation of $H(F_v)$. In particular, we may assume that $L$ is of the form $L_vL^v$ in which $L_v$ is hyperspecial maximal. By \cite{RSZ20}*{Theorem~4.5} (or a more closely related discussion after \cite{LL}*{Proposition~7.1}), $X_L$ admits a finite \'{e}tale cover that has smooth reduction at every place $u$ of $E$ above $v$. Thus, $\rV_{\pi,L}$ is crystalline at $u$.
\end{proof}

\begin{lem}\label{le:galois}
There is a unique up to isomorphism semisimple continuous representation $\rho_\pi$ of $\Gal(\ol{E}/E)$ of dimension $n$ with coefficients in $\ol\dQ_p$ such that for every place $u$ of $E$ not above $\lozenge$ that is split over $F$, $\rho_\pi$ is unramified at $u$ and a geometric Frobenius at $u$ acts with a characteristic polynomial that coincides with the Satake polynomial of $\pi_u$, regarded as an unramified representation of $\GL_n(E_u)$. Moreover, we have $\rho_{\hat\pi}\simeq\rho_\pi^\tc\simeq\rho_\pi^\vee(1-n)$.
\end{lem}

\begin{proof}
The uniqueness of $\rho_\pi$ follows from its property and the Chebotarev density theorem; and the last statement follows from the uniqueness. It remains to show the existence of $\rho_\pi$.

Choose an isomorphism $\iota\colon\ol\dQ_p\xrightarrow\sim\dC$. By \cite{Mok15}, the automorphic base change of  $\otimes_{v\in\tV_F^{(\infty)}}\pi^{[r]}_v\otimes\iota\pi$ is an isobaric sum of distinct unitary cuspidal automorphic representations $\Pi_j$ of $\GL_{n_j}(\dA_E)$ for some partition $n=n_1+\cdots+n_s$. By \cite{CH13}*{Theorem~3.2.3}, for each $1\leq j\leq s$, we have a semisimple representation $\rho_{\Pi_j}$ of $\Gal(\ol{E}/E)$ such that for every place $u$ of $E$ not above $\lozenge$ that is split over $F$, the restriction of $\rho_{\Pi_j}$ to the place $u$ is unramified and corresponds to the irreducible admissible representation $\(\Pi_{j,u}\otimes|\;|_{E_u}^{\frac{1-n}{2}}\)\otimes_{\dC,\iota^{-1}}\ol\dQ_p$ of $\GL_{n_j}(E_u)$ under the unramified local Langlands correspondence. Then $\rho_\pi\coloneqq\bigoplus_{j=1}^s\rho_{\Pi_j}$ does the job.
\end{proof}

\begin{hypothesis}\label{hy:galois}
For every homomorphism $\iota\colon\dL_\pi\to\ol\dQ_p$ over $\dL$ and every irreducible $\ol\dQ_p[\Gal(\ol{E}/E)]$-module $\rho$ that is a subquotient of $\rV_{\pi,L}\otimes_{\dL_\pi,\iota}\ol\dQ_p$, $\rho$ is a direct summand of $\rho_\pi(r)$.
\end{hypothesis}

\begin{remark}\label{re:hypothesis}
We have the following remarks concerning Hypothesis \ref{hy:galois}.
\begin{enumerate}
  \item Hypothesis \ref{hy:galois} is equivalent to the parallel statement for $\hat\pi$.

  \item We understand that Hypothesis \ref{hy:galois} will follow from a sequel of the work \cite{KSZ}.

  \item A precise prediction of the semisimplification of $\rV_{\pi,L}\otimes_{\dL_\pi,\iota}\ol\dQ_p$, if not zero, can be found in \cite{LL}*{Hypothesis~6.6}. Such prediction is independent of $\iota$.

  \item It is conjectured that $\rV_{\pi,L}\otimes_{\dL_\pi,\iota}\ol\dQ_p$ is irreducible as a module over the algebra $\ol\dQ_p[L\backslash H(\dA_F^\infty)/L][\Gal(\ol{E}/E)]$. However, this does not seem reachable at this moment.
\end{enumerate}
\end{remark}

From this moment, we will \emph{assume Hypothesis \ref{hy:galois} without further mentioning}.

\begin{lem}\label{le:hodge3}
For every finite place $u$ of $E$ not above $p$, we have
\[
\rH^i(E_u,\rV_{\pi,L})=\rH^i(E_u,\rV_{\hat\pi,L})=0
\]
for every $i\in\dZ$.
\end{lem}

\begin{proof}
By symmetry, we only need to consider $\rV_{\pi,L}$. By Hypothesis \ref{hy:galois}, it suffices to show that $\rH^1(E_u,\rho_\pi(r))=0$ for such $u$. By \cite{Car12}*{Theorem~1.1} and \cite{TY07}*{Lemma~1.4(3)}, we know that the associated Weil--Deligne representation of $\rho_\pi(r)\res_{E_u}$ is pure of weight $-1$, which implies that $\rH^i(E_u,\rho_\pi(r))=0$ by (the proof of) \cite{Nek00}*{Proposition~2.5}.
\end{proof}

\begin{lem}\label{le:hodge2}
Take $v\in\tV_F^{(p)}$. If $\pi_v$ is Panchishkin unramified (Definition \ref{de:satake}), then both $\rV_{\pi,L}\res_{E_u}$ and $\rV_{\hat\pi,L}\res_{E_u}$ satisfy the Panchishkin condition (Definition \ref{de:panchishkin}) and are pure of weight $-1$ for $u$ above $v$.
\end{lem}

\begin{proof}
By symmetry and Lemma \ref{le:satake2}, we only need to consider $\rV_{\pi,L}$. We will use the results and notational conventions introduced in \S\ref{ss:recollection}. Since $\rV_{\pi,L}$ is crystalline (Lemma \ref{le:hodge4}), by Lemma \ref{le:panchishkin} and Hypothesis \ref{hy:galois}, it suffices to show that $\rho_\pi(r)\res_{E_u}$ satisfies the Panchishkin condition and is pure of weight $-1$ for $u$ above $v$. By \cite{Car14}*{Theorem~1.1}, we know that for every embedding $\tau\colon E_u\to\ol\dQ_p$,
\begin{enumerate}
  \item $\rho_\pi(r)\res_{E_u}$ is crystalline and has Hodge--Tate weights $\{-r,-r+1,\dots,r-1\}$ at $\tau$;

  \item the associated Weil--Deligne representation $\r{WD}(\rho_\pi(r)\res_{E_u})_\tau$ (see \S\ref{ss:recollection}) is unramified and its multiset of generalized geometric Frobenius eigenvalues is $\{\alpha_{v,1}\sqrt{q_v}^{-1},\dots,\alpha_{v,n}\sqrt{q_v}^{-1}\}$.
\end{enumerate}
By (2), we know that $\rho_\pi(r)\res_{E_u}$ is pure of weight $-1$. Moreover, by Lemma \ref{le:hodge4} and Remark \ref{rem:fr}, the multiset of generalized $\varphi$-eigenvalues on $\dD\coloneqq\dD_{\cris}(\rho_{\pi}(r)\res_{E_{u}})$ is $\{\alpha_{v,1}\sqrt{q_v}^{-1},\dots,\alpha_{v,n}\sqrt{q_v}^{-1}\}$ as well.

For the Panchishkin condition, by Lemma \ref{le:satake2}, we may assume that the unique subset $J$ of $\{1,\dots,n\}$ with $|J|=r$ such that $\sqrt{q_v}^{r^2}\prod_{j\in J}\alpha_{v,j}\in O_\dL^\times$ is $\{1,\dots,r\}$ without loss of generality. Then $\alpha_{v,j}\sqrt{q_v}^{-1}$ belongs to $\ol\dZ_p$ if and only if $i\geq r+1$. Let $\dD^{+}\subset\dD$ be the $\dL\otimes_{\dQ_{p}}E_{u,0}$-submodule spanned by the generalized eigenspaces with respect to the crystalline Frobenius for the eigenvalues $\{\alpha_{v,j}\sqrt{q_v}^{-1}\res 1\leq j\leq r\}$, which is the negative-slope submodule defined in general in Lemma \ref{le:pan-eq}. By the weak admissibility of $\dD$ and rank counting, the map \eqref{eq:pan-ses} for $\dD^{+}$ is an isomorphism, and by inspection of the Newton and Hodge polygons, $\dD^{+}$ is  weakly admissible. It follows that the equivalent Panchishkin condition of Lemma \ref{le:pan-eq} is satisfied.
\end{proof}

If $\pi_v$ is Panchishkin unramified for every $v\in\tV_F^{(p)}$, then we may apply \S\ref{ss:decomposition} to the case where $K=E$, $X=X_L$, $d=d'=r$, $\dL=\dL_\pi$, $\rV=\rV_{\pi,L}$ and $\rV'=\rV_{\hat\pi,L}$. Indeed, (V1) is due to Lemma \ref{le:hodge3}; (V2) and (V3) are due to Lemma \ref{le:hodge4} and Lemma \ref{le:hodge2}. Consequently, we have a canonical $p$-adic height pairing
\begin{align}\label{eq:heightv}
\langle\;,\;\rangle_{(\rV_{\pi,L},\rV_{\hat\pi,L}),E}
\colon\rH^1_f(E,\rV_{\pi,L})\times\rH^1_f(E,\rV_{\hat\pi,L})\to\Gamma_{E,p}\otimes_{\dZ_p}\dL_\pi.
\end{align}

\subsection{Selmer theta lifts}
\label{ss:theta}

We take a finite set $\blacklozenge$ of places of $\dQ$ containing $\{\infty\}$ and a subfield $\dM$ of $\dC$ containing $\dQ\langle\prod_{w\in\blacklozenge\setminus\{\infty\}}w\rangle$ and $\dQ_\pi$.

Consider a neat open compact subgroup $L\subseteq H(\dA_F^\infty)$. Recall that for every element $x\in V^m\otimes_F\dA_F^\infty$, we have \emph{Kudla's special cycle} $Z(x)_L\in\rZ^m(X_L)$ if $T(x)\in\Herm_m^\circ(F)^+$ and $Z(x)_L\in\CH^m(X_L)_\dQ$ in general. See \cite{LL}*{Section~4} for more details in our setting. For every $\phi\in\sS(V^m\otimes_F\dA_F^\infty,\dM)^{K_m^{\blacklozenge}\times L}$ and every $T\in\Herm_m(F)$, we put
\[
Z_T(\phi)_L\coloneqq\sum_{\substack{x\in L\backslash V^m\otimes_F\dA_F^\infty\\ T(x)=T}}\phi(x) Z(x)_L
\]
as an element in $\rZ^m(X_L)\otimes\dM$ if $T\in\Herm_m^\circ(F)^+$ and in $\CH^m(X_L)\otimes\dM$ in general. Denote by
\[
Z_T^\clubsuit(\phi)_L\in\rH^{2m}_{\et}(X_L,\dQ_p(m))\otimes_\dQ\dM
\]
the image of $Z_T(\phi)_L$ under the (absolute) cycle class map
\[
\rZ^m(X_L)\to\CH^m(X_L)\to\rH^{2m}_{\et}(X_L,\dQ_p(m)).
\]

\begin{definition}\label{de:selmer}
Suppose that $m=r$. We define the \emph{$\pi$-Selmer generating function} to be
\begin{align*}
Z_\phi^\pi(g)_L &\coloneqq\sum_{T\in\Herm_r(F)^+}\wp_{\pi}\(Z_T^\clubsuit(\omega_r(g)\phi)_L\)\cdot q^T \\
&\in \rH^1(E,\rV_{\pi,L})\otimes_{\dL_\pi}\SF_r(\dL\otimes_\dQ\dC)
\end{align*}
for $g\in G_r(\dA_F^\infty)$. Here, $\omega_r$ is the restricted tensor product of $\omega_{r,v}$ (\S\ref{ss:setup}(H7)) over all $v\in\tV_F^\fin$; and $\wp_{\pi}$ is the map in Lemma \ref{le:theta2}.\footnote{Since $\SF_r(\dL\otimes_\dQ\dM)$ \emph{strictly} contains $\dL\otimes_\dQ\SF_r(\dM)$, \emph{a priori} we do not know whether $Z_\phi^\pi(g)_L$ belongs to $\rH^1(E,\rV_{\pi,L})\otimes_{\dQ_\pi}\SF_r(\dC)$.}
\end{definition}

\begin{hypothesis}[Modularity of $\pi$-Selmer generating functions]\label{hy:modularity}
For every element $\phi\in\sS(V^r\otimes_F\dA_F^\infty,\dC)^{K_r^{\blacklozenge}\times L}$, there exists an element
\[
Z^\pi_{\phi,L}\in\rH^1_f(E,\rV_{\pi,L})\otimes_{\dQ_\pi}\cA_{r,\hol}^{[r]}
\]
such that $(1\otimes\bbq_r^\an)(g\cdot Z^\pi_{\phi,L})=Z^\pi_\phi(g)_L$ holds in $\rH^1(E,\rV_{\pi,L})\otimes_{\dL_\pi}\SF_r(\dL\otimes_\dQ\dC)$ for every $g\in G_r(\dA_F^\infty)$, where $\bbq_r^\an$ is the analytic $q$-expansion map (Definition \ref{de:fourier}).
\end{hypothesis}

\begin{remark}\label{re:modularity}
We have the following remarks concerning Hypothesis \ref{hy:modularity}.
\begin{enumerate}
  \item This hypothesis is implied by \cite{LL}*{Hypothesis~4.5}.

  \item The natural map $\dL\otimes_\dQ\dC\to\prod_{\iota\colon\dL\to\dC}\dC$ is injective. Indeed, it suffices to show that for every finitely generated subfields $\dM$ and $\dM'$ of $\dC$, the map $\dM\otimes_\dQ\dM'\to\prod_{\iota\colon\dM\to\dC}\dC$ is injective. For this statement, we immediately reduce it to the case where $\dM$ is purely transcendental, say generated by $x_1,\dots,x_d$ that are algebraically independent over $\dQ$. Then we may find a single embedding $\iota\colon\dM\to\dC$ such that $\iota x_1,\dots,\iota x_d$ are algebraically independent over $\dM'$, so that the map $\dM\otimes_\dQ\dM'\to\dC$ induced by $\iota$ is already injective.

  \item We warn the readers that this hypothesis is stronger than the following statement: For every embedding $\iota\colon\dL\to\dC$, there exists an element $Z^{\pi,\iota}_{\phi,L}$ in the space $\rH^1_f(E,\rV_{\pi,L})\otimes_{\dL_\pi,\iota\times 1}\cA_{r,\hol}^{[r]}$ such that for every $g\in G_r(\dA_F^\infty)$, $(1\otimes\bbq_r^\an)(g\cdot Z^{\pi,\iota}_{\phi,L})$ coincides with the natural image of $Z^\pi_\phi(g)_L$ in $\rH^1(E,\rV_{\pi,L})\otimes_{\dL_\pi,\iota\times 1}\SF_r(\dC)$ induced by $\iota$.\footnote{In particular, Theorem \ref{th:modularity} below is stronger than Theorem \ref{th:main2}.} The stronger statement in Hypothesis \ref{hy:modularity} reflects, in some sense, the conjecture that the image of the absolute cycle class map $\CH^r(X_L)\to\rH^{2r}_{\et}(X_L,\dQ_p(r))$ is a finitely generated abelian group.
\end{enumerate}
\end{remark}

Recall from Definition \ref{de:holomorphic} the $\dQ_p$-vector space $\cH_r^{[r]}$ and the subspaces $\cV_\pi,\cV_{\hat\pi}$ of $\cH_r^{[r]}\otimes_{\dQ_p}\dL$ introduced after Lemma \ref{le:dual}.

\begin{proposition}\label{pr:generating}
Assume that Hypothesis \ref{hy:modularity} holds for $\pi$. Then for every $\phi\in\sS(V^r\otimes_F\dA_F^\infty,\dM)^{K_r^{\blacklozenge}\times L}$, there exists a unique element
\[
\cZ^\pi_{\phi,L}\in\rH^1_f(E,\rV_{\pi,L})\otimes_{\dL_\pi}\(\cV_\pi\otimes_\dQ\dM\)
\]
such that for every embedding $\iota\colon\dL\to\dC$, $(\cZ^\pi_{\phi,L})^\iota$, regarded as an element in $\rH^1_f(E,\rV_{\pi,L})\otimes_{\dQ_\pi}\cA_{r,\hol}^{[r]}$ via the inclusion $\dM\subseteq\dC$, coincides with $Z^\pi_{\phi,L}$.
\end{proposition}

\begin{proof}
We first explain that it suffices to find the element $\cZ^\pi_{\phi,L}$ in
\[
\rH^1_f(E,\rV_{\pi,L})\otimes_{\dQ_p\otimes_\dQ\dQ_\pi}\(\cH_r^{[r]}\otimes_\dQ\dM\).
\]
Indeed, if we can find such elements, then the assignment $\phi\mapsto\cZ^\pi_{\phi,L}$ defines a functional in
\[
\Hom_?\(\sS(V^r\otimes_F\dA_F^\infty,\dM)^{K_r^{\blacklozenge}\times L},\rH^1_f(E,\rV_{\pi,L})\otimes_{\dQ_p\otimes_\dQ\dQ_\pi}\(\cH_r^{[r]}\otimes_\dQ\dM\)\)
\]
with $?=\dQ_\pi[K_r^{\blacklozenge}\backslash G_r(\dA_F^\infty)/K_r^{\blacklozenge}][L\backslash H(\dA_F^\infty)/L]$.
By the definition of $\rV_{\pi,L}$ from Lemma \ref{le:theta1} and Lemma \ref{le:theta}, the functional $\cZ^\pi_{-,L}$ must take values in the (possibly zero) subspace
\[
\rH^1_f(E,\rV_{\pi,L})\otimes_{\dL_\pi}\(\cV_\pi\otimes_\dQ\dM\).
\]

Now we show the existence of $\cZ^\pi_{\phi,L}$ as an element in
\[
\rH^1_f(E,\rV_{\pi,L})\otimes_{\dQ_p\otimes_\dQ\dQ_\pi}\(\cH_r^{[r]}\otimes_\dQ\dM\).
\]
Put $G'_r\coloneqq\Res_{F/\dQ}G_r$, which has been regarded as a subgroup of $\widetilde{G}_r$ in Remark \ref{no:pel}. For every $w\not\in\blacklozenge$, choose a nonnegative power $\Delta_w$ of $w$ such that the intersection of
\[
\widetilde{K}_{r,w}\coloneqq
\widetilde\cG_r(\dZ_w)\times_{\widetilde\cG_r(\dZ_w/\Delta_w)}\widetilde\cP_r(\dZ_w/\Delta_w)
\]
with $G'_r(\dQ_w)$ is contained in $K_{r,w}$ (and we may take $\Delta_w=1$ when $w$ is unramified in $E$). Put $\widetilde{K}_r^\blacklozenge\coloneqq\prod_{w\not\in\blacklozenge}\widetilde{K}_{r,w}$ and $K_r^{\prime\blacklozenge}\coloneqq G'_r(\dA^\blacklozenge)\cap\widetilde{K}_r^\blacklozenge\subseteq K_r^\blacklozenge$.

We claim that for every open compact subgroup $K'$ of $\prod_{w\in\blacklozenge}G'_r(\dQ_w)$, there exists an open compact subgroup $\widetilde{K}$ of $\prod_{w\in\blacklozenge}\widetilde{G}_r(\dQ_w)$ containing $K'$ such that the natural map
\begin{align*}
G'_r(\dQ)\backslash G'_r(\dR)^\ad \times G'_r(\dA^\infty) / K'K_r^{\prime\blacklozenge}
\to \widetilde{G}_r(\dQ)\backslash\widetilde{G}_r(\dR)^\ad \times \widetilde{G}_r(\dA^\infty) / \widetilde{K}\widetilde{K}_r^{\blacklozenge}
\end{align*}
is injective, and hence an open and closed immersion. In fact, since $\widetilde{G}_r(\dQ)$ is discrete in $\widetilde{G}_r(\dA^\infty)$, we have
\[
\varprojlim_{K'\subseteq\widetilde{K}}\widetilde{G}_r(\dQ)\backslash\widetilde{G}_r(\dR)^\ad \times \widetilde{G}_r(\dA^\infty) / \widetilde{K}\widetilde{K}_r^{\blacklozenge}
=\widetilde{G}_r(\dQ)\backslash\widetilde{G}_r(\dR)^\ad \times \widetilde{G}_r(\dA^\infty) /K'\widetilde{K}_r^{\blacklozenge}
\]
Then the claim follows from the obvious injectivity of the map
\[
G'_r(\dQ)\backslash G'_r(\dR)^\ad \times G'_r(\dA^\infty) / K'K_r^{\prime\blacklozenge}
\to \widetilde{G}_r(\dQ)\backslash\widetilde{G}_r(\dR)^\ad \times \widetilde{G}_r(\dA^\infty) /K'\widetilde{K}_r^{\blacklozenge}.
\]

Choose a sufficiently small $K'$ as above such that
\[
Z^\pi_{\phi,L}\in\rH^1_f(E,\rV_{\pi,L})\otimes_{\dQ_\pi}\cA_{r,\hol}^{[r]}(K'K_r^{\blacklozenge}).
\]
By the above claim, we may extend $Z^\pi_{\phi,L}$ by zero to obtain an element
\[
\widetilde{Z}^\pi_{\phi,L}\in\rH^1_f(E,\rV_{\pi,L})\otimes_{\dQ_\pi}
\widetilde\cA_{r,\hol}^{[r]}(\widetilde{K}\widetilde{K}_r^{\blacklozenge})
\]
for some $\widetilde{K}$ as above.

Note that for every $\dL_\pi$-module $\rM$, the commutative diagram
\[
\xymatrix{
\rM\otimes_{\dQ_\pi}\SF_r(\dM) \ar[r]\ar[d]& \rM\otimes_{\dL_\pi}\SF_r(\dL\otimes_\dQ\dM) \ar[d] \\
\rM\otimes_{\dQ_\pi}\SF_r(\dC) \ar[r]& \rM\otimes_{\dL_\pi}\SF_r(\dL\otimes_\dQ\dC)
}
\]
in the category of $\dL_\pi$-modules, in which all arrows are natural inclusions, is Cartesian. Thus, by Lemma \ref{le:rational1}, we have
\[
\widetilde\bbh_r\(\widetilde{Z}^\pi_{\phi,L}\)\in\rH^1_f(E,\rV_{\pi,L})\otimes_{\dQ_p\otimes_\dQ\dQ_\pi}
\(\widetilde\cH_r^{[r]}\otimes_\dQ\dM\).
\]
It follows from the construction that, in view of \eqref{eq:holomorphic5}, the element
\[
\xi_{r*}\zeta_r^*\widetilde\bbh_r\(\widetilde{Z}^\pi_{\phi,L}\)\in\rH^1_f(E,\rV_{\pi,L})\otimes_{\dQ_p\otimes_\dQ\dQ_\pi}
\(\rH^0(\bsigma_{r}(K'K_r^{\prime\blacklozenge}),\xi_{r*}(\bomega_r^\delta)^{\otimes r})\otimes_\dQ\dM\)
\]
belongs to the subspace
\[
\rH^1_f(E,\rV_{\pi,L})\otimes_{\dQ_p\otimes_\dQ\dQ_\pi}\(\rH^0(\bsigma_{r}(K'K_r^{\blacklozenge}),\bomega_r^{\otimes r})\otimes_\dQ\dM\)
\]
(along the canonical subbundle $\bomega_r^{\otimes r}\subseteq\xi_{r*}(\bomega_r^\delta)^{\otimes r}$). Then we define $\cZ^\pi_{\phi,L}$ to be $\xi_{r*}\zeta_r^*\widetilde\bbh_r\(\widetilde{Z}^\pi_{\phi,L}\)$, which satisfies the requirement. The proposition is proved.
\end{proof}

\begin{definition}[Selmer theta lift]\label{de:theta}
Suppose that Hypothesis \ref{hy:modularity} holds for $\pi$. For every $\phi\in\sS(V^r\otimes_F\dA_F^\infty,\dM)^{K_r^{\blacklozenge}\times L}$ and every $\varphi\in\cV_{\hat\pi}$, we put
\begin{align*}
\Theta_\phi^\sel(\varphi)_L\coloneqq\langle\varphi^\dag,\cZ^\pi_{\phi,L}\rangle_\pi
\end{align*}
(see Notation \ref{no:period} for the pairing) as an element of $\rH^1_f(E,\rV_{\pi,L})\otimes_{\dQ_\pi}\dM$, called a \emph{Selmer theta lift} of $\pi$. It is clear from the construction that $\Theta_\phi^\sel(\varphi)_L$ is compatible under pullbacks with respect to $L$.
\end{definition}

At last, we state our theorem concerning Hypothesis \ref{hy:modularity}, whose proof will be given in \S\ref{ss:proof}.

\begin{theorem}\label{th:modularity}
Suppose that we are in the situation of Assumption \ref{st:main} and $n<p$. If the vanishing order of $\sL_p^\lozenge(\pi)$ at $\b1$ is one, then Hypothesis \ref{hy:modularity} holds for $\pi$.
\end{theorem}

\subsection{A $p$-adic arithmetic inner product formula}
\label{ss:aipf}

Recall from \cite{LL}*{Definition~3.8} that we have a canonical volume $\vol^\natural(L)\in\dQ_{>0}$, which in fact equals the product of the constant $W_{2r}$ in Lemma \ref{le:section6} and the volume of $L$ under the Siegel--Weil measure in \S\ref{ss:setup}(H9). If Hypothesis \ref{hy:modularity} holds for both $\pi$ and $\hat\pi$, then for every $\varphi_1\in\cV_{\hat\pi}$, every $\varphi_2\in\cV_\pi$ and every pair $\phi_1,\phi_2\in\sS(V^r\otimes_F\dA_F^\infty,\dM)^{K_r^{\blacklozenge}\times L}$, we have the height
\begin{align*}
\vol^\natural(L)\cdot
\langle\Theta_{\phi_1}^\sel(\varphi_1)_L,\Theta_{\phi_2}^\sel(\varphi_2)_L\rangle_{(\rV_{\pi,L},\rV_{\hat\pi,L}),E}
&\in\Gamma_{E,p}\otimes_{\dZ_p}\dL_\pi\otimes_{\dQ_\pi}\dM \\
&=\Gamma_{E,p}\otimes_{\dZ_p}\dL\otimes_\dQ\dM
\end{align*}
from \eqref{eq:heightv}, which is independent of $L$. Denote the above canonical value by
\[
\langle\Theta_{\phi_1}^\sel(\varphi_1),\Theta_{\phi_2}^\sel(\varphi_2)\rangle_{\pi,E}^\natural,
\]
and then put
\begin{align*}
\langle\Theta_{\phi_1}^\sel(\varphi_1),\Theta_{\phi_2}^\sel(\varphi_2)\rangle_{\pi,F}^\natural
&\coloneqq\Nm_{E/F}\langle\Theta_{\phi_1}^\sel(\varphi_1),\Theta_{\phi_2}^\sel(\varphi_2)\rangle_{\pi,E}^\natural \\
&\in\Gamma_{F,p}\otimes_{\dZ_p}\dL\otimes_\dQ\dM.
\end{align*}
Now we can state our \emph{$p$-adic arithmetic inner product formula}, whose proof will be given in \S\ref{ss:proof1}.

\begin{theorem}\label{th:aipf}
Suppose that we are in the situation of Assumption \ref{st:main} and $n<p$.
\begin{enumerate}
  \item If the vanishing order of $\sL_p^\lozenge(\pi)$ at $\b1$ is one (so that Hypothesis \ref{hy:modularity} holds for both $\pi$ and $\hat\pi$ by Theorem \ref{th:modularity} and Remark \ref{re:padic}(3)), then for every choice of elements
     \begin{itemize}
       \item $\varphi_1=\otimes_v\varphi_{1,v}\in\cV_{\hat\pi}$ and $\varphi_2=\otimes_v\varphi_{2,v}\in\cV_\pi$ both fixed by $K_r^\lozenge$ such that $\langle\varphi_{1,v},\varphi_{2,v}\rangle_{\pi_v}=1$ for every $v\in\tV_F\setminus\tV_F^{(\lozenge)}$,

       \item $\phi_1=\otimes_v\phi_{1,v},\phi_2=\otimes_v\phi_{2,v}\in\sS(V^r\otimes_F\dA_F^\infty,\dM)^{K_r^{\blacklozenge}}$ with $\phi_1^\lozenge=\phi_2^\lozenge$ good (\S\ref{ss:setup}(H5)),
     \end{itemize}
     the identity
     \begin{align}\label{eq:aipf}
     \langle\Theta_{\phi_1}^\sel(\varphi_1),\Theta_{\phi_2}^\sel(\varphi_2)\rangle_{\pi,F}^\natural&=
     \partial\sL_p^\lozenge(\pi)\cdot
     \prod_{v\in\tV_F^{(p)}}\prod_{u\in\tP_v}\gamma(\tfrac{1+r}{2},\ul{\pi_u},\psi_{F,v}) \\
     &\quad\times\prod_{v\in\tV_F^{(\lozenge\setminus\{\infty\})}}
     Z(\varphi_{1,v}^\dag\otimes\varphi_{2,v},f^{\r{SW}}_{\phi_{1,v}\otimes\phi_{2,v}}) \notag
     \end{align}
     holds in $\Gamma_{F,p}\otimes_{\dZ_p}\dL\otimes_\dQ\dM$, where
     \begin{itemize}
       \item $\gamma(\tfrac{1+r}{2},\ul{\pi_u},\psi_{F,v})$ is the unique element in $\dL^\times$ satisfying $\iota\gamma(\tfrac{1+r}{2},\ul{\pi_u},\psi_{F,v})=\gamma(\tfrac{1+r}{2},\iota\ul{\pi_u},\psi_{F,v})$ for every embedding $\iota\colon\dL\to\dC$;

       \item the term $Z(\varphi_{1,v}^\dag\otimes\varphi_{2,v},f^{\r{SW}}_{\phi_{1,v}\otimes\phi_{2,v}})\in\dL\otimes_\dQ\dM$ is from Lemma \ref{le:zeta1}.
     \end{itemize}

  \item If the vanishing order of $\sL_p^\lozenge(\pi)$ at $\b1$ is not one, then
      \[
      \Nm_{E/F}\left\langle\wp_{\pi}\(Z_{T_1}^\clubsuit(\phi_1)_L\),
      \wp_{\hat\pi}\(Z_{T_2}^\clubsuit(\phi_2)_L\)\right\rangle_{(\rV_{\pi,L},\rV_{\hat\pi,L}),E}
      =0
      \]
      for every $\phi_1,\phi_2\in\sS(V^r\otimes_F\dA_F^\infty,\dC)^L$ and $T_1,T_2\in\Herm_r(F)^+$.
\end{enumerate}
\end{theorem}

\begin{remark}\label{re:aipf}
We have the following remarks concerning Theorem \ref{th:aipf}.
\begin{enumerate}
  \item By the interpolation property of $\sL_p^\lozenge(\pi)$ and Lemma \ref{le:zeta2}, the right-hand side of \eqref{eq:aipf} does not change when enlarging $\lozenge$. In particular, we may enlarge $\lozenge$ to prove the theorem.

  \item Note that when we vary $\varphi_{1,v},\varphi_{2,v},\phi_{1,v},\phi_{2,v}$ for $v\in\tV_F^{(\lozenge\setminus\{\infty\})}$, both sides of \eqref{eq:aipf} define elements in the space
      \[
      \bigotimes_{v\in\tV_F^{(\lozenge\setminus\{\infty\})}}\Hom_{G_r(F_v)\times G_r(F_v)}\(\rI_{r,v}^\Box(\CF),\pi_v\boxtimes\hat\pi_v\),
      \]
      which is one-dimensional if $V$ is $\pi$-coherent \cite{LL2}*{Proposition~4.8(1)}\footnote{Here, we regard $G_r\times G_r$ as a subgroup of $G_{2r}$ through \eqref{eq:embedding} rather than $\imath$ from \S\ref{ss:degenerate}, which explains the change from $\pi_v^\vee$ to $\hat\pi_v$ (see Remark \ref{re:dagger}).} and vanishes if not. In particular, when $V$ is not $\pi$-coherent, all quantities in the theorem are trivially zero.

  \item In the situation of Assumption \ref{st:main}, we have $\epsilon(\pi_v)=-1$ (resp.\ $\epsilon(\pi_v)=1$) if $v\in\tS_\pi$ (resp.\ $v\in\tV_F^\fin\setminus\tS_\pi$), where $\epsilon(\pi_v)$ is introduced before Proposition \ref{pr:complete}. By \cite{Liu22}*{Theorem~1.2} and \cite{LL2}*{Proposition~3.9}, $V$ is $\pi$-coherent if and only if $\eta_{E/F}\((-1)^r\det V_v\)=\epsilon(\pi_v)$ for every $v\in\tV_F^\fin$. In particular, the theorem is trivial unless $r[F:\dQ]+|\tS_\pi|$ is odd by Remark \ref{re:coherent}(2).

  \item It is clear that Theorem \ref{th:aipf}(2) implies Theorem \ref{th:main1}(2).

  \item The role of the set $\blacklozenge$ is only to control the coefficient field $\dM$ (the smaller $\blacklozenge$ is, the smaller $\dM$ we can take). For the proof of the theorem, we may just take $\dM=\dC$ and ignore the choice of $\blacklozenge$.
\end{enumerate}
\end{remark}

\begin{proof}[Proof of Corollary \ref{th:main} assuming Theorem \ref{th:aipf}]
When the vanishing order of $\sL_p^\lozenge(\pi)$ at $\b1$ is one, we may choose $V$ that is $\pi$-coherent by Proposition \ref{pr:complete} and Remark \ref{re:aipf}(3). In particular, we may find data $\varphi_1,\varphi_2,\phi_1,\phi_2$ such that the right-hand side of \eqref{eq:aipf} is nonzero. Thus, $\Theta_{\phi_1}^\sel(\varphi_1)\neq0$, which implies that $\rH^1_f(E,\rV_{\pi,L})\neq 0$. Then the corollary follows from Hypothesis \ref{hy:galois}.
\end{proof}

The rest of this section is devoted to the proof of Theorem \ref{th:modularity} and Theorem \ref{th:aipf}. Once again, for the proof of these theorems, we may assume $\dM=\dC$ hence the choice of $\blacklozenge$ is irrelevant.

From now on, we will assume that we are in the situation of Assumption \ref{st:main}. In particular, $\dQ_\pi=\dQ$ (Remark \ref{re:weil}(3)). We may also that $V$ is $\pi$-coherent (Definition \ref{de:coherent}) hence the vanishing order of $\sL_p(\pi)$ at $\b1$ is at least one (Remark \ref{re:coherent}), since otherwise both theorems are trivial.

To shorten notation, we put
\begin{align*}
\tR&\coloneqq\tV_F^{(\lozenge\setminus\{p\})}\cap\tV_F^\spl\cap\tV_F^\heartsuit,\\
\tT&\coloneqq\tV_F^\fin\setminus(\tV_F^{(\lozenge)}\cap\tV_F^\spl\cap\tV_F^\heartsuit),
\end{align*}
so that $\tV_F^{(p)}\cup\tR\cup\tT$ is a partition of $\tV_F^\fin$. By enlarging $\lozenge$, we also assume the following assumption.

\begin{assumption}\label{eq:classes}
The set of primes of $E$ above $\tR$ is nonempty and generate the relative class group of $E/F$.
\end{assumption}

\subsection{Strategy for the modularity}
\label{ss:strategy}

We first reduce Hypothesis \ref{hy:modularity} to a problem about height pairing. The lemma below is the starting point.

\begin{lem}\label{le:hodge1}
If $\pi_v$ is unramified for every $v\in\tV_F^{(p)}$, then the image of the composite map
\[
\CH^r(X_L)\to\rH^{2r}_{\et}(X_L,\dQ_p(r))\xrightarrow{\wp_\pi}\rH^1(E,\rV_{\pi,L})
\]
is contained in $\rH^1_f(E,\rV_{\pi,L})$ for every neat open compact subgroup $L$ of $H(\dA_F^\infty)$.
\end{lem}

\begin{proof}
By Lemma \ref{le:hodge3} (which relies on Hypothesis \ref{hy:galois}), it suffices to show that the image of the above composite map is a crystalline class at every $p$-adic place of $E$. This then follows from Lemma \ref{le:hodge4} and \cite{Nek00}*{Theorem~3.1}.
\end{proof}

\begin{definition}\label{de:strong}
We say that an element $\varphi\in\dL\otimes_\dQ\cA_{r,\hol}^{[r]}$ is \emph{strongly nonzero} if for every embedding $\iota\colon\dL\to\dC$, the induced element $\iota\varphi\in\cA_{r,\hol}^{[r]}$ is nonzero.
\end{definition}

\begin{lem}\label{le:modularity}
Suppose that we can find
\begin{itemize}
  \item a neat open compact subgroup $L$ of $H(\dA_F^\infty)$,

  \item an element $\phi_1\in\sS(V^r\otimes_F\dA_F^\infty)^L$,

  \item an element $\zeta\in\rH^1_f(E,\rV_{\hat\pi,L})$,

  \item an element $\lambda\in\Hom_{\dZ_p}(\Gamma_{E,p},\dZ_p)$,

  \item a strongly nonzero element $\varphi_1\in\dL\otimes_\dQ\cA_{r,\hol}^{[r]}$,
\end{itemize}
such that
\begin{align*}
(1\times\bbq_r^\an)(g\cdot\varphi_1)=\sum_{T\in\Herm_r(F)^+}
\lambda\left\langle\wp_{\pi}\(Z_T^\clubsuit(\omega_r(g)\phi_1)_L\),
\zeta\right\rangle_{(\rV_{\pi,L},\rV_{\hat\pi,L}),E}\cdot q^T
\end{align*}
holds in $\SF_r(\dL\otimes_\dQ\dC)$ for every $g\in G_r(\dA_F^\infty)$ (in which the height pairing makes sense by Lemma \ref{le:hodge1}). Then Hypothesis \ref{hy:modularity} holds for every $\phi\in\sS(V^r\otimes_F\dA_F^\infty)$.
\end{lem}

\begin{proof}
First, we note that $\pi$ can actually be defined over a number field $\dE$ contained in $\dL$ and we will assume that $\pi$ has coefficients in $\dE$ in this proof. For every embedding $\varepsilon\colon\dE\to\dC$, put
\[
\rV_\pi^\varepsilon\coloneqq
\Hom_{H(\dA_F^\infty)}\(\theta(\varepsilon\pi),\varinjlim_{L}\rH^{2r-1}_\et(\ol{X}_L,\dL(r))\otimes_{\dE,\varepsilon}\dC\)
\]
as an $(\dL\otimes_{\dE,\varepsilon}\dC)[\Gal(\ol{E}/E)]$-module, where we recall $\theta(\varepsilon\pi)$ from \S\ref{ss:setup}(H10). Then for each individual $L$, $\rV_{\pi,L}\otimes_{\dE,\varepsilon}\dC=\theta(\varepsilon\pi)^L\otimes_\dC\rV_\pi^\varepsilon$, so that
\[
\rH^1_f(E,\rV_{\pi,L})\otimes_{\dE,\varepsilon}\dC=\theta(\varepsilon\pi)^L\otimes_\dC\rH^1_f(E,\rV_\pi^\varepsilon).
\]

For every neat open compact subgroup $L$ of $H(\dA_F^\infty)$, the assignment $\phi\mapsto Z^\pi_{\phi}(-)_L$
(Definition \ref{de:selmer}) defines a functional
\[
Z\in\Hom_{?}\(\sS(V^r\otimes_F\dA_F^\infty)^L,\rH^1_f(E,\rV_{\pi,L})\otimes_\dL\c{SF}_r(\dL\otimes_\dQ\dC)\)
\]
with $?=\dC[G_r(\dA_F^\infty)][L\backslash H(\dA_F^\infty)/L]$ (Definition \ref{de:sf1}). For the lemma, it suffices to show that for every embedding $\varepsilon\colon\dE\to\dC$, the functional $Z^\varepsilon$ in
\[
\Hom_{?}\(\sS(V^r\otimes_F\dA_F^\infty)^L,
\theta(\varepsilon\pi)^L\otimes_\dC\rH^1_f(E,\rV_\pi^\varepsilon)\otimes_{\dL\otimes_\dQ\dC}\c{SF}_r(\dL\otimes_\dQ\dC)\)
\]
factors through the subspace $\theta(\varepsilon\pi)^L\otimes_\dC\rH^1_f(E,\rV_\pi^\varepsilon)\otimes_\dC\(\bbq_r^\infty\cA_{r,\hol}^{[r]}\)$ (Definition \ref{de:sf1}) of the target. By Lemma \ref{le:theta}, there exists an irreducible $\dC[G_r(\dA_F^\infty)]$-submodule $\cM$ of $\rH^1_f(E,\rV_\pi^\varepsilon)\otimes_{\dL\otimes_\dQ\dC}\c{SF}_r(\dL\otimes_\dQ\dC)$ such that $Z^\varepsilon$ takes values in $\theta(\varepsilon\pi)^L\otimes_\dC\cM$. Thus, it suffices to show that $\cM$ and $\rH^1_f(E,\rV_\pi^\varepsilon)\otimes_\dC\(\bbq_r^\infty\cA_{r,\hol}^{[r]}\)$ have nonzero intersection, which is implied by the situation of the lemma.
\end{proof}

In practice, we are not able to study the height pairing
\[
\left\langle\wp_{\pi}\(Z_T^\clubsuit(\omega_r(g)\phi_1)_L\),
\zeta\right\rangle_{(\rV_{\pi,L},\rV_{\hat\pi,L}),E}
\]
for every $g\in G_r(\dA_F^\infty)$ for given $\phi_1$ and $\zeta$. However, the following lemma shows that it suffices to consider a much smaller set of $g$. Recall the subgroups $M_r\subseteq P_r\subseteq G_r$ from \S\ref{ss:run}(G3,G4).

\begin{lem}\label{le:modularity1}
Let $\ell\colon\CH^r(X_L)_\dQ\to\dL$ be a $\dQ$-linear map. For every $\phi\in\sS(V^r\otimes_F\dA_F^\infty)^L$ that is fixed under $\prod_{v\in\tV_F^\fin\setminus\tR}\(K_{r,v}\cap M_r(F_v)\)$ and every $\varphi\in\dL\otimes_\dQ\cA_{r,\hol}^{[r]}$, if
\begin{align*}
(1\times\bbq_r^\an)(g\cdot\varphi)=\sum_{T\in\Herm_r(F)^+}
\ell\(Z_T(\omega_r(g)\phi)_L\)\cdot q^T
\end{align*}
holds in $\SF_r(\dL\otimes_\dQ\dC)$ for every $g\in M_r(F_\tR)$, then it holds for all $g\in G_r(\dA_F^\infty)$.
\end{lem}

\begin{proof}
To prove the lemma, it suffices to show the identity for every embedding $\iota\colon\dL\to\dC$ by Remark \ref{re:modularity}(2). Thus, we may assume that the coefficients are in $\dC$ instead of $\dL\otimes_\dQ\dC$ and that $\ell\colon\CH^r(X_L)_\dC\to\dC$ is a complex linear functional.

We first turn the generating function into the automorphic setting. For every $v\in\tV_F^{(\infty)}$, let $V_v^\std\coloneqq(E_v)^{2r}$ be the standard positive definite hermitian space defined by the identity matrix $1_{2r}$, $\phi_v^\std$ the standard Gaussian function on $(V_v^\std)^r$, and $\omega_{r,v}$ the Weil representation of $G_r(F_v)$ generated by $\phi_v^\std$ in which every function factors through the moment map $T\colon(V_v^\std)^r\to\Herm_r(F_v)$. Put $\bomega_r\coloneqq\otimes_{v\in\tV_F}\omega_{v,r}$ and $\bphi\coloneqq(\otimes_{v\in\tV_F^{(\infty)}}\phi_v^\std)\otimes\phi$. For every $T\in\Herm_r(F)^+$ and $g\in G_r(\dA_F)$, put
\[
Z_T(\bomega_r(g)\bphi)_L\coloneqq\sum_{\substack{x\in L\backslash V^r\otimes_F\dA_F^\infty\\ T(x)=T}}(\bomega_r(g)\bphi)(T,x)Z(x)_L
\in\CH^r(X_L)_\dC.
\]
Denote $\sG$ the subset of $G_r(\dA_F^\infty)$ consisting of $g$ such that for every $g_\infty\in G_r(F_\infty)$, the sum
\[
\sum_{T\in\Herm_r(F)^+}\ell\(Z_T(\bomega_r(g_\infty g)\bphi)_L\)
\]
is absolutely convergent and equals $\varphi(g_\infty g)$. Thus, the lemma is equivalent to the following statement: If $M_r(F_\tR)\subseteq\sG$, then $\sG=G_r(\dA_F^\infty)$.

Choose an open compact subgroup $K$ of $G_r(\dA_F^\infty)$ that fixes $\phi$ and contains the subgroup
\[
\prod_{v\in\tV_F^\fin\setminus\tR}\(K_{r,v}\cap M_r(F_v)\).
\]
It is clear that $\sG$ is preserved under the right translation by $K$. On the other hand, Assumption \ref{eq:classes} implies that $M_r(F_\tR)$ maps surjectively to the double quotient $G_r(F)\backslash G_r(\dA_F^\infty)/K$. Thus, it suffices to show the following claim:
\begin{itemize}
  \item[$(*)$] If $g\in\sG$, then $\gamma g\in\sG$ for every $\gamma\in G_r(F)$.
\end{itemize}

The above claim is slightly stronger than the formal modularity property of Kudla's generating functions as proved in \cite{Liu11}*{Theorem~3.5} as we do not assume the absolute convergence \emph{a priori}. Nevertheless, it can be proved by essentially the same argument. First note that $G_r(F)$ is generated by $P_r(F)$ and the element
\[
w\coloneqq
\begin{pmatrix}
 1_{r-1} &  &  &  \\
     &  & & 1 \\
     &  &1_{r-1} &  \\
     & -1 & &  \\
\end{pmatrix}.
\]
Claim $(*)$ is obvious for $\gamma\in P_r(F)$. Thus, it remains to consider $\gamma=w$. Denote by $\partial\colon\Herm_r\to\Herm_{r-1}$ the map that sends $T$ to its upper-left block of size $r-1$. The proof of \cite{Liu11}*{Theorem~3.5(1)} indeed shows the following: If $g\in\sG$, then for every $g_\infty\in G_r(F_\infty)$, the sum
\[
\sum_{T'\in\Herm_{r-1}(F)^+}\(\sum_{\substack{T\in\Herm_r(F)^+ \\ \partial T=T'}}\ell\(Z_T(\bomega_r(w g_\infty g)\bphi)_L\)\)
\]
is absolutely convergent \emph{in order}, and equals $\varphi(w g_\infty g)$. It remains to show that the above sum is indeed absolutely convergent as a double sum. Since $\Herm_r(F)$ is dense in $\Herm_r(\dA_F^\infty)$, it follows that
\[
\varphi(n(b) w g_\infty g)=\sum_{T'\in\Herm_{r-1}(F)^+}\(\sum_{\substack{T\in\Herm_r(F)^+ \\ \partial T=T'}}\ell\(Z_T(\bomega_r(n(b)w g_\infty g)\bphi)_L\)\)
\]
for every $b\in\Herm_r(\dA_F)$, in which the right-hand side is again understood as a convergent sum in order. Then it is easy to see that for every $T\in\Herm_r^+(F)$,
\[
\int_{\Herm_r(F)\backslash\Herm_r(\dA_F)}\varphi(n(b) w g_\infty g)\psi_F^{-1}(\tr Tb)\rd b
=\ell\(Z_T(\bomega_r(w g_\infty g)\bphi)_L\).
\]
Thus, $\sum_{T\in\Herm_r(F)^+}\ell\(Z_T(\bomega_r(w g_\infty g)\bphi)_L\)$ is absolutely convergent and equals $\varphi(w g_\infty g)$. In other words, $wg\in\sG$. Claim $(*)$ hence the lemma are proved.
\end{proof}

The candidate $\zeta$ in Lemma \ref{le:modularity} will also be (the limit of) elements of the form $\wp_{\hat\pi}(Z_{T_2}^\clubsuit(\phi_2)_L)$ for some $T_2\in\Herm_r^\circ(F)^+$ and $\phi_2\in\sS(V^r\otimes_F\dA_F^\infty)^L$. Next, we construct some pairs of Schwartz functions in $\sS(V_v^r)$ for every $v\in\tV_F^\fin$ that will be candidates in Lemma \ref{le:modularity}.

\begin{notation}
For every $v\in\tV_F^{(p)}$, we denote by $\varepsilon_v\in\dN^{\tP_v}$ the element that takes value $1$ on $\tP_v\cap\tP_{\r{CM}}$ (\S\ref{ss:run}(F2)) and value $0$ on $\tP_v\setminus\tP_{\r{CM}}$. Put $\varepsilon\coloneqq(\varepsilon_v)_v\in\dN^\tP$.
\end{notation}

For $v\in\tR$, define
\begin{align}\label{eq:support}
\sR_v\coloneqq\left\{\left.(\phi_{v,1},\phi_{v,2})\in\sS(V_v^r,\dZ_{(p)})^2\right|
\supp(\phi_{v,1}\otimes\phi_{v,2})\subseteq(V_v^{2r})_\reg\right\}
\end{align}
(\S\ref{ss:setup}(H3)), which is stable under the action of $M_r(F_v)\times M_r(F_v)$.

\begin{enumerate}[label=(S\arabic*)]
  \item For $v\in\tR$, choose an arbitrary pair $(\phi_{v,1},\phi_{v,2})\in\sR_v$.

  \item For $v\in\tT$, choose a good lattice $\Lambda_v$ of $V_v$ and put $\phi_{v,1}=\phi_{v,2}\coloneqq\CF_{\Lambda_v^r}\in\sS(V_v^r,\dZ)$.

  \item For $v\in\tV_F^{(p)}$, choose a good lattice $\Lambda_v$ of $V_v$ and a polarization $\Lambda_v=\Lambda_{v,1}\oplus\fp_v^{-d_v}\Lambda_{v,2}$ of free $O_{E_v}$-modules, namely, $\Lambda_{v,1}$ and $\Lambda_{v,2}$ are free isotropic $O_{E_v}$-submodules of $\Lambda_v$ of rank $r$. For $e\in\dN^{\tP_v}$, define
      \begin{itemize}
        \item $\Lambda_{v,1}^{[e]}$ to be the subset of $\(\varpi_v^{-e-\varepsilon_v}\cdot\Lambda_{v,1}\oplus \varpi_v^{-e+\varepsilon_v^\tc}\cdot\Lambda_{v,2}\)^r$ consisting of $x$ such that $T(x)\in\Herm_r(O_{F_v})$ and $x\modulo\Lambda_{v,2}\otimes\dQ$ generates $\varpi_v^{-e-\varepsilon_v}\cdot\Lambda_{v,1}$;

        \item $\Lambda_{v,2}^{[e]}$ to be the subset of $\(\varpi_v^{-e}\cdot\Lambda_{v,1}\oplus\varpi_v^{-e}\cdot\Lambda_{v,2}\)^r$ consisting of $x$ such that $T(x)\in\Herm_r(O_{F_v})$ and $x\modulo\Lambda_{v,1}\otimes\dQ$ generates $\varpi_v^{-e}\cdot\Lambda_{v,2}$.
      \end{itemize}
      For $i=1,2$, let $\phi_{v,i}^{[e]}\in\sS(V_v^r,\dZ)$ be the characteristic function of $\Lambda_{v,i}^{[e]}$.
\end{enumerate}
For $i=1,2$ and $e\in\dN^\tP$, we put
\begin{align*}
\phi_i^{[e]}\coloneqq\(\bigotimes_{v\in\tV_F^{(p)}}\phi_{v,i}^{[e_v]}\)\otimes
\(\bigotimes_{v\in\tV_F^\fin\setminus\tV_F^{(p)}}\phi_{v,i}\)\in\sS(V^r\otimes_F\dA_F^\infty,\dZ_{(p)}).
\end{align*}

At last, we choose an open compact subgroup $L_v\subseteq H(F_v)$ for every $v\in\tV_F^\fin$.
\begin{itemize}
  \item For $v\in\tR$, we choose some $L_v$ that fixes $\phi_{v,i}$ for $i=1,2$.

  \item For $v\in\tT$, we define $L_v$ to be the stabilizer of $\Lambda_v$.

  \item For $v\in\tV_F^{(p)}$, we define $L_v$ to be the stabilizer of the lattice chain
      \[
      \Lambda_{v,1}\oplus\fp_v\Lambda_{v,2}\subseteq\Lambda_{v,1}\oplus\Lambda_{v,2}.
      \]
\end{itemize}
Put $L\coloneqq\prod_vL_v\subseteq H(\dA_F^\infty)$ so that $L^\lozenge$ is good. We may assume that $L$ is neat by shrinking $L_v$ for $v\in\tR(\neq\emptyset)$.

\begin{lem}\label{le:regular}
Take an element $v\in\tV_F^{(p)}$.
\begin{enumerate}
  \item For $i=1,2$ and $e,e'\in\dN^{\tP_v}$, we have $\rU_v^{e'}\phi_{v,i}^{[e]}=\phi_{v,i}^{[e+e']}$.

  \item For $i=1,2$ and $e\in\dN^{\tP_v}$, $\phi_{v,i}^{[e]}$ is fixed by $L_v$.

  \item For every $(e_1,e_2)\in\dN^{\tP_v}\times\dN^{\tP_v}$, the support of $\phi_{v,1}^{[e_1]}\otimes\phi_{v,2}^{[e_2]}$ is contained in $(V_v^{2r})_\reg$ (\S\ref{ss:setup}(H3)); and we have
      \begin{align}\label{eq:regular}
      f^{\r{SW}}_{\phi_{v,1}^{[e_1]}\otimes\phi_{v,2}^{[e_2]}}=b_{2r,v}(\CF)^{-1}\vol(L_v,\rd h_v)\cdot(\tf^{[e_1^\tc+\varepsilon_v^\tc+e_2]}_{\CF_v})^{\CF_v},
      \end{align}
      where $\vol(L_v,\rd h_v)$ denotes the volume of $L_v$ under the Siegel--Weil measure $\rd h_v$ in \S\ref{ss:setup}(H9).
\end{enumerate}
\end{lem}

\begin{proof}
For (1), by induction, it suffices to consider the case where $e'=1_u$ for some $u\in\tP_v$. We will prove the case where $i=1$ and leave the other similar case to the reader. By definition, we have
\begin{align*}
(\rU_v^{1_u}\phi_{v,1}^{[e]})(x)&=
\sum_{b\in\Herm_r(O_{F_{v}}/\varpi_{v})}(\omega_{r,v}(n(\varpi_v^{-d_v}b^\sharp)m(\varpi_v^{1_u}))\phi_{v,1}^{[e]})(x) \notag \\
&=(\omega_{r,v}(m(\varpi_v^{1_u}))\phi_{v,1}^{[e]})(x)\sum_{b\in\Herm_r(O_{F_{v}}/\varpi_{v})}\psi_{F,v}(\tr \varpi_v^{-d_v}b^\sharp T(x)) \notag \\
&=q_v^{-r^2}\phi_{v,1}^{[e]}(\varpi_v^{1_u}x)\sum_{b\in\Herm_r(O_{F_{v}}/\varpi_{v})}\psi_{F,v}(\tr \varpi_v^{-d_v}b^\sharp T(x)).
\end{align*}
Since
\begin{align*}
&\quad\sum_{b\in\Herm_r(O_{F_{v}}/\varpi_{v})}\psi_{F,v}(\varpi_v^{-d_v}\tr b^\sharp T(x)) \\
&=
\begin{dcases}
q_v^{r^2} & \text{if $T(x)\in\Herm_{2r}(O_{F_v})$,} \\
0 & \text{if $T(x)\in\varpi_v^{-1}\Herm_{2r}(O_{F_v})\setminus\Herm_{2r}(O_{F_v})$,}
\end{dcases}
\end{align*}
we have $(\rU_v^{1_u}\phi_{v,1}^{[e]})(x)=\phi_{v,1}^{[e+1_u]}(x)$.

For (2), by (1), it suffices to consider the case where $e=0$, for which the invariance under $L_v$ is obvious.

For (3), it is easy to see that the image of $\Lambda_{v,1}^{[e_1]}\times\Lambda_{v,2}^{[e_2]}$ under the moment map $T\colon V_v^{2r}\to\Herm_{2r}(F_v)$ is contained in the set $\fT_v^{[e_1^\tc+\varepsilon_v^\tc+e_2]}$ in Construction \ref{co:section}, which is contained in $\Herm_{2r}^\circ(F_v)$. For \eqref{eq:regular}, by (1) and Lemma \ref{le:section5}, it suffices to consider the case where $e_1=e_2=0$. In the definition of $\Lambda_{v,i}^{[0]}$, the condition that $T(x)\in\Herm_{2r}(O_{F_v})$ is automatic. Then it is a straightforward exercise in linear algebra that the image of $\Lambda_{v,1}^{[0]}\times\Lambda_{v,2}^{[0]}$ under the moment map $T$ is exactly $\fT_v^{[\varepsilon_v^\tc]}$; and that for every $x\in\Lambda_{v,1}^{[0]}\times\Lambda_{v,2}^{[0]}$, an element $h_v\in H(F_v)$ keeps $x$ in $\Lambda_{v,1}^{[0]}\times\Lambda_{v,2}^{[0]}$ if and only if $h_v\in L_v$. It follows from \S\ref{ss:setup}(H9) that
\[
W_{T^\Box}(f^{\r{SW}}_{\phi_{v,1}^{[0]}\otimes\phi_{v,2}^{[0]}})=b_{2r,v}(\CF)^{-1}\vol(L_v,\rd h_v)\cdot\CF_{\fT_v^{[\varepsilon_v^\tc]}}(T^\Box)
\]
for every $T^\Box\in\Herm_{2r}^\circ(F_v)$, which implies \eqref{eq:regular} (when $e_1=e_2=0$).

The lemma is proved.
\end{proof}

Recall the ideals $\dS^\lozenge_{\pi,L}$ and $\dS^\lozenge_{\hat\pi,L}$ of $\dS^\lozenge_\dL$ introduced in front of Lemma \ref{le:hodge}. For $(T_1,T_2)\in\Herm_r^\circ(F)^+\times\Herm_r^\circ(F)^+$, $(\rs_1,\rs_2)\in\dS^\lozenge_{\pi,L}\times\dS^\lozenge_{\hat\pi,L}$ and $(e_1,e_2)\in\dN^\tP\times\dN^\tP$, we have
$Z_{T_1}^\clubsuit(\rs_1\phi_1^{[e_1]})_L=\rs_1^*Z_{T_1}^\clubsuit(\phi_1^{[e_1]})_L$ and $Z_{T_2}^\clubsuit(\rs_2\phi_2^{[e_2]})_L=\rs_2^*Z_{T_2}^\clubsuit(\phi_2^{[e_2]})_L$ by \cite{LL}*{Lemma~4.4} (see Remark \ref{re:error1} below). In particular,
\begin{align*}
Z_{T_1}^\clubsuit(\rs_1\phi_1^{[e_1]})_L&=\wp_\pi\(\rs_1^*Z_{T_1}^\clubsuit(\phi_1^{[e_1]})\)
\in\rH^1_f(E,\rV_{\pi,L}),\\
Z_{T_2}^\clubsuit(\rs_2\phi_2^{[e_2]})_L&=\wp_{\hat\pi}\(\rs_2^*Z_{T_2}^\clubsuit(\phi_2^{[e_2]})_L\)
\in\rH^1_f(E,\rV_{\hat\pi,L})
\end{align*}
by Lemma \ref{le:hodge1}. By \cite{LL}*{Lemma~6.4} (in which we may take $\tR'$ as $\tR\cup\tV_F^{(p)}$ by Lemma \ref{le:regular}(3)), the algebraic cycles $Z_{T_1}(\rs_1\phi_1^{[e_1]})_L$ and $Z_{T_2}(\rs_2\phi_2^{[e_2]})_L$ do not intersect. Therefore, by the discussion in \S\ref{ss:decomposition}, we have a decomposition formula
\begin{align}\label{eq:decompose}
&\quad\langle Z_{T_1}^\clubsuit(\rs_1\phi_1^{[e_1]})_L,Z_{T_2}^\clubsuit(\rs_2\phi_2^{[e_2]})_L
\rangle_{(\rV_{\pi,L},\rV_{\hat\pi,L}),E} \\
&=\sum_{u\nmid\infty}\langle Z_{T_1}(\rs_1\phi_1^{[e_1]})_L,
Z_{T_2}(\rs_2\phi_2^{[e_2]})_L\rangle_{(\rV_{\pi,L},\rV_{\hat\pi,L}),E_u}\in\Gamma_{E,p}\otimes_{\dZ_p}\dL \notag
\end{align}
for our $p$-adic height pairing. In what follows, to shorten notation, we will suppress the part $(\rV_{\pi,L},\rV_{\hat\pi,L})$ in the subscript of height pairings.

\begin{notation}
For a finite place $u$ (resp.\ $v$) of $E$ (resp.\ $F$) not above $p$, we denote by $[u]$ (resp.\ $[v]$) the image of an arbitrary uniformizer at $u$ (resp.\ $v$) in $\Gamma_{E,p}$ (resp.\ $\Gamma_{F,p}$).
\end{notation}

\subsection{Local height away from $p$}
\label{ss:height1}

In this subsection and the next one, we study the local summands in \eqref{eq:decompose}.

\begin{lem}\label{le:section7}
For every $v\in\tT$ and every $T^\Box\in\Herm_{2r}^\circ(F_v)$, there exists a unique element $\sfW_{T^\Box,v}\in\dZ[X]$ such that
\[
\sfW_{T^\Box,v}(\chi_v(\varpi_v))=b_{2r,v}(\chi)\cdot W_{T^\Box}(f_{\chi_v})
\]
holds for every finite character $\chi\colon \Gamma_{F,p}\to\dC^\times$, where $f_{\chi_v}\in\rI^\Box_{r,v}(\chi_v)$ is the unique section that satisfies $f_{\chi_v}\res_{K_{2r,v}}=f^{\r{SW}}_{\CF_{\Lambda_v^{2r}}}\res_{K_{2r,v}}$ and $\varpi_v$ is an arbitrary uniformizer of $F_v$.
\end{lem}

\begin{proof}
When $v\in\tV_F^\ram$, this follows from \cite{LL2}*{Remark~2.18~\&~Lemma~2.19}. When $v\in\tS_\pi$, this follows from the discussion in \cite{LZ}*{Section~9}. The remaining cases have been settled in Lemma \ref{le:section4}(1) as in these cases $f_{\chi_v}=f_{\chi_v}^\sph$ (Notation \ref{no:section1}(2)).
\end{proof}

\begin{notation}\label{no:incoherent}
For every $T^\Box\in\Herm_{2r}^\circ(F)^+$, put
\[
\Diff(T^\Box,V)\coloneqq\{v\in\tV_F^\fin\res(V^{2r}_v)_{T^\Box}=\emptyset\},
\]
which is a finite subset of $\tV_F^\fin\setminus\tV_F^\spl$ of odd cardinality.\footnote{The fact that $\Diff(T^\Box,V)$ has odd cardinality follows from the Hasse principle for hermitian spaces, in view of the signature of $V$.} We define $\Herm_{2r}^\circ(F)^+_V$ to be the subset of $\Herm_{2r}^\circ(F)^+$ consisting of $T^\Box$ such that $\Diff(T^\Box,V)$ is a singleton, whose unique element we denote by $v_{T^\Box}$.
\end{notation}

\begin{proposition}\label{pr:height}
There exists a pair $(\rt_1,\rt_2)\in\dS^\lozenge_{\dL}\times\dS^\lozenge_{\dL}$ satisfying $\chi_{\hat\pi}^\lozenge(\rt_1)\chi_\pi^\lozenge(\rt_2)\neq 0$, such that for every $(T_1,T_2)\in\Herm_r^\circ(F)^+\times\Herm_r^\circ(F)^+$, every $(\rs_1,\rs_2)\in\dS^\lozenge_{\pi,L}\times\dS^\lozenge_{\hat\pi,L}$ and every $(e_1,e_2)\in\dN^\tP\times\dN^\tP$, we have
\begin{align*}
&\quad\Nm_{E/F}\(\vol^\natural(L)\sum_{u\nmid\infty p}
\langle Z_{T_1}(\rt_1\rs_1\phi_1^{[e_1]})_L,Z_{T_2}(\rt_2\rs_2\phi_2^{[e_2]})_L\rangle_{E_u}\) \\
&=W_{2r}\(\sum_{\substack{T^\Box\in\Herm_{2r}^\circ(F)^+_V \\ \partial_{r,r}T^\Box=(T_1,T_2)}}
\sfW'_{T^\Box,v_{T^\Box}}(1)\cdot I_{T^\Box}((\rt_1\rs_1\phi_1^{[e_1]}\otimes\rt_2\rs_2\phi_2^{[e_2]})^{v_{T^\Box}})\cdot[v_{T^\Box}]\) \\
&\quad+W_{2r}\sum_{v\in\tS_\pi}\frac{2}{q_v^{2r}-1}
\(\sum_{\substack{T^\Box\in\Herm_{2r}^\circ(F)^+ \\ \partial_{r,r}T^\Box=(T_1,T_2)}}
\sfW_{T^\Box,v}^\sph(1)\cdot I_{T^\Box}((\rt_1\rs_1\phi_1^{[e_1]}\otimes\rt_2\rs_2\phi_2^{[e_2]})^v)\)[v]
\end{align*}
in $\Gamma_{F,p}\otimes_{\dZ_p}\dL$, where $W_{2r}$ is the rational constant in Lemma \ref{le:section6}, $\sfW_{T^\Box,v}^\sph\in\dZ[X]$ is the polynomial in Lemma \ref{le:section4}(1), and $I_{T^\Box}$ is (the product of) the functional in \S\ref{ss:setup}(H9).
\end{proposition}

\begin{proof}
We first note that by Proposition \ref{th:index}, the local $p$-adic height at $u\nmid\infty p$ coincides with Beilinson's local index. To compute the local indices at different $u$, we have four cases:

Suppose that $u$ lies over $\tV_F^\spl$. By \cite{LL2}*{Proposition~4.20} in which we may take $\tR'$ to be $\tR\cup\tV_F^{(p)}$ which has cardinality at least $2$ (see Remark \ref{re:error} below), we can find a pair $(\rt^u_1,\rt^u_2)\in\dS^\lozenge_{\dL}\times\dS^\lozenge_{\dL}$ satisfying $\chi_{\hat\pi}^\lozenge(\rt^u_1)\chi_\pi^\lozenge(\rt^u_2)\neq 0$ such that
\[
\langle Z_{T_1}(\rt^u_1\rs_1\phi_1^{[e_1]})_L,Z_{T_2}(\rt^u_2\rs_2\phi_2^{[e_2]})_L\rangle_{E_u}=0.
\]
Moreover, we may take $\rt^u_1=\rt^u_2=1$ for all but finitely many $u$.

Suppose that $u$ lies over an element $v\in\tV_F^\unr\setminus\tS_\pi$. By \cite{LL}*{Proposition~8.1} and Remark \ref{re:beilinson}, we have
\begin{align*}
&\quad\vol^\natural(L)\cdot\langle Z_{T_1}(\rs_1\phi_1^{[e_1]})_L,Z_{T_2}(\rs_2\phi_2^{[e_2]})_L\rangle_{E_u} \\
&=-W_{2r}\(\sum_{\substack{T^\Box\in\Herm_{2r}^\circ(F)^+_V \\ \partial_{r,r}T^\Box=(T_1,T_2) \\ v_{T^\Box}=v}}
\frac{b_{2r,v}(\CF)}{\log q_v^2}\cdot W'_{T^\Box}(0,1_{4r},\CF_{\Lambda_v^{2r}})
\cdot I_{T^\Box}((\rs_1\phi_1^{[e_1]}\otimes\rs_2\phi_2^{[e_2]})^v)\)[u],
\end{align*}
where $W_{T^\Box}(s,1_{4r},\CF_{\Lambda_v^{2r}})$ denotes the usual Siegel--Whittaker function with complex variable $s$ (see \cite{LL}*{(3.3)} for example). In our case, the character $\chi_v$ plays the role as $|\;|_{F_v}^s$, which implies that
\begin{align}\label{eq:height3}
\sfW_{T^\Box,v}(q_v^{-s})=\prod_{i=1}^n L(s+i,\eta_{E/F,v}^{n-i})\cdot W_{T^\Box}(s,1_{4r},\CF_{\Lambda_v^{2r}}).
\end{align}
Together with the relation $\Nm_{E/F}[u]=2[v]$, we obtain
\begin{align}\label{eq:height}
&\quad\Nm_{E/F}\(\vol^\natural(L)\cdot\langle Z_{T_1}(\rs_1\phi_1^{[e_1]})_L,Z_{T_2}(\rs_2\phi_2^{[e_2]})_L\rangle_{E_u}\) \\
&=W_{2r}\(\sum_{\substack{T^\Box\in\Herm_{2r}^\circ(F)^+_V \\ \partial_{r,r}T^\Box=(T_1,T_2) \\ v_{T^\Box}=v}}
\sfW'_{T^\Box}(1)
\cdot I_{T^\Box}((\rs_1\phi_1^{[e_1]}\otimes\rs_2\phi_2^{[e_2]})^v)\)[v]. \notag
\end{align}

Suppose that $u$ lies over an element $v\in\tV_F^\ram$. By \cite{LL2}*{Proposition~4.28} and Remark \ref{re:beilinson}, we have
\begin{align*}
&\quad\vol^\natural(L)\cdot\langle Z_{T_1}(\rs_1\phi_1^{[e_1]})_L,Z_{T_2}(\rs_2\phi_2^{[e_2]})_L\rangle_{E_u} \\
&=-W_{2r}\(\sum_{\substack{T^\Box\in\Herm_{2r}^\circ(F)^+_V \\ \partial_{r,r}T^\Box=(T_1,T_2) \\ v_{T^\Box}=v}}
\frac{b_{2r,v}(\CF)}{\log q_v} W'_{T^\Box}(0,1_{4r},\CF_{\Lambda_v^{2r}})
\cdot I_{T^\Box}((\rs_1\phi_1^{[e_1]}\otimes\rs_2\phi_2^{[e_2]})^v)\)[u].
\end{align*}
Now we have \eqref{eq:height3} again but $\Nm_{E/F}[u]=[v]$, which imply \eqref{eq:height} as well.

Suppose that $u$ lies over an element $v\in\tS_\pi$. By \cite{LL}*{Proposition~9.1} (see Remark \ref{re:error} below) and Remark \ref{re:beilinson}, we can find a pair $(\rt^u_1,\rt^u_2)\in\dS^\lozenge_{\dL}\times\dS^\lozenge_{\dL}$ satisfying $\chi_{\hat\pi}^\lozenge(\rt^u_1)\chi_\pi^\lozenge(\rt^u_2)\neq 0$ such that
\begin{align*}
&\quad\vol^\natural(L)\cdot\langle Z_{T_1}(\rt^u_1\rs_1\phi_1^{[e_1]})_L,Z_{T_2}(\rt^u_2\rs_2\phi_2^{[e_2]})_L\rangle_{E_u} \\
&=-W_{2r}\(\sum_{\substack{T^\Box\in\Herm_{2r}^\circ(F)^+_V \\ \partial_{r,r}T^\Box=(T_1,T_2) \\ v_{T^\Box}=v}}
\frac{b_{2r,v}(\CF)}{\log q_v^2} W'_{T^\Box}(0,1_{4r},\CF_{\Lambda_v^{2r}})
\cdot I_{T^\Box}((\rt^u_1\rs_1\phi_1^{[e_1]}\otimes\rt^u_2\rs_2\phi_2^{[e_2]})^v)\)[u] \\
&\quad+W_{2r}\frac{1}{q_v^{2r}-1}
\(\sum_{\substack{T^\Box\in\Herm_{2r}^\circ(F)^+ \\ \partial_{r,r}T^\Box=(T_1,T_2)}}
\sfW_{T^\Box,v}^\sph(1)\cdot I_{T^\Box}((\rt^u_1\rs_1\phi_1^{[e_1]}\otimes\rt^u_2\rs_2\phi_2^{[e_2]})^v)\)[u].
\end{align*}
Now we have \eqref{eq:height3} and $\Nm_{E/F}[u]=2[v]$, which imply
\begin{align*}
&\quad\Nm_{E/F}\(\vol^\natural(L)\cdot\langle Z_{T_1}(\rt^u_1\rs_1\phi_1^{[e_1]})_L,Z_{T_2}(\rt^u_2\rs_2\phi_2^{[e_2]})_L\rangle_{E_u}\) \\
&=W_{2r}\(\sum_{\substack{T^\Box\in\Herm_{2r}^\circ(F)^+_V \\ \partial_{r,r}T^\Box=(T_1,T_2) \\ v_{T^\Box}=v}}
\sfW'_{T^\Box,v}(1) \cdot I_{T^\Box}((\rt^u_1\rs_1\phi_1^{[e_1]}\otimes\rt^u_2\rs_2\phi_2^{[e_2]})^v)\)[v] \\
&\quad+W_{2r}\frac{2}{q_v^{2r}-1}
\(\sum_{\substack{T^\Box\in\Herm_{2r}^\circ(F)^+ \\ \partial_{r,r}T^\Box=(T_1,T_2)}}
\sfW_{T^\Box,v}^\sph(1)\cdot I_{T^\Box}((\rt^u_1\rs_1\phi_1^{[e_1]}\otimes\rt^u_2\rs_2\phi_2^{[e_2]})^v)\)[v].
\end{align*}

Finally, for $i=1,2$, we take $\rt_i=\prod_u\rt^u_i$ to be the (finite) product of the above auxiliary Hecke operators. The proposition follows by taking the sum over all $u\nmid\infty p$, which is a finite sum.
\end{proof}

\subsection{Local height above $p$}
\label{ss:height2}

Take an element $u\in\tP$ with $v\in\tV_F^{(p)}$ its underlying place. For technical purposes, we fix an $E$-linear isomorphism $\ol{E_u}\xrightarrow\sim\dC$.

\begin{lem}\label{le:height}
Suppose that $n<p$. There exists a pair $(\rt_1,\rt_2)\in\dS^\lozenge_{\dL}\times\dS^\lozenge_{\dL}$ satisfying $\chi_{\hat\pi}^\lozenge(\rt_1)\chi_\pi^\lozenge(\rt_2)\neq 0$, such that for every $(T_1,T_2)\in\Herm_r^\circ(F)^+\times\Herm_r^\circ(F)^+$, every $(\rs_1,\rs_2)\in\dS^\lozenge_{\pi,L}\times\dS^\lozenge_{\hat\pi,L}$ and every $(e_1,e_2)\in\dN^\tP\times\dN^\tP$, we have
\[
\langle Z_{T_1}(\rt_1\rs_1\phi_1^{[e_1]})_L,Z_{T_2}(\rt_2\rs_2\phi_2^{[e_2]})_L\rangle_{E_u}
\in(O_{E_u}^\times)^\free\otimes_{\dZ_p}\dL.
\]
\end{lem}

\begin{proof}
In view of Remark \ref{re:height}, we would like to apply Theorem \ref{th:crystalline}, for which we need an integral model of $X_L$ over $O_{E_u}$. For this, we need an auxiliary Shimura variety that admits such a model via moduli interpretation. Choose a CM type $\bPhi$ of $E$ such that the $p$-adic places of $E$ induced by $\bPhi$ via the fixed isomorphism $\ol{E_u}\xrightarrow\sim\dC$ form a subset $\tP_{\bPhi}$ of $\tP$ of cardinality $[F:\dQ]$ that contains $u$. Then the reflex field $E_{\bPhi}\subseteq\dC$ of $\bPhi$ is contained in $E_u$. Recall that we have the $\dQ$-torus $T$ from \S\ref{ss:hermitian} and fix a neat open compact subgroup $K_T$ of $T(\dA^\infty)$ that is maximal at primes not in $\lozenge\setminus\{p\}$. We have the Shimura variety $Y_{K_T}$ of $T$ with respect to the CM type $\bPhi$ at level $K_T$, which is finite \'{e}tale over $\Spec E_{\bPhi}$. Put $X\coloneqq(X_L\otimes_EE_u)\otimes_{E_{\bPhi}}Y_{K_T}$, which is a finite \'{e}tale cover of $X_L\otimes_EE_u$ and hence a smooth projective scheme over $E_u$ of pure dimension $n-1$. The ring $\dS^\lozenge$ extends naturally to a ring of finite \'{e}tale correspondences (see \S\ref{ss:etale}) of $X$. For every $x\in V^r\otimes_F\dA_F^\infty$, we denote by $Z(x)$ the pullback of $Z(x)_L$ to $X$.

Now for the lemma, it suffices to find elements $(\rt_1,\rt_2)\in\dS^\lozenge_{\dL}\times\dS^\lozenge_{\dL}$ satisfying $\chi_{\hat\pi}^\lozenge(\rt_1)\chi_\pi^\lozenge(\rt_2)\neq 0$, such that for every $x_1,x_2\in V^r\otimes_F\dA_F^\infty$ satisfying
\begin{align}\label{eq:height1}
T(x_i)\in\Herm_r^\circ(F)^+,\quad x_{i,v}\in\bigcup_{e\in\dN^{\tP_v}}\Lambda_{v,i}^{[e]},\quad i=1,2,
\end{align}
we have $\langle \rt_1^*Z(x_1),\rt_2^*Z(x_2)\rangle_{X,E_u}\in O_{E_u}^\times\otimes_{\dZ_p}\dL$.

Put $K\coloneqq E_u$ with the residue field $\kappa$. The $K$-scheme $X$ admits an integral model $\cX$ over $O_K$ such that for every $S\in\Sch'_{/O_K}$, $\cX(S)$ is the set of equivalence classes (given by $p$-principal isogenies) of tuples $(A_0,\lambda_0,\eta_0;A,\lambda,\eta;G_{u^\tc}\to G_{A,u^\tc})$ where
\begin{itemize}
  \item $A_0$ is an abelian scheme over $S$ with an action of $O_E$ of signature type $\bPhi$, together with a compatible $p$-principal polarization $\lambda_0$ and a level structure $\eta_0$ away from $p$,

  \item $A$ is an abelian scheme over $S$ with an action of $O_E$ of signature type $n\bPhi-\inc+\inc^\tc$ ($\inc$ being the inclusion $E\hookrightarrow\dC$), together with a compatible $p$-principal polarization $\lambda$, so that $G_{A,u^\tc}\coloneqq A[(u^\tc)^\infty]$ is an $O_{F_v}$-divisible module of dimension $1$ and relative height $n$,

  \item $\eta$ is an $L^v$-level structure for the hermitian space $\Hom_{O_E}(A_0,A)\otimes_F\dA_F^{\infty v}$,

  \item $G_{u^\tc}\to G_{A,u^\tc}$ is an isogeny of $O_{F_v}$-divisible modules over $S$ whose kernel is contained in $G_{u^\tc}[\varpi_v]$ and has degree $q_v^r$.
\end{itemize}
The reader may consult \cite{LL}*{Section~7} for more details about the first three items, which are not so related to our argument below. By the same argument for \cite{TY07}*{Proposition~3.4}, we know that $\cX$ is a projective strictly semistable scheme over $O_K$ to which finite \'{e}tale correspondences in $\dS^\lozenge_{\dL}$ naturally extend. Moreover, if we put $\sfX\coloneqq\cX\otimes_{O_K}\kappa$ and let $\sfX_1$ (resp.\ $\sfX_2$) be the closed locus of $\sfX$ on which the kernel of $G_{u^\tc}\to G_{A,u^\tc}$ (resp.\ $G_{A,u^\tc}\to G_{u^\tc}/G_{u^\tc}[\varpi_v]$) is not \'{e}tale, then under the notation of \S\ref{ss:crystalline},
\[
\sfX^{(1)}=\sfX_1\coprod\sfX_2,\quad
\sfX^{(2)}=\sfX_1\bigcap\sfX_2,\quad
\sfX^{(3)}=\sfX^{(4)}=\cdots=\emptyset.
\]

We then would like to apply Theorem \ref{th:crystalline} with $\dT=\dS^\lozenge_{\dL}$, $\fm=\Ker\chi^\lozenge_{\hat\pi}$ and $\fm'=\Ker\chi^\lozenge_\pi$. To check \eqref{eq:vanishing}, we realize that both $\chi^\lozenge_{\hat\pi}$ and $\chi^\lozenge_\pi$ can be defined over a number field $\dE$ contained in $\dL$. Thus, by \cite{KM74}*{Theorem~2}, it suffices to show that
\begin{align}\label{eq:height2}
\bigoplus_{q\geq 0}\rH^q_\et(\sfX^{(2)}\otimes_\kappa\ol\kappa,\dE_\ell)_\fm=
\bigoplus_{q\geq 0}\rH^q_\et(\sfX^{(2)}\otimes_\kappa\ol\kappa,\dE_\ell)_{\fm'}=0
\end{align}
where $\ell$ is an arbitrary prime of $\dE$ not above $p$. Indeed, there is a finite flat morphism $\cX_1\to\cX$ to which finite \'{e}tale correspondences in $\dS^\lozenge_{\dL}$ naturally extend, in which $\cX_1$ is the integral model with a Drinfeld level-$1$ structure at $v$ as the one used in \cite{LL}*{Section~7}. Then \eqref{eq:height2} follows from claim (2) in the proof of \cite{LL}*{Lemma~7.3} with $m=j=1$.

Denote by $\cZ(x)$ the Zariski closure of $Z(x)$ in $\cX$. By Theorem \ref{th:crystalline} and Remark \ref{re:height}, it suffices to show the following two claims for the lemma.
\begin{enumerate}
  \item For every $x_1,x_2\in V^r\otimes_F\dA_F^\infty$ satisfying \eqref{eq:height1} and every $\rt_1,\rt_2\in\dS^\lozenge_{\dL}$, we have $\rt_1^*\cZ(x_1)\cap\rt_2^*\cZ(x_2)=\emptyset$.

  \item For every $x\in V^r\otimes_F\dA_F^\infty$ with $T(x)\in\Herm_r^\circ(F)^+$, the dimension of $\cZ(x)\cap\sfX^{(h)}$ is at most $r-h$ for $h=1,2$.
\end{enumerate}

Part (1) follows from (the same argument for) \cite{LL}*{Lemma~7.2}.

For (2), since $\cZ(x)$ remains the same if we scale $x$ by an element in $F^\times$, we may assume that $x_v\in(\Lambda_{v,1}\oplus\Lambda_{v,2})^r$ for every $v\in\tV_F^{(p)}$. Up to a Hecke translation away from $p$, which does not affect the conclusion of (2), we may also assume that $x\in V^r$. We have a moduli scheme $\cY(x)$ finite over $\cX$, such that for every object $S=(A_0,\lambda_0,\eta_0;A,\lambda,\eta;G_{u^\tc}\to G_{A,u^\tc})$ of $\cX$, the set $\cY(x)(S)$ consists of elements $y\in\Hom_{O_E}(A_0^{\oplus r},A)\otimes\dZ_{(p)}$ satisfying $T(y)=T(x)$ and $y^p\in \eta(L^p x)$. By \cite{LL}*{Lemma~5.4}, $\cZ(x)$ is contained in the image of $\cY(x)$ in $\cX$. Thus, it suffices to show that the dimension of $\sfY(x)^{(h)}$ is at most $r-h$ for $h=1,2$, where $\sfY(x)^{(h)}\coloneqq\cY(x)\times_\cX\sfX^{(h)}$.

Let $V_x\subseteq V$ be the hermitian subspace (of dimension $r$) that is the orthogonal complement of the subspace spanned by $x$. Put $H_x\coloneqq\rU(V_x)$ which is naturally a subgroup of $H$, and put $L_x\coloneqq L\cap H_x(\dA_F^\infty)$. We have a similar moduli scheme $\cX_x$ over $O_K$ for $V_x$ similar to the one for $V$ but with the hyperspecial level structure at $p$. More precisely, for every $S\in\Sch'_{/O_K}$, $\cX_x(S)$ is the set of equivalence classes (given by $p$-principal isogenies) of tuples $(A_0,\lambda_0,\eta_0;A_1,\lambda_1,\eta_1)$ where
\begin{itemize}
  \item $(A_0,\lambda_0,\eta_0)$ is like the one in the definition of $\cX$,

  \item $A_1$ is an abelian scheme over $S$ with an action of $O_E$ of signature type $r\bPhi-\inc+\inc^\tc$, together with a compatible $p$-principal polarization $\lambda_1$,

  \item $\eta_1$ is an $L_x^p$-level structure for the hermitian space $\Hom_{O_E}(A_0,A_1)\otimes_F\dA_F^{\infty p}$.
\end{itemize}
In particular, $\cX_x$ is a projective smooth scheme over $O_K$ of pure relative dimension $r-1$. Put $\sfX_x\coloneqq\cX_x\otimes_{O_K}\kappa$; and for $h\geq 1$, denote by $\sfX_x^{[h]}$ the closed locus of $\sfX_x$ where the height of the connected part of $G_{A_1,u^\tc}\coloneqq A_1[(u^\tc)^\infty]$ is at least $h$. It is known that $\sfX_x^{[h]}$ has pure dimension $r-h$. Claim (2) will follow if there is a finite morphism $f\colon\cY(x)\to\cX_x$ that sends $\sfY(x)^{(h)}$ into $\sfX_x^{[h]}$ for $h=1,2$, which we now construct.

Take a point $P=(A_0,\lambda_0,\eta_0;A,\lambda,\eta;G_{u^\tc}\to G_{A,u^\tc};y)$ of $\cY(x)(S)$. Put $A'\coloneqq(A^\vee/(\lambda\circ y)A_0^{\oplus r})^\vee$, which inherits an action of $O_E$ which has signature type $r\bPhi-\inc+\inc^\tc$ and admits a natural map to $A$. Since $T(x)\in\Herm_r^\circ(F)$, the induced map $\lambda'\colon A'\to A\xrightarrow{\lambda} A^\vee\to A'^\vee$ is a quasi-polarization such that $\lambda'[p^\infty]$ is an isogeny. For every $\tilde{u}\in\tP$, we have the induced isogeny $\lambda'_{\tilde{u}^\tc}\colon G_{A',\tilde{u}^\tc}\to G_{A',\tilde{u}}^\vee$. Put
\[
A_1\coloneqq A'\left/\bigoplus_{\tilde{u}\in\tP_{\bPhi}}\Ker\lambda'_{\tilde{u}^\tc}\right.
\]
and let $\lambda_1\colon A_1\to A_1^\vee$ be the quasi-polarization induced from $\lambda'$, which is in fact $p$-principal from the construction. We can also define a natural $L_x^p$-level structure $\eta_1$ for $A_1$ whose details we leave to the reader. Then we define $f(P)$ to be $(A_0,\lambda_0,\eta_0;A_1,\lambda_1,\eta_1)$. Since the $O_{F_v}$-divisible module $G_{A_0,u^\tc}$ is \'{e}tale, the height of the connected part of $G_{A_1,u^\tc}$ equals to that of $G_{A,u^\tc}$. In particular, $f$ sends $\sfY(x)^{(h)}$ into $\sfX_x^{[h]}$ for $h=1,2$. It remains to show that $f$ is finite. Since $\cY(x)$ is proper over $O_K$, it suffices to show that the fiber of $f$ over an arbitrary $\ol\kappa$-point is finite. Indeed, when $S=\Spec\ol\kappa$, $G_{A',\tilde{u}^\tc}$ has dimension $1$ (resp. is \'{e}tale) if $\tilde{u}=u$ (resp.\ $\tilde{u}\in\tP_{\bPhi}\setminus\{u\}$), and the degree of $\lambda'_{\tilde{u}^\tc}$ is bounded by the moment matrix $T(x)$. It follows that up to isomorphism, there are only finitely many such isogenies $G_{A',\tilde{u}^\tc}\to G_{A_1,\tilde{u}}$ with fixed $G_{A_1,\tilde{u}}$ for every $\tilde{u}\in\tP_{\bPhi}$. Thus, $f$ is finite and claim (2) is confirmed.

The lemma is finally proved.
\end{proof}

\begin{proposition}\label{pr:height1}
Suppose that $n<p$. There exist an integer $M\geq 0$ and a pair $(\rt_1,\rt_2)\in\dS^\lozenge_{O_\dL}\times\dS^\lozenge_{O_\dL}$ satisfying $\chi_{\hat\pi}^\lozenge(\rt_1)\chi_\pi^\lozenge(\rt_2)\neq 0$, such that for every $(T_1,T_2)\in\Herm_r^\circ(F)^+\times\Herm_r^\circ(F)^+$, every $(\rs_1,\rs_2)\in(\dS^\lozenge_{\pi,L}\cap\dS^\lozenge_{O_\dL})\times(\dS^\lozenge_{\hat\pi,L}\cap\dS^\lozenge_{O_\dL})$ and every $(e_1,e_2)\in\dN^\tP\times\dN^\tP$, we have
\[
\langle Z_{T_1}(\rt_1\rs_1\phi_1^{[e_1]})_L,Z_{T_2}(\rt_2\rs_2\phi_2^{[e_2]})_L\rangle_{E_u}\in (O_{E_u}^\times)^\free\otimes_{\dZ_p}(p^{|e_{2,v}|-M}O_\dL).
\]
\end{proposition}

The rest of this subsection is devoted to the proof of this proposition. We may assume that $\rV_{\pi,L}\neq 0$ and $\rV_{\hat\pi,L}\neq 0$ since otherwise the proposition is trivial.

Let $S$ be the kernel of the norm map $\Nm_{E/F}\colon\Res_{O_E/O_F}\bG\to\bG_{O_F}$. Consider the reciprocity map
\begin{align}\label{eq:reciprocity}
\rec\colon\Aut(\dC/E)\to E^\times\backslash\dA_E^{\infty,\times}\to S(F)\backslash S(\dA_F^\infty)
\end{align}
in which the first one is from the global class field theory and the second (surjective) one sends $a$ to $a^\tc/a$. For $d\in\dN$, we
\begin{itemize}
  \item put $L_{S,v}^{[d]}\coloneqq S(O_{F_v})\cap(1+\fp_v^d O_{E_v})$,

  \item let $E^{[d]}\subseteq\dC$ be the finite abelian extension of $E$ such that the map $\rec$ \eqref{eq:reciprocity} induces an isomorphism
     \[
     \Gal(E^{[d]}/E)\simeq S(F)\backslash S(\dA_F^\infty)\left/L_{S,v}^{[d]}\prod_{\tilde{v}\in\tV_F^\fin\setminus\{v\}}S(O_{F_{\tilde{v}}})\right.,
     \]

  \item denote by $\rZ^r_{[d]}(X_L)$ the image of the norm map
     \[
     \Nm_{E^{[d]}/E}\colon\rZ^r(X_L\otimes_EE^{[d]})\to\rZ^r(X_L).
     \]
\end{itemize}

\begin{lem}\label{le:height1}
For every $x\in V^r\otimes_F\dA_F^\infty$ satisfying $T(x)\in\Herm_r^\circ(F)^+$ and $x_v\in\Lambda_{v,2}^{[e]}$ for some $e\in\dN^{\tP_v}$, we have
\[
Z(x)_L\in\rZ^r_{[|e|]}(X_L).
\]
\end{lem}

\begin{proof}
Up to a Hecke translation away from $p$, which does not affect the conclusion of the lemma, we may assume that $x\in V^r$. Let $V_x\subseteq V$ be the hermitian subspace (of dimension $r$) that is the orthogonal complement of the subspace spanned by $x$. Put $H_x\coloneqq\rU(V_x)$ which is naturally a subgroup of $H$, and put $L_x\coloneqq L\cap H_x(\dA_F^\infty)$. We have the Shimura variety $X_{x,L_x}$ for $H_x$ with level $L_x$, similar to $X_L$. By definition, $Z(x)_L$ is the fundamental cycle of the finite unramified morphism $X_{x,L_x}\to X_L$ defined over $E$.

We have the determinant map $\dtm\colon H_x\to S\otimes_{O_F}F$ which identifies $S\otimes_{O_F}F$ with the maximal abelian quotient of $H_x$. Then the set of connected components of $X_{x,L_x}\otimes_E\dC$ is canonically parameterized by the set $S(F)\backslash S(\dA_F^\infty)/\dtm L_x$. For every $s\in S(F)\backslash S(\dA_F^\infty)/\dtm L_x$, we denote by $X_{x,L_x}^s$ the corresponding connected component. The definition of canonical models of Shimura varieties implies that $\gamma X_{x,L_x}^s=X_{x,L_x}^{\rec(\gamma)s}$ for every $\gamma\in\Aut(\dC/E)$, where $\rec$ is the map \eqref{eq:reciprocity}.

We claim that $\dtm L_{x,v}\subseteq L_{S,v}^{[|e|]}$. Then we have the quotient map
\[
S(F)\backslash S(\dA_F^\infty)/\dtm L_x \to S(F)\backslash S(\dA_F^\infty)\left/L_{S,v}^{[|e|]}\prod_{\tilde{v}\in\tV_F^\fin\setminus\{v\}}S(O_{F_{\tilde{v}}})\right.
\]
Let $\fS$ be the fiber of $1$ in the above map. Then $\sum_{s\in\fS}X_{x,L_x}^s$ is defined over $E^{[|e|]}$; and $\Nm_{E^{[|e|]}/E}\sum_{s\in\fS}X_{x,L_x}^s=X_{x,L_x}$. The lemma then follows.

It remains to show the claim, which is an exercise in linear algebra. We assume $e\neq 0$ as the case for $e=0$ is trivial. By definition, $L_{x,v}$ is simply the subgroup of $L_v$ that fixes $x_v$, or equivalently, $x'_v\coloneqq \varpi_v^e\cdot x_v$. By the definition of $\Lambda_{v,2}^{[e]}$ in \S\ref{ss:strategy}(S3), $x'_v$ belongs to $(\Lambda_{v,1}\oplus\Lambda_{v,2})^r$ such that $T(x'_v)\in\varpi_v^{|e|}\Herm_r(O_{F_v})$ and that $x'_v\modulo\Lambda_{v,1}$ generates $\Lambda_{v,2}$. It follows that the image of $x'_v$ in $(\Lambda_{v,1}\oplus\Lambda_{v,2})^r\otimes_{O_{F_v}}O_{F_v}/\fp_v^{|e|}$ generates a Lagrangian $O_{E_v}\otimes_{O_{F_v}}O_{F_v}/\fp_v^{|e|}$-submodule of $(\Lambda_{v,1}\oplus\Lambda_{v,2})\otimes_{O_{F_v}}O_{F_v}/\fp_v^{|e|}$. In particular, every element in $L_{x,v}$, which stabilizes $\Lambda_{v,1}\oplus\Lambda_{v,2}$, has determinant $1$ modulo $\fp_v^{|e|}$. The claim follows.
\end{proof}

For $d\in\dN$, let $u_d$ be the place of $E^{[d]}$ induced from the fixed isomorphism $\ol{E_u}\xrightarrow\sim\dC$, which is above $u$. Put $K\coloneqq E_u$, $K_d\coloneqq(E^{[d]})_{u_d}$ for $d\in\dN$ and $K_\infty\coloneqq\bigcup_{d\geq 0}K_d$. Then $K_0/K$ is unramified and $K_d/K_0$ is totally ramified of degree $(q_v-1)q_v^{d-1}/|U_E|$ for $d>0$, where $U_E$ is the torsion subgroup of $O_E^\times$.

Recall that $\rV_{\pi,L}$ and $\rV_{\hat\pi,L}$ are subspaces of $\rH^{2r-1}_{\et}(\ol{X}_L,\dL(r))$. Put
\begin{align*}
\rT_{\pi,L}&\coloneqq\rV_{\pi,L}\cap\rH^{2r-1}_{\et}(\ol{X}_L,O_\dL(r))^\free,\\
\rT_{\hat\pi,L}&\coloneqq\rV_{\hat\pi,L}\cap\rH^{2r-1}_{\et}(\ol{X}_L,O_\dL(r))^\free,
\end{align*}
both being $O_\dL[\Gal(\ol{E}/E)]$-modules. For $d\in\dN$, we put
\[
\rN_\infty\rH^1_f(K_d,\rT_{\pi,L})\coloneqq
\bigcap_{d'\geq d}\IM\(\Cor_{K_{d'}/K_d}\colon\rH^1_f(K_{d'},\rT_{\pi,L})\to\rH^1_f(K_d,\rT_{\pi,L})\),
\]
in which $\Cor_{K_{d'}/K_d}$ denotes the corresponding corestriction map.

\begin{lem}\label{le:height2}
There exists an integer $M\geq 0$ such that $p^M$ annihilates
\[
\rH^1_f(K_d,\rT_{\pi,L})/\rN_\infty\rH^1_f(K_d,\rT_{\pi,L})
\]
for every $d\in\dN$.
\end{lem}

\begin{proof}
By Lemma \ref{le:hodge2}, we know that $\rV_{\pi,L}\res_{K_d}$ satisfies the Panchishkin condition (Definition \ref{de:panchishkin}) for every $d\in\dN$. By Lemma \ref{le:hodge4}, we may apply \cite{Nek93}*{Theorem~6.9} to $\rV_{\pi,L}\res_{K_d}$.\footnote{Though our extension $K_\infty/K_d$ is in general not a $\dZ_p$-extension as assumed in \cite{Nek93}*{6.2}, the argument for \cite{Nek93}*{Theorem~6.9} works without change.} Thus, by the same argument at the end of the proof of \cite{Nek95}*{Proposition~II.5.10}, it suffices to show that $\rH^0(K_\infty,\rV_{\pi,L})=\rH^0(K_\infty,\rV_{\hat\pi,L})=0$.

We follow the strategy in \cite{Shn16}*{Section~8}. We may choose an element $\xi\in S(F)$ such that $\xi=(\varpi_v^{1_u-1_{u^\tc}})^{[K_0:K]}$ in $\Gal(E^{[0]}/E)$. Then by the same argument for \cite{Shn16}*{Proposition~8.3}, $K_\infty$ is contained in $K_\xi$ -- the field attached to the Lubin--Tate group relative to the extension $K_0/K$ with parameter $\xi$. Let $\chi_\xi\colon\Gal(K_\xi/K_0)\to O_K^\times$ be the character given by the Galois action on the torsion points of this relative Lubin--Tate group; and let $K(\chi_\xi)$ be the corresponding one-dimensional representation. Let $L$ be the maximal subfield of $K$ that is unramified over $\dQ_p$. By the same argument for \cite{Shn16}*{Proposition~8.4}, we see that $K(\chi_\xi)$ is crystalline, and that the $q_v$-Frobenius map (which is $L$-linear) acts on $\dD_\cris(K(\chi_\xi))$, which is a free $K\otimes_{\dQ_p}L$-module of rank $1$, by multiplication by $\xi^{-1}$.\footnote{Note that in this article, we always use the covariant version for $\dD_{\dr}$ and $\dD_\cris$.} Note that $\dL$ is a subfield of $\dC$ and hence $\ol{K}$ via the fixed isomorphism $\ol{K}=\ol{E_u}\xrightarrow\sim\dC$. Let $\rV$ be either $\rV_{\pi,L}\otimes_\dL\ol{K}$ or $\rV_{\hat\pi,L}\otimes_\dL\ol{K}$. Repeating the argument in \cite{Shn16}*{Proposition~8.9} (which followed an approach in \cite{Nek95}) to $\rV$, we obtain $\rH^0(K_\xi,\rV)=0$ since $\rV$ is crystalline of pure weight $-1$ by Lemma \ref{le:hodge2}.

The lemma is proved.
\end{proof}

\begin{proof}[Proof of Proposition \ref{pr:height1}]
Let $M\geq 0$ be the integer in Lemma \ref{le:height2} and $(\rt_1,\rt_2)\in\dS^\lozenge_{O_\dL}\times\dS^\lozenge_{O_\dL}$ the pair in Lemma \ref{le:height}. We first note that $Z_{T_i}(\rt_i\rs_i\phi_i^{[e_i]})_L\in\rZ^r(X_L)\otimes O_\dL$ for $i=1,2$. By Lemma \ref{le:height1}, we may find an element $Z\in\rZ^r(X_L\otimes_EE^{[|e_{2,v}|]})\otimes O_\dL$ such that $\Nm_{E^{[|e_{2,v}|]}/E}Z=Z_{T_2}(\rt_2\phi_2^{[e_2]})_L$. We may also assume that the support of $Z$ is contained in the support of $Z_{T_2}(\rt_2\phi_2^{[e_2]})_L$. Put
\[
Z_2\coloneqq\Nm_{E^{[|e_{2,v}|]}\otimes_EK/K_{|e_{2,v}|}}Z\otimes_EK\in\rZ^r(X_L\otimes_EK_{|e_{2,v}|})\otimes O_\dL,
\]
so that $\Nm_{K_{|e_{2,v}|}/K}Z_2=Z_{T_2}(\rt_2\phi_2^{[e_2]})_L\otimes_EK$. Since the natural map
\[
\rH^{2r-1}_{\et}(\ol{X}_L,O_\dL(r))^\free/\rT_{\pi,L}\to\rH^{2r-1}_{\et}(\ol{X}_L,\dL(r))/\rV_{\pi,L}
\]
is injective, the class $Z_{T_1}^\et(\rt_1\rs_1\phi_1^{[e_1]})_L=\rs_1^*Z^\et_{T_1}(\rs_1\phi_1^{[e_1]})_L$ sits in $\rH^1_f(K,\rT_{\pi,L})$. Similarly, the cycle class of $\rs_2^*Z_2$ sits in $\rH^1_f(K_{|e_{2,v}|},\rT_{\hat\pi,L})$. By \cite{Nek95}*{II.(1.9.1)}, we have
\begin{align*}
&\quad\langle Z_{T_1}(\rt_1\rs_1\phi_1^{[e_1]})_L,Z_{T_2}(\rt_2\rs_2\phi_2^{[e_2]})_L\rangle_K \\
&=\Nm_{K_{|e_{2,v}|}/K}\langle \rs_1^*Z_{T_1}(\rt_1\phi_1^{[e_1]})_L\otimes_EK_{|e_{2,v}|},\rs_2^*Z_2\rangle_{K_{|e_{2,v}|}}.
\end{align*}
By Lemma \ref{le:height}, we have
\[
\langle\rs_1^*Z_{T_1}(\rt_1\phi_1^{[e_1]})_L\otimes_EK_{|e_{2,v}|},\rs_2^*Z_2\rangle_{K_{|e_{2,v}|}}
\in(O_{K_{|e_{2,v}|}}^\times)^\free\otimes_{\dZ_p}\dL.
\]
In other words, the corresponding bi-extension is crystalline (Remark \ref{re:height}). By the argument for \cite{Nek95}*{Proposition~II.1.11}, we have
\[
\langle\rs_1^*Z_{T_1}(\rt_1\phi_1^{[e_1]})_L\otimes_EK_{|e_{2,v}|},\rs_2^*Z_2\rangle_{K_{|e_{2,v}|}}\in (O_{K_{|e_{2,v}|}}^\times)^\free\otimes_{\dZ_p}(p^{-M}O_\dL).
\]
Finally, since the image of the norm map $\Nm_{K_d/K}\colon(O_{K_d}^\times)^\free\to (O_K^\times)^\free$ is precisely $p^d (O_K^\times)^\free$ for $d\in\dN$, we have
\[
\langle Z_{T_1}(\rt_1\rs_1\phi_1^{[e_1]})_L,Z_{T_2}(\rt_2\rs_2\phi_2^{[e_2]})_L\rangle_K
\in(O_K^\times)^\free\otimes_{\dZ_p}(p^{|e_{2,v}|-M}O_\dL).
\]
The proposition is proved.
\end{proof}

\subsection{Proof of Theorem \ref{th:modularity}}
\label{ss:proof}

In this subsection, we prove Theorem \ref{th:modularity}. As we have already pointed out, the key strategy is to use Lemma \ref{le:modularity}, which, in some sense, transfers the modularity problem to the scalar case by using $p$-adic height pairing as a projector. To apply this lemma, the main work is to find the testing element $\zeta$. It turns out that we can take $\zeta$ to be special cycles themselves, or more precisely, their $p$-adic limits. The reason we need to consider limits is that for the computation of local $p$-adic heights above $p$ in Proposition \ref{pr:height1}, we only know their asymptotic behavior.

Take $e\in\dN$, which is regarded as a constant tuple according to the context.

For every $v\in\tR$, we
\begin{itemize}
  \item choose a pair $(\phi_{v,1},\phi_{v,2})\in\sR_v$ \eqref{eq:support} and put $\Phi_v\coloneqq\phi_{v,1}\otimes\phi_{v,2}$,

  \item let $\tf_v\in\sS(\Herm_{2r}(F_v),\dZ_{(p)})$ be the unique element such that $\tf_v^{\CF_v}=f^{\r{SW}}_{\Phi_v}$,

  \item put $f_{\chi_v}\coloneqq\tf_v^{\chi_v}\in\rI^\Box_{r,v}(\chi_v)$ \eqref{eq:section} for every finite character $\chi\colon\Gamma_{F,p}\to\dC^\times$.
\end{itemize}

For every $v\in\tT$, we
\begin{itemize}
  \item put $\Phi_v\coloneqq\phi_{v,1}\otimes\phi_{v,2}=\CF_{\Lambda_v^{2r}}$ (\S\ref{ss:strategy}(S2)),

  \item let $f_{\chi_v}\in\rI^\Box_{r,v}(\chi_v)$ be the section from Lemma \ref{le:section7} for every finite character $\chi\colon\Gamma_{F,p}\to\dC^\times$, that is, the standard section such that $f_{\b1_v}=f_{\Phi_v}^{\r{SW}}$.
\end{itemize}

For every $v\in\tV_F^{(p)}$, we
\begin{itemize}
  \item put $\Phi_v^{[e]}\coloneqq\phi_{v,1}^{[0]}\otimes\phi_{v,2}^{[e]}$ (\S\ref{ss:strategy}(S3)),

  \item put $f_{\chi_v}^{[e]}\coloneqq b_{2r,v}(\CF)^{-1}\cdot\vol(L_v,\rd h_v)\cdot(\tf_{\chi_v}^{[e+\varepsilon^\tc]})^{\chi_v}\in\rI^\Box_{r,v}(\chi_v)$ for every finite character $\chi\colon\Gamma_{F,p}\to\dC^\times$, so that $f_{\CF_v}^{[e]}=f^{\r{SW}}_{\Phi_v^{[e]}}$ by Lemma \ref{le:regular}(3).
\end{itemize}

Then we
\begin{itemize}
  \item put
     \begin{align}\label{eq:phi}
     \phi_1&\coloneqq\(\bigotimes_{v\in\tV_F^{(p)}}\phi_{v,1}^{[0]}\)\otimes
     \(\bigotimes_{v\in\tV_F^\fin\setminus\tV_F^{(p)}}\phi_{v,1}\),\\
     \phi_2^{[e]}&\coloneqq\(\bigotimes_{v\in\tV_F^{(p)}}\phi_{v,2}^{[e]}\)\otimes
     \(\bigotimes_{v\in\tV_F^\fin\setminus\tV_F^{(p)}}\phi_{v,2}\) \notag
     \end{align}
     in $\sS(V^r\otimes_F\dA_F^\infty,\dZ_{(p)})^{K_r^\lozenge\times L}$,

  \item put
     \[
     \Phi^{[e]}\coloneqq\(\bigotimes_{v\in\tV_F^{(p)}}\Phi_v^{[e]}\)\otimes
     \(\bigotimes_{v\in\tV_F^\fin\setminus\tV_F^{(p)}}\Phi_v\)
     \in\sS(V^r\otimes_F\dA_F^\infty,\dZ_{(p)})^{K_{2r}^\lozenge\times L},
     \]

  \item put
     \[
     f_{\chi^\infty}^{[e]}\coloneqq\(\bigotimes_{v\in\tV_F^{(p)}}f_{\chi_v}^{[e]}\)\otimes
     \(\bigotimes_{v\in\tV_F^\fin\setminus\tV_F^{(p)}}f_{\chi_v}\)\in\rI^\Box_r(\chi)^\infty
     \]
     for every finite character $\chi\colon\Gamma_{F,p}\to\dC^\times$.
\end{itemize}

Finally, we
\begin{itemize}
  \item let $M\in\dN$ be the smallest element such that all of the following elements
     \[
     p^M\vol^\natural(L), \quad \frac{p^M W_{2r}^{\lozenge\setminus\{p\}}}{q_{v'}^{2r}-1}\prod_{v\in\tV_F^{(p)}}\vol(L_v,\rd h_v),\forall v'\in\tS_\pi
     \]
     belong to $\dZ_{(p)}$, where $W_{2r}^{\lozenge\setminus\{p\}}$ is defined in \eqref{eq:eisenstein1};

  \item fix an open compact subgroup $K^\dag$ of $G_{r,r}(\dA_F^\infty)=G_r(\dA_F^\infty)\times G_r(\dA_F^\infty)$ of the form
     \[
     K^\dag_{\lozenge\setminus\{\infty,p\}}\times
     \(\prod_{v\in\tV_F^{(p)}}\cG_{r,r}(O_{F_v})\times_{\cG_{r,r}(O_{F_v}/\varpi_v)}\cP_{r,r}(O_{F_v}/\varpi_v)\)
     \times\(K_r^\lozenge\times K_r^\lozenge\)
     \]
     (Definition \ref{de:lattice}) in which $K^\dag_{\lozenge\setminus\{\infty,p\}}$ contains
     \[
     \prod_{v\in\tV_F^{(\lozenge\setminus\{\infty,p\})}\setminus\tR}\(K_{r,v}\cap M_r(F_v)\)\times \(K_{r,v}\cap M_r(F_v)\)
     \]
     and fixes $\prod_{v\in\tV_F^{(\lozenge\setminus\{\infty,p\})}}f_{\chi_v}$ for every finite character $\chi\colon\Gamma_{F,p}\to\dC^\times$;

  \item fix $(\rt_1,\rt_2)\in\dS^\lozenge_{O_\dL}\times\dS^\lozenge_{O_\dL}$ that is the product of those pairs from Proposition \ref{pr:height} and Proposition \ref{pr:height1} for every $u\in\tP$ (and a suitable scalar), which satisfies $\chi_{\hat\pi}^\lozenge(\rt_1)\chi_\pi^\lozenge(\rt_2)\neq 0$;

  \item fix $(\rs_1,\rs_2)\in(\dS^\lozenge_{\pi,L}\cap\dS^\lozenge_{O_\dL})\times(\dS^\lozenge_{\hat\pi,L}\cap\dS^\lozenge_{O_\dL})$ such that $\chi_{\hat\pi}^\lozenge(\rs_1)\chi_\pi^\lozenge(\rs_2)\neq 0$, which is possible by Lemma \ref{le:hodge}(2).
\end{itemize}

\begin{remark}\label{re:hecke}
For every $v\in\tV_F^\spl\setminus\tV_F^{(\lozenge)}$, we have a canonical isomorphism
\[
\dZ[L_v\backslash H(F_v)/L_v]\simeq\dZ[K_{r,v}\backslash G_r(F_v)/K_{r,v}]
\]
of rings via Satake isomorphisms. By \cite{Liu11}*{Proposition~A.5}, we know that the action of $\rs\in\dZ[L_v\backslash H(F_v)/L_v]$ on $\sS(V_v^r)^{K_{r,v}\times L_v}$ via the Weil representation $\omega_{r,v}$ coincides with that of $\hat\rs\in\dZ[K_{r,v}\backslash G_r(F_v)/K_{r,v}]$, where $\hat\rs$ denotes the adjoint of $\rs$.
\end{remark}

Motivated by Proposition \ref{pr:height}, for every $e\in\dN$ and every pair $(g_1,g_2)\in M_r(F_\tR)\times M_r(F_\tR)$, we define the following elements in $\Gamma_{F,p}\otimes_{\dZ_p}\dL$:
\begin{align*}
(\cE^{[e]}_{(g_1,g_2)})_{T^\Box}\coloneqq p^M W_{2r}\cdot\sfW'_{T^\Box,v_{T^\Box}}(1)\cdot I_{T^\Box}((\rt_1\rs_1g_1\phi_1\otimes\rt_2\rs_2g_2\phi_2^{[e]})^{v_{T^\Box}})\cdot[v_{T^\Box}]
\end{align*}
for every $T^\Box\in\Herm_{2r}^\circ(F)^+_V$ (Notation \ref{pr:height}), and
\begin{align*}
(\pres{v}\cE^{[e]}_{(g_1,g_2)})_{T^\Box}\coloneqq  p^M W_{2r}\frac{2}{q_v^{2r}-1}\cdot
\sfW_{T^\Box,v}^\sph(1)\cdot I_{T^\Box}((\rt_1\rs_1g_1\phi_1\otimes\rt_2\rs_2g_2\phi_2^{[e]})^v)\cdot[v]
\end{align*}
for every $T^\Box\in\Herm_{2r}^\circ(F)^+$ and every $v\in\tS_\pi$.

For every finite character $\chi\colon\Gamma_{F,p}\to\ol\dQ_p^\times$, we denote by $\dL_\chi$ the finite (normal) extension of $\dL$ generated by values of $\chi$, which is a subfield of $\ol\dQ_p$.

\begin{lem}\label{le:modularity2}
We have
\begin{enumerate}
  \item There exists a (module-)finite $\dZ_{(p)}$-ring $\dO$ contained in $\dC$ such that for every $e\in\dN$, every pair $(g_1,g_2)\in M_r(F_\tR)\times M_r(F_\tR)$, and every $T^\Box\in\Herm_{2r}^\circ(F)^+$, there is a unique integral $\dL\otimes_{\dZ_{(p)}}\dO$-valued $p$-adic measure (Definition \ref{de:integral}) $(\sE^{[e]}_{(g_1,g_2)})_{T^\Box}$ on $\Gamma_{F,p}$, such that for every finite character $\chi\colon\Gamma_{F,p}\to\ol\dQ_p^\times$ and every embedding $\iota\colon\ol\dQ_p\to\dC$,
      \begin{align}\label{eq:modularity1}
      \iota(\sE^{[e]}_{(g_1,g_2)})_{T^\Box}(\chi)=p^M W_{2r}^\lozenge\cdot b_{2r}^\lozenge(\iota\chi) \cdot W_{T^\Box}((g_1,g_2)\cdot\iota(\hat{\rs_1}\hat{\rt_1},\hat{\rs_2}\hat{\rt_2})f_{\iota\chi^\infty}^{[e]}),
      \end{align}
      where $W_{2r}^\lozenge$ is from \eqref{eq:eisenstein1}.

  \item The measure $(\sE^{[e]}_{(g_1,g_2)})_{T^\Box}$ in (1) satisfies $(\sE^{[e]}_{(g_1,g_2)})_{T^\Box}(\b1)=0$ and
      \[
      \partial(\sE^{[e]}_{(g_1,g_2)})_{T^\Box}=
      \begin{dcases}
      (\cE^{[e]}_{(g_1,g_2)})_{T^\Box} &\text{if $T^\Box\in\Herm_{2r}^\circ(F)^+_V$,} \\
      0 &\text{if $T^\Box\in\Herm_{2r}^\circ(F)^+\setminus\Herm_{2r}^\circ(F)^+_V$.}
      \end{dcases}
      \]

  \item For every finite character $\chi\colon\Gamma_{F,p}\to\ol\dQ_p^\times$, the assignment
      \[
      \sE^{[e]}_{-}(\chi)\colon(g_1,g_2)\mapsto \sE^{[e]}_{(g_1,g_2)}(\chi)\coloneqq\sum_{T^\Box\in\Herm_{2r}^\circ(F)^+}(\sE^{[e]}_{(g_1,g_2)})_{T^\Box}(\chi)q^{T^\Box}
      \]
      belongs to $O_{\dL_\chi}\otimes_{\dZ_{(p)}}\SF_{2r}(\dO)^{M_r(F_\tR)\times M_r(F_\tR)}$.\footnote{Once again, we have strict inclusions
      \[
      O_{\dL_\chi}\otimes_{\dZ_{(p)}}\SF_{2r}(\dO)^{M_r(F_\tR)\times M_r(F_\tR)}\subsetneq(O_{\dL_\chi}\otimes_{\dZ_{(p)}}\SF_{2r}(\dO))^{M_r(F_\tR)\times M_r(F_\tR)}\subsetneq\SF_{2r}(O_{\dL_\chi}\otimes_{\dZ_{(p)}}\dO)^{M_r(F_\tR)\times M_r(F_\tR)}.
      \]} For every $\sigma\in\Gal(\ol\dQ_p/\dL)$, we have $\sE^{[e]}_{(g_1,g_2)}(\sigma\chi)=\sigma\sE^{[e]}_{(g_1,g_2)}(\chi)$, where $\sigma$ acts on $O_{\dL_\chi}\otimes_{\dZ_{(p)}}\SF_{2r}(\dO)$ via the first factor.

  \item For every $v\in\tS_\pi$, the assignment
      \[
      \pres{v}\cE^{[e]}_{-}\colon(g_1,g_2)\mapsto
      \pres{v}\cE^{[e]}_{(g_1,g_2)}\coloneqq\sum_{T^\Box\in\Herm_{2r}^\circ(F)^+}(\pres{v}\cE^{[e]}_{(g_1,g_2)})_{T^\Box}q^{T^\Box}
      \]
      belongs to $O_\dL\otimes_{\dZ_{(p)}}\SF_{2r}(\dZ_{(p)})^{M_r(F_\tR)\times M_r(F_\tR)}$.
\end{enumerate}
\end{lem}

\begin{proof}
By Lemma \ref{le:regular}, the right-hand side of \eqref{eq:modularity1} equals the product of the following five terms
\begin{align*}
& p^M W_{2r}^\lozenge\prod_{v\in\tV_F^{(p)}}b_{2r,v}(\CF)^{-1}\vol(L_v,\rd h_v), \\
& \iota\chi_p(\Nm_{E_p/F_p}\det T^\Box_{12})\CF_{\fT^{[e+\varepsilon^\tc]}_p}(T^\Box), \\
& \prod_{v\in\tR}W_{T^\Box}((g_{1,v},g_{2,v})f_{\iota\chi_v}), \\
& \prod_{v\in\tV_F^{(\lozenge\setminus\{\infty,p\})}\setminus\tR}W_{T^\Box}(f_{\iota\chi_v}), \\
& b_{2r}^\lozenge(\iota\chi)\cdot W_{T^\Box}(\iota(\hat{\rs_1}\hat{\rt_1},\hat{\rs_2}\hat{\rt_2})f_{\iota\chi^\lozenge}).
\end{align*}
We have
\begin{itemize}
  \item the element
      \begin{align*}
      C &\coloneqq p^M W_{2r}^\lozenge\prod_{v\in\tV_F^{(p)}}b_{2r,v}(\CF)^{-1}\vol(L_v,\rd h_v) \\
      &=p^M W_{2r}^{\lozenge\setminus\{p\}}\prod_{v\in\tV_F^{(p)}}\vol(L_v,\rd h_v)
      \end{align*}
      belongs to $\dZ_{(p)}$;

  \item there is an element $\sW_{T^\Box,p}\in\dZ[\Gamma_{F,p}]$ satisfying
      \[
      \iota\sW_{T^\Box,p}(\chi)=\iota\chi_p(\Nm_{E_p/F_p}\det T^\Box_{12})\CF_{\fT^{[e+\varepsilon^\tc]}_p}(T^\Box)
      \]
      for every $\chi,\iota$ as above;

  \item for every $v\in\tR$, there is a constant $C_{T^\Box,(g_{1,v},g_{2,v})}\in\dZ_{(p)}$ that equals $W_{T^\Box}((g_{1,v},g_{2,v})f_{\iota\chi_v})$ for every $\chi$ as above;

  \item by Lemma \ref{le:section7}, for every $v\in\tT$, there is an element $\sW_{T^\Box,v}\in\dZ[\Gamma_{F,p}]$ satisfying $\iota\sW_{T^\Box,v}(\chi)=b_{2r,v}(\iota\chi)\cdot W_{T^\Box}(f_{\iota\chi_v})$ for every $\chi,\iota$ as above;

  \item it is easy to see that for every $v\in\tV_F^{(\lozenge\setminus\{\infty,p\})}\setminus\tR$, there is an element $\sB_v\in\dZ[\Gamma_{F,p}]$ satisfying $\iota\sB_v(\chi)=b_{2r,v}(\iota\chi)^{-1}$ for every $\chi,\iota$ as above;

  \item by Lemma \ref{le:section4}(1), there exist $\dO$ as in the statement of (1) and finitely many elements $c_1,\dots,c_t\in O_\dL$ such that for every $T^\Box\in\Herm_{2r}^\circ(F)^+$, there are elements $\sW^\lozenge_{T^\Box,1},\dots,\sW^\lozenge_{T^\Box,t}$ in $\dO[\Gamma_{F,p}]$ satisfying
      \[
      \iota\sum_{i=1}^t c_i\otimes\sW^\lozenge_{T^\Box,i}(\chi)=b_{2r}^\lozenge(\iota\chi)\cdot W_{T^\Box}(\iota(\hat{\rs_1}\hat{\rt_1},\hat{\rs_2}\hat{\rt_2})f_{\iota\chi^\lozenge})
      \]
      for every $\chi,\iota$ as above.
\end{itemize}
For each $i$, put
\[
\sW_{T^\Box,\tT,i}\coloneqq
\(\prod_{v\in\tV_F^{(\lozenge\setminus\{\infty,p\})}\setminus\tR}\sW_{T^\Box,v}\)\cdot\sW^\lozenge_{T^\Box,i}
\in\dO[\Gamma_{F,p}].
\]
Moreover, for every $v\in\tT\setminus\tV_F^\spl$, we can write $\sW_{T^\Box,\tT,i}=\sW_{T^\Box,v}\cdot\sW_{T^\Box,\tT,i}^v$ for a unique element $\sW_{T^\Box,\tT,i}^v\in\dO[\Gamma_{F,p}]$.

For (1), we may take $(\sE^{[e]}_{(g_1,g_2)})_{T^\Box}$ to be
\begin{align*}
\sum_{i=1}^t c_i\otimes\(C\prod_{v\in\tR}C_{T^\Box,(g_{1,v},g_{2,v})}\)\cdot
\(\prod_{v\in\tV_F^{(\lozenge\setminus\{\infty,p\})}\setminus\tR}\sB_v\)
\cdot\sW_{T^\Box,p}\cdot\sW_{T^\Box,\tT,i}
\end{align*}
in $(O_\dL\otimes_{\dZ_{(p)}}\dO)[\Gamma_{F,p}]$. The uniqueness is automatic.

Part (3) is obvious from the construction in (1).

The proof of (4) is similar by realizing that $\sum_{i=1}^t c_i\otimes\sW_{T^\Box,\tT,i}^v(\b1)\in\dZ_{(p)}$.

For (2), note that $W_{T^\Box}(f_{\CF_v})=0$ for every $v\in\Diff(T^\Box,V)$, which implies that $(\sE^{[e]}_{(g_1,g_2)})_{T^\Box}(\b1)=0$, and for $T^\Box\in\Herm_{2r}^\circ(F)^+\setminus\Herm_{2r}^\circ(F)^+_V$ that $\partial(\sE^{[e]}_{(g_1,g_2)})_{T^\Box}=0$. For $T^\Box\in \Herm_{2r}^\circ(F)^+_V$, we have $\sE_{T^\Box,v_{T^\Box}}(\b1)=0$. Then by the $p$-adic Leibniz rule, we have
\begin{align*}
&\quad\partial(\sE^{[e]}_{(g_1,g_2)})_{T^\Box} \\
&=p^M W_{2r} \cdot b_{2r}^\infty(\b1) b_{2r,v_{T^\Box}}(\b1)^{-1}\cdot
W_{T^\Box}((g_1,g_2)\cdot\iota(\hat{\rs_1}\hat{\rt_1},\hat{\rs_2}\hat{\rt_2})f_{\iota\chi^{\infty v_{T^\Box}}}^{[e]})
\cdot\partial\sW_{T^\Box,v_{T^\Box}}.
\end{align*}
By Remark \ref{re:hecke} and \S\ref{ss:setup}(H9), we have
\begin{align*}
&\quad b_{2r}^\infty(\b1) b_{2r,v_{T^\Box}}(\b1)^{-1}\cdot
W_{T^\Box}((g_1,g_2)\cdot\iota(\hat{\rs_1}\hat{\rt_1},\hat{\rs_2}\hat{\rt_2})f_{\iota\chi^{\infty v_{T^\Box}}}^{[e]}) \\
&=I_{T^\Box}((\rt_1\rs_1g_1\phi_1\otimes\rt_2\rs_2g_2\phi_2^{[e]})^{v_{T^\Box}}).
\end{align*}
This, it remains to show that
\[
\partial\sW_{T^\Box,v_{T^\Box}}=\sfW'_{T^\Box,v_{T^\Box}}(1)\cdot[v_{T^\Box}],
\]
which is tautological as $\sW_{T^\Box,v_{T^\Box}}=\sfW_{T^\Box,v_{T^\Box}}([v_{T^\Box}])$ from Lemma \ref{le:section7}.

The lemma is proved.
\end{proof}

We search for Eisenstein series whose $q$-expansions are given by $\sE^{[e]}_{-}(\chi)$ and $\pres{v}\cE^{[e]}_{-}$. We refer to \S\ref{ss:eisenstein} for the notation concerning Eisenstein series. For every $e\in\dN$, every finite character $\chi\colon\Gamma_{F,p}\to\ol\dQ_p^\times$ and every embedding $\iota\colon\ol\dQ_p\to\dC$, define an Eisenstein series
\begin{align*}
E^{[e]}_{\iota\chi}\coloneqq p^M\cdot b_{2r}^\lozenge(\CF)^{-1} \cdot b_{2r}^\lozenge(\iota\chi) \cdot E(-,f_\infty^{[r]}\otimes \iota(\hat{\rs_1}\hat{\rt_1},\hat{\rs_2}\hat{\rt_2})f_{\iota\chi^\infty}^{[e]}) \in \cA_{2r,\hol}^{[r]}.
\end{align*}
and, for every $v\in\tS_\pi$, an Eisenstein series
\begin{align*}
\pres{v}E^{[e]}_\iota\coloneqq p^M\frac{2}{q_v^{2r}-1}\cdot E(-,f_\infty^{[r]}\otimes \iota(\hat{\rs_1}\hat{\rt_1},\hat{\rs_2}\hat{\rt_2})\pres{v}f_{\b1^\infty}^{[e]})\in\cA_{2r,\hol}^{[r]},
\end{align*}
where $\pres{v}f_{\b1^\infty}^{[e]}$ is obtained from $f_{\b1^\infty}^{[e]}$ after replacing the component $f_{\b1_v}$ by $f_{\b1_v}^\sph$ from Notation \ref{no:section1}.

\begin{lem}\label{le:modularity3}
We have
\begin{enumerate}
  \item For every finite character $\chi\colon\Gamma_{F,p}\to\ol\dQ_p^\times$ and every embedding $\iota\colon\ol\dQ_p\to\dC$,
     \[
     \bbq_{2r}^\an((g_1,g_2)\cdot E^{[e]}_{\iota\chi})=\iota\sE^{[e]}_{(g_1,g_2)}(\chi)
     \]
     holds for every $(g_1,g_2)\in M_r(F_\tR)\times M_r(F_\tR)$.

  \item For every $v\in\tS_\pi$ and every embedding $\iota\colon\ol\dQ_p\to\dC$,
     \[
     \bbq_{2r}^\an((g_1,g_2)\cdot \pres{v}E^{[e]}_\iota)=\iota\pres{v}\cE^{[e]}_{(g_1,g_2)}
     \]
     holds for every $(g_1,g_2)\in M_r(F_\tR)\times M_r(F_\tR)$.
\end{enumerate}
\end{lem}

\begin{proof}
Since for $v\in\tR$ (which is nonempty), $g_{1,v}\phi_{1,v}\otimes g_{2,v}\phi_{2,v}$ again belongs to $\sR_v$, both cases follow from the discussion in \cite{Liu12}*{Section~2B} and Lemma \ref{le:section6}.
\end{proof}

\begin{definition}\label{de:sf}
For every open compact subgroup $K\subseteq G_r(\dA_F^\infty)$ and every subring $\dM$ of $\dC$, we define $\cA^K_\dM$ the $\dM$-submodule of $\cA_{r,r,\hol}^{[r]}$ consisting of all $\varphi$ that are fixed by $K$ and satisfy $\bbq_{r,r}^\an((g_1,g_2)\cdot\varphi)\in\SF_{r,r}(\dM)$ for every $(g_1,g_2)\in M_r(F_\tR)\times M_r(F_\tR)$.
\end{definition}

\begin{lem}\label{le:sf}
Suppose that $K$ contains
\[
\prod_{v\in\tV_F^{\lozenge\setminus\{\infty\}}\setminus\tR}\(K_{r,v}\cap M_r(F_v)\)\times\(K_{r,v}\cap M_r(F_v)\)
\]
so that the tautological map $\cA^K_\dM\to\SF_{r,r}(\dM)^{M_r(F_\tR)\times M_r(F_\tR)}$ sending $\varphi$ to the assignment $(g_1,g_2)\mapsto \bbq_{r,r}^\an((g_1,g_2)\cdot\varphi)$ is injective.
\begin{enumerate}
  \item For rings $\dZ_{(p)}\subseteq\dM\subseteq\dM'\subseteq\dC$, the natural diagram
     \[
     \xymatrix{
     O_\dL\otimes_{\dZ_{(p)}}\cA^K_\dM \ar[r]\ar[d] & O_\dL\otimes_{\dZ_{(p)}}\cA^K_{\dM'} \ar[d] \\
     O_\dL\otimes_{\dZ_{(p)}}\SF_{r,r}(\dM)^{M_r(F_\tR)\times M_r(F_\tR)}
     \ar[r] & O_\dL\otimes_{\dZ_{(p)}}\SF_{r,r}(\dM')^{M_r(F_\tR)\times M_r(F_\tR)}
     }
     \]
     is Cartesian.

  \item For a ring $\dZ_{(p)}\subseteq\dM\subseteq\dC$, the natural diagram
     \[
     \xymatrix{
     O_\dL\otimes_{\dZ_{(p)}}\cA^K_\dM \ar[r]\ar[d] & \prod\limits_{\iota\colon\dL\to\dC}\cA^K_\dC \ar[d] \\
     O_\dL\otimes_{\dZ_{(p)}}\SF_{r,r}(\dM)^{M_r(F_\tR)\times M_r(F_\tR)}
     \ar[r] & \prod\limits_{\iota\colon\dL\to\dC}\SF_{r,r}(\dC)^{M_r(F_\tR)\times M_r(F_\tR)}
     }
     \]
     is Cartesian.
\end{enumerate}
\end{lem}

\begin{proof}
Part (1) follows from Definition \ref{de:sf} and the fact that $O_\dL$ is flat over $\dZ_{(p)}$.

For (2), consider an element $x=\sum_{j=1}^s c_j\otimes x_j$ of $O_\dL\otimes_{\dZ_{(p)}}\SF_{r,r}(\dM)^{M_r(F_\tR)\times M_r(F_\tR)}$ in which $c_1,\dots,c_s$ are $\dZ_{(p)}$-linearly independent elements of $O_\dL$, satisfying that for every $\iota\colon\dL\to\dC$, its image in $\SF_{r,r}(\dC)^{M_r(F_\tR)\times M_r(F_\tR)}$ comes from $\cA^K_\dC$. Since $\dL$ has characteristic zero, we may find embeddings $\iota_1,\dots,\iota_s$ such that $A\coloneqq(\iota_i c_j)_{1\leq i,j\leq s}$ is invertible. If we write $\iota_i x=y_i$ for $y_i\in\cA^K_\dC$, then $\pres{\rt}(x_1,\dots,x_s)=A^{-1}\pres{\rt}(y_1,\dots, y_s)$. In particular, $x$ belongs to $O_\dL\otimes_{\dZ_{(p)}}\SF_{r,r}(\dC)^{M_r(F_\tR)\times M_r(F_\tR)}$. Applying (1) with $\dM'=\dC$, we obtain (2).
\end{proof}

\begin{lem}\label{le:modularity4}
Recall the map $\varrho_{r,r}$ from Definition \ref{de:fourier}.
\begin{enumerate}
  \item For every finite character $\chi\colon\Gamma_{F,p}\to\ol\dQ_p^\times$, there exists $e_\chi\in\dN$ such that for every $e\geq e_\chi$, there exists a (unique) element $D^{[e]}_\chi\in O_{\dL_\chi}\otimes_{\dZ_{(p)}}\cA^{K^\dag}_\dO$ satisfying
      \[
      (1\times\bbq_{r,r}^\an)((g_1,g_2)\cdot D^{[e]}_\chi)=\varrho_{r,r}\sE^{[e]}_{(g_1,g_2)}(\chi)
      \]
      for every $(g_1,g_2)\in M_r(F_\tR)\times M_r(F_\tR)$. Moreover, the sequence $\{D^{[N!]}_\chi\}$ converges in $O_{\dL_\chi}\otimes_{\dZ_{(p)}}\cA^{K^\dag}_\dO$ when $N\to\infty$.

  \item For every $v\in\tS_\pi$ and $e\in\dN$, there exists a (unique) element $\pres{v}D^{[e]}\in O_\dL\otimes_{\dZ_{(p)}}\cA^{K^\dag}_{\dZ_{(p)}}$ satisfying
      \[
      (1\times\bbq_{r,r}^\an)((g_1,g_2)\cdot\pres{v}D^{[e]})=\varrho_{r,r}\pres{v}\cE^{[e]}_{(g_1,g_2)}
      \]
      for every $(g_1,g_2)\in M_r(F_\tR)\times M_r(F_\tR)$. Moreover, the sequence $\{\pres{v}D^{[N!]}\}$ converges in $O_\dL\otimes_{\dZ_{(p)}}\cA^{K^\dag}_{\dZ_{(p)}}$ when $N\to\infty$.
\end{enumerate}
\end{lem}

\begin{proof}
For (1), note that by Lemma \ref{le:modularity2}(3),
\[
\varrho_{r,r}\sE^{[e]}_{(g_1,g_2)}(\chi)\in O_{\dL_\chi}\otimes_{\dZ_{(p)}}\SF_{r,r}(\dO)^{M_r(F_\tR)\times M_r(F_\tR)}.
\]
Then by Lemma \ref{le:modularity3}(1) and Lemma \ref{le:sf}(2), for every $e\in\dN$, there is a (unique) element $D^{[e]}_\chi\in O_{\dL_\chi}\otimes_{\dZ_{(p)}}\cA^{K^{[e]}}_{\dO}$ for some subgroup $K^{[e]}\subseteq K^\dag$ of finite index such that
\[
(1\times\bbq_{r,r}^\an)((g_1,g_2)\cdot D^{[e]}_\chi)=\varrho_{r,r}\sE^{[e]}_{(g_1,g_2)}(\chi)
\]
holds for every $(g_1,g_2)\in M_r(F_\tR)\times M_r(F_\tR)$. Now Lemma \ref{le:open} tells us that we may take $K^{[e]}=K^\dag$ when $e\geq e_\chi$ in that lemma.

For the convergence, we have a natural inclusion
\begin{align}\label{eq:modularity0}
O_{\dL_\chi}\otimes_{\dZ_{(p)}}\cA^{K^\dag}_\dO
\hookrightarrow\dL_\chi\otimes_{\dQ_p}\(\cH_{r,r}^{[r]}(K^\dag)\otimes_\dQ\dM_{K^\dag}\)
\end{align}
(Definition \ref{de:holomorphic}) for some number field $\dM_{K^\dag}\subseteq\dC$ containing $\dO$ depending on $K^\dag$. It is well-known that the limit of the operators $\{\rU_p^{N!/2}\times\rU_p^{N!/2}\}_{N\geq 2}$ exists in $\End_{\dQ_p}\(\cH_{r,r}^{[r]}(K^\dag)\)$, which is the projection to the (Siegel-)ordinary part (see, for example, \cite{Hid98}*{Page~685}). Thus, by Lemma \ref{le:section5}, $\{D^{[N!]}_\chi\}$ converges in
\[
\dL_\chi\otimes_{\dQ_p}\(\cH_{r,r}^{[r]}(K^\dag)\otimes_\dQ\dM_{K^\dag}\).
\]
Since the inclusion \eqref{eq:modularity0} is closed, the limit belongs to the source.

The proof for (2) is similar, by using Lemma \ref{le:modularity2}(4) and Lemma \ref{le:modularity3}(2).
\end{proof}

\begin{lem}\label{le:modularity5}
For every $g_2\in G_r(F_\tR)$ and every $T_2\in\Herm_r^\circ(F)^+$, the sequence $\{Z_{T_2}^\clubsuit(\rt_2\rs_2g_2\phi_2^{[N!]})_L\}_{N}$ converges in $\rH^1_f(E,\rV_{\hat\pi,L})$.
\end{lem}

\begin{proof}
Fix an embedding $\dQ\langle p\rangle\hookrightarrow\ol\dQ_p$ and all representations will have coefficients in $\ol\dQ_p$. The assignment
\[
\bphi_p\mapsto\sS(V^r\otimes_FF_p,\ol\dQ_p)^{L_p}\to Z_{T_2}^\clubsuit(\rt_2\rs_2g_2\phi_2^p\bphi_p)_L
\]
factors through $\theta(\hat\pi_p)^L$ (\S\ref{ss:setup}(H10) but with $\dC$ replaced by $\ol\dQ_p$) as a $\ol\dQ_p[L_p\backslash H(F_p)/L_p]$-module, by the influence of $\rs_2$. Write $(\ol\dQ_p)_{T_2}$ the character of $N_r(F_p)$ such that for every $b\in\Herm_r(F_p)$, $n(b)$ acts by $\psi_{F,p}(\tr T_2b)$. Then, by Lemma \ref{le:theta}, there exists an element $w\in\Hom_{N_r(F_p)}(\hat\pi_p,(\ol\dQ_p)_{T_2})$ such that the assignment $Z_{T_2}^\clubsuit(\rt_2\rs_2g_2\phi_2^p-)_L$ factors through the composition
\[
\sS(V^r\otimes_FF_p,\ol\dQ_p)^{L_p}\to\hat\pi_p\otimes_{\ol\dQ_p}\theta(\hat\pi_p)^L\xrightarrow{w\otimes 1}\theta(\hat\pi_p)^L.
\]
By Lemma \ref{le:regular}(1), $\phi_{p,2}^{[e]}=\rU_p^e\phi_{p,2}^{[0]}$ hence is invariant under $I_p\coloneqq\prod_{v\in\tV_F^{(p)}}I_v$ \eqref{eq:iwahori}. Since $\{\rU_p^{N!}\}$ is convergent as a sequence of endomorphisms on $\hat\pi_p^{I_p}$ and $w\res_{\hat\pi_p^{I_p}}$ is continuous, the lemma follows.
\end{proof}

In what follows, we put
\begin{align*}
D_\chi&\coloneqq\lim\limits_{N\to\infty}D^{[N!]}_\chi\in O_{\dL_\chi}\otimes_{\dZ_{(p)}}\cA^{K^\dag}_\dO,\\
\pres{v}D&\coloneqq\lim\limits_{N\to\infty}\pres{v}D^{[N!]} \in O_\dL\otimes_{\dZ_{(p)}}\cA^{K^\dag}_{\dZ_{(p)}}, \\
\zeta_{g_2,T_2}&\coloneqq p^M\lim\limits_{N\to\infty}Z_{T_2}^\clubsuit(\rt_2\rs_2g_2\phi_2^{[N!]})_L\in\rH^1_f(E,\rV_{\hat\pi,L}).
\end{align*}
For every integer $d\geq1$, we recall the subgroup $U_d$ of $\Gamma_{F,p}$ and the set of representative $\Gamma_d$ of $\Gamma_{F,p}/U_d$ fixed from \S\ref{ss:measure}. For every $\lambda\in\Hom_{\dZ_p}(\Gamma_{F,p},\dZ_p)$, define
\[
I_{\lambda,d}\coloneqq
\(\sum_{x\in\Gamma_d}\frac{\lambda(x)}{|\Gamma_{F,p}/U_d|}\sum_{\chi\colon\Gamma_{F,p}/U_d\to\ol\dQ_p^\times}\chi(x)^{-1}D_\chi\)
+\(\sum_{v\in\tS_\pi}\pres{v}D\).
\]
By Lemma \ref{le:modularity2}, we have $D^{[e]}_{\sigma\chi}=\sigma D^{[e]}_\chi$ hence $D_{\sigma\chi}=\sigma D_\chi$ for every $\sigma\in\Gal(\ol\dQ_p/\dL)$. Thus, $I_{\lambda,d}$ belongs to $O_\dL\otimes_{\dZ_{(p)}}\cA^{K^\dag}_\dO$.

\begin{proposition}\label{pr:modularity2}
Suppose that $n<p$. Take an element $\lambda\in\Hom_{\dZ_p}(\Gamma_{F,p},\dZ_p)$ and put $\lambda_E\coloneqq\lambda\circ\Nm_{E/F}$.
\begin{enumerate}
  \item The sequence $\{I_{\lambda,d}\}_{d\geq 1}$ is a convergent sequence in $\dL\otimes_{\dZ_{(p)}}\cA^{K^\dag}_\dO$.

  \item Put $I_\lambda\coloneqq\lim\limits_{d\to\infty}I_{\lambda,d}$. Then
      \begin{align*}
      &\quad(1\times\bbq_{r,r}^\an)((g_1,g_2)\cdot I_\lambda) \\
      &=\vol^\natural(L)\sum_{T_1,T_2\in\Herm_r^\circ(F)^+\times\Herm_r^\circ(F)^+}
      \lambda_E\left\langle Z_{T_1}^\clubsuit(\rt_1\rs_1g_1\phi_1)_L,\zeta_{g_2,T_2}\right\rangle_E\cdot q^{T_1,T_2}
      \end{align*}
      holds for every $(g_1,g_2)\in M_r(F_\tR)\times M_r(F_\tR)$.

  \item The limit $I_\lambda$ belongs to $O_\dL\otimes_{\dZ_{(p)}}\cA^{K^\dag}_{\dZ_{(p)}}$.
\end{enumerate}
\end{proposition}

\begin{proof}
To simplify the notation in the proof, we introduce the following.
\begin{itemize}
  \item For a function $\sF$ on the set of finite characters $\chi\colon\Gamma_{F,p}\to\ol\dQ_p^\times$, we put
      \[
      \sF(\lambda_d)\coloneqq\sum_{x\in\Gamma_d}\frac{\lambda(x)}{|\Gamma_{F,p}/U_d|}
      \sum_{\chi\colon\Gamma_{F,p}/U_d\to\ol\dQ_p^\times}\chi(x)^{-1}\sF(\chi)
      \]
      for every integer $d\geq 1$.\footnote{Here, we can think of $\lambda_d$ as the $U_d$-invariant function that agrees with $\lambda$ on $\Gamma_d$; though it is not a character.}

  \item For two elements $x,y\in\dL\otimes_{\dZ_{(p)}}\dO$ and an integer $d$, we write $x\equiv_d y$ if $x-y$ belongs to $p^d O_\dL\otimes_{\dZ_{(p)}}\dO$.
\end{itemize}

For every $(g_1,g_2)\in M_r(F_\tR)\times M_r(\tR)$ and $(T_1,T_2)\in\Herm_r^\circ(F)^+\times\Herm_r^\circ(F)^+$, denote by $(\sD_{(g_1,g_2)})_{T_1,T_2}(\chi)$ the $q^{T_1,T_2}$-th coefficient of $(1\times\bbq_{r,r}^\an)((g_1,g_2)\cdot D_\chi)$ and $\pres{v}(\cD_{(g_1,g_2)})_{T_1,T_2}$ the $q^{T_1,T_2}$-th coefficient of $(1\times\bbq_{r,r}^\an)((g_1,g_2)\cdot\pres{v}D)$. We claim that the two sequences
\begin{align}\label{eq:converge1}
(\sD_{(g_1,g_2)})_{T_1,T_2}(\lambda_d)+
\(\sum_{v\in\tS_\pi}\pres{v}(\cD_{(g_1,g_2)})_{T_1,T_2}\),
\end{align}
\begin{align}\label{eq:converge2}
\sum_{\substack{T^\Box\in\Herm_{2r}^\circ(F)^+ \\ \partial_{r,r}T^\Box=(T_1,T_2)}}
\(\partial_\lambda\sE^{[N!]}_{(g_1,g_2)})_{T^\Box}+\sum_{v\in\tS_\pi}\pres{v}(\cE^{[N!]}_{(g_1,g_2)})_{T^\Box}\),
\end{align}
in $\dL\otimes_{\dZ_{(p)}}\dO$, both converge to $\lambda_E\left\langle Z_{T_1}^\clubsuit(\rt_1\rs_1g_1\phi_1)_L,\zeta_{g_2,T_2}\right\rangle_E$ when $d\to\infty$ and $N\to\infty$, respectively.

Parts (1) and (2) already follow from the limit formula for \eqref{eq:converge1} since $\dL\otimes_{\dZ_{(p)}}\cA^{K^\dag}_\dO$ is a finite dimensional $\dL$-vectors space. For (3), by the convergence of \eqref{eq:converge2}, the assignment
\[
(g_1,g_2)\mapsto(1\times\bbq_{r,r}^\an)((g_1,g_2)\cdot I_\lambda)
\]
belongs to $\SF_{r,r}(O_\dL)^{M_r(F_\tR)\times M_r(F_\tR)}$. It is straightforward to check that
\begin{align*}
&\quad\dL\otimes_{\dZ_{(p)}}\SF_{r,r}(\dO)^{M_r(F_\tR)\times M_r(F_\tR)}\cap\SF_{r,r}(O_\dL)^{M_r(F_\tR)\times M_r(F_\tR)} \\
&=O_\dL\otimes_{\dZ_{(p)}}\SF_{r,r}(\dZ_{(p)})^{M_r(F_\tR)\times M_r(F_\tR)}.
\end{align*}
Then (3) follows from (1) and Lemma \ref{le:sf}(1) (with $\dM=\dZ_{(p)}$ and $\dM'=\dO$).

It remains to prove the claim. Without loss of generality, we assume $(g_1,g_2)=(1_r,1_r)$ and suppress it (together with redundant parentheses) from the notation (in particular, $\zeta_{T_2}$ means $\zeta_{1_r,T_2}$). For every $d\geq 1$, we may find an element $N_d\in\dN$ satisfying:
\begin{enumerate}[label=(\alph*)]
  \item $N_d!\geq e_\chi$ (Lemma \ref{le:modularity4}(1)) for every $\chi\colon\Gamma_{F,p}/U_d\to\ol\dQ_p^\times$;

  \item for every integer $N\geq N_d$,
     \[
     \sD_{T_1,T_2}(\lambda_d)\equiv_d
     \sum_{\substack{T^\Box\in\Herm_{2r}^\circ(F)^+ \\ \partial_{r,r}T^\Box=(T_1,T_2)}}\sE^{[N!]}_{T^\Box}(\lambda_d);
     \]

  \item for every integer $N\geq N_d$,
     \[
     \pres{v}\cD_{T_1,T_2}\equiv_d\sum_{\substack{T^\Box\in\Herm_{2r}^\circ(F)^+ \\ \partial_{r,r}T^\Box=(T_1,T_2)}}\pres{v}\cE^{[N!]}_{T^\Box}
     \]
     holds for every $v\in\tS_\pi$

  \item $N_d!\geq d+M_u$ for every $u\in\tP$, where $M_u$ is the integer from Proposition \ref{pr:height1} (for $u$);

  \item for every integer $N\geq N_d$,
     \begin{align*}
     &\quad\vol^\natural(L)\cdot\lambda_E\left\langle Z_{T_1}^\clubsuit(\rt_1\rs_1\phi_1)_L,
     p^MZ_{T_2}^\clubsuit(\rt_2\rs_2\phi_2^{[N!]})_L\right\rangle_E\\
     &\equiv_d
     \vol^\natural(L)\cdot\lambda_E\left\langle Z_{T_1}^\clubsuit(\rt_1\rs_1\phi_1)_L,\zeta_{T_2}\right\rangle_E.
     \end{align*}
\end{enumerate}
By (b) and (c), for the claim, it suffices to show that
\begin{align}\label{eq:modularity2}
&\quad\sum_{\substack{T^\Box\in\Herm_{2r}^\circ(F)^+ \\ \partial_{r,r}T^\Box=(T_1,T_2)}}\sE^{[N!]}_{T^\Box}(\lambda_d)
+\sum_{v\in\tS_\pi}\sum_{\substack{T^\Box\in\Herm_{2r}^\circ(F)^+ \\ \partial_{r,r}T^\Box=(T_1,T_2)}}\pres{v}\cE^{[N!]}_{T^\Box}\\
&\equiv_d
\vol^\natural(L)\cdot\lambda_E\left\langle Z_{T_1}^\clubsuit(\rt_1\rs_1\phi_1)_L,\zeta_{T_2}\right\rangle_E \notag
\end{align}
for every integer $N\geq N_d$. By Lemma \ref{le:modularity2}(1,2) and Lemma \ref{le:integral},
\[
\sum_{\substack{T^\Box\in\Herm_{2r}^\circ(F)^+ \\ \partial_{r,r}T^\Box=(T_1,T_2)}}\sE^{[N!]}_{T^\Box}(\lambda_d)
\equiv_d\sum_{\substack{T^\Box\in\Herm_{2r}^\circ(F)^+ \\ \partial_{r,r}T^\Box=(T_1,T_2)}}
\partial\sE^{[N!]}_{T^\Box}.
\]
By \eqref{eq:decompose}, Proposition \ref{pr:height} and Proposition \ref{pr:height1} (which is applicable by (d)),
\begin{align*}
&\quad\sum_{\substack{T^\Box\in\Herm_{2r}^\circ(F)^+ \\ \partial_{r,r}T^\Box=(T_1,T_2)}}
\(\partial\sE^{[N!]}_{T^\Box}+\sum_{v\in\tS_\pi}\pres{v}\cE^{[N!]}_{T^\Box}\)\\
&\equiv_d
\vol^\natural(L)\cdot\lambda_E\left\langle Z_{T_1}^\clubsuit(\rt_1\rs_1\phi_1)_L,p^MZ_{T_2}^\clubsuit(\rt_2\rs_2\phi_2^{[N!]})_L\right\rangle_E.
\end{align*}
Thus, \eqref{eq:modularity2} follows from the above two relations and (e).

The proposition is proved.
\end{proof}

\begin{proposition}\label{pr:modularity3}
Suppose that
\begin{enumerate}[label=(\alph*)]
  \item $n<p$;

  \item $\partial\sL_p^\lozenge(\pi)\neq 0$;

  \item for every $v\in\tR$, there exist $\varphi_v^\vee\in\pi_v^\vee$, $\varphi_v\in\pi_v$ such that $Z(\varphi_v^\vee\otimes\varphi_v,f_{\Phi_v}^\r{SW})\neq 0$ (Lemma \ref{le:zeta1}).
\end{enumerate}
Then there exists $\lambda\in\Hom_{\dZ_p}(\Gamma_{F,p},\dZ_p)$ such that $I_\lambda\neq 0$.
\end{proposition}

The proof of the above proposition will be given in the next subsection. Now we move to the proof of Theorem \ref{th:modularity}.

\begin{proof}[Proof of Theorem \ref{th:modularity}]
By \cite{LL}*{Proposition~3.13}, for every $v\in\tR$, we may choose a pair $(\phi_{v,1},\phi_{v,2})\in\sR_v$ (\S\ref{ss:strategy}(S1)) such that condition (c) in Proposition \ref{pr:modularity3} holds. Choose $\lambda\in\Hom_{\dZ_p}(\Gamma_{F,p},\dZ_p)$ such that $I_\lambda\neq 0$ by this proposition. In particular, we may choose some $g_2\in G_r(F_\tR)$ and $T_2\in\Herm_r^\circ(F)^+$, such that the $q^{T_2}$-th coefficient of $(1,g_2)\cdot I_\lambda$, which we denote by $\varphi_{g_2,T_2,\lambda}$, is nonzero. Since $I_\lambda$ belongs to $O_\dL\otimes_{\dZ_{(p)}}\cA^K_{\dZ_{(p)}}$ by Proposition \ref{pr:modularity2}(3), $\varphi_{g_2,T_2,\lambda}$ is a strongly nonzero element in $\dL\otimes_\dQ\cA_{r,\hol}^{[r]}$ (Definition \ref{de:strong}), which satisfies
\begin{align*}
&\quad(1\times\bbq_r^\an)(g_1\cdot \varphi_{g_2,T_2,\lambda}) \\
&=\sum_{T_1\in\Herm_r^\circ(F)^+}
\lambda_E\left\langle Z_{T_1}^\clubsuit(g_1\rt_1\rs_1\phi_1)_L,\zeta_{g_2,T_2}\right\rangle_E\cdot q^{T_1}\\
&=\sum_{T_1\in\Herm_r^\circ(F)^+}
\lambda_E\left\langle\wp_\pi(Z_{T_1}^\clubsuit(g_1\rt_1\rs_1\phi_1)_L),\zeta_{g_2,T_2}\right\rangle_E\cdot q^{T_1}
\end{align*}
for every $g_1\in G_r(F_\tR)$ by Proposition \ref{pr:modularity2}(2). By Lemma \ref{le:modularity1}, the above identity indeed holds for every $g_1\in G_r(\dA_F^\infty)$. Thus, we may apply Lemma \ref{le:modularity} with $L$, $\rt_1\rs_1\phi_1$, $\zeta_{g_2,T_2}$, $\lambda_E$ and $\varphi_{g_2,T_2,\lambda}$, hence Theorem \ref{th:modularity} follows.
\end{proof}

\begin{remark}
Unfortunately, the strategy for proving Theorem \ref{th:modularity} hence giving an unconditional construction of the Selmer theta lifts can not be applied to give an unconditional construction of the arithmetic theta lifts (on the level of Chow groups) appeared in \cites{LL,LL2}, since our strategy relies on the fact that $\rH^1(E,\rH^{2r-1}(\ol{X}_L,\dQ_p(r)))$ as a $\dQ_p[L\backslash H(\dA_F^\infty)/L]$-module is semisimple and automorphic -- this is not known for $\CH^r(X_L)$.
\end{remark}

\subsection{Proof of Theorem \ref{th:aipf}}
\label{ss:proof1}

In this subsection, we prove Proposition \ref{pr:modularity3} and Theorem \ref{th:aipf}. Both proofs require choices of vectors from $\hat\pi$ and $\pi$, which we do now. Choose decomposable elements $\varphi_1=\otimes_v\varphi_{1,v}\in\cV_{\hat\pi}$ and $\varphi_2=\otimes_v\varphi_{2,v}\in\cV_\pi$ satisfying
\begin{enumerate}[label=(T\arabic*)]\setcounter{enumi}{1}

  \item $\varphi^\dag_{1,v}\in(\pi_v^\vee)^-$, $\varphi_{2,v}\in\pi_v^-$ and $\langle\pi_v^\vee(\tw_r)\varphi^\dag_{1,v},\varphi_{2,v}\rangle_{\pi_v}=q_v^{-d_vr^2}$ for $v\in\tV_F^{(p)}$,

  \item $\varphi^\dag_{1,v}\in(\pi_v^\vee)^{K_{r,v}}$, $\varphi_{2,v}\in\pi_v^{K_{r,v}}$ and $\langle\varphi^\dag_{1,v},\varphi_{2,v}\rangle_{\pi_v}=1$ for $v\in\tT\setminus\tS_\pi$,

  \item $\varphi^\dag_{1,v},\varphi_{2,v}$ are new vectors\footnote{A new vector in an almost unramified representation of $G_r(F_v)$ is a vector in the (one-dimensional) space in \cite{Liu22}*{Definition~5.3(2)}.} with respect to $K_{r,v}$ and $\langle\varphi^\dag_{1,v},\varphi_{2,v}\rangle_{\pi_v}=1$ for $v\in\tS_\pi$.
\end{enumerate}

\begin{proposition}\label{pr:modularity1}
Suppose that $n<p$. Take an element $\lambda\in\Hom_{\dZ_p}(\Gamma_{F,p},\dZ_p)$ and regard $I_\lambda$ as an element of $\dL\otimes_{\dQ_p}\(\cH_{r,r}^{[r]}(K^\dag)\otimes_\dQ\dC\)$ (Definition \ref{de:holomorphic}). Then
\begin{align*}
&\quad\left\langle\varphi_1\otimes\varphi_2,I_\lambda\right\rangle_{\pi,\hat\pi} \\
&=p^M\cdot \chi^\lozenge_{\hat\pi}(\rt_1\rs_1)\chi^\lozenge_{\pi}(\rt_2\rs_2)\cdot\partial_\lambda\sL_p^\lozenge(\pi)
\cdot\prod_{v\in\tV_F^{(\lozenge\setminus\{\infty,p\})}}Z(\varphi_{1,v}^\dag\otimes\varphi_{2,v},f_{\CF_v}),
\end{align*}
where $\left\langle\;,\;\right\rangle_{\pi,\hat\pi}$ is introduced in Notation \ref{no:doubling} and $Z(\varphi_{1,v}^\dag\otimes\varphi_{2,v},f_{\CF_v})$ is from Lemma \ref{le:zeta1}.
\end{proposition}

\begin{proof}
We first compute $\left\langle\varphi_1\otimes\varphi_2,D_\chi\right\rangle_{\pi,\hat\pi}$ and $\left\langle\varphi_1\otimes\varphi_2,\pres{v}D\right\rangle_{\pi,\hat\pi}$.

Let $\chi\colon\Gamma_{F,p}\to\ol\dQ_p^\times$ be a finite character. By definition, we have
\[
\left\langle\varphi_1\otimes\varphi_2,D_\chi\right\rangle_{\pi,\hat\pi}=\lim_{N\to\infty}
\left\langle\varphi_1\otimes\varphi_2,D^{[N!]}_\chi\right\rangle_{\pi,\hat\pi}.
\]
For the right-hand side, we perform a computation similar to the one in the proof of Theorem \ref{th:measure}. For every embedding $\iota\colon\ol\dQ_p\to\dC$, we have
\begin{align*}
&\quad\iota\left\langle\varphi_1\otimes\varphi_2,D^{[N!]}_\chi\right\rangle_{\pi,\hat\pi} \\
&=\frac{1}{(\rP^\iota_\pi)^2}\iint\limits_{\(G_r(F)\backslash G_r(\dA_F)\)^2}\varphi_1^\iota(g_1^\dag)\varphi_2^\iota(g_2^\dag)E^{[N!]}_{\iota\chi}((g_1,g_2))\rd g_1\rd g_2 \notag\\
&=\frac{1}{(\rP^\iota_\pi)^2}\iint\limits_{\(G_r(F)\backslash G_r(\dA_F)\)^2}(\varphi_1^\dag)^\iota(g_1)\varphi_2^\iota(g_2)E^{[N!]}_{\iota\chi}(\imath(g_1,g_2))\rd g_1\rd g_2
\end{align*}
by Lemma \ref{le:modularity4}(1) and Lemma \ref{le:modularity3}(1). By the doubling integral expansion and Lemma \ref{le:zeta2},
\begin{align*}
&\quad\iota\left\langle\varphi_1\otimes\varphi_2,D^{[N!]}_\chi\right\rangle_{\pi,\hat\pi} \\
&=p^M\cdot \iota\chi^\lozenge_{\hat\pi}(\rt_1\rs_1) \iota\chi^\lozenge_{\pi}(\rt_2\rs_2)
\cdot\frac{1}{\rP^\iota_\pi}\cdot\frac{Z_r^{[F:\dQ]}}{b_{2r}^\lozenge(\CF)}
\cdot L(\tfrac{1}{2},\BC(\iota\pi^\lozenge)\otimes(\iota\chi^\lozenge\circ\Nm_{E/F}))\\
&\quad\times\prod_{v\in\tV_F^{(p)}}Z^\iota(\varphi_{1,v}^\dag\otimes\varphi_{2,v},(\tf_{\iota\chi_v}^{[N!]})^{\iota\chi_v})
\prod_{v\in\tV_F^{(\lozenge\setminus\{\infty,p\})}}Z^\iota(\varphi_{1,v}^\dag\otimes\varphi_{2,v},f_{\iota\chi_v}).
\end{align*}
By (T2) and Lemma \ref{le:section5}, for every $v\in\tV_F^{(p)}$,
\[
Z^\iota(\varphi_{1,v}^\dag\otimes\varphi_{2,v},(\tf_{\iota\chi_v}^{[N!]})^{\iota\chi_v})=\(\iota\prod_{u\in\tP_v}\alpha(\pi_u)\)^{-N!}
Z^\iota(\varphi_{1,v}^\dag\otimes\varphi_{2,v},(\tf_{\iota\chi_v}^{[0]})^{\iota\chi_v}).
\]
By Proposition \ref{pr:zeta} and (T2), for every $v\in\tV_F^{(p)}$,
\[
Z^\iota(\varphi_{1,v}^\dag\otimes\varphi_{2,v},(\tf_{\iota\chi_v}^{[0]})^{\iota\chi_v})=
\prod_{u\in\tP_v}\gamma(\tfrac{1+r}{2},\iota\ul{\pi_u}\otimes\chi_v,\psi_{F,v})^{-1}.
\]
Together, we have
\begin{align*}
&\quad\iota\left\langle\varphi_1\otimes\varphi_2,D^{[N!]}_\chi\right\rangle_{\pi,\hat\pi} \\
&=\frac{1}{\rP^\iota_\pi}\cdot\frac{Z_r^{[F:\dQ]}}{b_{2r}^\lozenge(\CF)}\cdot
\prod_{v\in\tV_F^{(p)}}\prod_{u\in\tP_v}\gamma(\tfrac{1+r}{2},\iota(\ul{\pi_u}\otimes\chi_v),\psi_{F,v})^{-1} \\
&\quad\times L(\tfrac{1}{2},\BC(\iota\pi^\lozenge)\otimes(\iota\chi^\lozenge\circ\Nm_{E/F})) \\
&\quad\times p^M\cdot \iota\chi^\lozenge_{\hat\pi}(\rt_1\rs_1) \iota\chi^\lozenge_{\pi}(\rt_2\rs_2)  \\ &\quad\times\(\iota\prod_{u\in\tP}\alpha(\pi_u)\)^{-N!}\cdot
\prod_{v\in\tV_F^{(\lozenge\setminus\{\infty,p\})}}Z^\iota(\varphi_{1,v}^\dag\otimes\varphi_{2,v},f_{\iota\chi_v}),
\end{align*}
which, by Theorem \ref{th:measure} and Lemma \ref{le:zeta1}, equals
\begin{align*}
& p^M\cdot \iota\(\chi^\lozenge_{\hat\pi}(\rt_1\rs_1)\chi^\lozenge_{\pi}(\rt_2\rs_2)\)\cdot
\iota\sL_p^\lozenge(\pi)(\chi)\\
&\times \iota\(\prod_{u\in\tP}\alpha(\pi_u)^{-N!}\)\cdot\iota
\(\prod_{v\in\tV_F^{(\lozenge\setminus\{\infty,p\})}}Z(\varphi_{1,v}^\dag\otimes\varphi_{2,v},f_{\chi_v})\).
\end{align*}
As a consequence, we have
\begin{align*}
&\quad\left\langle\varphi_1\otimes\varphi_2,D^{[N!]}_\chi\right\rangle_{\pi,\hat\pi} \\
&= p^M\cdot\chi^\lozenge_{\hat\pi}(\rt_1\rs_1)\chi^\lozenge_{\pi}(\rt_2\rs_2)\cdot \sL_p^\lozenge(\pi)(\chi)\\
&\quad\times\(\prod_{u\in\tP}\alpha(\pi_u)^{-N!}\)\cdot
\prod_{v\in\tV_F^{(\lozenge\setminus\{\infty,p\})}}Z(\varphi_{1,v}^\dag\otimes\varphi_{2,v},f_{\chi_v}),
\end{align*}
hence
\begin{align}\label{eq:modularity9}
&\quad\left\langle\varphi_1\otimes\varphi_2,D_\chi\right\rangle_{\pi,\hat\pi} \\
&=p^M\cdot\chi^\lozenge_{\hat\pi}(\rt_1\rs_1)\chi^\lozenge_{\pi}(\rt_2\rs_2)\cdot
\sL_p^\lozenge(\pi)(\chi)\cdot\prod_{v\in\tV_F^{(\lozenge\setminus\{\infty,p\})}}
Z(\varphi_{1,v}^\dag\otimes\varphi_{2,v},f_{\chi_v}). \notag
\end{align}

By a similar argument, for every $v\in\tS_\pi$, we have
\[
\left\langle\varphi_1\otimes\varphi_2,\pres{v}D^{[N!]}\right\rangle_{\pi,\hat\pi}=0
\]
since $Z(\varphi_{1,v}^\dag\otimes\varphi_{2,v},f_{\b1_v}^\sph)=0$. Thus, $\left\langle\varphi_1\otimes\varphi_2,\pres{v}D\right\rangle_{\pi,\hat\pi}=0$.

Now the proposition follows from \eqref{eq:modularity9}, \eqref{eq:derivative}, and the $p$-adic Leibniz rule.
\end{proof}

\begin{lem}\label{le:modularity6}
For every $v\in\tV_F^{(\lozenge\setminus\{\infty,p\})}\setminus\tR$, we have $Z(\varphi_{1,v}^\dag\otimes\varphi_{2,v},f_{\CF_v})\neq 0$.
\end{lem}

\begin{proof}
By \cite{Liu22}*{Proposition~5.6} and (T4) when $v\in\tS_\pi$, \cite{LL2}*{Proposition~3.6} and (T3) when $v\in\tV_F^\ram$, Lemma \ref{le:zeta2} and (T3) when $v\in\tT\setminus(\tS_\pi\cup\tV_F^\ram)$, we have
\[
Z(\varphi_{1,v}^\dag\otimes\varphi_{2,v},f_{\CF_v})=C_v\cdot\frac{L(\tfrac{1}{2},\BC(\pi_v))}{b_{2r,v}(\CF)}
\]
for a constant $C_v\in\dQ^\times$. Then the nonvanishing is clear.
\end{proof}

\begin{proof}[Proof of Proposition \ref{pr:modularity3}]
We would like to apply Proposition \ref{pr:modularity1}. By condition (c), for every $v\in\tR$, we may find $\varphi_v^\vee\in\pi_v^\vee$, $\varphi_v\in\pi_v$ such that $Z(\varphi_v^\vee\otimes\varphi_v,f_{\Phi_v}^\r{SW})\neq 0$, that is, $Z(\varphi_{1,v}^\dag\otimes\varphi_{2,v},f_{\CF_v})\neq 0$. Together with Lemma \ref{le:modularity6}, we have
\[
\prod_{v\in\tV_F^{(\lozenge\setminus\{\infty,p\})}}Z(\varphi_{1,v}^\dag\otimes\varphi_{2,v},f_{\CF_v})\neq 0.
\]
By condition (b), there exists $\lambda\in\Hom_{\dZ_p}(\Gamma_{F,p},\dZ_p)$ such that $\partial_\lambda\sL_p^\lozenge(\pi)\neq 0$. Thus, by Proposition \ref{pr:modularity1}, $I_\lambda\neq 0$. The proposition is proved.
\end{proof}

\begin{proof}[Proof of Theorem \ref{th:aipf}]
For (1), we may apply Theorem \ref{th:modularity} so that we have elements $\cZ^\pi_{\phi_1,L}$ and $\cZ^{\hat\pi}_{\phi_2,L}$ from Proposition \ref{pr:generating}. By Remark \ref{re:aipf}(2), it suffices to show \eqref{eq:aipf} for a single choice of data $(\varphi_1,\varphi_2,\phi_1,\phi_2)$ (as in the statement of Theorem \ref{th:aipf}) satisfying
\[
\prod_{v\in\tV_F^{(\lozenge\setminus\{\infty\})}}
Z(\varphi_{1,v}^\dag\otimes\varphi_{2,v},f^{\r{SW}}_{\phi_{1,v}\otimes\phi_{2,v}})\neq 0.
\]
Thus, by Lemma \ref{le:modularity6}, it suffices to show \eqref{eq:aipf} for our particular choices of $(\phi_1,\phi_2\coloneqq\phi_2^{[0]})$ as in \eqref{eq:phi} and $(\varphi_1,\varphi_2)$ from (T2--T4), together satisfying the following extra requirement
\begin{enumerate}[label=(T\arabic*)]
  \item $Z(\varphi_{1,v}^\dag\otimes\varphi_{2,v},f_{\phi_{1,v}\otimes\phi_{2,v}}^{\r{SW}})\neq0$ for $v\in\tR$.
\end{enumerate}
This is possible by \cite{LL}*{Proposition~3.13}.

By Remark \ref{re:hecke} and Lemma \ref{le:regular}(1),
\begin{align*}
\Theta_{\phi_1}^\sel(\varphi_1)_L&= \chi^\lozenge_{\hat\pi}(\rt_1\rs_1)^{-1}\Theta_{\rt_1\rs_1\phi_1}^\sel(\varphi_1)_L, \\
\Theta_{\phi_2}^\sel(\varphi_2)_L&= \chi^\lozenge_{\pi}(\rt_2\rs_2)^{-1}\(\prod_{u\in\tP}\alpha(\pi_u)^{-e}\)
\Theta_{\rt_2\rs_2\phi_2^{[e]}}^\sel(\varphi_2)_L
\end{align*}
hold for every $e\in\dN$. By Definition \ref{de:theta},
\begin{align*}
\Theta_{\rt_1\rs_1\phi_1}^\sel(\varphi_1)_L &= \langle\varphi_1^\dag,\cZ^\pi_{\rt_1\rs_1\phi_1,L}\rangle_{\pi},\\
\Theta_{\rt_2\rs_2\phi_2^{[e]}}^\sel(\varphi_2)_L &=
\langle\varphi_2^\dag,\cZ^{\hat\pi}_{\rt_2\rs_2\phi_2^{[e]},L}\rangle_{\hat\pi},
\end{align*}
in which
\begin{align*}
\cZ^\pi_{\rt_1\rs_1\phi_1,L}&\in\rH^1_f(E,\rV_{\hat\pi,L})\otimes_\dL\(\cV_{\hat\pi}\otimes_\dQ\dM\),\\ \cZ^{\hat\pi}_{\rt_2\rs_2\phi_2^{[e]},L}&\in\rH^1_f(E,\rV_{\hat\pi,L})\otimes_\dL\(\cV_{\hat\pi}\otimes_\dQ\dM\)
\end{align*}
for some field $\dM\subseteq\dC$. Indeed, for given $\tR$-components of $\phi_1$ and $\phi_2$, we can shrink $\dM$ to a number field, which is in particular independent of $e$. Again by Lemma \ref{le:regular}(1), the sequence $\{\cZ^{\hat\pi}_{\rt_2\rs_2\phi_2^{[N!]},L}\}$ converges when $N\to\infty$, whose limit we simply denote by $\cZ_2$. Then
\[
\Theta_{\phi_2}^\sel(\varphi_2)_L =\chi^\lozenge_{\pi}(\rt_2\rs_2)^{-1}\langle\varphi_2^\dag,\cZ_2\rangle_{\hat\pi}.
\]
Therefore, for every element $\lambda\in\Hom_{\dZ_p}(\Gamma_{F,p},\dZ_p)$ with $\lambda_E\coloneqq\lambda\circ\Nm_{F/F}$, we have
\begin{align*}
&\quad\lambda\langle\Theta_{\phi_1}^\sel(\varphi_1),\Theta_{\phi_2}^\sel(\varphi_2)\rangle_{\pi,F}^\natural \\
&=\vol^\natural(L)\cdot\lambda_E\langle\Theta_{\phi_1}^\sel(\varphi_1)_L,\Theta_{\phi_2}^\sel(\varphi_2)_L\rangle_E \notag\\
&=\chi^\lozenge_{\hat\pi}(\rt_1\rs_1)^{-1}\chi^\lozenge_{\pi}(\rt_2\rs_2)^{-1}\cdot
\left\langle\varphi_1\otimes\varphi_2,
\vol^\natural(L)\cdot\lambda_E\left\langle\cZ^\pi_{\rt_1\rs_1\phi_1,L},\cZ_2\right\rangle_E\right\rangle_{\pi,\hat\pi}.
\end{align*}

By Proposition \ref{pr:generating},
\begin{align*}
\bbq_{r,r}(g_1\cdot\cZ^\pi_{\rt_1\rs_1\phi_1,L})&=\sum_{T_1\in\Herm_r^\circ(F)^+}Z_{T_1}^\clubsuit(\rt_1\rs_1g_1\phi_1)_L q^{T_1},\\
\bbq_{r,r}(g_2\cdot\cZ_2)&=p^{-M}\sum_{T_2\in\Herm_r^\circ(F)^+}\zeta_{g_2,T_2}
\end{align*}
hold for every pair $(g_1,g_2)\in M_r(F_\tR)\times M_r(F_\tR)$, where we recall that
\[
\zeta_{g_2,T_2}=p^M\lim\limits_{N\to\infty}Z_{T_2}^\clubsuit(\rt_2\rs_2g_2\phi_2^{[N!]})_L.
\]
Thus, we have, by Proposition \ref{pr:modularity2}(2),
\[
\vol^\natural(L)\cdot\lambda_E\left\langle\cZ^\pi_{\rt_1\rs_1\phi_1,L},\cZ_2\right\rangle_E=p^{-M}I_\lambda,
\]
and by Proposition \ref{pr:modularity1},
\begin{align*}
&\quad\lambda\langle\Theta_{\phi_1}^\sel(\varphi_1),\Theta_{\phi_2}^\sel(\varphi_2)\rangle_{\pi,F}^\natural \\
&=\chi^\lozenge_{\hat\pi}(\rt_1\rs_1)^{-1}\chi^\lozenge_{\pi}(\rt_2\rs_2)^{-1}\cdot
\left\langle\varphi_1\otimes\varphi_2,p^{-M}I_\lambda\right\rangle_{\pi,\hat\pi} \\
&=\partial_\lambda\sL_p^\lozenge(\pi)
\cdot\prod_{v\in\tV_F^{(\lozenge\setminus\{\infty,p\})}}Z(\varphi_{1,v}^\dag\otimes\varphi_{2,v},f_{\CF_v}).
\end{align*}
In other words,
\begin{align*}
\langle\Theta_{\phi_1}^\sel(\varphi_1),\Theta_{\phi_2}^\sel(\varphi_2)\rangle_{\pi,F}^\natural &=
\partial\sL_p^\lozenge(\pi)
\cdot\prod_{v\in\tV_F^{(\lozenge\setminus\{\infty,p\})}}Z(\varphi_{1,v}^\dag\otimes\varphi_{2,v},f_{\CF_v}) \\
&=\partial\sL_p^\lozenge(\pi)\cdot
\prod_{v\in\tV_F^{(\lozenge\setminus\{\infty,p\})}}Z(\varphi_{1,v}^\dag\otimes\varphi_{2,v},f_{\phi_{1,v}\otimes\phi_{2,v}}^{\r{SW}}).
\end{align*}
Finally, by Proposition \ref{pr:zeta} and (T2),
\[
Z(\varphi_{1,v}^\dag\otimes\varphi_{2,v},f_{\phi_{1,v}\otimes\phi_{2,v}}^{\r{SW}})
=\prod_{u\in\tP_v}\gamma(\tfrac{1+r}{2},\ul{\pi_u},\psi_{F,v})^{-1}
\]
for every $v\in\tV_F^{(p)}$. Together, we obtain \eqref{eq:aipf}. Part (1) is proved.

For (2), it suffices to show the vanishing under every embedding $\iota\colon\dL\to\dC$. Thus, we may regard $\dL$ as a subfield of $\dC$ and $\pi$ as defined over $\dC$. For every $\lambda\in\Hom_{\dZ_p}(\Gamma_{F,p},\dZ_p)$ with $\lambda_E\coloneqq\lambda\circ\Nm_{E/F}$, we have a map
\[
\bbi_{r,r}^\lambda\colon
\sS(V^{2r}\otimes_F\dA_F^\infty)=\sS(V^r\otimes_F\dA_F^\infty)\otimes_\dC\sS(V^r\otimes_F\dA_F^\infty)\to\c{SF}_{r,r}(\dC)
\]
(Definition \ref{de:sf1}) of $\dC[G_{r,r}(\dA_F^\infty)]$-modules sending $(\phi_1,\phi_2)$ to the assignment
\[
(g_1,g_2)\mapsto\sum_{(T_1,T_2)}\lambda_E\left\langle\wp_{\pi}\(Z_{T_1}^\clubsuit(\omega_r(g_1)\phi_1)_L\),
\wp_{\hat\pi}\(Z_{T_2}^\clubsuit(\omega_r(g_2)\phi_2)_L\)\right\rangle_E\cdot q^{T_1,T_2}
\]
in which the sum is taken over $\Herm_r(F)^+\times\Herm_r(F)^+$.

We prove (2) by contradiction. Assume the opposite hence $\bbi_{r,r}^\lambda$ is nontrivial for some $\lambda$. Then it is clear from the construction that $\bbi_{r,r}^\lambda$ factors through successive $G_{r,r}(\dA_F^\infty)$-equivariant quotient maps
\[
\sS(V^{2r}\otimes_F\dA_F^\infty)\to\rI_r^\Box(\CF)=\prod_{v\in\tV_F^{\infty}}\rI_{r,v}^\Box(\CF)
\to\pi\boxtimes\hat\pi.
\]
We claim that the image of $\bbi_{r,r}^\lambda$ is contained in $\bbq_{r,r}^\infty\cA_{r,r,\hol}^{[r]}$ (Definition \ref{de:sf1}). By \cite{LL2}*{Proposition~4.8(1)}, it suffices to show that $\bbi_{r,r}^\lambda(\phi_1,\phi_2)\in\bbq_{r,r}^\infty\cA_{r,r,\hol}^{[r]}$ for one choice of pair $(\phi_1,\phi_2)$ such that $\phi_1\otimes\phi_2$ has nonzero image under the unique nontrivial map in
$\Hom_{G_r(\dA_F^\infty)\times G_r(\dA_F^\infty)}\(\rI_r^\Box(\CF),\pi\boxtimes\hat\pi\)$. Indeed, we choose the pair to be $(\rt_1\rs_1\phi_1,\rt_2\rs_2\phi'_2)$ in which $\phi_1$ and $(\phi'_2)^p$ (away-from-$p$ part) are from the proof of (1), and $\phi'_{2,v}$ for $v\in\tV_F^{(p)}$ is an arbitrary element in $\sS(V_v^r)$ whose image in the quotient $\hat\pi_v\boxtimes\theta(\hat\pi_v)$ (Lemma \ref{le:theta}) is the limit of the images of $\phi_{v,2}^{[N!]}$ in that quotient when $N\to\infty$ (which exists by Lemma \ref{le:regular}(1)). By Lemma \ref{le:modularity1} (applied to both variables), it suffices to show that there exists $J_\lambda\in\cA_{r,r,\hol}^{[r]}$ such that
\begin{align*}
&\quad \bbq_{r,r}^\an((g_1,g_2)\cdot J_\lambda)\\
&=
\sum_{(T_1,T_2)}
\lambda_E\left\langle\wp_{\pi}\(Z_{T_1}^\clubsuit(\omega_r(g_1)\rt_1\rs_1\phi_1)_L\),
\wp_{\hat\pi}\(Z_{T_2}^\clubsuit(\omega_r(g_2)\rt_2\rs_2\phi'_2)_L\)\right\rangle_E\cdot q^{T_1,T_2} \\
&=
\sum_{(T_1,T_2)}
\lambda_E\left\langle Z_{T_1}^\clubsuit(\omega_r(g_1)\rt_1\rs_1\phi_1)_L,
Z_{T_2}^\clubsuit(\omega_r(g_2)\rt_2\rs_2\phi'_2)_L\right\rangle_E\cdot q^{T_1,T_2}
\end{align*}
for every pair $(g_1,g_2)\in M_r(F_\tR)\times M_r(F_\tR)$, in which the sums are taken over $\Herm_r(F)^+\times\Herm_r(F)^+$. Then by Proposition \ref{pr:modularity2}(2), we may take $J_\lambda$ to be $\vol^\natural(L)^{-1}p^{-M}\cdot I_\lambda$ (regarded as an element of $\cA_{r,r,\hol}^{[r]}$). It remains to show that $I_\lambda$ vanishes hence the map $\bbi_{r,r}^\lambda$ vanishes, resulting in a contradiction. Once again, since $\bbi_{r,r}^\lambda$ factors through $\pi\boxtimes\hat\pi$, it suffices to show that $\left\langle\varphi_1\otimes\varphi_2,I_\lambda\right\rangle_{\pi,\hat\pi}=0$ for a single (decomposable) pair $(\varphi_1,\varphi_2)$ such that $Z(\varphi_{1,v}^\dag\otimes\varphi_{2,v},f_{\CF_v})\neq 0$ for every $v\in\tV_F^\fin$. Indeed, we can just take $(\varphi_1,\varphi_2)$ to be the pair from the proof of (1). Then the vanishing of $\left\langle\varphi_1\otimes\varphi_2,I_\lambda\right\rangle_{\pi,\hat\pi}$ follows from Proposition \ref{pr:modularity1}, since we have assumed that the vanishing order of $\sL_p^\lozenge(\pi)$ at $\b1$ is at least one. Part (2) is proved.
\end{proof}

\subsection{Errata for [LL21, LL22]}

In this subsection, we correct two errors in the two preceding articles \cites{LL,LL2} in two remarks.

\begin{remark}\label{re:error}
In both \cite{LL} and \cite{LL2}, the authors mistakenly identified $\chi^\tR_{\pi^\vee}$ with $(\chi^\tR_\pi)^\tc$, where $\chi^\tR_\pi\colon\dT^\tR_{\dQ^\ac}\to\dQ^\ac$ is the Hecke character in \cite{LL}*{Definition~6.8} (and similarly for $\chi^\tR_{\pi^\vee}$); in fact, they only coincide when restricted to $\dT^\tR_{\dQ^\ac\cap\dR}$. As a consequence, one should replace $\chi^\tR_\pi(\rs)^\tc$ by $\chi^\tR_{\pi^\vee}(\rs)$ in \cite{LL}*{Proposition~6.10(1)}; and whenever one asks for two elements in $\dS_{\dQ^\ac}^\tR\setminus\fm_\pi^\tR$, they should actually be in $\dS_{\dQ^\ac}^\tR\setminus\fm_{\pi^\vee}^\tR$. Such modifications do not affect the proof of the results.
\end{remark}

\begin{remark}\label{re:error1}
The result in a recent preprint \cite{Sha} indicated that two formulae we claimed in the proof of \cite{LL}*{Lemma~4.4} regarding the Hecke correspondences of special cycles and Schwartz functions are both incorrect in general (though the statement of the lemma does hold). Here, we give a correct (and alternative) argument. For every moment matrix $T\in\Herm_m(F)^+$, by the same argument for \cite{Kud97}*{Proposition~5.10}, the assignment $\phi^\infty\mapsto Z_T(\phi^\infty)_{L'}$ is compatible with changing $L'$ that are neat open compact subgroups of $H(\dA_F^\infty)$ fixing $\phi$ under pullbacks, hence defines a functional
\[
\sS(V^m\otimes_{\dA_F}\dA_F^\infty)\to\varinjlim_{L'}\CH^m(X_{L'})_\dC.
\]
It is a map of smooth representations of $H(\dA_F^\infty)$, where $H(\dA_F^\infty)$ on the right-hand side by Hecke translations. Thus, for every neat open compact subgroup $L$ of $H(\dA_F^\infty)$, the above map restricts to a map
\[
\sS(V^m\otimes_{\dA_F}\dA_F^\infty)^L\to\(\varinjlim_{L'}\CH^m(X_{L'})_\dC\)^L
\]
of $\dC[L\backslash H(\dA_F^\infty)/L]$-modules. However, by definition, the above map factors through the submodule $\CH^m(X_L)_\dC$, which confirms \cite{LL}*{Lemma~4.4}.
\end{remark}

\appendix

\section{Bi-extensions and $p$-adic height pairings}
\label{ss:a}

In this appendix, we develop further the theory of $p$-adic heights on general varieties. We fix a prime number $p$ and an integer $n\geq 2$. Moreover, $\dW$ denotes a finite flat local extension of $\dZ_p$ and $\dL$ denotes a finite product of finite extensions of $\Frac(\dW)$.

\subsection{\'{E}tale correspondences}
\label{ss:etale}

Let $X$ be a scheme. An \emph{\'{e}tale correspondence} on $X$ is a diagram
\[
t\colon X \xleftarrow{f} X' \xrightarrow{g} X
\]
in which both $f$ and $g$ are finite \'{e}tale morphisms.

The collection of all \'{e}tale correspondences on $X$
forms a monoidal category $\text{\'{E}tCor}(X)$. See \cite{Liu19}*{Definition~2.11} for more details. We denote by $\r{EC}(X)$ the (unital) $\dZ$-algebra generated by the underlying monoid of isomorphism classes of objects of $\text{\'{E}tCor}(X)$. For every ring $R$, an \emph{$R$-ring of \'{e}tale correspondences} on $X$ is an $R$-ring $\dT$ together with a homomorphism $\dT\to\r{EC}(X)_R$ that is $R$-linear and unital. Usually, we only write $\dT$ for the notation when the homomorphism $\dT\to\r{EC}(X)_R$ is clear.

\begin{notation}\label{no:correspondence}
Let $S$ be a subset of $X$. For an \'{e}tale correspondence $t$ as above, we put $S^t\coloneqq f(g^{-1}(S))$. For a finite linear combination $t=\sum c_it_i$ with $c_i\neq 0$, we put $S^t\coloneqq\bigcup_i S^{t_i}$.
\end{notation}

\subsection{Remarks on sheaves}
\label{ss:sheaves}

Let $S$ be a site. For a diagram of rings $R_\bullet$, we denote by
\begin{itemize}
  \item $\bM(S,R_\bullet)$ the abelian category of sheaves of $R_\bullet$-modules on $S$,

  \item $\bC^+(S,R_\bullet)=\rC^+(\bM(S,R_\bullet))$ the abelian category of bounded below complexes in $\bM(S,R_\bullet)$,

  \item $\bD^+(S,R_\bullet)=\rD^+(\bM(S,R_\bullet))$ the derived category of $\bM(S,R_\bullet)$ with bounded below cohomology.
\end{itemize}
In this article, we mainly use two kinds of diagrams of rings. The first is a singleton valued in a ring $R$ so that $\bM(S,R_\bullet)=\bM(S,R)$ is the usual category of sheaves of $R$-modules on $S$. The second is the diagram
\[
\dW_\bullet\colon \cdots\to\dW/p^3 \to \dW/p^2 \to \dW/p
\]
so that $\bM(S,\dW_\bullet)$ consists of sequences $F_\bullet=(\cdots\to F_3\to F_2\to F_1)$ of sheaves of $\dW$-modules on $S$ in which $F_l$ is annihilated by $p^l$. Then we have the left-exact functor
\[
\varprojlim\colon\bM(S,\dW_\bullet)\to\bM(S,\dW)
\]
and the exact restriction functor
\[
-_l\colon\bM(S,\dW_\bullet)\to\bM(S,\dW/p^l)
\]
for every $l\geq 1$.

In what follows, we suppress $S$ (together with the comma after it) in the notation when $S$ is a point. We have the global section functor $\Gamma(S,-)\colon\bM(S,R_\bullet)\to\bM(R_\bullet)$. For every $q\in\dZ$, denote by $\bH^q(S,-)$ the $q$-th cohomology of $\rR\Gamma(S,-)$.

\begin{example}\label{ex:jannsen}
For a scheme $X$ and an integer $r$, we have the object $\mu_{p^\bullet}^{\otimes r}\in\bM(X_{\et},\dZ_{p\bullet})$ and put
\begin{align*}
\dL(r)_X\coloneqq\(\rR\varprojlim\(\mu_{p^\bullet}^{\otimes r}\)_X\)\otimes_{\dZ_p}\dL\in\bD^+(X_{\et},\dL).
\end{align*}
Then $\bH^q(X_{\et},\dL(r)_X)$ coincides with Jannsen's continuous \'{e}tale cohomology \cite{Jan88} of $X$ (with coefficients in $\dL$), usually denoted by $\rH^q(X,\dL(r))$. For a relative scheme $\pi\colon X\to S$ of finite type, $\bH^q(S_\et,\rR\pi_!\dL(r)_X)$ coincides with the continuous cohomology of $X/S$ with proper support, usually denoted by $\rH^q_\rc(X,\dL(r))$ when the base $S$ is clear.
\end{example}

Now let $\rG$ be a \emph{profinite} group that acts on $S$. We similarly define three categories $\bM_\rG(S,R_\bullet)$, $\bC^+_\rG(S,R_\bullet)$, and $\bD^+_\rG(S,R_\bullet)$ for compatible $\rG$-equivariant $R_\bullet$-modules on $S$. We have the global section functor $\Gamma(S,-)\colon\bM_\rG(S,R_\bullet)\to\bM_\rG(R_\bullet)$ and the $\rG$-invariants functor $\Gamma_\rG\colon\bM_\rG(R_\bullet)\to\bM(R_\bullet)$. For every $q\in\dZ$, denote by $\bH^q_\rG(S,-)$ the $q$-th cohomology of $\rR\Gamma_\rG\circ\rR\Gamma(S,-)$. The natural transformation from $\Gamma_\rG$ to the forgetful functor (forgetting the $\rG$-action) induces, for every $q\in\dZ$, a functor $\bH^q_\rG(S,-)\to\bH^q(S,-)$, whose kernel we denote by $\bH^q_\rG(S,-)^0$.

For an object $\sC\in\bD^+_\rG(S,\dW_\bullet)$, we put
\begin{align}\label{eq:adic}
\sC_\dL\coloneqq\(\rR\varprojlim\sC\)\otimes_{\dW}\dL\in\bD^+_\rG(S,\dL).
\end{align}

\begin{definition}\label{de:admissible}
We say that an object $\sC\in\bD^+_\rG(S,\dW_\bullet)$ is \emph{admissible} if
\begin{enumerate}
  \item $\bH^q(S,\sC_l)$ is a discrete $\rG$-module for every $q\in\dZ$ and $l\geq 1$, that is, every element of $\bH^q(S,\sC_l)$ has an open stabilizer;

  \item $\rR^1\varprojlim_l\bH^q(S,\sC_l)$ has finite exponent for every $q\in\dZ$.
\end{enumerate}
\end{definition}

It is clear that if $\sC$ is admissible, then the natural map
\[
\bH^q(S,\sC_\dL)\to\(\varprojlim_l\bH^q(S,\sC_l)\)\otimes_{\dW}\dL
\]
is an isomorphism of $\rG$-modules, through which we view $\bH^q(S,\sC_\dL)$ as a topological $\rG$-module. The following lemma slightly generalizes \cite{Jan88}*{Corollary~3.4} in the case of rational coefficients.

\begin{lem}\label{le:spectral}
Let $\sC\in\bD^+_\rG(S,\dZ_{p\bullet})$ be an admissible object. Then there is a spectral sequence
\[
\rH^p_\cont(\rG,\bH^q(S,\sC_\dL))\Rightarrow\bH^{p+q}_\rG(S,\sC_\dL).
\]
In particular, we have the edge map
\[
\bH^q_\rG(S,\sC_\dL)^0\to\rH^1_\cont(\rG,\bH^{q-1}(S,\sC_\dL)).
\]
\end{lem}

\begin{proof}
We first note that $(\Gamma_\rG\circ\Gamma(S,-))\circ\varprojlim=(\Gamma_\rG\circ\varprojlim)\circ\Gamma(S,-)$ and all of the three functors preserve injectives. By the same argument for \cite{Jan88}*{Theorem~3.3}, there is a spectral sequence
\[
\rH^p(\rG,(\bH^q(S,\sC_l)))\Rightarrow\bH^{p+q}_\rG(S,\rR\varprojlim\sC)
\]
which is simply the Grothendieck spectral sequence for the composition $\rR(\Gamma_\rG\circ\varprojlim)\circ\rR\Gamma(S,-)$. Here, we recall that $\sC_l\in\bD^+_\rG(S,\dZ_p/p^l)$ is the restriction of $\sC$. By the similar argument for \cite{Jan88}*{Theorem~2.2} and Definition \ref{de:admissible}, there is a canonical isomorphism
\[
\rH^p(\rG,(\bH^q(S,\sC_l)))\otimes_{\dW}\dL\simeq\rH^p_\cont(\rG,\bH^q(S,\sC_\dL)).
\]
The lemma then follows as $\bH^{p+q}_\rG(S,\rR\varprojlim\sC)\otimes_{\dW}\dL=\bH^{p+q}_\rG(S,\sC_\dL)$.
\end{proof}

The lemma below will only be used in \S\ref{ss:alteration}.

\begin{lem}\label{le:hochschild}
Let $S$ be a site with an action of $\rG$ and $q$ an integer. Consider a distinguished triangle
\[
\sA\xrightarrow{\mu}\sB\xrightarrow{\nu}\sC\xrightarrow{+1}
\]
of admissible objects in $\bD^+_\rG(S,\dW_\bullet)$, inducing the following commutative diagram
\[
\xymatrix{
\bH^{q-1}_\rG(S,\sC_\dL) \ar[r]^-{\lambda}\ar[d] & \bH^q_\rG(S,\sA_\dL) \ar[r]^-{\mu}\ar[d]^-{\alpha} & \bH^q_\rG(S,\sB_\dL) \ar[r]^-{\nu}\ar[d]^-{\beta} &  \bH^q_\rG(S,\sC_\dL) \ar[d]^-{\gamma} \\
\bH^{q-1}(S,\sC_\dL) \ar[r]^-{\ol{\lambda}} & \bH^q(S,\sA_\dL) \ar[r]^-{\ol{\mu}} & \bH^q(S,\sB_\dL) \ar[r]^-{\ol{\nu}} &  \bH^q(S,\sC_\dL)
}
\]
in which all vertical maps are induced from forgetting $\rG$-actions. Then for every element $b\in\bH^q_\rG(S,\sB_\dL)$ satisfying $\gamma(\nu(b))=0$, the image of $\nu(b)$ under the composite map
\[
\bH^q_\rG(S,\sC_\dL)^0\to\rH^1_\cont(\rG,\bH^{q-1}(S,\sC_\dL))\to
\rH^1_\cont\(\rG,\frac{\bH^{q-1}(S,\sC_\dL)}{\bH^{q-1}(S,\sB_\dL)}\)
\]
can be represented by the (continuous) $1$-cocycle $g\mapsto g\ol{a}-\ol{a}$ for $g\in\rG$, where $\ol{a}$ is an arbitrary element in $\bH^q(S,\sA_\dL)$ satisfying $\ol{\mu}(\ol{a})=\beta(b)$.
\end{lem}

\begin{proof}
It is easy to check that the $1$-cocycle does not depend on the choice of $\ol{x}$. Thus, it suffices to check the statement for one such element.

Take an injective resolution $0\to A_\bullet^\bullet\xrightarrow{\mu^\bullet}B_\bullet^\bullet\xrightarrow{\nu^\bullet}C_\bullet^\bullet\to 0$ in $\bC^+_\rG(\dW_\bullet)$ (the bullet in the superscript denotes the cohomological degree) of the exact triangle
\[
\rR\Gamma(S,\sA)\to\rR\Gamma(S,\sB)\to\rR\Gamma(S,\sC)\xrightarrow{+1}.
\]
Put $X_\dL^\bullet\coloneqq\(\varprojlim_lX_l^\bullet\)\otimes_{\dW}\dL\in\bC^+_\rG(\dL)$ for $X=A,B,C$. Then the diagram in the statement can be replaced by
\[
\xymatrix{
\rH^{q-1}\((C^\bullet_\dL)^\rG\) \ar[r]^-{\lambda}\ar[d] & \rH^q\((A^\bullet_\dL)^\rG\) \ar[r]^-{\mu}\ar[d]^-{\alpha} & \rH^q\((B^\bullet_\dL)^\rG\) \ar[r]^-{\nu}\ar[d]^-{\beta} &  \rH^q\((C^\bullet_\dL)^\rG\) \ar[d]^-{\gamma} \\
\rH^{q-1}\(C^\bullet_\dL\) \ar[r]^-{\ol{\lambda}} & \rH^q\(A^\bullet_\dL\) \ar[r]^-{\ol{\mu}} & \rH^q\(B^\bullet_\dL\) \ar[r]^-{\ol{\nu}} &  \rH^q\(C^\bullet_\dL\)
}
\]
in which all vertical arrows are induced by natural inclusions. By definition, we have
\[
\bH^q_\rG(S,\sC_\dL)^0=\rH^q\((C^\bullet_\dL)^\rG\)^0=
\frac{\Ker\((E^q_\dL)^\rG\xrightarrow{\rd}(C^{q+1}_\dL)^\rG\)\cap\IM\(C^{q-1}_\dL\xrightarrow{\rd}E^q_\dL\)}
{\IM\((C^{q-1}_\dL)^\rG\xrightarrow{\rd}(E^q_\dL)^\rG\)}.
\]
It follows from the formation of the spectral sequence in Lemma \ref{le:spectral} that for $c\in\rH^q\((C^\bullet_\dL)^\rG\)^0$, its image in $\rH^1_\cont(\rG,\rH^{q-1}(C^\bullet_\dL))$ under the edge map can be represented by the $1$-cocycle that sends $g\in\rG$ to the (cohomology class in $\rH^{q-1}(C^\bullet_\dL)$) of $g c^{q-1}-c^{q-1}$, where $c^{q-1}\in C^{q-1}_\dL$ is an arbitrary element whose differential $\r{d}c^{q-1}$ in $C^q_\dL$ represents $c$ (so that $g c^{q-1}-c^{q-1}$ is closed).

To prove the lemma, let $b^q\in(B^q_\dL)^\rG$ be a closed element that represents $b\in\rH^q\((B^\bullet_\dL)^\rG\)$. Since $\nu^q(y^q)$ induces an element in $\rH^q\((C^\bullet_\dL)^\rG\)^0$, we may choose an element $c^{q-1}\in C^{q-1}_\dL$ whose differential $\r{d}c^{q-1}$ in $C^q_\dL$ equals $\nu^q(b^q)$. To construct an element $\ol{a}$ that lifts $\beta(b)$, we take an element $b^{q-1}\in B^{q-1}_\dL$ such that $\mu^{q-1}(b^{q-1})=c^{q-1}$. Then $\nu^q(b^q-\r{d}b^{q-1})=0$, which implies that $a^q\coloneqq b^q-\r{d}b^{q-1}$ belongs to $A^q_\dL$ and is closed. Then we may take $\ol{a}\in\rH^q\(A^\bullet_\dL\)$ to be the cohomology class of $a^q$. Since $b^q$ is fixed by $\rG$, we have for every $g\in\rG$, $g a^q-a^q=g\r{d}b^{q-1}-\r{d}b^{q-1}=\r{d}(g b^{q-1}-b^{q-1})$, which is a closed element in $A^q_\dL$. However, $\r{d}(g b^{q-1}-b^{q-1})$ exactly represents the class of $g c^{q-1}-c^{q-1}$ under the coboundary map $\ol{\lambda}\colon\rH^{q-1}\(C^\bullet_\dL\)\to\rH^q\(A^\bullet_\dL\)$. The lemma is proved.
\end{proof}

In what follows, $\rG$ will be the absolute Galois group of a field $K$ of characteristic different from $p$, which we denote by $\rG_K$. As for common practice, we simply write $\rH^q(K,-)$ for $\rH^q_\cont(\rG_K,-)$. We also introduce the additive category $\bM_K(\dL)$ of topological $\dL$-vector spaces with continuous actions by $\rG_K$.

\subsection{Bi-extensions of cycles}
\label{ss:biextension}

Let $K$ be a field of characteristic different from $p$ with a fixed algebraic closure $\ol{K}$. Let $X$ be a projective smooth scheme over $K$ of pure dimension $n-1$. For every integer $d$, put
\[
\rZ^d(X)_\dL^0\coloneqq\Ker\(\rZ^d(X)_\dL\to\rH^{2d}(X_{\ol{K}},\dL(d))\).
\]

Now we consider two elements $c\in\rZ^d(X)_\dL^0$ and $c'\in\rZ^{d'}(X)_\dL^0$ with $d+d'=n$, such that $c$ and $c'$ have disjoint supports. Choose disjoint nonempty closed subsets $Z$ and $Z'$ of $X$ of pure codimension $d$ and $d'$ containing the supports of $c$ and $c'$, respectively. Denote by ${i}\colon Z\to X$ and ${i}'\colon Z'\to X$ the closed immersions and put $U\coloneqq X\setminus Z$ and $U'\coloneqq X\setminus Z'$. We have the following diagram
\[
\xymatrix{
U\cap U' \ar[r]^-{\jmath'}\ar[d]_-{\jmath} & U \ar[d]^-{{j}} \\
U' \ar[r]^-{{j}'} & X
}
\]
of open immersions. We have the following induced diagram
\begin{align}\label{eq:biextension}
\resizebox{\hsize}{!}{
\xymatrix{
& & 0 \ar[d] & 0 \ar[d] & \\
\rH^{2d-2}(X_{\ol{K}},\dL(d)) \ar[r]\ar[d]^-{\sim} & \rH^{2d-2}(Z'_{\ol{K}},\dL(d)) \ar[r]\ar@{=}[d]
& \rH^{2d-1}(X_{\ol{K}},{j}'_!\dL(d)) \ar[r]\ar[d] & \rH^{2d-1}(X_{\ol{K}},\dL(d)) \ar[r]\ar[d] & 0 \\
\rH^{2d-2}(U_{\ol{K}},\dL(d)) \ar[r] & \rH^{2d-2}(Z'_{\ol{K}},\dL(d)) \ar[r]
& \rH^{2d-1}(U_{\ol{K}},\jmath'_!\dL(d)) \ar[r]\ar[d] & \rH^{2d-1}(U_{\ol{K}},\dL(d)) \ar[r]\ar[d] & 0 \\
& & \rH^{2d}(Z_{\ol{K}},{i}^!\dL(d)) \ar@{=}[r]\ar[d] & \rH^{2d}(Z_{\ol{K}},{i}^!\dL(d)) \ar[d] \\
& & \rH^{2d}(X_{\ol{K}},{j}'_!\dL(d)) \ar[r]^-{\sim} & \rH^{2d}(X_{\ol{K}},\dL(d))
}
}
\end{align}
in $\bM_K(\dL)$.

The element $c$ gives rise to a map $\kappa^c\colon\dL\to\rH^{2d}(Z_{\ol{K}},{i}^!\dL(d))$ whose image is contained in the kernel of the map $\rH^{2d}(Z_{\ol{K}},{i}^!\dL(d))\to\rH^{2d}(X_{\ol{K}},\dL(d))$. The element $c'$ gives rise to a map $\kappa_{c'}\colon\rH^{2d-2}(Z'_{\ol{K}},\dL(d-1))\to\dL$ that vanishes on the image of the map $\rH^{2d-2}(X_{\ol{K}},\dL(d))\to\rH^{2d-2}(Z'_{\ol{K}},\dL(d))$. Applying the pullback along $\kappa^c$ and the pushforward along $\kappa_{c'}(1)$ to the diagram \eqref{eq:biextension}, we obtain the following bi-extension diagram
\begin{align*}
\xymatrix{
& & 0 \ar[d] & 0 \ar[d] & \\
0 \ar[r] & \dL(1) \ar[r]\ar@{=}[d]
& \rE_{c'} \ar[r]\ar[d] & \rH^{2d-1}(X_{\ol{K}},\dL(d)) \ar[r]\ar[d] & 0 \\
0 \ar[r] & \dL(1) \ar[r]
& \rE^c_{c'} \ar[r]\ar[d] & \rE^c \ar[r]\ar[d] & 0 \\
& & \dL \ar@{=}[r]\ar[d] & \dL \ar[d] \\
& & 0  & 0
}
\end{align*}
in $\bM_K(\dL)$. It is easy to see that the above diagram does not depend on the choices of $Z$ and $Z'$. We have three induced extension classes
\begin{itemize}
  \item $[\rE^c]\in\rH^1(K,\rH^{2d-1}(X_{\ol{K}},\dL(d)))$,

  \item $[\rE_{c'}^\vee(1)]=[\rE^{c'}]\in\rH^1(K,\rH^{2d'-1}(X_{\ol{K}},\dL(d')))$,

  \item $[\rE^c_{c'}]\in\rH^1(K,\rE_{c'})$.
\end{itemize}

\subsection{Relation with Beilinson's local index}
\label{ss:beilinson}

In this subsection, we assume that $K$ is a non-archimedean local field whose residue field has characteristic different from $p$. We study the relation between Nekov\'{a}\v{r}'s local $p$-adic height and Beilinson's local index.

Assume that the cycle classes of $c$ and $c'$ in $\rH^{2d}(X,\dL(d))$ and $\rH^{2d'}(X,\dL(d'))$ vanish, respectively.\footnote{This is automatic if the monodromy-weight conjecture holds for $\rH^{2d-1}(X_{\ol{K}},\dL(d))$.} Then the image of $[\rE^c_{c'}]$ in $\rH^1(K,\rH^{2d-1}(X_{\ol{K}},\dL(d)))$ vanishes, hence $[\rE^c_{c'}]$ belongs to the image of the map $\rH^1(K,\dL(1))\to\rH^1(K,\rE_{c'})$ which is injective. Following Nekov\'{a}\v{r} \cite{Nek93}, we denote by $\langle c,c'\rangle_{X,K}^\rN$ the image of $[\rE^c_{c'}]$ under the natural isomorphism $\rH^1(K,\dL(1))=\widehat{K^\times}\otimes_{\widehat\dZ}\dL$ given by the Kummer maps.

We recall the definition of Beilinson's local index \cite{Bei87}*{Section~2} (see also \cite{LL}*{Appendix~B}). We have the refined cycle class $[c]\in\rH^{2d}_Z(X,\dL(d))=\rH^{2d}(Z,\jmath^!\dL(d))$, which is contained in the kernel of the map $\rH^{2d}_Z(X,\dL(d))\to\rH^{2d}(X,\dL(d))$, hence we may choose an element $\gamma\in\rH^{2d-1}(U,\dL(d))$ that maps to $[c]$ under the coboundary map $\rH^{2d-1}(U,\dL(d))\to\rH^{2d}_Z(X,\dL(d))$. Similarly, we can choose an element $\gamma'\in\rH^{2d'-1}(U',\dL(d'))$ for $c'$. Beilinson's local index, which we denote by $\langle c,c'\rangle_{X,K}^\rB$, is defined as the image of $\gamma\cup\gamma'$ under the composite map
\begin{align*}
\rH^{2n-2}(U\cap U',\dL(n))&\to\rH^{2n-1}(X,\dL(n)) \\
&\xrightarrow{\Tr_{X/K}}\rH^1(\Spec K,\dL(1))=
\rH^1(K,\dL(1))=\widehat{K^\times}\otimes_{\widehat\dZ}\dL,
\end{align*}
in which the first map is the coboundary map in the Mayer--Vietoris exact sequence for the open covering $X=U\cup U'$.

\begin{remark}\label{re:beilinson}
In fact, in \cite{Bei87}*{Section~2} and \cite{LL}*{Appendix~B}, the local index $\langle c,c'\rangle_K^\rB$ takes value in $\dL$ via the canonical isomorphism $\rH^1(\Spec K,\dQ_p(1))\simeq\dQ_p$ that is the composition
\[
\rH^1(\Spec K,\dQ_p(1))\to\rH^2_{\Spec\kappa}(\Spec O_K,\dQ_p(1))\xrightarrow\sim
\rH^0(\Spec\kappa,\dQ_p)\simeq\dQ_p
\]
in which $\kappa$ is the residue field of $K$. By \cite{GD77}*{2.1.3} or \cite{Nek95}*{II.(2.16.1)}, the induced isomorphism $\widehat{K^\times}\otimes_{\widehat\dZ}\dQ_p\to\dQ_p$ sends a uniformizer of $K$ to $-1$, rather than $1$.
\end{remark}

The following proposition is simply \cite{Sch94}*{Theorem~5.3}.

\begin{proposition}\label{th:index}
Let the situation be as above. Then
\[
\langle c,c'\rangle_{X,K}^\rN=\langle c,c'\rangle_{X,K}^\rB.
\]
\end{proposition}

\subsection{Crystalline property of bi-extensions}
\label{ss:crystalline}

In this subsection, we assume that $K$ is a finite extension of $\dQ_p$ with residue field $\kappa$. Denote by $W$ the Witt ring of $\kappa$, and by $K_0$ the fraction field of $W$, which is canonically a subfield of $K$.

We assume that $X$ admits a proper strictly semistable model $\cX$ over $O_K$. Put $\sfX\coloneqq\cX\otimes_{O_K}\kappa$. For every integer $h\geq 1$, denote by $\sfX^{(h)}$ the disjoint union of intersections of $h$ different irreducible components of $\sfX$, which is either empty or a proper smooth scheme over $\kappa$.

\begin{theorem}\label{th:crystalline}
Suppose that $n<p$. Let $\dT$ be an $\dL$-ring of \'{e}tale correspondences on $\cX$ and $\fm,\fm'$ two maximal ideals of $\dT$ satisfying that
\begin{align}\label{eq:vanishing}
\bigoplus_{h>1,q\geq 0}\(\rH^q_\cris(\sfX^{(h)}/W)\otimes_{\dZ_p}\dL\)_\fm=
\bigoplus_{h>1,q\geq 0}\(\rH^q_\cris(\sfX^{(h)}/W)\otimes_{\dZ_p}\dL\)_{\fm'}=0.
\end{align}
Then there exist elements $t\in\dT\setminus\fm$ and $t'\in\dT\setminus\fm'$ depending only on $\cX$ such that the following holds: For two arbitrary elements $c\in\rZ^d(X)_\dL^0$ and $c'\in\rZ^{d'}(X)_\dL^0$ with $d+d'=n$ satisfying that
\begin{enumerate}
  \item $(\supp\cC)^t\cap(\supp\cC')^{t'}=\emptyset$ for every $t,t'\in\dT$, where $\cC$ and $\cC'$ denote the Zariski closures of $c$ and $c'$ in $\cX$, respectively;

  \item the codimension of $\supp\cC'$ in $\sfX^{(h)}$ is at least $d'$ for every $h\geq 1$,
\end{enumerate}
the following properties
\begin{itemize}
  \item $[\rE^{t^*c}]\in\rH^1_f(K,\rH^{2d-1}(X_{\ol{K}},\dL(d)))$,

  \item $[\rE_{t'^*c'}^\vee(1)]\in\rH^1_f(K,\rH^{2d'-1}(X_{\ol{K}},\dL(d')))$,

  \item $[\rE^{t^*c}_{t'^*c'}]\in\rH^1_f(K,\rE_{t'^*c'})$,
\end{itemize}
hold simultaneously. Here, the bi-extension $\rE^{t^*c}_{t'^*c'}$ exists by (1).
\end{theorem}

\begin{remark}
By taking $\dT=\dL$, Theorem \ref{th:crystalline} asserts that $\rE^c_{c'}$ is crystalline as long as $\supp\cC\cap\supp\cC'=\emptyset$ when $\cX$ is a proper smooth model of $X$ over $O_K$.\footnote{We warn the readers that this assertion is wrong if one replaces the word \emph{smooth} by \emph{strictly semistable}.} This confirms the (equivalent) conjecture in the remark after \cite{Shn16}*{Theorem~8.7} when $n<p$.\footnote{However, our strategy for the proof of Theorem \ref{th:crystalline} is different from the case of (local systems over) curves in \cite{Shn16}*{Theorem~8.7}.}
\end{remark}

Our main strategy of proving Theorem \ref{th:crystalline} is similar to \cite{Sat13}. The main difficulty is to show that the bi-extension $[\rE^{t^*c}_{t'^*c'}]$ is a crystalline class. We consider the Abel--Jacobi map from the homologically trivial part of the degree $2d$ syntomic cohomology of $\cX\setminus(\supp\cC')$ with proper support to $\rH^1(K,\rH^{2d-1}_\rc(U'_{\ol{K}},\dQ_p(d)))$. After we show that $\rE_{t'^*c'}$ is crystalline for suitable $t'$, it suffices to show that if a syntomic class comes from a (homologically trivial) cycle, then its Abel--Jacobi image vanishes in $\rH^1(K,\rH^{2d-1}_\rc(U'_{\ol{K}},\dQ_p)\otimes_{\dQ_p}\dB_\cris)$ ``after localization at $\fm'$'' (and replace $K$ by a finite extension in fact), under the conditions in Theorem \ref{th:crystalline}. However, the main challenge for us is that unlike the situation in \cite{Sat13} where $U'=X$, we do not have a comparison theorem for $\rH^{2d-1}_\rc(U'_{\ol{K}},\dQ_p)$ with some cohomology on the special fiber in general. Also, due to the constraint of the conditions in the theorem, we can not reduce the theorem to an alteration. We solve this problem in the following way. First, we show that the kernel of the Abel--Jacobi map from the above syntomic cohomology to $\rH^1(K,\rH^{2d-1}_\rc(U'_{\ol{K}},\dQ_p)\otimes_{\dQ_p}\dB_\cris)$ contains the kernel of another map whose range is a certain space defined by log rigid cohomology (Proposition \ref{pr:comparison}). Second, we show that a cycle class will vanish in this space ``after localization at $\fm'$''. For the first step, which shall hold more generally without the conditions in the theorem, we pass to a strict semistable alteration.

The full proof of Theorem \ref{th:crystalline} will occupy the entire Appendix \ref{ss:b}.

\subsection{Recollection on $p$-adic Galois representation}
\label{ss:recollection}

Let $K$ be as in the previous subsection. Let $\rV$ be a finite-dimensional continuous representation of $\rG_K$ with coefficients in $\dL$, which we assume to be \emph{de Rham}, that is
\[
\dD_{\dr}(\rV)\coloneqq\(\rV\otimes_{\dQ_p}\dB_{\dr}\)^{\rG_K}
\]
is a free $\dL\otimes_{\dQ_p}K$-module of rank $\dim_\dL\rV$. In the rest of this subsection, we assume that $\dL$ is a subfield of $\ol\dQ_p$; and one can generalize the discussion to a finite product of finite extensions of $\dQ_p$ by considering all homomorphisms $\dL\to\ol\dQ_p$ over $\dQ_p$.

We have a decreasing filtration $\rF^i\dD_{\dr}(\rV)$ of $\dL\otimes_{\dQ_p}K$-submodules of $\dD_{\dr}(\rV)$, known as the de Rham filtration. Moreover, for every embedding $\tau\colon K\to\ol\dQ_p$, Fontaine constructed a Weil--Deligne representation $\r{WD}(\rV)_\tau$ of the Weil group of $K$ with coefficients in $\ol\dQ_p$, with underlying $\ol\dQ_p$-vector space $\dD_{\dr}(\rV)\otimes_{\dL\otimes_{\dQ_p}K,1\otimes\tau}\ol\dQ_p$. The isomorphism class of $\r{WD}(\rV)_\tau$ is independent of $\tau$. See, for example, \cite{TY07}*{Section~1} for more details. We also recall that $\rV$ is crystalline if
\[
\dD_\cris(\rV)\coloneqq\(\rV\otimes_{\dQ_p}\dB_{\cris}\)^{\rG_K}
\]
is a free $\dL\otimes_{\dQ_p}K_0$-module of rank $\dim_\dL\rV$.

\begin{remark}\label{rem:fr}
It follows from the construction of $\r{WD}(\rV)_\tau$ that if $\rV$ is crystalline, then the following polynomials coincide:
\begin{itemize}
  \item the characteristic polynomial of a geometric Frobenius on $\r{WD}(\rV)_\tau$, for any $\tau$;

  \item the characteristic polynomial of $\varphi$ on $\dD_\cris(\rV)$, where $\varphi$ denotes the $[K_{0}:\dQ_p]$-th power of the crystalline Frobenius.
\end{itemize}
\end{remark}

\begin{definition}\label{de:weight}
Let $\mu$ be a real number. We say that $\rV$ is \emph{pure of weight $\mu$} if for some (hence every) $\tau\colon K\to\ol\dQ_p$, all geometric Frobenius eigenvalues of $\gr_i\r{WD}(\rV)_\tau$ are Weil $|\kappa|^{\mu+i}$-numbers for every $i\in\dZ$, where $\gr_i\r{WD}(\rV)_\tau$ denotes the $i$-th graded piece of the monodromy filtration on $\r{WD}(\rV)_\tau$.\footnote{In particular, $\dL(1)$ is pure of weight $-2$.}
\end{definition}

We make the following definition, after \cite{Nek93}*{6.7}.

\begin{definition}\label{de:panchishkin}
We say that $\rV$ satisfies the \emph{Panchishkin condition} if there exists a necessarily unique $\dL[\rG_K]$-submodule $\rV^+\subseteq\rV$ (with $\rV^-\coloneqq\rV/\rV^+$) such that
\[
\rF^{0}\dD_{\dr}(\rV^+)= \dD_{\dr}(\rV^-)/\rF^{0} \dD_{\dr}(\rV^-)=0.
\]
\end{definition}

\begin{lem}\label{le:pan-eq}
For a crystalline representation $\rV$ of $\rG_K$ over $\dL$, the following are equivalent:
\begin{enumerate}
  \item $\rV$ satisfies the Panchishkin condition;

  \item the $\dL\otimes_{\dQ_p}K_0$-submodule $\dD_\cris^+(\rV)\subset\dD_\cris(\rV)$ on which the crystalline Frobenius acts with negative slopes\footnote{In other words, the submodule $\dD^+_\cris(\rV)$ generated by the generalized $\varphi$-eigenspaces relative to the eigenvalues of negative valuation.} is weakly admissible, and the natural map
      \begin{align}\label{eq:pan-ses}
      (\dD^+_\cris(\rV)\otimes_{K_0}K)\oplus\rF^0\dD_{\dr}(V)\xrightarrow{\sim}\dD_\dr(\rV)
      \end{align}
      is a splitting of the Hodge filtration on $\dD_\dr(\rV)$.
\end{enumerate}
\end{lem}

\begin{proof}
Assume that $\rV$ is Panchishkin with a subrepresentation $\rV^{+}$ as in the definition. By the weak admissibility of $\dD_{\cris}(\rV^{\pm})$, the crystalline Frobenius acts on $\dD_{\cris}(\rV^{+})$ with negative slopes and on $\dD_{\cris}(\rV^{-})$ with non-negative slopes; it follows that $\dD_{\cris}^{+}(\rV)=\dD_{\cris}(\rV^{+})$ satisfies the second condition. Conversely, that condition implies that $\dD^{+}_{\cris}(\rV)\subset \dD_{\cris}(\rV)$ is weakly admissible, hence by \cite{CF00} it arises as $\dD^{+}_{\cris}(\rV)=\dD_{\cris}(\rV^{+})$ from a subrepresentation $\rV^{+}\subset\rV$, which witnesses the Panchishkin condition.
\end{proof}

\begin{lem}\label{le:panchishkin}
Let $\mu$ be a nonzero integer and $0\to\rV_1\to\rV\to\rV_2\to 0$ a short exact sequence of crystalline representations of $\rG_K$. If $\rV_1$ and $\rV_2$ are pure of weight~$\mu$, then so is $\rV$. If $\rV_1$ and $\rV_2$ satisfy the Panchishkin condition, then so does $\rV$.
\end{lem}

\begin{proof}
The first statement is clear. We prove the second one. The sequence
\[
0 \to \rF^0\dD_{\dr}(\rV_1) \to \rF^0\dD_{\dr}(\rV) \to \rF^0\dD_{\dr}(\rV_2) \to 0
\]
is exact, and by definition, so is
\[
0\to\dD_{\cris}^+(\rV_1)\to\dD_{\cris}^+(\rV)\to\dD_{\cris}^+(\rV_2)\to 0.
\]
This implies that \eqref{eq:pan-ses} is an isomorphism. By \cite{CF00}*{Proposition~3.4}, the module $\dD_{\cris}^+(\rV)$ is weakly admissible too, so that the equivalent Panchishkin condition of Lemma \ref{le:pan-eq} is satisfied.
\end{proof}

\subsection{Decomposition of $p$-adic height pairing}
\label{ss:decomposition}

In this subsection, we take $K$ to be a number field. Let $X$ be a proper smooth scheme over $K$ of pure dimension $n-1$ and take two positive integers $d,d'$ satisfying $d+d'=n$.

Consider $\dL[\rG_K]$-submodules $\rV$ and $\rV'$ of $\rH^{2d-1}(X_{\ol{K}},\dL(d))$ and $\rH^{2d'-1}(X_{\ol{K}},\dL(d'))$, respectively, satisfying
\begin{enumerate}[label=(V\arabic*)]
  \item For every nonarchimedean place $u$ of $K$ not above $p$, $\rH^i(K_u,\rV)=\rH^i(K_u,\rV')=0$ for $i\in\dZ$.

  \item For every place $u$ of $K$ above $p$, both $\rV\res_{K_u}$ and $\rV'\res_{K_u}$ are semistable and pure of weight $-1$ (Definition \ref{de:weight}).

  \item For every place $u$ of $K$ above $p$, both $\rV\res_{K_u}$ and $\rV'\res_{K_u}$ satisfy the Panchishkin condition (Definition \ref{de:panchishkin}).
\end{enumerate}
We have the canonical $p$-adic height pairing
\[
\langle\;,\;\rangle_{(\rV,\rV'),K}\colon\rH^1_f(K,\rV)\times\rH^1_f(K,\rV')\to\Gamma_{K,p}\otimes_{\dZ_p}\dL
\]
constructed in \cite{Nek93}, using the Hodge splitting map \eqref{eq:pan-ses}. The pairing is $\dL$-bilinear.

Recall that we have the Abel--Jacobi map
\[
\AJ\colon\rZ^d(X)^0_\dL\to\rH^1(K,\rH^{2d-1}(X_{\ol{K}},\dL(d))).
\]
We denote by $\rZ^d_\rV(X)^0_\dL$ the subspace of $\rZ^d(X)^0_\dL$ that is the inverse image of $\rH^1(K,\rV)$. By \cite{Nek00}*{Theorem~3.1}, the image of $\AJ$ is contained in $\rH^1_{\r{st}}(K,\rH^{2d-1}(X_{\ol{K}},\dL(d)))$. Moreover, (V2) implies that $\rH^1_{\r{st}}(K,\rV)=\rH^1_f(K,\rV)$. Thus, we have the Abel--Jacobi map
\[
\AJ\colon\rZ^d_\rV(X)^0_\dL\to\rH^1_f(K,\rV).
\]
Similarly, we have $\rZ^{d'}_{\rV'}(X)^0_\dL$ and the corresponding Abel--Jacobi map
\[
\AJ\colon\rZ^{d'}_{\rV'}(X)^0_\dL\to\rH^1_f(K,\rV').
\]
Combining the $p$-adic height pairing with the two Abel--Jacobi maps, we obtain a pairing
\begin{align}\label{eq:decomposition}
\langle\;,\;\rangle_{(\rV,\rV'),K}\colon\rZ^d_\rV(X)^0_\dL\times\rZ^{d'}_{\rV'}(X)^0_\dL\to\Gamma_{K,p}\otimes_{\dZ_p}\dL.
\end{align}

Take two elements $c\in\rZ^d_\rV(X)_\dL^0$ and $c'\in\rZ^{d'}_{\rV'}(X)_\dL^0$ with disjoint supports. Then according to \cite{Nek93}*{Section~4}, we have a decomposition
\begin{align*}
\langle c,c'\rangle_{(\rV,\rV'),K}=\sum_{u\nmid\infty}\langle c,c'\rangle_{(\rV,\rV'),K_u}
\end{align*}
of the pairing \eqref{eq:decomposition} into local ones $\langle c,c'\rangle_{(\rV,\rV'),K_u}\in \widehat{K_u^\times}\otimes_{\widehat\dZ}\dL$ over all \emph{nonarchimedean} places $u$ of $K$, in which $\langle c,c'\rangle_{(\rV,\rV'),K_u}=\langle c,c'\rangle_{X_u,K_u}^\rN$ for $u$ not above $p$.

\begin{remark}\label{re:height}
For a place $u$ of $K$ above $p$, the bi-extension class $[\rE^c_{c'}]\in\rH^1(K_u,\rE_{c'})$ belongs to $\rH^1_f(K_u,\rE_{c'})$ if and only if $\langle c,c'\rangle_{(\rV,\rV'),K_u}\in O_{K_u}^\times\otimes_{\dZ_p}\dL$.
\end{remark}

\section{Proof of Theorem \ref{th:crystalline}}
\label{ss:b}

Let $K$ be a finite extension of $\dQ_p$ with residue field $\kappa$. Denote by $W$ the Witt ring of $\kappa$, and by $K_0$ the fraction field of $W$, which is canonically a subfield of $K$. We fix an algebraic closure $\ol{K}$ of $K$ with the residue field $\ol\kappa$. For every finite extension $K'$ of $K_0$ contained in $\ol{K}$, put $\rG_{K'}\coloneqq\Gal(\ol{K}/K')$ as a profinite group.

\subsection{Preparation}

For a scheme $\cZ$ of finite type over $O_{K'}$ with $K'$ a finite extension of $K_0$ contained in $\ol{K}$, we
\begin{itemize}
  \item put $Z\coloneqq\cZ\otimes_{O_{K'}}K'$ for the generic fiber,

  \item put $\cZ_l\coloneqq\cZ\otimes\dZ/p^l$ for every integer $l\geq 1$,

  \item put $\sfZ\coloneqq\cZ\otimes_{O_{K'}}\kappa'$ for its special fiber, where $\kappa'$ is the residue field of $K'$,

  \item denote by $\fZ$ the formal completion of $\cZ$ along $\sfZ$,

  \item denote by $\fZ_\eta$ the generic fiber of $\fZ$, regarded as an analytic space over $K'$ in the sense of Berkovich,

  \item put $\ol\cZ\coloneqq\cZ\otimes_{O_{K'}}O_{\ol{K}}$, $\ol{Z}\coloneqq Z\otimes_{K'}\ol{K}$, and $\ol\sfZ\coloneqq\sfZ\otimes_{\kappa'}\ol\kappa$.
\end{itemize}
We apply the similar notational convention to morphisms over $O_{K'}$ as well.\footnote{Later, we will see $\cC,\cD,\cE,\cT,\cU,\cV,\cX,\cY,\cZ$ for schemes and $\cA,\cF,\cG,\cH$ for morphisms over $O_{K'}$.}

Suppose that $\sfX$ is a subscheme of $\sfZ$, we denote by $]\sfX[_{\fZ_\eta}$ its tubular neighbourhood in $\fZ_\eta$. We have the quasi-\'{e}tale site $]\sfX[_{\fZ_\eta,\qet}$ \cite{Ber94}*{\S3} with the natural map $]\sfX[_{\fZ_\eta,\qet}\to(\widehat{\fZ_{/\sfX}})_{\et}$, where $\widehat{\fZ_{/\sfX}}$ denotes the formal completion of $\fZ$ along $\sfX$. On the other hand, the natural map $(\widehat{\fZ_{/\sfX}})_{\et}\to\sfX_{\et}$ is an equivalence of sites \cite{Ber96}*{Proposition~2.1}. Together, we obtain the \emph{specialization} map $\ts_{(\sfX,\fZ)}\colon]\sfX[_{\fZ_\eta,\qet}\to\sfX_{\et}$, and will simply write $\ts$ when no confusion arises.

\begin{definition}\label{de:restriction}
Let $R$ be a ring. In the situation above, suppose that $\sfX$ is a closed subscheme of $\sfZ$ and $\sfU$ an open subscheme of $\sfX$, we define a functor
\[
\tf^!_{(\sfU,\sfX)}\colon\bM(]\sfX[_{\fZ_\eta,\qet},R)\to\bM(]\sfX[_{\fZ_\eta,\qet},R)
\]
to be the kernel of the unit transform $\id\to g_*\circ g^*$, where $g$ denotes the open immersion $]\sfX\setminus\sfU[_{\fZ_\eta}\to]\sfX[_{\fZ_\eta}$.
\end{definition}

\begin{remark}
The functors $g^*$, $g_*$, and $\tf^!_{(\sfU,\sfX)}$ are all exact. Moreover, there is in general no functor $\tf\colon\bD^+(\sfX_{\et},R)\to\bD^+(\sfX_{\et},R)$ such that $\tf\circ\rR\ts_{(\sfX,\fZ)*}\simeq\rR\ts_{(\sfX,\fZ)*}\circ\tf^!_{(\sfU,\sfX)}$, even when $\sfX=\sfZ$.
\end{remark}

\begin{lem}\label{le:restriction}
Let the situation be as in Definition \ref{de:restriction}. The diagram
\[
\xymatrix{
\sfF_!\circ\sfF^*\circ\rR\ts_{(\sfX,\fZ)*} \ar[r]\ar[d] & 0 \ar[d] \\
\rR\ts_{(\sfX,\fZ)*} \ar[r] & \rR\ts_{(\sfX,\fZ)*}\circ g_* \circ g^*
}
\]
of functors from $\bD^+(]\sfX[_{\fZ_\eta,\qet},R)$ to $\bD^+(\sfX_{\et},R)$ commutes, where $\sfF\colon\sfU\to\sfX$ denotes the open immersion. In particular, there is a canonical natural transform
\[
\sfF_!\circ\sfF^*\circ\rR\ts_{(\sfX,\fZ)*} \to \rR\ts_{(\sfX,\fZ)*}\circ\tf^!_{(\sfU,\sfX)}\colon\bD^+(]\sfX[_{\fZ_\eta,\qet},R)\to\bD^+(\sfX_{\et},R).
\]
\end{lem}

\begin{proof}
It suffices to notice that the unit transform $\rR\ts_{(\sfX,\fZ)*} \to \rR\ts_{(\sfX,\fZ)*}\circ g_*\circ g^*$ factors through
\[
\rR\ts_{(\sfX,\fZ)*} \to \sfG_*\circ\sfG^*\circ\rR\ts_{(\sfX,\fZ)*} \to
\sfG_*\circ\rR\ts_{(\sfX\setminus\sfU,\fZ)*}\circ g^* \xrightarrow\sim
\rR\ts_{(\sfX,\fZ)*}\circ g_*\circ g^*,
\]
where $\sfG\colon\sfX\setminus\sfU\to\sfX$ denotes the closed immersion.
\end{proof}

In what follows, we will work with log-schemes, written as $(X,L)$ with the first variable the underlying scheme and the second variable the log structure. Since the integral model in Theorem \ref{th:crystalline} is strictly semistable, we assume that the log structures are defined in the Zariski topology.

For a scheme $X$ and a closed subset $Y$, we denote by $L_X^Y$ the log structure $\sO_X\cap j_*\sO_{X\setminus Y}^\times\to\sO_X$, where $j\colon X\setminus Y\to X$ denotes the open immersion. For a log-scheme $(X,L)$ and a morphism $f\colon X'\to X$, we write $f^*L$ for the pullback log structure or simply $L\res_{X'}$ when $f$ is clear from the context.

We write $W^\triv$ for $(\Spec W,W^\times)$, $W[t]^\circ$ for $(\Spec W[t],L)$ where $L$ is the log structure associated with $1\mapsto t$, $W^\circ$ for the fiber of $W[t]^\circ$ at $t=0$, and $\kappa^\circ$ for the fiber of $W^\circ$ at $p=0$. Note that the natural morphism $W[t]^\circ\to W^\triv$ is log-smooth. For every extension $K'/K_0$ contained in $\ol{K}$ with the residue field $\kappa'$, we put $O_{K'}^\can\coloneqq(\Spec O_{K'},L_{\Spec O_{K'}}^{\Spec\kappa'})$.

Let $(X,M)$ be a fine log-scheme over a fine base log-scheme $(S,L)$ of finite type. Recall that an \emph{embedding system} for $(X,M)/(S,L)$ is a projective system  $\{(X^\star,M^\star)\hookrightarrow(Z^\star,N^\star)\}_{\star=0,1,\dots}$ of exact closed immersions of log-schemes over $(S,L)$ in which $X^\star$ is a Zariski hypercovering of $X$, $M^\star=M\res_{X^\star}$, and $(Z^\star,N^\star)$ is a fine log-scheme log-smooth over $(S,L)$ of finite type. Note that embedding system always exists.

In the case where $(S,L)=W[t]^\circ$ and $(X,M)$ is a strictly semistable log-scheme over $\kappa^\circ$ \cite{GK05}*{\S2.1} of finite type, we say that an embedding system $\{(X^\star,M^\star)\hookrightarrow(Z^\star,N^\star)\}$ for $(X,M)/(S,L)$ is \emph{admissible} if
\begin{itemize}
  \item $(X^0,M^0)\to(Z^0,N^0)$ induces an isomorphism $(X^0,M^0)\simeq(Z^0,N^0)\times_{W[t]^\circ}\kappa^\circ$;

  \item $Z^0$ is flat and generically smooth over $W[t]$, and is smooth over $W$;

  \item $Y^0\coloneqq Z^0\otimes_{W[t]}W$ is a relative strict normal crossings divisor of $Z^0$ over $W$;

  \item $N^0=L^{Y^0}_{Z^0}$;

  \item $(X^\star,M^\star)\to(Z^\star,N^\star)$ is (isomorphic to the one) induced from $(X^0,M^0)\to(Z^0,N^0)$ in the process described in \cite{GK05}*{\S5.1}.
\end{itemize}

\subsection{Rigid de Rham--Witt complexes}
\label{ss:witt}

Let $(\sfX,L)$ be a fine log-scheme over $\kappa^\circ$ of finite type. Let $\sfF\colon\sfU\to\sfX$ be an open subscheme. For every Zariski hypercovering $\sfX^\star$ of $\sfX$, we put $\sfF^\star\colon\sfU^\star\coloneqq\sfU\times_\sfX\sfX^\star\to\sfX^\star$.

Choose an embedding system $\{(\sfX^\star,L^\star)\hookrightarrow(\cY^\star,M^\star)\}$ for $(\sfX,L)/W^\circ$. Denote by $\sfu\colon\sfX^\star\to\sfX$ the augmentation map for the hypercovering, and put $\ts^\star\coloneqq\ts_{(\sfX^\star,\fY^\star)}$. We define the rigid de Rham--Witt complex of $(\sfX,L)$ to be\footnote{More precisely, it should be called \emph{convergent de Rham--Witt complex} since it gives the log convergent cohomology in general. However, later we will take $(\sfX,L)$ to be strictly semistable and proper.}
\[
\bomega_{\sfX}\coloneqq\rR\sfu_*
\(\rR\ts^\star_*\(\Omega^\bullet_{(\cY^\star,M^\star)/W^\circ}\otimes_{\sO_{\cY^\star}}\sO_{]\sfX^\star[_{\fY^\star}}\)\)
\in\bD^+(\sfX_{\et},K_0),
\]
where $\Omega^\bullet_{(\cY^\star,M^\star)/W^\circ}$ denotes the complex on $\cY^\star_{\et}$ of relative logarithmic differentials of the log-smooth morphism $(\cY^\star,M^\star)/W^\circ$. By (the same argument in) \cite{GK05}*{Lemma~1.4}, the complex $\bomega_{\sfX}$ does not depend on the choice of the embedding system for $(\sfX,M)/W^\circ$.\footnote{Of course, $\bomega_{\sfX}$ depends on the log structure $L$. However, as a common practice for de Rham--Witt complexes, we will not include $L$ in the notation.}

On the other hand, choose an embedding system $\{(\sfX^\star,L^\star)\hookrightarrow(\cZ^\star,N^\star)\}$ for $(\sfX,L)/W[t]^\circ$, hence for $(\sfX,L)/W^\triv$. Then
\[
\left\{(\sfX^\star,L^\star)\hookrightarrow(\cY^\star,M^\star)\coloneqq(\cZ^\star,N^\star)\times_{W[t]^\circ}W^\circ\right\}
\]
is an embedding system for $(\sfX,L)/W^\circ$. We have the short exact sequence
\[
0 \to \Omega^{q-1}_{(\cY^\star,M^\star)/W^\circ} \to \Omega^q_{(\cZ^\star,N^\star)/W^\triv}\otimes_{\sO_{\cZ^\star}}\sO_{\cY^\star} \to \Omega^q_{(\cY^\star,M^\star)/W^\circ}
\to 0
\]
of coherent sheaves on $\cY^\star_{\et}$, in which the first map is given by $\wedge\rd\log t$, for every $q\geq 0$ compatible with differentials. We put
\[
\widetilde\bomega_{\sfX}\coloneqq\rR\sfu_*
\(\rR\ts^\star_*\(\Omega^\bullet_{(\cZ^\star,N^\star)/W^\triv}\otimes_{\sO_{\cZ^\star}}\sO_{]\sfX^\star[_{\fY^\star}}\)\)
\in\bD^+(\sfX_{\et},K_0).
\]
Then there is a distinguished triangle
\begin{align}\label{eq:monodromy1}
\bomega^\triangle_\sfX\colon\quad
\bomega_{\sfX}[-1] \to \widetilde\bomega_{\sfX} \to \bomega_{\sfX} \xrightarrow{N}\bomega_{\sfX}
\end{align}
in $\bD^+(\sfX_{\et},K_0)$, where $N$ denotes the connecting map; it is independent of the choice of the embedding system for $(\sfX,L)/W[t]^\circ$. Moreover, for a morphism $f\colon (\sfX',L')\to(\sfX,L)$ of fine log-schemes over $\kappa^\circ$ of finite type, we have an induced map $\bomega^\triangle_\sfX\to\rR f_*\widetilde\bomega^\triangle_{\sfX'}$ of distinguished triangles in $\bD^+(\sfX_{\et},K_0)$.

We put
\begin{align*}
\bomega_{(\sfU,\sfX)}&\coloneqq\rR\sfu_*
\(\rR\ts^\star_*\tf^!_{(\sfU^\star,\sfX^\star)}
\(\Omega^\bullet_{(\cY^\star,M^\star)/W^\circ}\otimes_{\sO_{\cY^\star}}\sO_{]\sfX^\star[_{\fY^\star}}\)\), \\
\widetilde\bomega_{(\sfU,\sfX)}&\coloneqq\rR\sfu_*
\(\rR\ts^\star_*\tf^!_{(\sfU^\star,\sfX^\star)}
\(\Omega^\bullet_{(\cZ^\star,N^\star)/W^\triv}\otimes_{\sO_{\cZ^\star}}\sO_{]\sfX^\star[_{\fY^\star}}\)\),
\end{align*}
both in $\bD^+(\sfX_{\et},K_0)$ (see Definition \ref{de:restriction}). Then by definition, we have a distinguished triangle
\begin{align}\label{eq:monodromy3}
\bomega^\triangle_{(\sfU,\sfX)}\colon\quad
\bomega_{(\sfU,\sfX)}[-1] \to \widetilde\bomega_{(\sfU,\sfX)} \to \bomega_{(\sfU,\sfX)} \xrightarrow{N} \bomega_{(\sfU,\sfX)}
\end{align}
in $\bD^+(\sfX_{\et},K_0)$, and a distinguished triangle of distinguished triangles
\begin{align}\label{eq:derhamwitt}
\bomega^\triangle_{(\sfU,\sfX)} \to \bomega^\triangle_\sfX \to \sfG_*\bomega^\triangle_{\sfX\setminus\sfU} \xrightarrow{+1}
\end{align}
where $\sfG\colon\sfX\setminus\sfU\to\sfX$ denotes the closed immersion.

From now on, we assume that $(\sfX,L)$ is strictly semistable over $\kappa^\circ$. We recall the construction of several crystalline complexes of $(\sfX,L)$. For every $l\geq 1$, let $\cD^\star_l$ and $\cE^\star_l$ be the (scheme part of the) PD envelopes of $\sfX^\star$ in $\cY^\star_l$ and $\cZ^\star_l$ (over the base $W$ equipped with the usual PD structure), respectively.\footnote{The natural morphism $\cD^\star_l\to\cE^\star_l\otimes_{W\langle t\rangle}W$ is an isomorphism, where $W\langle t\rangle$ denotes the PD envelope of $(W[t],(t))$ over $W$.} We have complexes
\begin{align*}
&\Omega^\bullet_{(\cY^\star,M^\star)/W^\circ}\otimes_{\sO_{\cY^\star}}\sO_{\cD^\star_\bullet},\\
&\Omega^\bullet_{(\cZ^\star,N^\star)/W^\triv}\otimes_{\sO_{\cZ^\star}}\sO_{\cD^\star_\bullet},\\
&\Omega^\bullet_{(\cZ^\star,N^\star)/W^\triv}\otimes_{\sO_{\cZ^\star}}\sO_{\cE^\star_\bullet},\\
&\Omega^\bullet_{(\cZ^\star,N^\star)/W[t]^\circ}\otimes_{\sO_{\cZ^\star}}\sO_{\cE^\star_l},
\end{align*}
in $\bC^+(\sfX^\star_{\et},W_\bullet)$. Put
\begin{align*}
\sC_{(\sfX,L)/W^\circ}&\coloneqq
\rR\sfu_*\(\Omega^\bullet_{(\cY^\star,M^\star)/W^\circ}\otimes_{\sO_{\cY^\star}}\sO_{\cD^\star_l}\), \\
\widetilde\sC_{(\sfX,L)/W^\circ}&\coloneqq
\rR\sfu_*\(\Omega^\bullet_{(\cZ^\star,N^\star)/W^\triv}\otimes_{\sO_{\cZ^\star}}\sO_{\cD^\star_l}\), \\
\sC_{(\sfX,L)/W^\triv}&\coloneqq
\rR\sfu_*\(\Omega^\bullet_{(\cZ^\star,N^\star)/W^\triv}\otimes_{\sO_{\cZ^\star}}\sO_{\cE^\star_l}\), \\
\sC_{(\sfX,L)/W[t]^\circ}&\coloneqq
\rR\sfu_*\(\Omega^\bullet_{(\cZ^\star,N^\star)/W[t]^\circ}\otimes_{\sO_{\cZ^\star}}\sO_{\cE^\star_l}\),
\end{align*}
all in $\bD^+(\sfX_{\et},W_\bullet)$. It is well-known that the above objects do not depend on the choice of the embedding system for $(\sfX,L)/W[t]^\circ$. In fact, $\sC_{(\sfX,L)/W^\circ}$ is nothing but the modified de Rham--Witt complex $W\omega^\bullet_\sfX$ \cites{Hyo91,HK94}, and $\widetilde\sC_{(\sfX,L)/W^\circ}$ is simply $W\widetilde\omega^\bullet_\sfX$.

Applying the notation \eqref{eq:adic} (with $\dW=W$ and $\dL=K_0$), we obtain a distinguished triangle
\begin{align}\label{eq:monodromy2}
W\omega^\bullet_{\sfX,K_0}[-1] \to W\widetilde\omega^\bullet_{\sfX,K_0} \to W\omega^\bullet_{\sfX,K_0} \xrightarrow{N}W\omega^\bullet_{\sfX,K_0}
\end{align}
in $\bD^+(\sfX_{\et},K_0)$, similar to \eqref{eq:monodromy1}, in which the first arrow is given by $\wedge\rd\log t$, the second arrow is the natural one, and the last arrow is the connecting map.

We would like to compare \eqref{eq:monodromy1} and \eqref{eq:monodromy2}. We have a canonical map
\[
\sO_{\widehat{\fY^\star_{/\sfX^\star}}}\to\varprojlim_l\sO_{\cD^\star_l}=\rR\varprojlim_l\sO_{\cD^\star_l}
\]
as in \cite{Ber97}*{(1.9.2)}, which induces maps
\begin{align*}
&\quad\rR\ts^\star_*\(\Omega^\bullet_{(\cY^\star,M^\star)/W^\circ}\otimes_{\sO_{\cY^\star}}\sO_{]\sfX^\star[_{\fY^\star}}\) \\
&\simeq
\(\Omega^\bullet_{(\cY^\star,M^\star)/W^\circ}\otimes_{\sO_{\cY^\star}}\sO_{\widehat{\fY^\star_{/\sfX^\star}}}\)\otimes_{W}K_0 \\
&\to\(\rR\varprojlim_l\Omega^\bullet_{(\cY^\star,M^\star)/W^\circ}\otimes_{\sO_{\cY^\star}}\sO_{\cD^\star_l}\)\otimes_{W}K_0 \end{align*}
and
\begin{align*}
&\quad\rR\ts^\star_*\(\Omega^\bullet_{(\cZ^\star,N^\star)/W^\triv}\otimes_{\sO_{\cZ^\star}}\sO_{]\sfX^\star[_{\fY^\star}}\) \\
&\simeq
\(\Omega^\bullet_{(\cZ^\star,N^\star)/W^\triv}\otimes_{\sO_{\cZ^\star}}\sO_{\widehat{\fY^\star_{/\sfX^\star}}}\)\otimes_{W}K_0  \\ &\to\(\rR\varprojlim_l\Omega^\bullet_{(\cZ^\star,N^\star)/W^\triv}\otimes_{\sO_{\cZ^\star}}\sO_{\cD^\star_l}\)\otimes_{W}K_0
\end{align*}
in $\bD^+(\sfX^\star_{\et},K_0)$. These maps are in fact equivalences since $(\sfX^\star,L^\star)$ is strictly semistable over $\kappa^\circ$ by an argument similar to \cite{Ber97}*{\S1.9}. Applying $\rR\sfu_*$, we obtain equivalences
\begin{align}\label{eq:witt1}
\bomega_{\sfX}\xrightarrow\sim W\omega^\bullet_{\sfX,K_0},\quad
\widetilde\bomega_{\sfX}\xrightarrow\sim W\widetilde\omega^\bullet_{\sfX,K_0},
\end{align}
under which \eqref{eq:monodromy1} is equivalent to \eqref{eq:monodromy2}. On the other hand, by Lemma \ref{le:restriction}, we have a natural map
\begin{align}\label{eq:witt2}
\sfF_!\sfF^*W\widetilde\omega^\bullet_{\sfX,K_0}&\to\rR\sfu_*\sfF^\star_!(\sfF^\star)^*
\(\(\rR\varprojlim_l\Omega^\bullet_{(\cZ^\star,N^\star)/W^\triv}\otimes_{\sO_{\cZ^\star}}\sO_{\cD^\star_l}\)\otimes_{W}K_0\) \\
&\xrightarrow{\sim}\rR\sfu_*\(\sfF^\star_!(\sfF^\star)^*
\rR\ts^\star_*\(\Omega^\bullet_{(\cZ^\star,N^\star)/W^\triv}\otimes_{\sO_{\cZ^\star}}\sO_{]\sfX^\star[_{\fY^\star}}\)\) \notag \\
&\to\rR\sfu_*\(\rR\ts^\star_*\tf^!_{(\sfU^\star,\sfX^\star)}
\(\Omega^\bullet_{(\cZ^\star,N^\star)/W^\triv}\otimes_{\sO_{\cZ^\star}}\sO_{]\sfX^\star[_{\fY^\star}}\)\) \notag \\
&=\widetilde\bomega_{(\sfU,\sfX)} \notag
\end{align}
in $\bD^+(\sfX_{\et},K_0)$.

When the model in Theorem \ref{th:crystalline} is not smooth, we also need a cohomological variant of the rigid de Rham--Witt complex, which we now introduce. We choose an admissible embedding system $\{(\sfX^\star,L^\star)\hookrightarrow(\cZ^\star,N^\star)\}$ for $(\sfX,L)/W[t]^\circ$.

For every $q\geq 0$, we have a natural subsheaf $\Omega^q_{\cZ^\star/W}\subseteq\Omega^q_{(\cZ^\star,N^\star)/W^\triv}$. Put
\[
\Xi^q_{\cZ^\star}\coloneqq\frac{\Omega^{q+1}_{(\cZ^\star,N^\star)/W^\triv}}{\Omega^{q+1}_{\cZ^\star/W}},
\]
which is an $\sO_{\cY^\star}$-module. The map $\wedge\rd\log t$ induces a diagram
\begin{align}\label{eq:witt6}
\xymatrix{
\Omega^q_{(\cY^\star,M^\star)/W^\circ} \ar[d]\ar[dr] \\
\Omega^{q+1}_{(\cZ^\star,N^\star)/W^\triv}\otimes_{\sO_{\cZ^\star}}\sO_{\cY^\star} \ar[r]& \Xi^q_{\cZ^\star}
}
\end{align}
of coherent sheaves on $\cY^\star_{\et}$, compatible with differentials. Define
\begin{align*}
\bomega^+_{\sfX}&\coloneqq\rR\sfu_*
\(\rR\ts^\star_*\(\Xi^\bullet_{\cZ^\star}\otimes_{\sO_{\cY^\star}}\sO_{]\sfX^\star[_{\fY^\star}}\)\), \\
\bomega^+_{(\sfU,\sfX)}&\coloneqq\rR\sfu_*\(\rR\ts^\star_*\tf^!_{(\sfU^\star,\sfX^\star)}
\(\Xi^\bullet_{\cZ^\star}\otimes_{\sO_{\cY^\star}}\sO_{]\sfX^\star[_{\fY^\star}}\)\),
\end{align*}
and we have a natural diagram
\begin{align}\label{eq:witt4}
\xymatrix{
\bomega_{(\sfU,\sfX)}[-1] \ar[d]\ar[dr] \\
\widetilde\bomega_{(\sfU,\sfX)} \ar[r]& \bomega^+_{(\sfU,\sfX)}[-1]
}
\end{align}
in $\bD^+(\sfX_{\et},K_0)$, in which the vertical arrow is same as the first arrow in the first line in \ref{eq:derhamwitt}. Let $W\Xi^\bullet_\sfX\in\bD^+(\sfX_{\et},W_\bullet)$ be the cohomological de Rham--Witt complex defined in \cite{Sat13}*{Definition~8.3}. Similar to \eqref{eq:witt1}, we have a natural equivalence
\[
\bomega^+_{\sfX}\simeq W\Xi^\bullet_{\sfX,K_0}
\]
in $\bD^+(\sfX_{\et},K_0)$ by \cite{Sat13}*{Proposition~8.4}. Similar to \eqref{eq:witt2}, we have a natural map
\begin{align}\label{eq:witt3}
\sfF_!\sfF^*W\Xi^\bullet_{\sfX,K_0}\to\bomega^+_{(\sfU,\sfX)}
\end{align}
in $\bD^+(\sfX_{\et},K_0)$, fitting into the following diagram
\begin{align}\label{eq:witt5}
\xymatrix{
\sfF_!\sfF^*W\widetilde\omega^\bullet_{\sfX,K_0} \ar[r]\ar[d]_-{\eqref{eq:witt2}} & \sfF_!\sfF^*W\Xi^\bullet_{\sfX,K_0}[-1] \ar[d]^-{\eqref{eq:witt3}} \\
\widetilde\bomega_{(\sfU,\sfX)} \ar[r] & \bomega^+_{(\sfU,\sfX)}[-1]
}
\end{align}
in which the upper horizontal arrow is induced by the one in \cite{Sat13}*{Proposition~8.10}.

\subsection{Log rigid cohomology with proper support}

Let the situation be as in the previous subsection with $(\sfX,L)$ strictly semistable over $\kappa^\circ$. We also assume that $\sfX$ is proper of pure dimension $n-1$. For every $h\geq 1$, let $\sfX^{(h)}$ be the disjoint union of intersections of $h$ different irreducible components of $\sfX$, and put $\sfU^{(h)}\coloneqq\sfU\times_{\sfX}\sfX^{(h)}$.

Recall from \cite{GK05}*{\S1.5} that for a scheme $\sfY$ over $\sfX$ of finite type, we have the log rigid cohomology $\rH^\bullet_\rig(\sfY/W^\circ)$ and the log convergent cohomology $\rH^\bullet_\conv(\sfY/W^\circ)$ for the log-scheme $(\sfY,L\res_\sfY)$, with a natural map $\rH^\bullet_\rig(\sfY/W^\circ)\to\rH^\bullet_\conv(\sfY/W^\circ)$. In particular, we have
\begin{align*}
\bH^q(\sfX_{\et},\bomega_\sfX)&=\rH^q_\conv(\sfX/W^\circ),\\
\bH^q(\sfX_{\et},\sfG_*\bomega_{\sfX\setminus\sfU})&=\rH^q_\conv(\sfX\setminus\sfU/W^\circ),
\end{align*}
for every $q\geq 0$.\footnote{Here, we use the fact that computing cohomology of coherent sheaves in the quasi-\'{e}tale topology of analytic spaces is the same as in the G-topology.} Moreover, the natural map $\rH^\bullet_\rig(\sfX/W^\circ)\to\rH^\bullet_\conv(\sfX/W^\circ)$ is an isomorphism by \cite{GK05}*{Theorem~5.3(ii)}.

\begin{definition}\label{de:rigid}
We define, in an \emph{ad hoc} way, the \emph{log rigid cohomology of $\sfU$ with proper support} to be
\[
\rH^q_\rig((\sfU,\sfX)/W^\circ)\coloneqq\bH^q(\sfX_{\et},\bomega_{(\sfU,\sfX)}),
\]
which \emph{a priori} depends on the embedding $\sfU\hookrightarrow\sfX$. For $\sfU\subseteq\sfU'\subseteq\sfX$, we have a natural pushforward map $\rH^q_\rig((\sfU,\sfX)/W^\circ)\to\rH^q_\rig((\sfU',\sfX)/W^\circ)$ by construction.
\end{definition}

The distinguished triangle \eqref{eq:derhamwitt} induces a long exact sequence
\begin{align}\label{eq:rigid}
\cdots &\to \rH^{q-1}_\conv(\sfX\setminus\sfU/W^\circ)
\to \rH^q_\rig((\sfU,\sfX)/W^\circ)  \\
&\to \rH^q_\rig(\sfX/W^\circ)
\to \rH^q_\conv(\sfX\setminus\sfU/W^\circ) \to \cdots \notag
\end{align}
of $K_0$-vector spaces.

We now review the weight spectral sequence for log rigid cohomology from \cite{GK05}*{\S5}, which is the rigid analogue of Mokrane's spectral sequence for log crystalline cohomology \cite{Mok93}. Take an admissible embedding system $\{(\sfX^\star,L^\star)\hookrightarrow(\cZ^\star,N^\star)\}$ for $(\sfX,L)/W[t]^\circ$. For $j\geq 0$, put
\[
P_j\Omega^q_{(\cZ^\star,N^\star)/W^\triv}\coloneqq
\IM\(\Omega^j_{(\cZ^\star,N^\star)/W^\triv}\otimes\Omega^{q-j}_{\cZ^\star/W}\to\Omega^q_{(\cZ^\star,N^\star)/W^\triv}\).
\]
We have the double complex
\[
A^{ij}_{\cZ^\star}\coloneqq
\frac{\Omega^{i+j+1}_{(\cZ^\star,N^\star)/W^\triv}}{P_j\Omega^{i+j+1}_{(\cZ^\star,N^\star)/W^\triv}}
\]
of $\sO_{\cY^\star}$-modules, in which the differential $A^{ij}_{\cZ^\star}\to A^{(i+1)j}_{\cZ^\star}$ is given by $(-1)^j\rd$ and the differential $A^{ij}_{\cZ^\star}\to A^{i(j+1)}_{\cZ^\star}$ is given by $\wedge\rd\log t$, with the filtration
\[
P_kA^{ij}_{\cZ^\star}\coloneqq
\frac{P_{2j+k+1}\Omega^{i+j+1}_{(\cZ^\star,N^\star)/W^\triv}}{P_j\Omega^{i+j+1}_{(\cZ^\star,N^\star)/W^\triv}}
\]
for $k\geq -j$. In particular, $A^{\bullet 0}_{\cZ^\star}$ is nothing but the complex $\Xi^\bullet_{\cZ^\star}$ from the previous subsection. Let $A^\bullet_{\cZ^\star}$ be the total complex of $A^{\bullet\bullet}_{\cZ^\star}$. It is shown in \cite{GK05}*{\S5.2} that the augmentation map $\Omega^\bullet_{(\cY^\star,M^\star)/W^\circ}\to \Xi^\bullet_{\cZ^\star}=A^{\bullet 0}_{\cZ^\star}$ in \eqref{eq:witt6} induces an equivalence $\Omega^\bullet_{(\cY^\star,M^\star)/W^\circ}\xrightarrow\sim A^\bullet_{\cZ^\star}$ in $\bD^+(\cY^\star_{\et},K_0)$. Then the total filtration on $A^\bullet_{\cZ^\star}$ induces spectral sequences
\begin{align*}
\rE(\sfX)^{-k,q+k}_1&=\bigoplus_{j\geq\max\{0,-k\}}\rH^{q-2j-k}_\rig(\sfX^{(2j+k+1)}/K_0) \Rightarrow \rH^q_\rig(\sfX/W^\circ), \\
\rE(\sfU)^{-k,q+k}_1&=\bigoplus_{j\geq\max\{0,-k\}}\rH^{q-2j-k}_\rig(\sfU^{(2j+k+1)}/K_0) \Rightarrow \rH^q_\rig(\sfU/W^\circ),
\end{align*}
which already appeared in \cite{GK05}*{(4)}, and
\begin{align}\label{eq:duality1}
&\quad\rE(\sfX\setminus\sfU)^{-k,q+k}_1 \\
&=\bigoplus_{j\geq\max\{0,-k\}}\rH^{q-2j-k}_\rig(\sfX^{(2j+k+1)}\setminus\sfU^{(2j+k+1)}/K_0) \Rightarrow \rH^q_\conv(\sfX\setminus\sfU/W^\circ),\notag
\end{align}
\begin{align}\label{eq:duality2}
&\quad\rE(\sfU,\sfX)^{-k,q+k}_1 \\
&=\bigoplus_{j\geq\max\{0,-k\}}\rH^{q-2j-k}_\rig((\sfU^{(2j+k+1)},\sfX^{(2j+k+1)})/K_0) \Rightarrow \rH^q_\rig((\sfU,\sfX)/W^\circ).\notag
\end{align}
Here, $\rH^\bullet_\rig((\sfU^{(2j+k+1)},\sfX^{(2j+k+1)})/K_0)$ is defined similarly as for $\rH^\bullet_\rig((\sfU,\sfX)/W^\circ)$ but without the log structure, which in fact coincides with the rigid cohomology with proper support $\rH^\bullet_{\rc,\rig}(\sfU^{(2j+k+1)}/K_0)$ defined by Berthelot since $\sfX^{(2j+k+1)}$ is proper. In particular, the spectral sequence \eqref{eq:duality2} can also be written as
\[
\rE_\rc(\sfU)^{-k,q+k}_1=\bigoplus_{j\geq\max\{0,-k\}}\rH^{q-2j-k}_{\rc,\rig}(\sfU^{(2j+k+1)}/K_0) \Rightarrow \rH^q_\rig((\sfU,\sfX)/W^\circ).
\]

The following two lemmas will be used later.

\begin{lem}\label{le:duality}
Let $d\geq 1$ be an integer. Suppose that $\dim(\sfX^{(h)}\setminus\sfU^{(h)})\leq d-h$ for every $h\geq 1$.
\begin{enumerate}
  \item The natural map $\rH^q_\rig((\sfU,\sfX)/W^\circ)\to\rH^q_\rig(\sfX/W^\circ)$ is an isomorphism for $q\geq 2d$.

  \item The natural map $\rE_\rc(\sfU)^{-k,q+k}_1\to\rE(\sfX)^{-k,q+k}_1$ is an isomorphism for $q\geq 2d-|k|$.

  \item For the map $\rE_\rc(\sfU)^{0,2d-1}_1\to\rE(\sfX)^{0,2d-1}_1$, the direct summand
      \[
      \bigoplus_{j\geq 1}\rH^{2d-1-2j}_{\rc,\rig}(\sfU^{(2j+1)}/K_0)\to\bigoplus_{j\geq 1}\rH^{2d-1-2j}_\rig(\sfX^{(2j+1)}/K_0)
      \]
      is an isomorphism.
\end{enumerate}
\end{lem}

\begin{proof}
For (1), it suffices to show that $\rH^q_\conv(\sfX\setminus\sfU/W^\circ)=0$ for $q\geq 2d-1$. By the spectral sequence \eqref{eq:duality1}, it suffices to show that
\[
\rH^{q-2j-k}_\rig(\sfX^{(2j+k+1)}\setminus\sfU^{(2j+k+1)}/K_0)=0
\]
for every $j$, $k$, and $q\geq 2d-1$, which follows from the fact that
\[
2d-1-2j-k>2(d-2j-k-1)\geq 2\dim(\sfX^{(2j+k+1)}\setminus\sfU^{(2j+k+1)}).
\]

For (2), it follows from the fact that
\[
\rH^{q-2j-k-1}_\rig(\sfX^{(2j+k+1)}\setminus\sfU^{(2j+k+1)}/K_0)
=\rH^{q-2j-k}_\rig(\sfX^{(2j+k+1)}\setminus\sfU^{(2j+k+1)}/K_0)=0
\]
for every $j\geq\max\{0,-k\}$ when $q\geq 2d-|k|$.

For (3), it follows from the fact that
\[
\rH^{2d-1-2j-1}_\rig(\sfX^{(2j+1)}\setminus\sfU^{(2j+1)}/K_0)=\rH^{2d-1-2j}_\rig(\sfX^{(2j+1)}\setminus\sfU^{(2j+1)}/K_0)=0
\]
for $j\geq 1$.
\end{proof}

\begin{lem}\label{le:duality1}
There is a spectral sequence $\rE^+_\rc(\sfU)^{-k,q+k}_1\Rightarrow\bH^q(\sfX_{\et},\bomega^+_{(\sfU,\sfX)})$ with
\[
\rE^+_\rc(\sfU)^{-k,q+k}_1=
\begin{dcases}
\rH^{q-k}_{\rc,\rig}(\sfU^{(k+1)}/K_0), & k\geq 0; \\
0 ,& k<0.
\end{dcases}
\]
Moreover, the map $\rH^q_\rig((\sfU,\sfX)/W^\circ)=\bH^q(\sfX_{\et},\bomega_{(\sfU,\sfX)})\to\bH^q(\sfX_{\et},\bomega^+_{(\sfU,\sfX)})$ is abutted by the map $\rE_\rc(\sfU)^{-k,q+k}_1\to\rE^+_\rc(\sfU)^{-k,q+k}_1$ given by the obvious projections.
\end{lem}

\begin{proof}
The spectral sequence follows from the filtration $P_kA^{i0}_{\cZ^\star}$ of $A^{i0}_{\cZ^\star}=\Xi^i_{\cZ^\star}$. Recall that the map $\bomega_{(\sfU,\sfX)}\to\bomega^+_{(\sfU,\sfX)}$ is induced by the natural projection map $A^\bullet_{\cZ^\star}\to A^{\bullet 0}_{\cZ^\star}$. The lemma follows since $P_kA^{i0}_{\cZ^\star}$ is the image of $P_kA^i_{\cZ^\star}$ under this map.
\end{proof}

\subsection{Abel--Jacobi map via rigid cohomology}
\label{ss:comparison}

We start to prove Theorem \ref{th:crystalline}. We fix a uniformizer $\varpi$ of $K$, and regard $O_K$ as an $W[t]$-ring via $t\mapsto\varpi$, making $O_K^\can$ an exact closed log-subscheme of $W[t]^\circ$. Let $\dB_\cris$ be the crystalline period ring and $\dB_\st$ the semistable period ring with respect to $\varpi$.

We fix a proper strictly semistable scheme $\cX$ over $O_K$ of pure (absolute) dimension $n\geq 2$. Then $(\cX,L_\cX^\sfX)$ is log-smooth over $O_K^\can$, and $(\sfX,L\coloneqq L_\cX^\sfX\res_\sfX)$ is strictly semistable over $\kappa^\circ$ of finite type. Consider an open immersion $\cF\colon\cU\to\cX$. The following definition will be frequently used later, which is related to condition (2) in Theorem \ref{th:crystalline}.

\begin{definition}
For an integer $d\geq 1$, we say that $\cU$ is \emph{$d$-dense} if $\dim(\sfX^{(h)}\setminus\sfU^{(h)})\leq d-h$ for every $h\geq 1$.
\end{definition}

Put $i\colon\sfX\to\cX$ and $j\colon X\to\cX$ for the special fiber and the generic fiber of $\cX$, respectively. Take an integer $d$ satisfying $0\leq d<p-1$. Let $\sS(d)_\cX\in\bD^+(\sfX_\et,\dZ_{p\bullet})$ be Kato's (log) syntomic complexes for $(\cX,L_\cX^\sfX)$.\footnote{We will recall the construction of many syntomic and crystalline complexes in \S\ref{ss:alteration} in which $\sS(d)_\cX$ is a special case.} We have the period map $\sS(d)_\cX\to i^*\rR j_*(\mu_{p^\bullet}^{\otimes d})_X$ which induces equivalences $\sS(d)_\cX\xrightarrow\sim\tau^{\leq d}i^*\rR j_*(\mu_{p^\bullet}^{\otimes d})_X$ (\cites{Kat94,Tsu00}). Put $i_\cU\colon\sfU\to\cU$ and $j_\cU\colon U\to\cU$ for the special fiber and the generic fiber of $\cU$, respectively. We have a sequence of maps
\[
\sfF_!\sfF^*i^*\rR j_*(\mu_{p^\bullet}^{\otimes d})_X
\xrightarrow{\sim}\sfF_! i_\cU^* \rR j_{\cU*}(\mu_{p^\bullet}^{\otimes d})_U
\xrightarrow{\sim}i^*\cF_!\rR j_{\cU*}(\mu_{p^\bullet}^{\otimes d})_U
\to i^*\rR j_* F_!(\mu_{p^\bullet}^{\otimes d})_U
\]
in $\bD^+(\sfX_\et,\dZ_{p\bullet})$, in which the last one is given by adjunction. Then we obtain the maps
\begin{align*}
&\quad\rR\Gamma\(\sfX_{\et},\sfF_!\sfF^*\rR\varprojlim\sS(d)_\cX\) \\
&\to\rR\Gamma\(\sfX_{\et},\rR\varprojlim\sfF_!\sfF^*\sS_l(d)_\cX\) \\
&\to\rR\Gamma\(\sfX_{\et},\rR\varprojlim\sfF_!\sfF^*i^*\rR j_*(\mu_{p^\bullet}^{\otimes d})_X\) \\
&\to\rR\Gamma\(\sfX_{\et},\rR\varprojlim i^*\rR j_* F_!(\mu_{p^\bullet}^{\otimes d})_U\) \\
&\xrightarrow\sim\rR\Gamma\(X_{\et},\rR\varprojlim F_!(\mu_{p^\bullet}^{\otimes d})_U\),
\end{align*}
where we have used the proper base change for the last equivalence. Put
\[
\rR\Gamma_\rc(U,\dQ_p(d))\coloneqq\rR\Gamma\(X_{\et},\rR\varprojlim F_!(\mu_{p^\bullet}^{\otimes d})_U\)\otimes_{\dZ_p}\dQ_p
\simeq\rR\Gamma\(X_{\et},F_!\dQ_p(d)_U\),
\]
whose $q$-th cohomology gives $\rH^q_\rc(U,\dQ_p(d))$. Then the composition of the previous sequence gives a map
\begin{align}\label{eq:syntomic}
\rR\Gamma\(\sfX_{\et},\sfF_!\sfF^*\sS(d)_{\cX,\dQ_p}\)\to\rR\Gamma_\rc(U,\dQ_p(d))
\end{align}
in $\bD^+(\sfX_{\et},\dQ_p)$, where we have again applied the notation \eqref{eq:adic} (with $\dW=\dZ_p$ and $\dL=\dQ_p$).

The Hochschild--Serre spectral sequence (Lemma \ref{le:spectral}) induces the edge map
\[
\rH^q_\rc(U,\dQ_p(d))^0\to\rH^1(K,\rH^{q-1}_\rc(\ol{U},\dQ_p(d))),
\]
where
\[
\rH^q_\rc(U,\dQ_p(d))^0\coloneqq\Ker\(\rH^q_\rc(U,\dQ_p(d)))\to\rH^q_\rc(\ol{U},\dQ_p(d))\).
\]
Let $\bH^q(\sfX_{\et},\sfF_!\sfF^*\sS(d)_{\cX,\dQ_p})^\heartsuit$ be the inverse image of $\rH^q_\rc(U,\dQ_p(d))^0$ under the map \eqref{eq:syntomic} (after taking $q$-th cohomology). Then we have the induced composite map
\begin{align}\label{eq:comparison}
\alpha_q\colon\bH^q(\sfX_{\et},\sfF_!\sfF^*\sS(d)_{\cX,\dQ_p})^\heartsuit
&\to\rH^1(K,\rH^{q-1}_\rc(\ol{U},\dQ_p(d))) \\
&\to\rH^1(K,\rH^{q-1}_\rc(\ol{U},\dQ_p)\otimes_{\dQ_p}\dB_\cris), \notag
\end{align}
where in the last map we use the canonical embedding $\dQ_p(d)\hookrightarrow\dB_\cris$ in the category $\bM_K(\dQ_p)$.

To study the kernel of $\alpha_{2d}$, we need to use crystalline complexes for $\cX$ rather than its special fiber. Let $\sC_{\cX/W^\triv}$ and $\sC_{\cX/W[t]^\circ}$ be the objects in $\bD^+(\sfX_{\et},W_\bullet)$ defined similarly as $\sC_{(\sfX,L)/W^\triv}$ and $\sC_{(\sfX,L)/W[t]^\circ}$ in \S\ref{ss:witt}, respectively, for which we use an embedding system for $(\cX,L_\cX^\sfX)/W[t]^\circ$ and just replace $\cE^\star_l$ by the PD envelope of $\cX_l$ in $\cZ^\star_l$. We have the following commutative diagram
\begin{align}\label{eq:comparison6}
\xymatrix{
\sS(d)_{\cX,\dQ_p} \ar[r]& \sC_{\cX/W^\triv,K_0} \ar[r]\ar[d]& \sC_{\cX/W[t]^\circ,K_0} \ar[d] \\
& W\widetilde\omega^\bullet_{\sfX,K_0} \ar[r] & W\omega^\bullet_{\sfX,K_0}
}
\end{align}
in $\bD^+(\sfX_{\et},\dQ_p)$, in which the first arrow is the natural map from the syntomic complex to the crystalline complex over $W^\triv$, and the vertical maps are induced by the specialization at $t=0$.

Using \eqref{eq:witt2}, we obtain a map
\begin{align}\label{eq:comparison5}
\sfF_!\sfF^*\sS(d)_{\cX,\dQ_p}\to\widetilde\bomega_{(\sfU,\sfX)}
\end{align}
in $\bD^+(\sfX_{\et},\dQ_p)$.

\begin{lem}\label{le:comparison}
Suppose that $\cU$ is $d$-dense if $d\geq 1$. Then the composite map
\begin{align*}
\bH^{2d}(\sfX_{\et},\sfF_!\sfF^*\sS(d)_{\cX,\dQ_p})
&\to\bH^{2d}(\sfX_{\et},\widetilde\bomega_{(\sfU,\sfX)}) \\
&\to\bH^{2d}(\sfX_{\et},\bomega_{(\sfU,\sfX)})=\rH^{2d}_\rig((\sfU,\sfX)/W^\circ)
\end{align*}
(Definition \ref{de:rigid}) vanishes on $\bH^{2d}(\sfX_{\et},\sfF_!\sfF^*\sS(d)_{\cX,\dQ_p})^\heartsuit$.
\end{lem}

\begin{proof}
When $d=0$, the natural map $\rH^0_\rig((\sfU,\sfX)/W^\circ)\to\rH^0_\rig(\sfX/W^\circ)$ is injective. When $d\geq 1$, since $\cU$ is $d$-dense, by Lemma \ref{le:duality}(1) and the long exact sequence \eqref{eq:rigid}, the natural map $\rH^{2d}_\rig((\sfU,\sfX)/W^\circ)\to\rH^{2d}_\rig(\sfX/W^\circ)$ is an isomorphism, in particular, injective as well. Thus, in both cases, we may assume $\cU=\cX$. Then the map $\bH^{2d}(\sfX_{\et},\sS(d)_{\cX,\dQ_p})\to\rH^{2d}_\rig(\sfX/W^\circ)$ factors through the map $\bH^{2d}(\sfX_{\et},\sS(d)_{\cX,\dQ_p})\to\bH^{2d}(\sfX_{\et},\sC_{\cX/W[t]^\circ,K_0})$ by \eqref{eq:comparison6}. We have the commutative diagram
\[
\xymatrix{
\bH^{2d}(\sfX_{\et},\sS(d)_{\cX,\dQ_p}) \ar[r]\ar@{-->}[dr]& \bH^{2d}(\sfX_{\et},\sC_{\cX/W^\triv,K_0}) \ar[r]\ar[d] & \bH^{2d}(\sfX_{\et},\sC_{\cX/W[t]^\circ,K_0}) \ar@{^(->}[d] \\
& \bH^{2d}(\ol\sfX_{\et},\sC_{\ol\cX/W^\triv,K_0}) \ar@{^(->}[r]& \bH^{2d}(\ol\sfX_{\et},\sC_{\ol\cX/W[t]^\circ,K_0})
}
\]
in which the injectivity of the two arrows follows from \cite{Sat13}*{Proposition~A.3.1}.\footnote{We will in fact review the definition of the objects $\sC_{\ol\cX/W^\triv}$ and $\sC_{\ol\cX/W[t]^\circ}$ of $\bD^+_{\rG_K}(\ol\sfX_\et,W_\bullet)$ in a more general setup in \S\ref{ss:alteration}. Meanwhile, it suffices to note that in terms of the notation, our map $\bH^{2d}(\ol\sfX_{\et},\sC_{\ol\cX/W^\triv,K_0})\to\bH^{2d}(\ol\sfX_{\et},\sC_{\ol\cX/W[t]^\circ,K_0})$ is parallel to the map $H^{2d}_{\r{crys}}((\ol{X},\ol{M})/W)_{\dQ_p}\to H^{2d}_{\r{crys}}((\ol{X},\ol{M})/(\sE,M_\sE))_{\dQ_p}$ in \cite{Sat13}.} By \cite{Sat13}*{Lemma~9.5}, $\bH^{2d}(\sfX_{\et},\sS(d)_{\cX,\dQ_p})^\heartsuit$ is contained in the kernel of the dashed arrow hence contained in the kernel of the map $\bH^{2d}(\sfX_{\et},\sS(d)_{\cX,\dQ_p})\to\bH^{2d}(\sfX_{\et},\sC_{\cX/W[t]^\circ,K_0})$. The lemma follows.
\end{proof}

The long exact sequence induced by \eqref{eq:monodromy3} gives an isomorphism
\begin{align*}
\frac{\rH^{q-1}_\rig((\sfU,\sfX)/W^\circ)}{N\rH^{q-1}_\rig((\sfU,\sfX)/W^\circ)}
\xrightarrow\sim\Ker\(\bH^q(\sfX_{\et},\widetilde\bomega_{(\sfU,\sfX)})
\to \bH^q(\sfX_{\et},\bomega_{(\sfU,\sfX)})\).
\end{align*}
Thus, by Lemma \ref{le:comparison}, we obtain a map
\begin{align}\label{eq:comparison4}
\rho_{2d}\colon\bH^{2d}(\sfX_{\et},\sfF_!\sfF^*\sS(d)_{\cX,\dQ_p})^\heartsuit
\to\frac{\rH^{2d-1}_\rig((\sfU,\sfX)/W^\circ)}{N\rH^{2d-1}_\rig((\sfU,\sfX)/W^\circ)}.
\end{align}
For every finite extension $K^\dag/K$ contained in $\ol{K}$ and every object $V\in\bM_K(\dQ_p)$, we denote by
\[
\r{res}_{K^\dag}\colon\rH^1(K,V)\to\rH^1(K^\dag,V)
\]
the restriction map.

The following proposition is the key to the proof of Theorem \ref{th:crystalline}, whose proof will be given in \S\ref{ss:alteration}.

\begin{proposition}\label{pr:comparison}
Suppose that $n<p$, $1\leq d< p-1$, and $\cU$ is $d$-dense. There exists a finite extension $K_U/K$ (depending on $U$) contained in $\ol{K}$ such that
\[
\Ker(\rho_{2d}) \subseteq \Ker(\r{res}_{K_U}\circ\alpha_{2d})
\]
holds. Moreover, we may take $K_X$ to be $K$.
\end{proposition}

Now we bring the $\dL$-ring $\dT$ of \'{e}tale correspondences on $\cX$ in Theorem \ref{th:crystalline}. In what follows, we write $V_\dL\coloneqq V\otimes_{\dQ_p}\dL$ for a $\dQ_p$-vector space $V$. In the situation of Theorem \ref{th:crystalline}, we may assume $1\leq d\leq n-1$ without loss of generality.

For every $t\in\dT$, put
\[
\cU_t\coloneqq\cX\setminus\((\cX\setminus\cU)^{t^\vee}\),\quad \cF_t\colon\cU_t\to\cX
\]
(see Notation \ref{no:correspondence}), where $t^\vee$ denotes the transpose of $t$. Then we have $(\cU_t)^t\subseteq\cU$. The element $t$ acts on various cohomology and spectral sequences compatibly,\footnote{We note that for a finite \'{e}tale morphism $f\colon\cX_0\to\cX$, one can choose admissible embedding systems $\{(\sfX^\star_0,L^\star_0)\hookrightarrow(\cZ^\star_0,N^\star_0)\}$ and $\{(\sfX^\star,L^\star)\hookrightarrow(\cZ^\star,N^\star)\}$ for $(\sfX_0,L_0)/W[t]^\circ$ and $(\sfX,L)/W[t]^\circ$, respectively, with an \'{e}tale morphism $(\cZ^\star_0,N^\star_0)\to(\cZ^\star,N^\star)$ that is compatible with $f$. We have the similar statement for embedding systems for $(\cX_0,L_{\cX_0}^{\sfX_0})/W[t]^\circ$ and $(\cX,L_\cX^\sfX)/W[t]^\circ$.} giving a commutative diagram
\begin{align}\label{eq:comparison7}
\resizebox{\hsize}{!}{
\xymatrix{
\rH^{2d-1}_\rig((\sfU_t,\sfX)/W^\circ)_\dL \ar[rrr]^-{t^*}\ar[d]^-{N} &&& \rH^{2d-1}_\rig((\sfU,\sfX)/W^\circ)_\dL \ar[d]^-{N} \\
\rH^{2d-1}_\rig((\sfU_t,\sfX)/W^\circ)_\dL \ar[rrr]^-{t^*}\ar[rd]\ar[dd] &&& \rH^{2d-1}_\rig((\sfU,\sfX)/W^\circ)_\dL \ar[ld]\ar[dd] \\
& \bH^{2d-1}(\sfX_{\et},\bomega^+_{(\sfU_t,\sfX)})_\dL \ar[r]^-{t^*} & \bH^{2d-1}(\sfX_{\et},\bomega^+_{(\sfU,\sfX)})_\dL \\
\bH^{2d}(\sfX_{\et},\widetilde\bomega_{(\sfU_t,\sfX)})_\dL \ar[rrr]^-{t^*}\ar[ur]\ar[d] &&& \bH^{2d}(\sfX_{\et},\widetilde\bomega_{(\sfU,\sfX)})_\dL \ar[ul]\ar[d] \\
\rH^{2d}_\rig((\sfU_t,\sfX)/W^\circ)_\dL \ar[rrr]^-{t^*} &&& \rH^{2d}_\rig((\sfU,\sfX)/W^\circ)_\dL \\
\bH^{2d}(\sfX_{\et},\sfF_{t!}\sfF_t^*\sS(d)_{\cX,\dQ_p})_\dL \ar[rrr]^-{t^*}
\ar@/^3.8pc/[uu]^-{\eqref{eq:comparison5}} &&&
\bH^{2d}(\sfX_{\et},\sfF_!\sfF^*\sS(d)_{\cX,\dQ_p})_\dL \ar@/_3.8pc/[uu]_-{\eqref{eq:comparison5}}
}
}
\end{align}
in which the two triangles are induced from \eqref{eq:witt4}.

Let $\dI\subseteq\dT$ be the annihilator of
\begin{align*}
\bigoplus_{h>1,q\geq 0}\rH^q_\cris(\sfX^{(h)}/W)\otimes_{\dZ_p}\dL=
\bigoplus_{h>1,q\geq 0}\rH^q_\rig(\sfX^{(h)}/K_0)_\dL.
\end{align*}

\begin{lem}\label{le:crystalline}
Suppose that $\cU$ is $d$-dense. Then for every $t\in\dI^{4n-5}$, the kernel of the map
\[
\frac{\rH^{2d-1}_\rig((\sfU_t,\sfX)/W^\circ)_\dL}{N\rH^{2d-1}_\rig((\sfU_t,\sfX)/W^\circ)_\dL}
\to\bH^{2d-1}(\sfX_{\et},\bomega^+_{(\sfU_t,\sfX)})_\dL
\]
is annihilated by $t^*$.
\end{lem}

\begin{proof}
Since $\cU$ is $d$-dense, $\cU_t$ is $d$-dense as well. Let
\[
0=F^{-1}\subseteq F^0\subseteq\cdots \subseteq F^{4d-2}=\rH^{2d-1}_\rig((\sfU_t,\sfX)/W^\circ)
\]
be the filtration induced by the spectral sequence $\rE_\rc(\sfU_t)^{-k,q+k}_1$. Let $V$ be the kernel of the map in the lemma. For every $i$, we regard $F^i_\dL\cap V$ as the intersection of $V$ and the image of $F^i_\dL$ in the target of the map in the lemma.

By Lemma \ref{le:duality}(2,3) and Lemma \ref{le:duality1}, we know that $F^i_\dL\cap V/F^{i-1}_\dL\cap V$ is a subquotient of $\bigoplus_{h>1,q}\rH^q_\rig(\sfX^{(h)}/K_0)_\dL$ for every $0\leq i\leq 4d-2$. Thus, every element in $\dI$ annihilates $F^i_\dL\cap V/F^{i-1}_\dL\cap V$ for $0\leq i\leq 4d-2$, which implies that $t^*$ annihilates $V$.
\end{proof}

\begin{proposition}\label{pr:crystalline}
Suppose that $n<p$ and $\cU$ is $d$-dense. Let $t$ be an element in $\dI^{4n-5}$. Then for every $c\in\rZ^d(X)_\dL^0$ such that the Zariski closure of its support in $\cX$ is contained in $\cU_t$, we have
\[
\r{res}_{K_U}(t^*\beta_c)\in\rH^1_f(K_U,\rH^{2d-1}_\rc(\ol{U},\dL(d))),
\]
where $\beta_c\in\rH^1(K,\rH^{2d-1}_\rc(\ol{U_t},\dL(d)))$ is the image of the cycle class of $c$ in $\rH^{2d}_\rc(U_t,\dL(d))^0$ under the edge map, and $K_U/K$ is the finite extension in Proposition \ref{pr:comparison}.
\end{proposition}

\begin{proof}
Let $\sT(d)_\cX$ be the object in $\bD^+(\sfX_\et,\dZ_{p\bullet})$ defined in \cite{Sat07}*{Definition~4.2.4}, which fits into a distinguished triangle
\begin{align}\label{eq:tate1}
\nu_\sfX^{d-1}[-d-1]\to\sT(d)_\cX\to \tau^{\leq d}i^*\rR j_*(\mu_{p^\bullet}^{\otimes d})_X\xrightarrow{+1}\nu_\sfX^{d-1}[-d]
\end{align}
where $\nu_\sfX^{d-1}\in\bM(\sfX_\et,\dZ_{p\bullet})$ is the (projective system of) logarithmic Hodge--Witt sheaves on $\sfX_\et$ defined in \cite{Sat07}*{\S2.2}.

Let $\cC$ be the Zariski closure of $c$ in $\cX$. Let $\{\cH_i\colon\cC_i\to\cX\}$ be the (finite) set of irreducible components of $\supp\cC$, which are projective flat schemes over $O_K$ of pure (absolute) dimension $n-d$. The construction of the refined cycle class of $\cC_i$ in \cite{Sat07}*{Definition~5.1.2} induces a map
\[
\bigoplus_i(\dZ_p/p^\bullet)_{\sfC_i}\to\sfH_i^!\sT(d)_\cX[2d]
\]
in $\bD^+(\sfX_\et,\dZ_{p\bullet})$ and hence a map
\[
\bigoplus_i(\dQ_p)_{\sfC_i}\to\sfH_i^!\sT(d)_{\cX,\dQ_p}[2d]
\]
in $\bD^+(\sfX_\et,\dZ_p)$. As $\supp\cC$ is contained in $\cU_t$, the natural map $\sfH_i^!\sfF_{t!}\sfF_t^*\sT(d)_{\cX,\dQ_p}\to\sfH_i^!\sT(d)_{\cX,\dQ_p}$ is an equivalence. Thus, we obtain a Gysin map
\[
\bigoplus_i\rH^0(\cC_i,\dL)\to\bH^{2d}(\sfX_{\et},\sfF_{t!}\sfF_t^*\sT(d)_{\cX,\dQ_p})_\dL.
\]
Let $\tau_c\in\bH^{2d}(\sfX_{\et},\sfF_{t!}\sfF_t^*\sT(d)_{\cX,\dQ_p})_\dL$ be the image of the cycle $\cC$ under the above map.

Since $d<n<p$, the period map induces equivalences $\sS(d)_\cX\xrightarrow\sim \tau^{\leq d}i^*\rR j_*(\mu_{p^\bullet}^{\otimes d})_X$. Replacing $\tau^{\leq d}i^*\rR j_*(\mu_{p^\bullet}^{\otimes d})_X$ by $\sS(d)_\cX$ in \eqref{eq:tate1} and applying \eqref{eq:adic}, we obtain a distinguished triangle
\[
\nu_{\sfX,\dQ_p}^{d-1}[-d-1]\to\sT(d)_{\cX,\dQ_p}\to \sS(d)_{\cX,\dQ_p}\xrightarrow{+1}\nu_{\sfX,\dQ_p}^{d-1}[-d]
\]
in $\bD^+(\sfX_{\et},\dQ_p)$. Denote by $\sigma_c\in\bH^{2d}(\sfX_{\et},\sfF_{t!}\sfF_t^*\sS(d)_{\cX,\dQ_p})_\dL$ the image of $\tau_c$ under the above sequence, which then belongs to $\bH^{2d}(\sfX_{\et},\sfF_{t!}\sfF_t^*\sS(d)_{\cX,\dQ_p})^\heartsuit_\dL$ since the cycle class of $c$ vanishes in $\rH^{2d}(\ol{X},\dL(d))$ and hence in $\rH^{2d}_\rc(\ol{U_t},\dL(d))$. Now we compute $\rho_{2d}(\sigma_c)$, where $\rho_{2d}$ is defined in \eqref{eq:comparison4}.

In the following diagram
\[
\xymatrix{
\sfF_{t!}\sfF_t^*\sS(d)_{\cX,\dQ_p} \ar[rr]^-{\eqref{eq:comparison6}}\ar[d] &&
\sfF_{t!}\sfF_t^*W\widetilde\omega^\bullet_{\sfX,\dQ_p} \ar[rr]^-{\eqref{eq:witt2}}\ar[d] &&
\widetilde\bomega_{(\sfU,\sfX)} \ar[d]^-{\eqref{eq:witt4}} \\
\sfF_{t!}\sfF_t^*\nu_{\sfX,\dQ_p}^{d-1}[-d] \ar[rr] && \sfF_{t!}\sfF_t^*W\Xi^\bullet_{\sfX,\dQ_p}[-1] \ar[rr]^-{\eqref{eq:witt3}} &&
\bomega^+_{(\sfU,\sfX)}[-1]
}
\]
in $\bD^+(\sfX_{\et},\dQ_p)$, the left square commutes as shown in the proof of \cite{Sat13}*{Proposition~9.10} (followed by taking limit and applying $\sfF_{t!}\sfF_t^*$), and the right square is \eqref{eq:witt5}. It follows that
\[
\rho_{2d}(\sigma_c)\in\ker\(\frac{\rH^{2d-1}_\rig((\sfU_t,\sfX)/W^\circ)_\dL}{N\rH^{2d-1}_\rig((\sfU_t,\sfX)/W^\circ)_\dL}
\to\bH^{2d-1}(\sfX_{\et},\bomega^+_{(\sfU_t,\sfX)})_\dL\).
\]
By Lemma \ref{le:crystalline} and the diagram \eqref{eq:comparison7}, we have $\rho_{2d}(t^*\sigma_c)=t^*\rho_{2d}(\sigma_c)=0$. Finally, by Proposition \ref{pr:comparison}, we have $\r{res}_{K_U}(\alpha_{2d}(t^*\sigma_c))=0$. Since $\r{res}_{K_U}(\alpha_{2d}(t^*\sigma_c))$ coincides with the image of $\r{res}_{K_U}(t^*\beta_c)$ under the map
\[
\rH^1(K_U,\rH^{2d-1}_\rc(\ol{U},\dL(d)))\to\rH^1(K_U,\rH^{2d-1}_\rc(\ol{U},\dQ_p)\otimes_{\dQ_p}\dB_\cris)_\dL,
\]
$\r{res}_{K_U}(t^*\beta_c)$ belongs to $\rH^1_f(K_U,\rH^{2d-1}_\rc(\ol{U},\dL(d)))$. The proposition is proved.
\end{proof}

\begin{proof}[Proof of Theorem \ref{th:crystalline}]
We consider the localized cohomology $\rH^{2d'-1}(\ol{X},\dL(d'))_{\fm'}$, which is a direct summand of $\rH^{2d'-1}(\ol{X},\dL(d'))$ in the category $\bM_K(\dL)$. By the $C_\st$-comparison theorem for $\cX$, Mokrane's weight spectral sequence \cite{Mok93} and \eqref{eq:vanishing}, we know that $\rH^{2d'-1}(\ol{X},\dL(d'))_{\fm'}$ is either zero or a semistable representation of $\rG_K$ pure of weight $-1$ (Definition \ref{de:weight}). In particular, \cite{Nek93}*{Proposition~1.25} implies the following
\begin{itemize}
  \item[$(*)$] For every short exact sequence
      \[
      0\to\dL(1)\to\rE\to\rH^{2d-1}(\ol{X},\dL(d))\to0
      \]
      in $\bM_K(\dL)$ such that $[\rE^\vee(1)]$ belongs to
      \[
      \rH^1_f(K,\rH^{2d'-1}(\ol{X},\dL(d')))\cap\rH^1(K,\rH^{2d'-1}(\ol{X},\dL(d'))_{\fm'}),
      \]
      we have a short exact sequence
      \[
      0 \to \rH^1_f(K^\dag,\dL(1)) \to \rH^1_f(K^\dag,\rE) \to \rH^1_f(K^\dag,\rH^{2d-1}(\ol{X},\dL(d))) \to 0
      \]
      for every finite extension $K^\dag/K$ contained in $\ol{K}$.
\end{itemize}

Let $\dJ'\subseteq\dT$ be the annihilator of
\[
\bigoplus_{d'=1}^{n-1}\Ker\(\rH^1(K,\rH^{2d'-1}(\ol{X},\dL(d')))\to\rH^1(K,\rH^{2d'-1}(\ol{X},\dL(d'))_{\fm'})\).
\]
In particular, we have $[\rE^{t'^*_1c'}]\in\rH^1(K,\rH^{2d'-1}(\ol{X},\dL(d'))_{\fm'})$ for every $t'_1\in\dJ'$. We need to apply Proposition \ref{pr:crystalline} three times.

First, we apply Proposition \ref{pr:crystalline} to the case $c=t'^*_1c'$, $d=d'$, and $\cU=\cX$. Then for every $t'_2\in\dI^{4n-5}$, the class $\beta_{t'^*_2t'^*_1c'}=t'^*_2\beta_{t'^*_1c'}$ belongs to $\rH^1_f(K,\rH^{2d'-1}(\ol{X},\dL(d')))$. In other words, we have
\[
[\rE_{t'^*c'}^\vee(1)]\in\rH^1_f(K,\rH^{2d'-1}(\ol{X},\dL(d')))\cap\rH^1(K,\rH^{2d'-1}(\ol{X},\dL(d'))_{\fm'}),
\]
where $t'\coloneqq t'_1 t'_2 $. By $(*)$, we have a short exact sequence
\begin{align*}
0 \to \rH^1_f(K^\dag,\dL(1)) \to \rH^1_f(K^\dag,\rE_{t'^*c'}) \to \rH^1_f(K^\dag,\rH^{2d-1}(\ol{X},\dL(d))) \to 0
\end{align*}
for every finite extension $K^\dag/K$ contained in $\ol{K}$. Now we denote by $\rH^1_\sharp(K^\dag,\rE_{t'^*c'})$ the inverse image of the subspace $\rH^1_f(K^\dag,\rH^{2d-1}(\ol{X},\dL(d)))$ under the map $\rH^1(K^\dag,\rE_{t'^*c'})\to\rH^1(K^\dag,\rH^{2d-1}(\ol{X},\dL(d)))$.\footnote{In fact, $\rH^1_\sharp(K^\dag,\rE_{t'^*c'})=\rH^1_\st(K^\dag,\rE_{t'^*c'})$.} Then the diagram
\begin{align}\label{eq:comparison8}
\xymatrix{
0  \ar[r] & \rH^1_f(K^\dag,\dL(1)) \ar[r]\ar[d] & \rH^1_f(K^\dag,\rE_{t'^*c'}) \ar[r]\ar[d] & \rH^1_f(K^\dag,\rH^{2d-1}(\ol{X},\dL(d))) \ar[r]\ar@{=}[d] & 0 \\
0  \ar[r] & \rH^1(K^\dag,\dL(1)) \ar[r] & \rH^1_\sharp(K^\dag,\rE_{t'^*c'}) \ar[r] & \rH^1_f(K^\dag,\rH^{2d-1}(\ol{X},\dL(d))) \ar[r] & 0
}
\end{align}
is a pushout of extensions.

Second, we apply Proposition \ref{pr:crystalline} to the case $c=c$, $d=d$, and $\cU=\cX$. Then for every $t_1\in\dI^{4n-5}$, the class $\beta_{t^*_1c}=t^*_1\beta_{c}$ belongs to $\rH^1_f(K,\rH^{2d-1}(\ol{X},\dL(d)))$. In other words, $[\rE^{t^*_1c}]\in\rH^1_f(K,\rH^{2d-1}(\ol{X},\dL(d)))$ and hence $[\rE^{t^*_1c}_{t'^*c'}]\in\rH^1_\sharp(K,\rE_{t'^*c'})$.

Third, we apply Proposition \ref{pr:crystalline} to the case $c=t^*_1c$, $d=d$, and $\cU=\cX\setminus(\supp\cC')^{t'}$. Then for every $t_2\in\dI^{4n-5}$, the class $\r{res}_{K_U}(\beta_{t^*_2t^*_1c})=\r{res}_{K_U}(t^*_2\beta_{t^*_1c})$ belongs to $\rH^1_f(K_U,\rH^{2d-1}_\rc(\ol{U_t},\dL(d)))$, where $t\coloneqq t_1t_2$. Note that the fact that $\cU$ is $d$-dense follows from condition (2), and that $\supp\cC\subseteq\cU_t$ follows from condition (1). Since $[\rE^{t^*c}_{t'^*c'}]$ is the image of $\beta_{t^*c}$ under the pushout map
\[
\rH^1(K,\rH^{2d-1}_\rc(\ol{U_t},\dL(d)))\to\rH^1(K,\rE_{t'^*c'}),
\]
we have $\r{res}_{K_U}([\rE^{t^*c}_{t'^*c'}])\in\rH^1_f(K_U,\rE_{t'^*c'})$. Since the inverse image of $\rH^1_f(K_U,\dL(1))$ under $\r{res}_{K_U}$ coincides with $\rH^1_f(K,\dL(1))$, we conclude that $[\rE^{t^*c}_{t'^*c'}]\in\rH^1_f(K,\rE_{t'^*c'})$ by the diagram \eqref{eq:comparison8} (for $K^\dag=K,K_U$).

From the above discussion, the conclusion of the theorem holds for every pair of elements $t\in(\dI^{4n-5}\setminus\fm)^2$ and $t'\in(\dJ'\setminus\fm')\cdot(\dI^{4n-5}\setminus\fm')$. It is clear that $\dJ'\setminus\fm'\neq\emptyset$. By \eqref{eq:vanishing}, we also have $\dI\setminus\fm\neq\emptyset$ and $\dI\setminus\fm'\neq\emptyset$.

The theorem is proved.
\end{proof}

\subsection{Further preparation}

Let $\bA$ be a $W$-linear additive category. Following Fontaine, we say that a \emph{$(\varphi,N)$-module} in $\bA$ is an object $C$ in $\bA$ with a $W$-semi-linear endomorphism $\varphi_C\colon C\to C$ (called the Frobenius operator) and a $W$-linear endomorphism $N_C\colon C\to C$ (called the monodromy operator) satisfying that $N_C\circ\varphi_C=p\cdot \varphi_C\circ N_C$. A map between $(\varphi,N)$-modules is a map that commutes with both Frobenius operators and monodromy operators.

Suppose that $\bA$ admits $W$-linear tensors. For two $(\varphi,N)$-modules $C$ and $D$ in $\bA$, we equip $C\otimes_WD$ with a $(\varphi,N)$-module structure with the obvious $W$-action, together with
\[
\varphi_{C\otimes_WD}\coloneqq\varphi_C\otimes\varphi_D,\quad
N_{C\otimes_WD}\coloneqq N_C\otimes 1+1\otimes N_D.
\]

When $\bA$ is an abelian category, we say that a $(\varphi,N)$-module is \emph{nilpotent} if it has finite length and the monodromy operator is nilpotent.

\begin{definition}\label{de:descent}
Let $\bD$ be a triangulated category. For a finite complex
\[
\sC_0 \to \sC_1 \to \cdots \to \sC_m
\]
in $\bD$, we have a canonical diagram
\[
\xymatrix{
\sC_{-1} \ar[r]\ar[d] & 0 \ar[d] \\
\sC_0 \ar[r]& \bullet \ar[r]\ar[d]& 0 \ar[d] \\
& \sC_1 \ar[r]& \bullet \ar[r]\ar[d] & \ddots \ar[d] \\
&& \ddots \ar[r] & \bullet \ar[d]\ar[r]& 0 \ar[d] \\
&&& \sC_{m-2} \ar[r]& \bullet \ar[d]\ar[r] & 0 \ar[d] \\
&&&& \sC_{m-1} \ar[r] & \sC_m
}
\]
in $\bD$, in which every square is a homotopy fiber. We call $\sC_{-1}$ the \emph{successive homotopy fiber} of the complex $\sC_0 \to \sC_1 \to \cdots \to \sC_m$.

If $\sC_0 \to \sC_1 \to \cdots \to \sC_m$ is a complex of $(\varphi,N)$-modules in a $W$-linear triangulated category, then $\sC_{-1}$ is canonically a $(\varphi,N)$-module.
\end{definition}

Next, we review some period rings. For every $l\geq 1$, let $R_l$ be the PD envelope of $O_K/p^l$ in $W[t]/p^l$, and let $P_l$ be the $R_l$-ring wit a discrete action by $\rG_K$ defined in \cite{Kat94}*{(3.2)} or \cite{Tsu99}*{\S1.6}. Then $R_\bullet$ and $P_\bullet$ are $(\varphi,N)$-modules in $\bM(W_\bullet)$. Put
\begin{align*}
\dK&\coloneqq\(\varprojlim_lR_l\)\otimes_{W}K_0\in\bM(K_0),\\
\widehat{\dB_\st^+}&\coloneqq\(\varprojlim_lP_l\)\otimes_{W}K_0\in\bM_K(K_0),
\end{align*}
which are again $(\varphi,N)$-modules in relevant categories. See \cite{Bre97}*{\S2~\&~\S4} for an explicit description.\footnote{Our $\dK$ and $\widehat{\dB_\st^+}$ are denoted as $S_{min}$ and $\widehat{B_\st^+}$ in \cite{Bre97}, respectively.}

\begin{lem}\label{le:descent7}
The following holds.
\begin{enumerate}
  \item There is a canonical isomorphism $\dB_\st^+\simeq(\widehat{\dB_\st^+})^{N\mathrm{-nilp}}$.

  \item The ring $\widehat{\dB_\st^+}$ is flat over $\dK$.

  \item The monodromy operator $N\colon \widehat{\dB_\st^+}\to \widehat{\dB_\st^+}$ is surjective.
\end{enumerate}
\end{lem}

\begin{proof}
For (1), this is \cite{Tsu99}*{Proposition~4.1.3}.

For (2), this is \cite{Tsu99}*{Proposition~4.1.5}.

For (3), we have $N=\(\varprojlim_l N_l\)\otimes_{\dZ_p}\dQ_p$, where $N_l$ is the monodromy operator on $P_l$ \cite{Kat94}*{Definition~3.4}. Then (3) follows from the fact that $N_l$ is surjective and $\rR^1\varprojlim_l\Ker N_l=0$ \cite{Kat94}*{Corollary~3.6}.
\end{proof}

For two $(\varphi,N)$-modules $C$ and $D$ in $\bM(K_0)$, the $(\varphi,N)$-module structure on $C\otimes_WD$ clearly descends to $C\otimes_{K_0}D$. Now we consider a nilpotent $(\varphi,N)$-module $D$ in $\bM(K_0)$. By Lemma \ref{le:descent7}(1,3), we have a Frobenius equivariant diagram
\begin{align*}
\xymatrix{
0 \ar[r]& \(D\otimes_{K_0}\dB_\st^+\)^{N=0} \ar@{=}[d]\ar[r]& D\otimes_{K_0}\dB_\st^+ \ar@{^(->}[d] \ar[r]^-{N}& D\otimes_{K_0}\dB_\st^+ \ar@{^(->}[d]\ar[r]& 0 \\
0 \ar[r]& \(D\otimes_{K_0}\dB_\st^+\)^{N=0} \ar[r]& D\otimes_{K_0}\widehat{\dB_\st^+} \ar[r]^-{N}& D\otimes_{K_0}\widehat{\dB_\st^+} \ar[r]& 0
}
\end{align*}
in $\bM_K(K_0)$, in which the two rows are short exact sequences. It induces a diagram
\[
\xymatrix{
\dfrac{D}{N D} \ar@{^(->}[r]\ar[d]& \rH^1\(K,\(D\otimes_{K_0}\dB_\st^+\)^{N=0}\) \ar@{=}[d] \\
\dfrac{D\otimes_{K_0}\dK}{N(D\otimes_{K_0}\dK)} \ar@{^(->}[r]& \rH^1\(K,\(D\otimes_{K_0}\dB_\st^+\)^{N=0}\)
}
\]
of edge maps as $(\dB_\st^+)^{\rG_K}=K_0$ and $(\widehat{\dB_\st^+})^{\rG_K}=\dK$ \cite{Bre97}*{Corollaire~4.1.3}.

\begin{lem}\label{le:preparation}
For every integer $r$, the restricted edge map
\[
\(\frac{D\otimes_{K_0}\dK}{N(D\otimes_{K_0}\dK)}\)^{\varphi=p^r}\to\rH^1\(K,\(D\otimes_{K_0}\dB_\st^+\)^{N=0}\)
\]
factors through the map
\[
\frac{D\otimes_{K_0}\dK}{N(D\otimes_{K_0}\dK)}\to\frac{D}{N D}
\]
induced by the specialization map $\dK\to K_0$ at $t=0$.
\end{lem}

\begin{proof}
This follows from the same proof of \cite{Lan99}*{Lemma~2.6}.\footnote{In \cite{Lan99}, the author works with $\widehat{K\langle t\rangle}$, which is different from our $\dK$ when $K/\dQ_p$ is ramified. However, such difference will not affect the proof in view of the explicit description of $\dK$ in \cite{Bre97}*{\S4}.}
\end{proof}

\begin{remark}
Lemma \ref{le:preparation} is certainly wrong without the restriction to the part $\varphi=p^r$ since the specialization map $\dK/N\dK\to K_0$ has a large kernel. However, we do not know whether the edge map $D\otimes_{K_0}\dK\to\rH^1\(K,(D\otimes_{K_0}\dB_\st)^{N=0}\)$ factors though the quotient $D/ND$, which is equivalent to the inclusion $\Ker(\dK\to K_0)\subseteq N(\widehat{\dB_\st})^{\rG_K}$ where $\widehat{\dB_\st}\coloneqq\widehat{\dB_\st^+}\otimes_{\dB_\cris^+}\dB_\cris$. See \cite{Bre97}*{\S5} for the mystery of $(\widehat{\dB_\st})^{\rG_K}$.
\end{remark}

\subsection{Proof of Proposition \ref{pr:comparison}}
\label{ss:alteration}

This subsection is devoted to the proof of Proposition \ref{pr:comparison}, for which we use de Jong's alterations. We may assume $1\leq d\leq n-1$.

By \cite{deJ96}*{Theorem~6.5}, we may find a finite extension $K'$ of $K$ (depending on $U$) contained in $\ol{K}$, a projective strictly semistable scheme $\cX'$ over $O_{K'}$ and a generically finite morphisms $\cA\colon\cX'\to\cX$ over $O_K$,\footnote{The letter $\cA$ stands for \emph{alteration}.} such that $(\cX',A^{-1}U)$ is a strict semistable pair \cite{deJ96}*{\S6.3}. Let $\cF'\colon\cU'\to\cX'$ be the open subscheme with $U'=A^{-1}U$ and such that $\cX'\setminus\cU'$ is flat over $O_{K'}$. Note that $\cU'$ may strictly contain $\cA^{-1}\cU$. We have objects with respect to $(\cX',\cU')$ and we put a \emph{prime} for their notations.

\begin{lem}\label{le:comparison1}
Suppose that $n<p$. Then for every $q\geq0$, the composite map
\begin{align*}
\bH^q(\sfX'_{\et},\sfF'_!\sfF'^*\sS(d)_{\cX',\dQ_p})
&\to\bH^q(\sfX'_{\et},\widetilde\bomega_{(\sfU',\sfX')}) \\
&\to\bH^q(\sfX'_{\et},\bomega_{(\sfU',\sfX')})=\rH^q_\rig((\sfU',\sfX')/W'^\circ)
\end{align*}
(Definition \ref{de:rigid}) vanishes on $\bH^q(\sfX'_{\et},\sfF'_!\sfF'^*\sS(d)_{\cX',\dQ_p})^\heartsuit$.
\end{lem}

Note that this lemma does not follow from Lemma \ref{le:comparison} even for $q=2d$, since $\cU'$ is not necessarily $d$-dense anymore. By Lemma \ref{le:comparison1}, we have the map
\begin{align}\label{eq:comparison2}
\rho'_q\colon\bH^q(\sfX'_{\et},\sfF'_!\sfF'^*\sS(d)_{\cX',\dQ_p})^\heartsuit
\to\frac{\rH^{q-1}_\rig((\sfU',\sfX')/W'^\circ)}{N\rH^{q-1}_\rig((\sfU',\sfX')/W'^\circ)}
\end{align}
similar to \eqref{eq:comparison4} for every $q\geq 0$.

\begin{lem}\label{le:comparison2}
Suppose that $n<p$. Then
\[
\Ker(\rho'_q)\subseteq\Ker(\alpha'_q)
\]
holds for every $q\geq 0$.
\end{lem}

We now prove Proposition \ref{pr:comparison} assuming the above two lemmas, whose proofs are postponed later.

\begin{proof}[Proof of Proposition \ref{pr:comparison}]
Take $K_U=K'$. It is also clear that we may take $K_X$ to be $K$. We have the commutative diagram
\[
\xymatrix{
\dfrac{\rH^{2d-1}_\rig((\sfU,\sfX)/W^\circ)}{N\rH^{2d-1}_\rig((\sfU,\sfX)/W^\circ)}
\ar[r]^-{\sfA^*} & \dfrac{\rH^{2d-1}_\rig((\sfU',\sfX')/W'^\circ)}{N\rH^{2d-1}_\rig((\sfU',\sfX')/W'^\circ)} \\
\bH^{2d}(\sfX_{\et},\sfF_!\sfF^*\sS(d)_{\cX,\dQ_p})^\heartsuit \ar[r]^-{\sfA^*} \ar[u]\ar[d] &
\bH^{2d}(\sfX'_{\et},\sfF'_!\sfF'^*\sS(d)_{\cX',\dQ_p})^\heartsuit \ar[u]\ar[d] \\
\rH^1(K,\rH^{2d-1}_\rc(\ol{U},\dQ_p)\otimes_{\dQ_p}\dB_\cris) \ar[r]^-{A^*} &
\rH^1(K',\rH^{2d-1}_\rc(\ol{U'},\dQ_p)\otimes_{\dQ_p}\dB_\cris)
}
\]
of $\dQ_p$-vector spaces.\footnote{Note that $\sfU'$ may properly contain $\sfA^{-1}\sfU$. By $\sfA^*$, we mean the restriction map with a possible composition of the pushforward map along the inclusion $\sfA^{-1}\sfU\subseteq\sfU'(\subseteq\sfX')$ of log rigid cohomology (Definition \ref{de:rigid}) or \'{e}tale cohomology with proper support.} By Lemma \ref{le:comparison2}, to prove the proposition, it suffices to show that the map
\[
\rH^1(K',\rH^{2d-1}_\rc(\ol{U},\dQ_p)\otimes_{\dQ_p}\dB_\cris)
\to\rH^1(K',\rH^{2d-1}_\rc(\ol{U'},\dQ_p)\otimes_{\dQ_p}\dB_\cris)
\]
is injective. However, this follows from the fact that the map $\rH^q_\rc(\ol{U},\dQ_p)\to\rH^q_\rc(\ol{U'},\dQ_p)$ in the category $\bM_{K'}(\dQ_p)$ admits a section, which is a consequence of the usual Poincar\'{e} duality for \'{e}tale cohomology of $\ol{U}$ and $\ol{U'}$. The proposition is proved.
\end{proof}

It remains to show Lemma \ref{le:comparison1} and Lemma \ref{le:comparison2}. Since we will only study $(\cX',\cU')$ from now on, we will suppress the \emph{prime} from all notation to release some burden. Put $\cV\coloneqq\cX\setminus\cU$, and for every $h\geq 1$, let $\cV^{(h)}$ be the disjoint union of intersections of $h$ different irreducible components of $\cV$, which is either empty or a strictly semistable scheme over $O_K$ of pure (absolute) dimension $n-h$. For notational convenience, we also put $\cV^{(0)}\coloneqq\cX$ and $\cV^{(-1)}\coloneqq(\cU,\cX)$. Denote by $\cG^{(h)}\colon\cV^{(h)}\to\cX$ the obvious morphism for $h\geq 0$.

\begin{lem}\label{le:log}
For every $h\geq 0$, the pullback of the log structure $L_\cX^\sfX$ for $\cX$ to $\cV^{(h)}$ coincides with $L_{\cV^{(h)}}^{\sfV^{(h)}}$.
\end{lem}

\begin{proof}
The question is local in the Zariski topology. By \cite{deJ96}*{\S6.4}, Zariski locally $\cX$ is smooth over $O_K[t_1,\dots,t_i,s_1,\dots,s_j]/(t_1\cdots t_i-\varpi)$. We may assume $j\geq h$ since otherwise $\cV^{(h)}$ is empty in this chart. It suffices to consider the open and closed subscheme $\cT$ of $\cV^{(h)}$ defined by $s_{j-h+1}=\cdots=s_j=0$. Now locally $L_\cX^\sfX$ and $L_\cT^\sfT$ are the log structures associated with the pre-log structures $\dN^i\to\sO_\cX$ and $\dN^i\to\sO_\cT$ sending $1$ in the $i'$-th factor to the pullback of $t_{i'}$ for $1\leq i'\leq i$, respectively. The lemma follows immediately.
\end{proof}

Our first step is to construct, for every $r\geq 0$, a syntomic complex $\sS(r)_{\ol{\cV^{(-1)}}}\in\bD^+_{\rG_K}(\ol\sfX_{\et},\dZ_{p\bullet})$, together with a period map
\begin{align}\label{eq:descent1}
\sS(r)_{\ol{\cV^{(-1)}}}\to \ol{i}^*\rR\ol{j}_*\ol{F}_!(\mu_{p^\bullet}^{\otimes r})_{\ol{U}}
\end{align}
when $0\leq r<p-1$, that becomes an equivalence when $n-1\leq r< p-1$. The construction is inspired by the observation in the following remark.

\begin{remark}\label{re:descent}
After choosing an order on the (finite) set of irreducible components of $\cV$, we have an exact sequence
\begin{align*}
0 \to \ol{F}_!(\mu_{p^\bullet}^{\otimes r})_{\ol{U}} &\to \ol{G^{(0)}}_*(\mu_{p^\bullet}^{\otimes r})_{\ol{V^{(0)}}}
\to \ol{G^{(1)}}_*(\mu_{p^\bullet}^{\otimes r})_{\ol{V^{(1)}}} \\
&\to \ol{G^{(2)}}_*(\mu_{p^\bullet}^{\otimes r})_{\ol{V^{(2)}}} \to\cdots \to \ol{G^{(n-1)}}_*(\mu_{p^\bullet}^{\otimes r})_{\ol{V^{(n-1)}}}
\to 0
\end{align*}
in $\bM_{\rG_K}(\ol{X}_{\et},\dZ_{p\bullet})$. In particular, $\ol{F}_!(\mu_{p^\bullet}^{\otimes r})_{\ol{U}}$ is canonically equivalent to the successive homotopy fiber of the complex
\[
\ol{G^{(0)}}_*(\mu_{p^\bullet}^{\otimes r})_{\ol{V^{(0)}}}
\to \ol{G^{(1)}}_*(\mu_{p^\bullet}^{\otimes r})_{\ol{V^{(1)}}} \to \ol{G^{(2)}}_*(\mu_{p^\bullet}^{\otimes r})_{\ol{V^{(2)}}} \to\cdots \to \ol{G^{(n-1)}}_*(\mu_{p^\bullet}^{\otimes r})_{\ol{V^{(n-1)}}}.
\]
\end{remark}

In order to unify the notation, we put
\begin{align*}
\sN(r)_{\ol{\cV^{(-1)}}}&\coloneqq\ol{i}^*\rR\ol{j}_*\ol{F}_!(\mu_{p^\bullet}^{\otimes r})_{\ol{U}}\in\bD^+_{\rG_K}(\ol\sfX_{\et},\dZ_{p\bullet}),\\
\sN(r)_{\cV^{(-1)}}&\coloneqq i^*\rR j_*F_!(\mu_{p^\bullet}^{\otimes r})_{U}\in\bD^+(\sfX_{\et},\dZ_{p\bullet}),
\end{align*}
and
\begin{align*}
\sN(r)_{\ol{\cV^{(h)}}}&\coloneqq\ol{i}^*\rR\ol{j}_*\ol{G^{(h)}}_*(\mu_{p^\bullet}^{\otimes r})_{\ol{V^{(h)}}}\in\bD^+_{\rG_K}(\ol\sfX_{\et},\dZ_{p\bullet}),\\
\sN(r)_{\cV^{(h)}}&\coloneqq i^*\rR j_*G^{(h)}_*(\mu_{p^\bullet}^{\otimes r})_{V^{(h)}}\in\bD^+(\sfX_{\et},\dZ_{p\bullet}),
\end{align*}
for $h\geq 0$.\footnote{The letter $\sN$ stands for \emph{nearby}.}

To define $\sS(r)_{\ol{\cV^{(-1)}}}$, we need to consider all extensions of $K$ in $\ol{K}$. A \emph{Galois embedding system} for $O_{\ol{K}}^\can/W[t]^\circ$ consists of
\begin{itemize}
  \item an increasing tower $K=K_1\subseteq K_2\subseteq K_3\subseteq\cdots$ of finite Galois extensions of $K$ with $\bigcup_m K_m=\ol{K}$ (and we regard $O_{K_m}^\can$ as a log-scheme over $O_K^\can$),

  \item for every $m\geq 1$, an embedding system $\{O_{K_m}^\can\hookrightarrow(\cZ_m^\flat,N_m^\flat)\}$ for $O_{K_m}^\can/W[t]^\circ$ with a compatible action of $\Gal(K_m/K)$ that fits into a commutative diagram
      \[
      \xymatrix{
      \vdots\ar[d] & \vdots\ar[d] \\
      O_{K_{m+1}}^\can \ar[d]\ar@{^(->}[r] & (\cZ_{m+1}^\flat,N_{m+1}^\flat) \ar[d] \\
      O_{K_m}^\can \ar[d]\ar@{^(->}[r] & (\cZ_m^\flat,N_m^\flat) \ar[d] \\
      \vdots & \vdots
      }
      \]
      of log-schemes over $W[t]^\circ$ that is $\rG_K$-equivariant.
\end{itemize}
It is clear that Galois embedding systems for $O_{\ol{K}}^\can/W[t]^\circ$ exist.

We now choose a Galois embedding system for $O_{\ol{K}}^\can/W[t]^\circ$ as above, and write $\kappa=\kappa_1\subseteq \kappa_2\subseteq \kappa_3\subseteq\cdots$ for the induced tower of residue fields. We also choose an embedding system $\{(\cX^\star,L^\star)\hookrightarrow(\cZ^\star,N^\star)\}$ for $(\cX,L_\cX^\sfX)/W[t]^\circ$. For $m\geq 1$, put
\[
(\cZ^\star_m,N^\star_m)\coloneqq(\cZ^\star,N^\star)\times_{W[t]^\circ}(\cZ_m^\flat,N_m^\flat).
\]

Let $\cT$ be an irreducible component of $\cV^{(h)}$ for some $h\geq 0$. Lemma \ref{le:log} implies that
\[
\left\{(\cT,L_\cT^\sfT)\times_{(\cX,L_\cX^\sfX)}(\cX^\star,L^\star)\hookrightarrow
(\cZ^\star,N^\star)
\right\}
\]
is an embedding system for $(\cT,L_\cT^\sfT)/W[t]^\circ$. For every $m\geq 1$, let $\cE^\star_{m,l}$ be the PD envelope of $\cT^\star\otimes_{O_K}O_{K_m}/p^l$ in $\cZ^\star_m\otimes\dZ/p^l$. For $i\geq 1$, let $\sJ_{m,l}^{[i]}\subseteq\sO_{\cE^\star_{m,l}}$ be the $i$-th divided power of the ideal $\sJ_{m,l}\coloneqq\Ker\(\sO_{\cE^\star_{m,l}}\to\sO_{\cT^\star\otimes_{O_K}O_{K_m}/p^l}\)$. For $i\leq 0$, we put $\sJ_{m,l}^{[i]}\coloneqq\sO_{\cE^\star_{m,l}}$. We have the complex
\begin{align}\label{eq:descent5}
&B(r)_{\cT,m}\colon \\
&\sJ_{m,\bullet}^{[r]}\to\sJ_{m,\bullet}^{[r-1]}\otimes_{\sO_{\cZ^\star_m}}\Omega^1_{(\cZ^\star_m,N^\star_m)/W^\triv}
\to\sJ_{m,\bullet}^{[r-2]}\otimes_{\sO_{\cZ^\star_m}}\Omega^2_{(\cZ^\star_m,N^\star_m)/W^\triv}
\to\cdots\notag
\end{align}
regarded as an object in $\bC^+_{\rG_K}((\sfT\times_\sfX\sfX^\star\otimes_\kappa\kappa_m)_{\et},W_\bullet)$, where $\sJ_{m,\bullet}^{[r]}$ is placed in degree $0$. Put
\begin{align}\label{eq:descent6}
B(r)_{\ol\cT}\coloneqq\varinjlim_m B(r)_{\cT,m}\res_{\ol\sfT\times_\sfX\sfX^\star},
\end{align}
where the colimit is taken in the abelian category $\bC^+_{\rG_K}((\ol\sfT\times_\sfX\sfX^\star)_{\et},W_\bullet)$. For every $h\geq0$, we put
\[
B(r)_{\ol{\cV^{(h)}}}\coloneqq\bigoplus_{\cT}B(r)_{\ol\cT},
\]
where the direct sum is taken over all irreducible components of $\cV^{(h)}$, regarded as an element in $\bC^+_{\rG_K}(\ol{\sfX^\star}_{\et},W_\bullet)$ via pushforward along closed immersions $\sfT\to\sfX$. Then parallel to Remark \ref{re:descent}, we have a complex
\[
B(r)_{\ol{\cV^{(0)}}}\to B(r)_{\ol{\cV^{(1)}}} \to B(r)_{\ol{\cV^{(2)}}} \to \cdots \to B(r)_{\ol{\cV^{(n-1)}}}
\]
in $\bC^+_{\rG_K}(\ol{\sfX^\star}_{\et},W_\bullet)$.

Take $h\geq 0$. Put $C_{\ol{\cV^{(h)}}/W^\triv}\coloneqq B(0)_{\ol{\cV^{(h)}}}$. Then we have the canonical map $B(r)_{\ol{\cV^{(h)}}}\to C_{\ol{\cV^{(h)}}/W^\triv}$ given by the inclusions $\sJ_{m,l}^{[i]}\to\sJ_{m,l}^{[0]}$. We also have the crystalline complex $C_{\ol{\cV^{(h)}}/W[t]^\circ}$, which is obtained in the same way as $C_{\ol{\cV^{(h)}}/W^\triv}$ except that in the definition of $B(0)_{\cT,m}$ \eqref{eq:descent5}, we replace $\Omega^q_{(\cZ^\star_m,N^\star_m)/W^\triv}$ by $\Omega^q_{(\cZ^\star_m,N^\star_m)/W[t]^\circ}$. We have natural Frobenius operators on both $C_{\ol{\cV^{(h)}}/W^\triv}$ and $C_{\ol{\cV^{(h)}}/W[t]^\circ}$, and a distinguished triangle
\begin{align}\label{eq:triangle1}
C^\triangle_{\ol{\cV^{(h)}}/W[t]^\circ}&\colon\quad
C_{\ol{\cV^{(h)}}/W[t]^\circ}[-1] \to  C_{\ol{\cV^{(h)}}/W^\triv} \to C_{\ol{\cV^{(h)}}/W[t]^\circ} \xrightarrow{N} C_{\ol{\cV^{(h)}}/W[t]^\circ}
\end{align}
in $\bD^+_{\rG_K}(\ol{\sfX^\star}_{\et},W_\bullet)$, where the first arrow is given by $\wedge\rd\log t$, the second arrow is the canonical one (which is Frobenius equivariant), and the third arrow is the connecting map, such that $C_{\ol{\cV^{(h)}}/W[t]^\circ}$ becomes a $(\varphi,N)$-module in $\bD^+_{\rG_K}(\ol{\sfX^\star}_{\et},W_\bullet)$. We define $S(r)_{\ol{\cV^{(h)}}}$ to be the homotopy fiber of the map
\[
1-p^{-r}\varphi_{r+\bullet}\colon B(r)_{\ol{\cV^{(h)}}}\to C_{\ol{\cV^{(h)}}/W^\triv}
\]
(see \cite{Tsu00}*{Page~540} for more details) in the category $\bD^+_{\rG_K}(\ol{\sfX^\star}_{\et},\dZ_{p\bullet})$.

Let $B(r)_{\ol{\cV^{(-1)}}}$ be the successive homotopy fiber (Definition \ref{de:descent}) of the complex
\[
B(r)_{\ol{\cV^{(0)}}}\to B(r)_{\ol{\cV^{(1)}}} \to B(r)_{\ol{\cV^{(2)}}} \to \cdots \to B(r)_{\ol{\cV^{(n-1)}}}
\]
in $\bD^+_{\rG_K}(\ol{\sfX^\star}_{\et},W_\bullet)$, and similarly for $C_{\ol{\cV^{(-1)}}/W^\triv}$, $C_{\ol{\cV^{(-1)}}/W[t]^\circ}$, and $S(r)_{\ol{\cV^{(-1)}}}$. Then $C_{\ol{\cV^{(-1)}}/W^\triv}$ is a $(\varphi,N)$-module and we have a similar distinguished triangle \eqref{eq:triangle1} for $h=-1$.

Finally, for $h\geq-1$, put
\begin{align*}
\sS(r)_{\ol{\cV^{(h)}}}&\coloneqq\rR\ol\sfu_*S(r)_{\ol{\cV^{(h)}}},\\
\sC_{\ol{\cV^{(h)}}/W^\triv}&\coloneqq\rR\ol\sfu_*C_{\ol{\cV^{(h)}}/W^\triv},\\
\sC_{\ol{\cV^{(h)}}/W[t]^\circ}&\coloneqq\rR\ol\sfu_*C_{\ol{\cV^{(h)}}/W[t]^\circ},
\end{align*}
in $\bD^+_{\rG_K}(\ol\sfX_{\et},\dZ_{p\bullet})$, $\bD^+_{\rG_K}(\ol\sfX_{\et},W_\bullet)$, and $\bD^+_{\rG_K}(\ol\sfX_{\et},W_\bullet)$, respectively. In particular, we have a canonical map
\begin{align}\label{eq:descent7}
\ol\xi_r\colon\sS(r)_{\ol{\cV^{(h)}}}\to\sC_{\ol{\cV^{(h)}}/W^\triv}
\end{align}
in $\bD^+_{\rG_K}(\ol\sfX_{\et},\dZ_{p\bullet})$, and a distinguished triangle
\begin{align}\label{eq:triangle2}
\sC^\triangle_{\ol{\cV^{(h)}}/W[t]^\circ}&\colon\quad
\sC_{\ol{\cV^{(h)}}/W[t]^\circ}[-1] \to  \sC_{\ol{\cV^{(h)}}/W^\triv} \to \sC_{\ol{\cV^{(h)}}/W[t]^\circ} \xrightarrow{N} \sC_{\ol{\cV^{(h)}}/W[t]^\circ}
\end{align}
in $\bD^+_{\rG_K}(\ol\sfX_{\et},W_\bullet)$.

When $0\leq r<p-1$, the usual period maps for $\cV^{(h)}$ with $h\geq 0$ give a commutative diagram
\[
\xymatrix{
\sS(r)_{\ol{\cV^{(0)}}} \ar[r]\ar[d] &
\ol{\sfG^{(0)}}_*\ol{i^{(0)}}^*\rR\ol{j^{(0)}}_*(\mu_{p^\bullet}^{\otimes r})_{\ol{V^{(0)}}} \ar[d] \\
\sS(r)_{\ol{\cV^{(1)}}} \ar[r]\ar[d] &
\ol{\sfG^{(1)}}_*\ol{i^{(1)}}^*\rR\ol{j^{(1)}}_*(\mu_{p^\bullet}^{\otimes r})_{\ol{V^{(1)}}} \ar[d] \\
\vdots \ar[d] & \vdots \ar[d] \\
\sS(r)_{\ol{\cV^{(n-1)}}}  \ar[r]&
\ol{\sfG^{(n-1)}}_*\ol{i^{(n-1)}}^*\rR\ol{j^{(n-1)}}_*(\mu_{p^\bullet}^{\otimes r})_{\ol{V^{(n-1)}}}
}
\]
in $\bD^+_{\rG_K}(\ol\sfX_{\et},\dZ_{p\bullet})$, where $i^{(h)}\colon\sfV^{(h)}\to\cV^{(h)}$ and $j^{(h)}\colon V^{(h)}\to\cV^{(h)}$ denote the special and generic fibers, respectively, for $h\geq 0$. However, since the natural map
\[
\sN(r)_{\ol{\cV^{(h)}}}=\ol{i}^*\rR\ol{j}_*\ol{G^{(h)}}_*(\mu_{p^\bullet}^{\otimes r})_{\ol{V^{(h)}}} \to \ol{\sfG^{(h)}}_*\ol{i^{(h)}}^*\rR\ol{j^{(h)}}_*(\mu_{p^\bullet}^{\otimes r})_{\ol{V^{(h)}}}
\]
is an equivalence for every $h\geq 0$, we obtain the period map
\begin{align}\label{eq:descent2}
\ol\pi_r\colon\sS(r)_{\ol{\cV^{(h)}}}\to\sN(r)_{\ol{\cV^{(h)}}}
\end{align}
in $\bD^+_{\rG_K}(\ol\sfX_{\et},\dZ_{p\bullet})$ for every $h\geq -1$ by Remark \ref{re:descent} and the process of taking successive homotopy fibers. If $n-1\leq r< p-1$, then \eqref{eq:descent2} is an equivalence. The desired map \eqref{eq:descent1} is simply \eqref{eq:descent2} for $h=-1$.

To proceed, we need versions of syntomic and crystalline complexes for $\cV^{(h)}$ rather than $\ol{\cV^{(h)}}$. The construction is similar to $\sS(r)_{\ol{\cV^{(h)}}}$ and $\sC_{\ol{\cV^{(h)}}/W^\triv}$ but only taking $m=1$ without the colimit \eqref{eq:descent6}. More precisely, for $h\geq -1$, we have
\begin{itemize}
  \item $\sS(r)_{\cV^{(h)}}$ in $\bD^+(\sfX_\et,\dZ_{p\bullet})$, which is obtained in the same way as $\sS(r)_{\ol{\cV^{(h)}}}$ but only taking $m=1$,\footnote{In particular, $\sS(d)_{\cV^{(0)}}$ coincides with $\sS(d)_\cX$ from \S\ref{ss:comparison}.}

  \item $\sC_{\cV^{(h)}/W^\triv}$ in $\bD^+(\sfX_\et,W_\bullet)$, which is obtained in the same way as $\sC_{\ol{\cV^{(h)}}/W^\triv}$ but only taking $m=1$,\footnote{In particular, $\sC_{\cV^{(0)}/W^\triv}$ coincides with $\sC_{l,\cX/W^\triv}$ from \S\ref{ss:comparison}.}

  \item $\sC_{\cV^{(h)}/W[t]^\circ}$ in $\bD^+(\sfX_\et,W_\bullet)$, which is obtained in the same way as $\sC_{\ol{\cV^{(h)}}/W[t]^\circ}$ but only taking $m=1$;\footnote{In particular, $\sC_{\cV^{(0)}/W[t]^\circ}$ coincides with $\sC_{l,\cX/W[t]^\circ}$ from \S\ref{ss:comparison}.} it is a $(\varphi,N)$-module,

  \item $\widetilde\sC_{\sfV^{(h)}/W^\circ}$ and $\sC_{\sfV^{(h)}/W^\circ}$ in $\in\bD^+(\sfX_\et,W_\bullet)$, which are obtained after we replace $B(0)_{\cT,1}$ \eqref{eq:descent5} by the following complexes
      \[
      C_{\sfT/W^\circ}\colon
      \sO_{\cD_\bullet^\star}\to\sO_{\cD_\bullet^\star}\otimes_{\sO_{\cZ^\star}}\Omega^1_{(\cZ^\star,N^\star)/W^\circ}
      \to\sO_{\cD_\bullet^\star}\otimes_{\sO_{\cZ^\star}}\Omega^2_{(\cZ^\star,N^\star)/W^\circ}\to\cdots
      \]
      and
      \[
      C_{\sfT/W^\circ}\colon
      \sO_{\cD_\bullet^\star}\to\sO_{\cD_\bullet^\star}\otimes_{\sO_{\cY^\star}}\Omega^1_{(\cY^\star,M^\star)/W^\circ}
      \to\sO_{\cD_\bullet^\star}\otimes_{\sO_{\cY^\star}}\Omega^2_{(\cY^\star,M^\star)/W^\circ}\to\cdots
      \]
      respectively, where $(\cY^\star,M^\star)\coloneqq(\cZ^\star,N^\star)\times_{W[t]^\circ}W^\circ$ as in \S\ref{ss:witt} and $\cD^\star_l$ denotes the PD envelope of $\sfT$ in $\cY^\star_l$ for $l\geq 1$.\footnote{In particular, $\widetilde\sC_{\sfV^{(0)}/W^\circ}$ and $\sC_{\sfV^{(0)}/W^\circ}$ coincide with $\widetilde\sC_{(\sfX,L_\cX^\sfX\res_\sfX)/W^\circ}$ and $\sC_{(\sfX,L_\cX^\sfX\res_\sfX)/W^\circ}$ from \S\ref{ss:witt}, respectively.}
\end{itemize}
By construction, we have maps
\begin{align*}
\xi_r&\colon\sS(r)_{\cV^{(h)}}\to\sC_{\cV^{(h)}/W^\triv},\\
\pi_r&\colon\sS(r)_{\cV^{(h)}}\to\sN(r)_{\cV^{(h)}},
\end{align*}
in $\bD^+(\sfX_{\et},\dZ_{p\bullet})$ similar to \eqref{eq:descent7} and \eqref{eq:descent2}, and a distinguished triangle
\begin{align}\label{eq:triangle3}
\sC^\triangle_{\cV^{(h)}/W[t]^\circ}&\colon\quad
\sC_{\cV^{(h)}/W[t]^\circ}[-1] \to  \sC_{\cV^{(h)}/W^\triv} \to \sC_{\cV^{(h)}/W[t]^\circ} \xrightarrow{N} \sC_{\cV^{(h)}/W[t]^\circ}
\end{align}
in $\bD^+(\sfX_{\et},\dZ_{p\bullet})$ similar to \eqref{eq:triangle2}, without \emph{bar}.

In order to prove Lemma \ref{le:comparison1} and Lemma \ref{le:comparison2}, we need to connect the syntomic cohomology to the log rigid cohomology via crystalline complexes we have just constructed. By construction, we have a commutative diagram
\begin{align}\label{eq:comparison9}
\xymatrix{
\sS(r)_{\cV^{(h)},\dQ_p} \ar[r]^-{\xi_r}& \sC_{\cV^{(h)}/W^\triv,K_0} \ar[r]\ar[d]& \sC_{\cV^{(h)}/W[t]^\circ,K_0} \ar[d] \\
& \widetilde\sC_{\sfV^{(h)}/W^\circ,K_0} \ar[r] & \sC_{\sfV^{(h)}/W^\circ,K_0}
}
\end{align}
in $\bD^+(\sfX_{\et},\dQ_p)$ for every $r\geq 0$ and every $h\geq -1$, similar to \eqref{eq:comparison6}. Now we study the cohomology of various crystalline complexes.

We start from $\sC_{\sfV^{(h)}/W^\circ}$. Note that, as in \S\ref{ss:witt}, for $h\geq0$, the object $\sC_{\sfV^{(h)}/W^\circ}$ is equivalent to the modified de Rham--Witt complex $W\omega^\bullet_{\sfV^{(h)}}$. By the construction and \eqref{eq:witt1}, we have
\begin{itemize}
  \item a distinguished triangle
     \[
     \sC^\triangle_{\sfV^{(h)}/W^\circ}\colon\quad
     \sC_{\sfV^{(h)}/W^\circ}[-1]\to\widetilde\sC_{\sfV^{(h)}/W^\circ}
     \to\sC_{\sfV^{(h)}/W^\circ}\xrightarrow{N}\sC_{\sfV^{(h)}/W^\circ}
     \]
     in $\bD^+(\sfX_{\et},W_\bullet)$ for every $h\geq0$ hence also $h=-1$, similar to \eqref{eq:triangle3}, so that $\sC_{\sfV^{(h)}/W^\circ}$ is a $(\varphi,N)$-module,

  \item a commutative diagram
     \begin{align}\label{eq:descent8}
     \xymatrix{
     \bomega^\triangle_{\sfV^{(-1)}} \ar[r]\ar[d]& \sfG^{(0)}_*\bomega^\triangle_{\sfV^{(0)}}  \ar[r]\ar[d]^-{\simeq}& \sfG^{(1)}_*\bomega^\triangle_{\sfV^{(1)}} \ar[r]\ar[d]^-{\simeq}& \cdots \ar[r]& \sfG^{(n-1)}_*\bomega^\triangle_{\sfV^{(n-1)}} \ar[d]^-{\simeq} \\
     \sC^\triangle_{\sfV^{(-1)}/W^\circ,K_0} \ar[r]& \sC^\triangle_{\sfV^{(0)}/W^\circ,K_0} \ar[r]& \sC^\triangle_{\sfV^{(1)}/W^\circ,K_0} \ar[r]& \cdots \ar[r]& \sC^\triangle_{\sfV^{(n-1)}/W^\circ,K_0}
     }
     \end{align}
     of distinguished triangles in $\bD^+(\sfX_{\et},K_0)$, in which terms in the top row are from \eqref{eq:monodromy1}; note that all vertical arrows starting from the second are equivalences.
\end{itemize}

\begin{lem}\label{le:descent1}
The first vertical arrow in \eqref{eq:descent8} is also an equivalence. In particular, for every $h\geq -1$ and $q\geq 0$, we have a canonical isomorphism
\[
\rH^q_\rig(\sfV^{(h)},W^\circ)\simeq\bH^q(\sfX_\et,\sC_{\sfV^{(h)}/W^\circ,K_0})
\]
of $K_0$-vector spaces which commutes with monodromy operators.
\end{lem}

\begin{proof}
It suffices to show that the map $\bomega_{\sfV^{(-1)}}=\bomega_{(\sfU,\sfX)}$ is the successive homotopy fiber of the complex
\[
\sfG^{(0)}_*\bomega_{\sfV^{(0)}} \to \sfG^{(1)}_*\bomega_{\sfV^{(1)}} \to \cdots \to \sfG^{(n-1)}_*\bomega_{\sfV^{(n-1)}}.
\]
However, this follows from the easy fact that for an embedding system $\{(\sfX^\star,L^\star)\hookrightarrow(\cY^\star,M^\star)\}$ for $(\sfX,L_\cX^\sfX\res_\sfX)/W^\circ$, the complex
\begin{align*}
0 & \to \tf^!_{(\sfU^\star,\sfX^\star)}
\(\Omega^\bullet_{(\cY^\star,M^\star)/W^\circ}\otimes_{\sO_{\cY^\star}}\sO_{]\sfX^\star[_{\fY^\star}}\)
\to \Omega^\bullet_{(\cY^\star,M^\star)/W^\circ}\otimes_{\sO_{\cY^\star}}\sO_{]\sfX^\star[_{\fY^\star}} \\
& \to \Omega^\bullet_{(\cY^\star,M^\star)/W^\circ}\otimes_{\sO_{\cY^\star}}\sO_{]\sfV^{(1)\star}[_{\fY^\star}}
\to \cdots \to \Omega^\bullet_{(\cY^\star,M^\star)/W^\circ}\otimes_{\sO_{\cY^\star}}\sO_{]\sfV^{(n-1)\star}[_{\fY^\star}} \to 0
\end{align*}
in $\bC^+(\fY^\star_{\eta,\qet},K_0)$ is exact. Here, for $h\geq 1$,
\[
\sO_{]\sfV^{(h)\star}[_{\fY^\star}}\coloneqq\bigoplus_\cT\sO_{]\sfT^\star[_{\fY^\star}}
\]
where the direct sum is taken over all irreducible components of $\cV^{(h)}$. The lemma is proved.
\end{proof}

Recall that by \cite{Tsu99}*{(4.5.1)}, we have the canonical identification
\begin{align*}
P_\bullet&=\Gamma(\Spec\ol\kappa,\sC_{\ol{\Spec O_K}/W[t]^\circ}),\\
\widehat{\dB_\st^+}&=\Gamma(\Spec\ol\kappa,\sC_{\ol{\Spec O_K}/W[t]^\circ,K_0})
\end{align*}
of $(\varphi,N)$-modules in $\bM_K(W_\bullet)$ and $\bM_K(K_0)$, respectively.

\begin{lem}\label{le:descent2}
The following holds for every $h\geq -1$.
\begin{enumerate}
  \item The object $\sC_{\cV^{(h)}/W[t]^\circ}$ of $\bD^+(\sfX_{\et},W_\bullet)$ is admissible (Definition \ref{de:admissible}).

  \item There is a canonical isomorphism
      \[
      \bH^q(\sfX_{\et},\sC_{\cV^{(h)}/W[t]^\circ,K_0})\simeq\bH^q(\sfX_{\et},\sC_{\sfV^{(h)}/W^\circ,K_0})\otimes_{K_0}\dK
      \]
      of $(\varphi,N)$-modules in $\bM(K_0)$ for every $q\geq 0$.

  \item The object $\sC_{\ol{\cV^{(h)}}/W[t]^\circ}$ of $\bD^+_{\rG_K}(\ol\sfX_{\et},W_\bullet)$ is admissible.

  \item The natural map
      \[
      \bH^q(\sfX_{\et},\sC_{\cV^{(h)}/W[t]^\circ,K_0})\otimes_{K_0}\widehat{\dB_\st^+}\to
      \bH^q(\ol\sfX_{\et},\sC_{\ol{\cV^{(h)}}/W[t]^\circ,K_0})
      \]
      of $(\varphi,N)$-modules induced by functoriality and cup product descends to an isomorphism
      \[
      \bH^q(\sfX_{\et},\sC_{\cV^{(h)}/W[t]^\circ,K_0})\otimes_{\dK}\widehat{\dB_\st^+}\xrightarrow{\sim}
      \bH^q(\ol\sfX_{\et},\sC_{\ol{\cV^{(h)}}/W[t]^\circ,K_0})
      \]
      of $(\varphi,N)$-modules in $\bM_K(K_0)$ for $q\geq 0$. In particular, the natural map
      \[
      \bH^q(\sfX_{\et},\sC_{\cV^{(h)}/W[t]^\circ,K_0})\to\bH^q(\ol\sfX_{\et},\sC_{\ol{\cV^{(h)}}/W[t]^\circ,K_0})
      \]
      is injective.

  \item The object $\sC_{\ol{\cV^{(h)}}/W^\triv}$ of $\bD^+_{\rG_K}(\ol\sfX_{\et},W_\bullet)$ is admissible.\footnote{However, the object $\sC_{\cV^{(h)}/W^\triv}$ of $\bD^+(\sfX_{\et},W_\bullet)$ is in general \emph{not} admissible.}

  \item The distinguished triangle $\sC^\triangle_{\ol{\cV^{(h)}}/W[t]^\circ}$ \eqref{eq:triangle2} induces a Frobenius equivariant isomorphism
      \[
      \bH^q(\ol\sfX_{\et},\sC_{\ol{\cV^{(h)}}/W^\triv,K_0})\simeq
      \(\bH^q(\sfX_{\et},\sC_{\sfV^{(h)}/W^\circ,K_0})\otimes_{K_0}\dB_\st^+\)^{N=0}
      \]
      in $\bM_K(K_0)$ for every $q\geq 0$.
\end{enumerate}
Moreover, the isomorphisms in (2,4,6) are compatible with $h$ in the obvious sense.
\end{lem}

\begin{proof}
It is clear that condition (1) of Definition \ref{de:admissible} holds for $\sC_{\cV^{(h)}/W[t]^\circ}$ trivially as $\rG$ is the trivial group and holds for $\sC_{\ol{\cV^{(h)}}/W[t]^\circ}$ and $\sC_{\ol{\cV^{(h)}}/W^\triv}$ since they are defined as injective limits over terms fixed by open subgroups of $\rG_K$. Thus, it remains to check condition (2) of Definition \ref{de:admissible} for the three objects. Below in the proof, we say that two objects $\sC$ and $\sC'$ in either $\bD^+(\sfX_{\et},W_\bullet)$ or $\bD^+_{\rG_K}(\ol\sfX_{\et},W_\bullet)$ are almost equivalent if there exists a map $\sC\to\sC'$ whose fiber is annihilated by some power of $p$. It is clear that condition (2) of Definition \ref{de:admissible} is preserved under almost equivalence.

First, we prove (1) and (2). Note that \cite{HK94}*{Lemma~5.2} is applicable to $\cV^{(h)}$ with $h\geq 0$, and hence also to the case for $h=-1$ by long exact sequences induced by taking successive homotopy fibers. In other words, for every $h\geq -1$, $\sC_{\cV^{(h)}/W[t]^\circ}$ and $\sC_{\sfV^{(h)}/W^\circ}\otimes_{W_\bullet}^\rL R_\bullet$ are almost equivalent. It is already known by the proof of \cite{Sat13}*{Proposition~A.3.1(1)} that, when $h\geq 0$, $\sC_{\sfV^{(h)}/W^\circ}\otimes_{W_\bullet}^\rL R_\bullet$ (hence $\sC_{\cV^{(h)}/W[t]^\circ}$) satisfies condition (2) of Definition \ref{de:admissible} and we have the canonical isomorphism in (2). To pass the same argument for $h=-1$, it suffices to show that $\bH^q(\sfX_\et,\sC_{\sfV^{(h)}/W^\circ,l})$ is finite for every $q\geq 0$ and $l\geq 1$, which follows from long exact sequences induced by taking successive homotopy fibers.

Second, we prove (3) and (4). It is known by \cite{Tsu99}*{Proposition~4.5.4} that for $h\geq 0$, the natural map
\[
\bH^q(\sfX_{\et},\sC_{\cV^{(h)}/W[t]^\circ})\otimes_{W_\bullet}\Gamma(\Spec\ol\kappa,\sC_{\ol{\Spec O_K}/W[t]^\circ})\to
\bH^q(\ol\sfX_{\et},\sC_{\ol{\cV^{(h)}}/W[t]^\circ})
\]
of $(\varphi,N)$-modules descends to an isomorphism
\begin{align}\label{eq:descent13}
\bH^q(\sfX_{\et},\sC_{\cV^{(h)}/W[t]^\circ})\otimes_{R_\bullet}\Gamma(\Spec\ol\kappa,\sC_{\ol{\Spec O_K}/W[t]^\circ})\xrightarrow\sim
\bH^q(\ol\sfX_{\et},\sC_{\ol{\cV^{(h)}}/W[t]^\circ})
\end{align}
of $(\varphi_N)$-modules in $\bM_{\rG_K}(W_\bullet)$ for every $q\geq 0$. The same holds for $h=-1$ by the long exact sequences induced by taking successive homotopy fibers and the fact that $\Gamma(\Spec\ol\kappa,\sC_{\ol{\Spec O_K}/W[t]^\circ})=P_\bullet$ is flat over $R_\bullet$ \cite{Tsu99}*{Proposition~4.1.5}.

To show (3), since $\sC_{\cV^{(h)}/W[t]^\circ}$ and $\sC_{\sfV^{(h)}/W^\circ}\otimes_{W_\bullet}^\rL R_\bullet$ are almost equivalent for $h\geq -1$, it suffices to show that $\rR^1\varprojlim_l\bH^q(\sfX_{\et},\sC_{\sfV^{(h)}/W^\circ,l})\otimes_{W/p^l}P_l=0$ which has been argued in the proof of \cite{Sat13}*{Proposition~A.3.1(2)}.

To show (4), by \eqref{eq:descent13} (for $h\geq-1$), we have
\[
\varprojlim_l\bH^q(\sfX_{\et},\sC_{\cV^{(h)}/W[t]^\circ,l})\otimes_{R_l}P_l\xrightarrow\sim
\varprojlim_l\bH^q(\ol\sfX_{\et},\sC_{\ol{\cV^{(h)}}/W[t]^\circ,l}).
\]
Since both $R_\bullet$ and $P_\bullet$ have surjective transition maps, we have
\begin{align*}
&\quad\(\varprojlim_l\bH^q(\sfX_{\et},\sC_{\cV^{(h)}/W[t]^\circ,l})\otimes_{R_l}P_l\)\otimes_WK_0 \\
&=\(\varprojlim_l\bH^q(\sfX_{\et},\sC_{\cV^{(h)}/W[t]^\circ,l})\otimes_WK_0\)\otimes_\dK\widehat{\dB_\st^+}.
\end{align*}
By (1) and (3), we have
\begin{align*}
\(\varprojlim_l\bH^q(\sfX_{\et},\sC_{\cV^{(h)}/W[t]^\circ,l})\)\otimes_WK_0
&=\bH^q(\sfX_{\et},\sC_{\cV^{(h)}/W[t]^\circ,K_0}),\\
\(\varprojlim_l\bH^q(\ol\sfX_{\et},\sC_{\ol{\cV^{(h)}}/W[t]^\circ,l})\)\otimes_WK_0
&=\bH^q(\ol\sfX_{\et},\sC_{\ol{\cV^{(h)}}/W[t]^\circ,K_0}),
\end{align*}
respectively. Together we obtain (4); and the injectivity follows from Lemma \ref{le:descent7}(2).

Finally, we prove (5) and (6). It is well-known that when $h\geq 0$, the monodromy operator on $\bH^q(\sfX_{\et},\sC_{\sfV^{(h)}/W^\circ})$ is nilpotent for every $q\geq 0$. The same holds for $h=-1$ by long exact sequences induced by taking successive homotopy fibers. Since the monodromy map on $P_\bullet$ is surjective \cite{Kat94}*{Corollary~3.6}, the monodromy map on $\bH^q(\sfX_{\et},\sC_{\sfV^{(h)}/W^\circ})\otimes_{W_\bullet}P_\bullet$ is surjective for every $h\geq -1$.

To show (5), since $\sC_{\cV^{(h)}/W[t]^\circ}$ and $\sC_{\sfV^{(h)}/W^\circ}\otimes_{W_\bullet}^\rL R_\bullet$ are almost equivalent and by \eqref{eq:descent13} for every $h\geq -1$, it suffices to show that
\[
\rR^1\varprojlim_l\(\bH^q(\sfX_{\et},\sC_{\sfV^{(h)}/W^\circ})\otimes_{W_\bullet}P_\bullet\)^{N=0}=0,
\]
which follows from the argument in the proof of \cite{Sat13}*{Proposition~A.3.1(3)}.

To show (6), we first note that
\[
\(\bH^q(\sfX_{\et},\sC_{\sfV^{(h)}/W^\circ,K_0})\otimes_{K_0}\dB_\st^+\)^{N=0}=
\(\bH^q(\sfX_{\et},\sC_{\sfV^{(h)}/W^\circ,K_0})\otimes_{K_0}\widehat\dB_\st^+\)^{N=0}
\]
by Lemma \ref{le:descent7}(1). Then by (2) and (4), it suffices to show that the monodromy map
\[
N\colon\bH^q(\sfX_{\et},\sC_{\sfV^{(h)}/W^\circ,K_0})\otimes_{K_0}\widehat\dB_\st^+
\to\bH^q(\sfX_{\et},\sC_{\sfV^{(h)}/W^\circ,K_0})\otimes_{K_0}\widehat\dB_\st^+
\]
is surjective for every $q\geq 0$. However, this follows from Lemma \ref{le:descent7}(3) and the fact that the monodromy operator on $\bH^q(\sfX_{\et},\sC_{\sfV^{(h)}/W^\circ,K_0})$ is nilpotent.
\end{proof}

\begin{lem}\label{le:descent3}
Consider integers $r$ satisfying $n-1\leq r<p-1$. For every $h\geq-1$ and every $q\geq 0$, the $\dB_\st$-linear extension of the composite map
\begin{align}\label{eq:descent9}
&\quad\bH^q(\ol\sfX_{\et},\sN(0)_{\ol{\cV^{(h)}},\dQ_p}) \\
& \xrightarrow{\sim} \bH^q(\ol\sfX_{\et},\sN(r)_{\ol{\cV^{(h)}},\dQ_p})\otimes_{\dQ_p}\dQ_p(-r) \notag \\
&\xrightarrow{\sim} \bH^q(\ol\sfX_{\et},\sS(r)_{\ol{\cV^{(h)}},\dQ_p})\otimes_{\dQ_p}\dQ_p(-r) \notag \\
& \to \bH^q(\ol\sfX_{\et},\sC_{\ol{\cV^{(h)}}/W^\triv,K_0})\otimes_{\dQ_p}\dQ_p(-r) \notag\\
& \xrightarrow{\sim} \(\bH^q(\sfX_{\et},\sC_{\sfV^{(h)}/W^\circ,K_0})\otimes_{K_0}\dB_\st^+\)^{N=0}\otimes_{\dQ_p}\dQ_p(-r) \notag\\
& \to \bH^q(\sfX_{\et},\sC_{\sfV^{(h)}/W^\circ,K_0})\otimes_{K_0}\dB_\st \notag
\end{align}
is independent of $r$ and induces an isomorphism
\begin{align*}
\bH^q(\ol\sfX_{\et},\sN(0)_{\ol{\cV^{(h)}},\dQ_p})\otimes_{\dQ_p}\dB_\st
\simeq\bH^q(\sfX_{\et},\sC_{\sfV^{(h)}/W^\circ,K_0})\otimes_{K_0}\dB_\st
\end{align*}
of $(\varphi,N)$-modules in $\bM_K(K_0)$. Here in \eqref{eq:descent9}, the second arrow is induced by the inverse of the period map $\ol\pi_r$ \eqref{eq:descent2}, the third arrow is induced by the map $\ol\xi_r$ \eqref{eq:descent7}, the fourth arrow is the isomorphism from Lemma \ref{le:descent2}(6), and the last arrow is induced by the canonical map $\dQ_p(-r)\hookrightarrow\dB_\st$.

In particular, the above isomorphism induces a Frobenius equivariant isomorphism
\begin{align*}
\bH^q(\ol\sfX_{\et},\sN(0)_{\ol{\cV^{(h)}},\dQ_p})\otimes_{\dQ_p}\dB_\cris
&\simeq\(\bH^q(\sfX_{\et},\sC_{\sfV^{(h)}/W^\circ,K_0})\otimes_{K_0}\dB_\st\)^{N=0}
\end{align*}
in $\bM_K(K_0)$.
\end{lem}

\begin{proof}
For $h\geq0$, the statement follows from \cite{Tsu99}*{Theorem~4.10.2} (the usual $C_\st$-comparison theorem for proper strictly semistable schemes) together with the compatibility properties \cite{Tsu99}*{Corollaries~4.8.8~\&~4.9.2} for the independence of $r$ (which is at least $\dim V^{(h)}$) of the map. The case for $h=-1$ follows from the series of long exact sequences induced by taking successive homotopy fibers.\footnote{Such an isomorphism for $h=-1$ has already been obtained in \cite{Yam11}. The results there are much stronger than ours and in particular they contain a $C_\dr$-comparison isomorphism. Thus, the log structure on $\cX$ used there is $L_\cX^{\sfX\cup\cV}$, which makes things more complicated.}
\end{proof}

\begin{lem}\label{le:descent5}
Suppose that $n<p$. For every $h\geq-1$, every $q\geq 0$, and every $0\leq r<p-1$, the following diagram
\[
\xymatrix{
\bH^q(\ol\sfX_{\et},\sS(r)_{\ol{\cV^{(h)}},\dQ_p})  \ar[r]^-{\ol\pi_r}_-{\eqref{eq:descent2}}
\ar[d]_-{\ol\xi_r}^-{\eqref{eq:descent7}}
& \bH^q(\ol\sfX_{\et},\sN(r)_{\ol{\cV^{(h)}},\dQ_p}) \ar[d]^-{\dQ_p(r)\hookrightarrow\dB_\cris} \\
\bH^q(\ol\sfX_{\et},\sC_{\ol{\cV^{(h)}}/W^\triv,K_0}) \ar[d]_-{\simeq}
& \bH^q(\ol\sfX_{\et},\sN(0)_{\ol{\cV^{(h)}},\dQ_p})\otimes_{\dQ_p}\dB_\cris \ar[d]^-{\simeq} \\
\(\bH^q(\sfX_{\et},\sC_{\sfV^{(h)}/W^\circ,K_0})\otimes_{K_0}\dB_\st^+\)^{N=0} \ar@{^(->}[r]
& \(\bH^q(\sfX_{\et},\sC_{\sfV^{(h)}/W^\circ,K_0})\otimes_{K_0}\dB_\st\)^{N=0}
}
\]
in $\bM_K(\dQ_p)$ commutes, in which the isomorphism in the left column is from Lemma \ref{le:descent2}(6), and the isomorphism in the right column is from Lemma \ref{le:descent3}.
\end{lem}

\begin{proof}
When $h\geq 0$, the commutativity follows from the compatibility properties \cite{Tsu99}*{Corollaries~4.8.8~\&~4.9.2}. The case for $h=-1$ follows from the series of long exact sequences induced by taking successive homotopy fibers.
\end{proof}

\begin{lem}\label{le:descent4}
Suppose that $n<p$. The map
\[
\sfF_!\sfF^*\sS(d)_{\cX,\dQ_p}\to\sS(d)_{\cX,\dQ_p}=\sS(d)_{\cV^{(0)},\dQ_p}
\]
factors though $\sS(d)_{\cV^{(-1)},\dQ_p}$. In particular,
\begin{enumerate}
  \item the map \eqref{eq:syntomic} factors as
     \begin{align*}
     \rR\Gamma(\sfX_{\et},\sfF_!\sfF^*\sS(d)_{\cX,\dQ_p})
     &\to\rR\Gamma(\sfX_{\et},\sS(d)_{\cV^{(-1)},\dQ_p}) \\
     &\to\rR\Gamma(\sfX_{\et},\sN(d)_{\cV^{(-1)},\dQ_p})
     =\rR\Gamma_\rc(U,\dQ_p(d))
     \end{align*}
     when $d<p-1$;

  \item the map \eqref{eq:comparison5} factors as
     \[
     \sfF_!\sfF^*\sS(d)_{\cX,\dQ_p}\to\sS(d)_{\cV^{(-1)},\dQ_p}
     \to\widetilde\sC_{\sfV^{(-1)}/W^\circ,K_0}\simeq\widetilde\bomega_{(\sfU,\sfX)}
     \]
     in which the last equivalence comes from Lemma \ref{le:descent1}.
\end{enumerate}
\end{lem}

\begin{proof}
Since the complex $\sS(d)_{\cV^{(h)},\dQ_p}$ is supported on $\sfV$ for every $h\geq 1$, the factorization follows from the construction.
\end{proof}

For every $h\geq -1$, every $q\geq 0$, and every $0\leq r< p-1$, put
\begin{align}\label{eq:comparison12}
\bH^q(\sfX_{\et},\sS(r)_{\cV^{(h)},\dQ_p})^\heartsuit\coloneqq
\Ker\(\bH^q(\sfX_{\et},\sS(r)_{\cV^{(h)},\dQ_p})\to\bH^q(\ol\sfX_{\et},\sN(r)_{\ol{\cV^{(h)}},\dQ_p})\).
\end{align}

Now we can give a proof of Lemma \ref{le:comparison1}.

\begin{proof}[Proof of Lemma \ref{le:comparison1}]
By Lemma \ref{le:descent4}(1), $\bH^q(\sfX_{\et},\sfF_!\sfF^*\sS(d)_{\cX,\dQ_p})^\heartsuit$ coincides with the kernel of the composite map
\begin{align*}
\bH^q(\sfX_{\et},\sfF_!\sfF^*\sS(d)_{\cX,\dQ_p})
&\to\bH^q(\sfX_{\et},\sS(d)_{\cV^{(-1)},\dQ_p}) \\
&\to\bH^q(\ol\sfX_{\et},\sS(d)_{\ol{\cV^{(-1)}},\dQ_p}) \\
&\xrightarrow{\ol\pi_r}\bH^q(\ol\sfX_{\et},\sN(d)_{\ol{\cV^{(-1)}},\dQ_p}).
\end{align*}
By Lemma \ref{le:descent5}, it is also the kernel of the composite map
\begin{align*}
\bH^q(\sfX_{\et},\sfF_!\sfF^*\sS(d)_{\cX,\dQ_p})
&\to\bH^q(\sfX_{\et},\sS(d)_{\cV^{(-1)},\dQ_p}) \\
&\to\bH^q(\ol\sfX_{\et},\sS(d)_{\ol{\cV^{(-1)}},\dQ_p}) \\
&\xrightarrow{\ol\xi_r}\bH^q(\ol\sfX_{\et},\sC_{\ol{\cV^{(-1)}}/W^\triv,K_0}).
\end{align*}
By Lemma \ref{le:descent4}(2) and \eqref{eq:comparison9}, we have the following commutative diagram
\[
\xymatrix{
\bH^q(\ol\sfX_{\et},\sS(d)_{\ol{\cV^{(-1)}},\dQ_p})\ar[r]^-{\ol\xi_r}
&  \bH^q(\ol\sfX_{\et},\sC_{\ol{\cV^{(-1)}}/W^\triv,K_0}) \ar[r]^-{\nu} & \bH^q(\ol\sfX_{\et},\sC_{\ol{\cV^{(-1)}}/W[t]^\circ,K_0})  \\
\bH^q(\sfX_{\et},\sS(d)_{\cV^{(-1)},\dQ_p})  \ar[r]^-{\xi_r}\ar[u] & \bH^q(\sfX_{\et},\sC_{\cV^{(-1)}/W^\triv,K_0}) \ar[r]\ar[u]\ar[d] &
\bH^q(\sfX_{\et},\sC_{\cV^{(-1)}/W[t]^\circ,K_0}) \ar[u]_-{\mu}\ar[d]\\
& \bH^q(\sfX_{\et},\widetilde\sC_{\cV^{(-1)}/W^\circ,K_0}) \ar[r] & \bH^q(\sfX_{\et},\sC_{\cV^{(-1)}/W^\circ,K_0})
}
\]
in which $\mu$ and $\nu$ are injective by Lemma \ref{le:descent2}(4) and (6), respectively. Thus, $\bH^q(\sfX_{\et},\sfF_!\sfF^*\sS(d)_{\cX,\dQ_p})^\heartsuit$ maps to zero all the way to the term at the lower-right corner $\bH^q(\sfX_{\et},\sC_{\cV^{(-1)}/W^\circ,K_0})=\rH^q_\rig((\sfU,\sfX)/W^\circ)$. The lemma is proved.
\end{proof}

In order to prove Lemma \ref{le:comparison2}, we need to compare edge maps in both \'{e}tale and crystalline settings. For every $h\geq-1$, every $q\geq 0$, and every $r\geq 0$, put
\begin{align*}
\bH^q(\sfX_{\et},\sN(r)_{\cV^{(h)},\dQ_p})^0&\coloneqq
\Ker\(\bH^q(\sfX_{\et},\sN(r)_{\cV^{(h)},\dQ_p})\to\bH^q(\ol\sfX_{\et},\sN(r)_{\ol{\cV^{(h)}},\dQ_p})\), \\
\bH^q(\sfX_{\et},\sC_{\cV^{(h)}/W^\triv,K_0})^0&\coloneqq
\Ker\(\bH^q(\sfX_{\et},\sC_{\cV^{(h)}/W^\triv,K_0})\to\bH^q(\ol\sfX_{\et},\sC_{\ol{\cV^{(h)}}/W^\triv,K_0})\).
\end{align*}
Suppose that $n<p$. Then by definition, $\pi_d$ induces a map
\begin{align}\label{eq:descent11}
\bH^q(\sfX_{\et},\sS(d)_{\cV^{(h)},\dQ_p})^\heartsuit \to \bH^q(\sfX_{\et},\sN(d)_{\cV^{(h)},\dQ_p})^0.
\end{align}
By Lemma \ref{le:descent5} and the definition of the syntomic complex, $\xi_d$ induces a map
\begin{align}\label{eq:descent12}
\bH^q(\sfX_{\et},\sS(d)_{\cV^{(h)},\dQ_p})^\heartsuit \to \(\bH^q(\sfX_{\et},\sC_{\cV^{(h)}/W^\triv,K_0})^0\)^{\varphi=p^d}.
\end{align}
As $\sN(d)_{\ol{\cV^{(h)}}}$ is clearly an admissible object of $\bD^+_{\rG_K}(\ol\sfX_\et,\dZ_{p\bullet})$ and $\sC_{\ol{\cV^{(h)}}/W^\triv}$ is an admissible object of $\bD^+_{\rG_K}(\ol\sfX_\et,W_\bullet)$ by Lemma \ref{le:descent2}(5), we have the spectral sequences and hence the corresponding edges maps for them from Lemma \ref{le:spectral}. Composing \eqref{eq:descent11} and \eqref{eq:descent12} with the corresponding edge maps, we obtain maps
\begin{align*}
\beta_q&\colon\bH^q(\sfX_{\et},\sS(d)_{\cV^{(h)},\dQ_p})^\heartsuit \to \rH^1\(K,\bH^{q-1}(\ol\sfX_{\et},\sN(d)_{\ol{\cV^{(h)}},\dQ_p})\), \\
\gamma_q&\colon\bH^q(\sfX_{\et},\sS(d)_{\cV^{(h)},\dQ_p})^\heartsuit \to \rH^1\(K,\bH^{q-1}(\ol\sfX_{\et},\sC_{\ol{\cV^{(h)}}/W^\triv,K_0})\),
\end{align*}
respectively.

The following lemma is an ``Abel--Jacobi'' version of Lemma \ref{le:descent5}.

\begin{lem}\label{le:descent6}
Suppose that $n<p$. For every $h\geq-1$ and every $q\geq 0$, the following diagram
\[
\resizebox{\hsize}{!}{
\xymatrix{
\bH^q(\sfX_{\et},\sS(d)_{\cV^{(h)},\dQ_p})^\heartsuit \ar[d]_-{\gamma_q}\ar[r]^-{\beta_q}
& \rH^1\(K,\bH^{q-1}(\ol\sfX_{\et},\sN(d)_{\ol{\cV^{(h)}},\dQ_p})\) \ar[d]^-{\dQ_p(d)\hookrightarrow\dB_\cris} \\
\rH^1\(K,\bH^{q-1}(\ol\sfX_{\et},\sC_{\ol{\cV^{(h)}}/W^\triv,K_0})\) \ar[d]_-{\simeq}
& \rH^1\(K,\bH^{q-1}(\ol\sfX_{\et},\sN(0)_{\ol{\cV^{(h)}},\dQ_p})\otimes_{\dQ_p}\dB_\cris\) \ar[d]^-{\simeq} \\
\rH^1\(K,\(\bH^{q-1}(\sfX_{\et},\sC_{\sfV^{(h)}/W^\circ,K_0})\otimes_{K_0}\dB_\st^+\)^{N=0}\)  \ar[r]
& \rH^1\(K,\(\bH^{q-1}(\sfX_{\et},\sC_{\sfV^{(h)}/W^\circ,K_0})\otimes_{K_0}\dB_\st\)^{N=0}\)
}
}
\]
of $\dQ_p$-vector spaces commutes, in which the isomorphism in the left column is from Lemma \ref{le:descent2}(6), and the isomorphism in the right column is from Lemma \ref{le:descent3}.
\end{lem}

This is \emph{not} an immediate consequence of Lemma \ref{le:descent5} since $\bH^q(\sfX_{\et},\sS(d)_{\cV^{(h)},\dQ_p})^\heartsuit$, by our definition \eqref{eq:comparison12}, is in general larger than
\[
\Ker\(\bH^q(\sfX_{\et},\sS(d)_{\cV^{(h)},\dQ_p})\to\bH^q(\ol\sfX_{\et},\sS(d)_{\ol{\cV^{(h)}},\dQ_p})\).
\]

\begin{proof}
Take an integer $r$ satisfying $(1\leq d\leq)n-1\leq r<p-1$ (which exists as $n<p$). We have
\begin{itemize}
  \item the map $\ol\pi_{r-d}\colon\sS(r-d)_{\ol{\Spec O_K}}\xrightarrow{\sim}\sN(r-d)_{\ol{\Spec O_K}}$ in $\bD^+_{\rG_K}(\dZ_{p\bullet})$, both equivalent to $\dZ_{p\bullet}(r-d)$,

  \item the object $\sC_{\ol{\Spec O_K}/W^\triv}$ in $\bD^+_{\rG_K}(W_\bullet)$, which is equivalent to $P_\bullet^{N=0}$,

  \item the map $\ol\xi_{r-d}\colon\sS(r-d)_{\ol{\Spec O_K}}\to\sC_{\ol{\Spec O_K}/W^\triv}$ in $\bD^+_{\rG_K}(\dZ_{p\bullet})$, which is equivalent to the natural map $\dZ_{p\bullet}(r-d)\hookrightarrow P_\bullet^{N=0}$.
\end{itemize}

For every $h\geq -1$, consider the following diagram
\begin{align}\label{eq:comparison13}
\xymatrix{
\rR\Gamma(\ol\sfX_{\et},\sC_{\ol{\cV^{(h)}}/W^\triv})\otimes_{\dZ_{p\bullet}}^\rL\sC_{\ol{\Spec O_K}/W^\triv}
\ar[r]& \rR\Gamma(\ol\sfX_{\et},\sC_{\ol{\cV^{(h)}}/W^\triv}) \\
\rR\Gamma(\ol\sfX_{\et},\sS(d)_{\ol{\cV^{(h)}}})\otimes_{\dZ_{p\bullet}}^\rL\sS(r-d)_{\ol{\Spec O_K}}
\ar[r]\ar[u]^-{\ol\xi_d\otimes\ol\xi_{r-d}}\ar[d]_-{\ol\pi_d\otimes\ol\pi_{r-d}} & \rR\Gamma(\ol\sfX_{\et},\sS(r)_{\ol{\cV^{(h)}}}) \ar[u]_-{\ol\xi_r}\ar[d]^-{\ol\pi_r} \\
\rR\Gamma(\ol\sfX_{\et},\sN(d)_{\ol{\cV^{(h)}}})\otimes_{\dZ_{p\bullet}}^\rL\sN(r-d)_{\ol{\Spec O_K}}
\ar[r]& \rR\Gamma(\ol\sfX_{\et},\sN(r)_{\ol{\cV^{(h)}}})
}
\end{align}
in $\bD^+_{\rG_K}(\dZ_{p\bullet})$, in which all horizontal maps are induced by cup products. We claim that \eqref{eq:comparison13} commutes. For $h\geq 0$, the upper square commutes by definition, and the lower square commutes by the compatibility of period maps with cup products (see \cite{Tsu99}*{\S3.1} in a more general context). The case for $h=-1$ follows from the process of taking successive homotopy fibers.

Using the equivalences
\begin{align*}
\rR\Gamma(\ol\sfX_{\et},\sS(d)_{\ol{\cV^{(h)}}})
&\simeq\rR\Gamma(\ol\sfX_{\et},\sS(d)_{\ol{\cV^{(h)}}})\otimes_{\dZ_{p\bullet}}^\rL\sS(r-d)_{\ol{\Spec O_K}}\otimes_{\dZ_{p\bullet}}^\rL\dZ_{p\bullet}(d-r), \\
\rR\Gamma(\ol\sfX_{\et},\sN(d)_{\ol{\cV^{(h)}}})
&\simeq\rR\Gamma(\ol\sfX_{\et},\sN(d)_{\ol{\cV^{(h)}}})\otimes_{\dZ_{p\bullet}}^\rL\sN(r-d)_{\ol{\Spec O_K}}\otimes_{\dZ_{p\bullet}}^\rL\dZ_{p\bullet}(d-r),
\end{align*}
and the fact that $\ol\pi_r$ (for $\ol\sfX_\et$) is an equivalence, we obtain the following commutative diagram
\[
\xymatrix{
\rR\Gamma(\ol\sfX_{\et},\sS(d)_{\ol{\cV^{(h)}}}) \ar[d]\ar[rrrd] \\
\rR\Gamma(\ol\sfX_{\et},\sN(d)_{\ol{\cV^{(h)}}}) \ar[rrr]^-{(\ol\xi_r\circ\ol\pi_r^{-1})\otimes\dZ_{p\bullet}(d-r)} &&&
\rR\Gamma(\ol\sfX_{\et},\sC_{\ol{\cV^{(h)}}/W^\triv})\otimes_{\dZ_{p\bullet}}^\rL\dZ_{p\bullet}(d-r)
}
\]
in $\bD^+_{\rG_K}(\dZ_{p\bullet})$ from \eqref{eq:comparison13}.\footnote{The above commutative diagram is parallel to \cite{Sat13}*{(A.6.4)}. However, somehow unfortunately, the roles of $d$ and $r$ here are switched from those there.} Composing with $\rR(\Gamma_{\rG_K}\circ\varprojlim)\otimes_{\dZ_p}\dQ_p$ and taking edge maps, we obtain a commutative diagram
\[
\resizebox{\hsize}{!}{
\xymatrix{
\bH^q(\sfX_{\et},\sS(d)_{\cV^{(h)},\dQ_p})^\heartsuit \ar[d]\ar[rrrd] \\
\rH^1\(K,\bH^{q-1}(\ol\sfX_{\et},\sN(d)_{\ol{\cV^{(h)}},\dQ_p})\) \ar[rrr]^-{(\ol\xi_r\circ\ol\pi_r^{-1})\otimes\dQ_p(d-r)} &&&
\rH^1\(K,\bH^{q-1}(\ol\sfX_{\et},\sC_{\ol{\cV^{(h)}}/W^\triv,K_0})\otimes_{\dQ_p}\dQ_p(d-r)\).
}
}
\]
The lemma follows since we have the commutative diagrams
\[
\resizebox{\hsize}{!}{
\xymatrix{
\bH^{q-1}(\ol\sfX_{\et},\sN(d)_{\ol{\cV^{(h)}},\dQ_p}) \ar[rr]^-{(\ol\xi_r\circ\ol\pi_r^{-1})\otimes\dQ_p(d-r)}\ar[dd]^-{\dQ_p(d)\hookrightarrow\dB_\cris}&& \bH^{q-1}(\ol\sfX_{\et},\sC_{\ol{\cV^{(h)}}/W^\triv,K_0})\otimes_{\dQ_p}\dQ_p(d-r) \ar[d]_-{\simeq}^-{\text{Lemma \ref{le:descent2}(6)}} \\
&& \(\bH^{q-1}(\sfX_{\et},\sC_{\sfV^{(h)}/W^\circ,K_0})\otimes_{K_0}\dB_\st^+\)^{N=0}\otimes_{\dQ_p}\dQ_p(d-r) \ar[d]^-{\dQ_p(d-r)\hookrightarrow\dB_\st} \\
\bH^{q-1}(\ol\sfX_{\et},\sN(0)_{\ol{\cV^{(h)}},\dQ_p})\otimes_{\dQ_p}\dB_\cris \ar[rr]^-{\sim}_-{\text{Lemma \ref{le:descent3}}}&& \(\bH^{q-1}(\sfX_{\et},\sC_{\sfV^{(h)}/W^\circ,K_0})\otimes_{K_0}\dB_\st\)^{N=0}
}
}
\]
and
\[
\resizebox{\hsize}{!}{
\xymatrix{
\bH^{q-1}(\ol\sfX_{\et},\sC_{\ol{\cV^{(h)}}/W^\triv,K_0}) \ar[r]\ar[d]_-{\simeq}^-{\text{Lemma \ref{le:descent2}(6)}} & \bH^{q-1}(\ol\sfX_{\et},\sC_{\ol{\cV^{(h)}}/W^\triv,K_0})\otimes_{\dQ_p}\dQ_p(d-r) \ar[d] \\
\(\bH^{q-1}(\sfX_{\et},\sC_{\sfV^{(h)}/W^\circ,K_0})\otimes_{K_0}\dB_\st^+\)^{N=0}
\ar[r]& \(\bH^{q-1}(\sfX_{\et},\sC_{\sfV^{(h)}/W^\circ,K_0})\otimes_{K_0}\dB_\st\)^{N=0}
}
}
\]
in which the upper horizontal arrow is induced by the canonical map
\[
\dQ_p(r-d)\hookrightarrow(\dB_\st^+)^{N=0}\simeq\sC_{\ol{\Spec O_K}/W^\triv,K_0}
\]
and the cup product, and the right vertical arrow is the composition of two right vertical arrows in the previous diagram.
\end{proof}

\begin{proof}[Proof of Lemma \ref{le:comparison2}]
By Lemma \ref{le:descent4}, we may replace the source of both $\alpha_q$ and $\rho_q$, which is originally $\bH^q(\sfX_{\et},\sfF_!\sfF^*\sS(d)_{\cX,\dQ_p})^\heartsuit$, by $\bH^q(\sfX_{\et},\sS(d)_{\cV^{(-1)},\dQ_p})^\heartsuit$. Then $\alpha_q$ \eqref{eq:comparison} coincides with the composite map
\[
\bH^q(\sfX_{\et},\sS(d)_{\cV^{(-1)},\dQ_p})^\heartsuit
\to\rH^1\(K,\bH^{q-1}(\ol\sfX_{\et},\sN(0)_{\ol{\cV^{(-1)}},\dQ_p})\otimes_{\dQ_p}\dB_\cris\)
\]
from the diagram in Lemma \ref{le:descent6} (with $h=-1$). By Lemma \ref{le:descent6},
\[
\Ker\(\bH^q(\sfX_{\et},\sS(d)_{\cV^{(-1)},\dQ_p})^\heartsuit\to
\rH^1\(K,\(\bH^{q-1}(\sfX_{\et},\sC_{\sfV^{(-1)}/W^\circ,K_0})\otimes_{K_0}\dB_\st^+\)^{N=0}\)\)
\]
is contained in $\Ker(\alpha_q)$. Thus, the lemma will follow if we can show
\begin{align}\label{eq:comparison1}
&\quad\Ker(\rho_q) \\
&=\Ker\(\bH^q(\sfX_{\et},\sS(d)_{\cV^{(-1)},\dQ_p})^\heartsuit\to
\rH^1\(K,\(\bH^{q-1}(\sfX_{\et},\sC_{\sfV^{(-1)}/W^\circ,K_0})\otimes_{K_0}\dB_\st^+\)^{N=0}\)\). \notag
\end{align}

Lemma \ref{le:descent2}(4,6) implies that
\begin{align*}
&\quad\bH^q(\sfX_{\et},\sC_{\cV^{(-1)}/W^\triv,K_0})^0 \\
&=\Ker\(\bH^q(\sfX_{\et},\sC_{\cV^{(-1)}/W^\triv,K_0})\to\bH^q(\sfX_{\et},\sC_{\cV^{(-1)}/W[t]^\circ,K_0})\),
\end{align*}
which induces an isomorphism
\[
\frac{\bH^{q-1}(\sfX_{\et},\sC_{\cV^{(-1)}/W[t]^\circ,K_0})}{N\bH^{q-1}(\sfX_{\et},\sC_{\cV^{(-1)}/W[t]^\circ,K_0})}
\xrightarrow\sim\bH^q(\sfX_{\et},\sC_{\cV^{(-1)}/W^\triv,K_0})^0
\]
by the distinguished triangle $\sC^\triangle_{\cV^{(-1)}/W[t]^\circ,K_0}$ \eqref{eq:triangle3}, under which the Frobenius operator on the right side corresponds to $p$ times the one on the left side.

Consider the following diagram
\begin{align}\label{eq:comparison14}
\resizebox{\hsize}{!}{
\xymatrix{
\(\dfrac{\bH^{q-1}(\sfX_{\et},\sC_{\cV^{(-1)}/W[t]^\circ,K_0})}
{N\bH^{q-1}(\sfX_{\et},\sC_{\cV^{(-1)}/W[t]^\circ,K_0})}\)^{\varphi=p^{d-1}}
\ar[r]^-{\sim} \ar[d] &
\(\bH^q(\sfX_{\et},\sC_{\cV^{(-1)}/W^\triv,K_0})^0\)^{\varphi=p^d} \ar[d] \\
\dfrac{\bH^{q-1}(\sfX_{\et},\sC_{\sfV^{(-1)}/W^\circ,K_0})}{N\bH^{q-1}(\sfX_{\et},\sC_{\sfV^{(-1)}/W^\circ,K_0})} \ar@{^(->}[r]^-{-\delta} &
\rH^1\(K,\(\bH^{q-1}(\sfX_{\et},\sC_{\sfV^{(-1)}/W^\circ,K_0})\otimes_{K_0}\dB_\st^+\)^{N=0}\)
}
}
\end{align}
in which
\begin{itemize}
  \item the map $\delta$ is the edge map induced from the short exact sequence
     \begin{align*}
     0 &\to \(\bH^{q-1}(\sfX_{\et},\sC_{\sfV^{(-1)}/W^\circ,K_0})\otimes_{K_0}\dB_\st^+\)^{N=0}
      \to \bH^{q-1}(\sfX_{\et},\sC_{\sfV^{(-1)}/W^\circ,K_0})\otimes_{K_0}\dB_\st^+ \\
     &\xrightarrow{N}\bH^{q-1}(\sfX_{\et},\sC_{\sfV^{(-1)}/W^\circ,K_0})\otimes_{K_0}\dB_\st^+
     \to 0
     \end{align*}
     in $\bM_K(K_0)$,

  \item the left vertical arrow is the specialization map at $t=0$, and

  \item the right vertical arrow is the composite map
     \begin{align*}
     \bH^q(\sfX_{\et},\sC_{\cV^{(-1)}/W^\triv,K_0})^0
     &\to \rH^1\(K,\bH^{q-1}(\ol\sfX_{\et},\sC_{\ol{\cV^{(-1)}}/W^\triv,K_0})\) \\
     &\xrightarrow{\sim} \rH^1\(K,\(\bH^{q-1}(\sfX_{\et},\sC_{\sfV^{(-1)}/W^\circ,K_0})\otimes_{K_0}\dB_\st^+\)^{N=0}\)
     \end{align*}
     in which the isomorphism is from Lemma \ref{le:descent2}(6).
\end{itemize}
We show that \eqref{eq:comparison14} commutes. Applying Lemma \ref{le:hochschild} to $S=\ol\sfX_{\et}$, the distinguished triangle
\[
\sC_{\ol{\cV^{(-1)}}/W[t]^\circ}[-1] \xrightarrow{-N[-1]}
\sC_{\ol{\cV^{(-1)}}/W[t]^\circ}[-1] \to  \sC_{\ol{\cV^{(-1)}}/W^\triv} \xrightarrow{+1} \sC_{\ol{\cV^{(-1)}}/W[t]^\circ}
\]
in $\bD^+_{\rG_K}(\ol\sfX_\et,W_\bullet)$ (which is a shift of the distinguished triangle $\sC^\triangle_{\ol{\cV^{(-1)}}/W[t]^\circ}$ \eqref{eq:triangle2}) in which all objects are admissible by Lemma \ref{le:descent2}(3,5), we know that the image of an element $c_1\in\bH^{q-1}(\sfX_{\et},\sC_{\cV^{(-1)}/W[t]^\circ,K_0})$ under the composite map
\begin{align*}
\bH^{q-1}(\sfX_{\et},\sC_{\cV^{(-1)}/W[t]^\circ,K_0})
&\to\bH^q(\sfX_{\et},\sC_{\cV^{(-1)}/W^\triv,K_0})^0 \\
&\to\bH^q_K(\ol\sfX_{\et},\sC_{\ol{\cV^{(-1)}}/W^\triv,K_0})^0 \\
&\to\rH^1\(K,\bH^{q-1}(\ol\sfX_{\et},\sC_{\ol{\cV^{(-1)}}/W^\triv,K_0})\)
\end{align*}
can be represented by the (continuous) $1$-cocycle $g\mapsto g\ol{c}_0-\ol{c}_0$ for $g\in\rG_K$, where $\ol{c}_0$ is an arbitrary element in $\bH^{q-1}(\ol\sfX_{\et},\sC_{\ol{\cV^{(-1)}}/W[t]^\circ,K_0})$ satisfying that $-N(\ol{c}_0)$ coincides with
\[
c_1\in\bH^{q-1}(\sfX_{\et},\sC_{\cV^{(-1)}/W[t]^\circ,K_0})
\subseteq\bH^{q-1}(\ol\sfX_{\et},\sC_{\ol{\cV^{(-1)}}/W[t]^\circ,K_0}).
\]
Then the commutativity of \eqref{eq:comparison14} follows from Lemma \ref{le:preparation} (with $D=\bH^{q-1}(\sfX_{\et},\sC_{\sfV^{(-1)}/W^\circ,K_0})$) and Lemma \ref{le:descent2}(2).

Now \eqref{eq:comparison1} follows since the composite map
\begin{align*}
\bH^q(\sfX_{\et},\sS(d)_{\cV^{(-1)},\dQ_p})^\heartsuit
&\xrightarrow{\eqref{eq:descent12}}\(\bH^q(\sfX_{\et},\sC_{\cV^{(-1)}/W^\triv,K_0})^0\)^{\varphi=p^d} \\
&\xrightarrow\sim \(\frac{\bH^{q-1}(\sfX_{\et},\sC_{\cV^{(-1)}/W[t]^\circ,K_0})}
{N\bH^{q-1}(\sfX_{\et},\sC_{\cV^{(-1)}/W[t]^\circ,K_0})}\)^{\varphi=p^{d-1}} \\
&\to\frac{\bH^{q-1}(\sfX_{\et},\sC_{\sfV^{(-1)}/W^\circ,K_0})}{N\bH^{q-1}(\sfX_{\et},\sC_{\sfV^{(-1)}/W^\circ,K_0})}
\end{align*}
is nothing but $\rho_q$ \eqref{eq:comparison2} composed with the isomorphism from Lemma \ref{le:descent1}.

Lemma \ref{le:comparison2} is proved.
\end{proof}


\begin{bibdiv}
\begin{biblist}


\bib{Bei87}{article}{
   author={Be\u{\i}linson, A.},
   title={Height pairing between algebraic cycles},
   conference={
      title={Current trends in arithmetical algebraic geometry},
      address={Arcata, Calif.},
      date={1985},
   },
   book={
      series={Contemp. Math.},
      volume={67},
      publisher={Amer. Math. Soc., Providence, RI},
   },
   date={1987},
   pages={1--24},
   review={\MR{902590}},
   doi={10.1090/conm/067/902590},
}



\bib{Ber94}{article}{
   author={Berkovich, Vladimir G.},
   title={Vanishing cycles for formal schemes},
   journal={Invent. Math.},
   volume={115},
   date={1994},
   number={3},
   pages={539--571},
   issn={0020-9910},
   review={\MR{1262943}},
   doi={10.1007/BF01231772},
}


\bib{Ber96}{article}{
   author={Berkovich, Vladimir G.},
   title={Vanishing cycles for formal schemes. II},
   journal={Invent. Math.},
   volume={125},
   date={1996},
   number={2},
   pages={367--390},
   issn={0020-9910},
   review={\MR{1395723}},
   doi={10.1007/s002220050078},
}


\bib{Ber97}{article}{
   author={Berthelot, Pierre},
   title={Finitude et puret\'{e} cohomologique en cohomologie rigide},
   language={French},
   note={With an appendix in English by Aise Johan de Jong},
   journal={Invent. Math.},
   volume={128},
   date={1997},
   number={2},
   pages={329--377},
   issn={0020-9910},
   review={\MR{1440308}},
   doi={10.1007/s002220050143},
}



\bib{BK90}{article}{
   author={Bloch, Spencer},
   author={Kato, Kazuya},
   title={$L$-functions and Tamagawa numbers of motives},
   conference={
      title={The Grothendieck Festschrift, Vol.\ I},
   },
   book={
      series={Progr. Math.},
      volume={86},
      publisher={Birkh\"auser Boston},
      place={Boston, MA},
   },
   date={1990},
   pages={333--400},
   review={\MR{1086888 (92g:11063)}},
}




\bib{Bre97}{article}{
   author={Breuil, Christophe},
   title={Repr\'{e}sentations $p$-adiques semi-stables et transversalit\'{e} de
   Griffiths},
   language={French},
   journal={Math. Ann.},
   volume={307},
   date={1997},
   number={2},
   pages={191--224},
   issn={0025-5831},
   review={\MR{1428871}},
   doi={10.1007/s002080050031},
}



\bib{BJ79}{article}{
   author={Borel, A.},
   author={Jacquet, H.},
   title={Automorphic forms and automorphic representations},
   note={With a supplement ``On the notion of an automorphic
   representation''\ by R. P. Langlands},
   conference={
      title={Automorphic forms, representations and $L$-functions},
      address={Proc. Sympos. Pure Math., Oregon State Univ., Corvallis,
      Ore.},
      date={1977},
   },
   book={
      series={Proc. Sympos. Pure Math., XXXIII},
      publisher={Amer. Math. Soc., Providence, R.I.},
   },
   date={1979},
   pages={189--207},
   review={\MR{546598}},
}


\bib{Car12}{article}{
   author={Caraiani, Ana},
   title={Local-global compatibility and the action of monodromy on nearby
   cycles},
   journal={Duke Math. J.},
   volume={161},
   date={2012},
   number={12},
   pages={2311--2413},
   issn={0012-7094},
   review={\MR{2972460}},
   doi={10.1215/00127094-1723706},
}


\bib{Car14}{article}{
   author={Caraiani, Ana},
   title={Monodromy and local-global compatibility for $l=p$},
   journal={Algebra Number Theory},
   volume={8},
   date={2014},
   number={7},
   pages={1597--1646},
   issn={1937-0652},
   review={\MR{3272276}},
   doi={10.2140/ant.2014.8.1597},
}



\bib{CH13}{article}{
   author={Chenevier, Ga\"{e}tan},
   author={Harris, Michael},
   title={Construction of automorphic Galois representations, II},
   journal={Camb. J. Math.},
   volume={1},
   date={2013},
   number={1},
   pages={53--73},
   issn={2168-0930},
   review={\MR{3272052}},
   doi={10.4310/CJM.2013.v1.n1.a2},
}



\bib{CF00}{article}{
   author={Colmez, Pierre},
   author={Fontaine, Jean-Marc},
   title={Construction des repr\'{e}sentations $p$-adiques semi-stables},
   language={French},
   journal={Invent. Math.},
   volume={140},
   date={2000},
   number={1},
   pages={1--43},
   issn={0020-9910},
   review={\MR{1779803}},
   doi={10.1007/s002220000042},
}



\bib{deJ96}{article}{
   author={de Jong, A. J.},
   title={Smoothness, semi-stability and alterations},
   journal={Inst. Hautes \'{E}tudes Sci. Publ. Math.},
   number={83},
   date={1996},
   pages={51--93},
   issn={0073-8301},
   review={\MR{1423020}},
}



\bib{Dis17}{article}{
   author={Disegni, Daniel},
   title={The $p$-adic Gross-Zagier formula on Shimura curves},
   journal={Compos. Math.},
   volume={153},
   date={2017},
   number={10},
   pages={1987--2074},
   issn={0010-437X},
   review={\MR{3692745}},
   doi={10.1112/S0010437X17007308},
}




\bib{Dis23}{article}{
   author={Disegni, Daniel},
   title={The universal $p$-adic Gross-Zagier formula},
   journal={Invent. Math.},
   volume={230},
   date={2022},
   number={2},
   pages={509--649},
   issn={0020-9910},
   review={\MR{4493324}},
   doi={10.1007/s00222-022-01133-w},
}



\bib{Dis22}{article}{
   author={Disegni, Daniel},
   title={The $p$-adic Gross-Zagier formula on Shimura curves, II: nonsplit
   primes},
   journal={J. Inst. Math. Jussieu},
   volume={22},
   date={2023},
   number={5},
   pages={2199--2240},
   issn={1474-7480},
   review={\MR{4624961}},
   doi={10.1017/S1474748021000608},
}





\bib{EL}{article}{
   author={Eischen, Ellen},
   author={Liu, Zheng},
   title={Archimedean zeta integrals for unitary groups},
   note={\href{https://arxiv.org/abs/2006.04302}{arXiv:2006.04302}},
}


\bib{EHLS}{article}{
   author={Eischen, Ellen},
   author={Harris, Michael},
   author={Li, Jianshu},
   author={Skinner, Christopher},
   title={$p$-adic $L$-functions for unitary groups},
   journal={Forum Math. Pi},
   volume={8},
   date={2020},
   pages={e9, 160},
   review={\MR{4096618}},
   doi={10.1017/fmp.2020.4},
}


\bib{GI16}{article}{
   author={Gan, Wee Teck},
   author={Ichino, Atsushi},
   title={The Gross-Prasad conjecture and local theta correspondence},
   journal={Invent. Math.},
   volume={206},
   date={2016},
   number={3},
   pages={705--799},
   issn={0020-9910},
   review={\MR{3573972}},
   doi={10.1007/s00222-016-0662-8},
}



\bib{GS12}{article}{
   author={Gan, Wee Teck},
   author={Savin, Gordan},
   title={Representations of metaplectic groups I: epsilon dichotomy and
   local Langlands correspondence},
   journal={Compos. Math.},
   volume={148},
   date={2012},
   number={6},
   pages={1655--1694},
   issn={0010-437X},
   review={\MR{2999299}},
   doi={10.1112/S0010437X12000486},
}


\bib{GG11}{article}{
   author={Gong, Z.},
   author={Greni\'{e}, L.},
   title={An inequality for local unitary theta correspondence},
   language={English, with English and French summaries},
   journal={Ann. Fac. Sci. Toulouse Math. (6)},
   volume={20},
   date={2011},
   number={1},
   pages={167--202},
   issn={0240-2963},
   review={\MR{2830396}},
}


\bib{GZ86}{article}{
   author={Gross, Benedict H.},
   author={Zagier, Don B.},
   title={Heegner points and derivatives of $L$-series},
   journal={Invent. Math.},
   volume={84},
   date={1986},
   number={2},
   pages={225--320},
   issn={0020-9910},
   review={\MR{833192}},
   doi={10.1007/BF01388809},
}


\bib{GK05}{article}{
   author={Grosse-Kl\"{o}nne, Elmar},
   title={Frobenius and monodromy operators in rigid analysis, and Drinfeld's symmetric space},
   journal={J. Algebraic Geom.},
   volume={14},
   date={2005},
   number={3},
   pages={391--437},
   issn={1056-3911},
   review={\MR{2129006}},
   doi={10.1090/S1056-3911-05-00402-9},
}


\bib{GD77}{article}{
   author={Grothendieck, A.},
   author={Deligne, P.},
   title={La classe de cohomologie associ\'{e}e \`a un cycle},
   language={French},
   conference={
      title={Cohomologie \'{e}tale},
   },
   book={
      series={Lecture Notes in Math.},
      volume={569},
      publisher={Springer, Berlin},
   },
   date={1977},
   pages={129--153},
   review={\MR{3727436}},
}


\bib{Har07}{article}{
   author={Harris, Michael},
   title={Cohomological automorphic forms on unitary groups. II. Period relations and values of $L$-functions},
   conference={
      title={Harmonic analysis, group representations, automorphic forms and invariant theory},
   },
   book={
      series={Lect. Notes Ser. Inst. Math. Sci. Natl. Univ. Singap.},
      volume={12},
      publisher={World Sci. Publ., Hackensack, NJ},
   },
   date={2007},
   pages={89--149},
   review={\MR{2401812}},
}


\bib{HKS96}{article}{
   author={Harris, Michael},
   author={Kudla, Stephen S.},
   author={Sweet, William J.},
   title={Theta dichotomy for unitary groups},
   journal={J. Amer. Math. Soc.},
   volume={9},
   date={1996},
   number={4},
   pages={941--1004},
   issn={0894-0347},
   review={\MR{1327161}},
   doi={10.1090/S0894-0347-96-00198-1},
}




\bib{Hid98}{article}{
   author={Hida, Haruzo},
   title={Automorphic induction and Leopoldt type conjectures for $\mathrm{GL}(n)$},
   note={Mikio Sato: a great Japanese mathematician of the twentieth
   century},
   journal={Asian J. Math.},
   volume={2},
   date={1998},
   number={4},
   pages={667--710},
   issn={1093-6106},
   review={\MR{1734126}},
   doi={10.4310/AJM.1998.v2.n4.a5},
}





\bib{Hyo91}{article}{
   author={Hyodo, Osamu},
   title={On the de Rham-Witt complex attached to a semi-stable family},
   journal={Compositio Math.},
   volume={78},
   date={1991},
   number={3},
   pages={241--260},
   issn={0010-437X},
   review={\MR{1106296}},
}



\bib{HK94}{article}{
   author={Hyodo, Osamu},
   author={Kato, Kazuya},
   title={Semi-stable reduction and crystalline cohomology with logarithmic
   poles},
   note={P\'{e}riodes $p$-adiques (Bures-sur-Yvette, 1988)},
   journal={Ast\'{e}risque},
   number={223},
   date={1994},
   pages={221--268},
   issn={0303-1179},
   review={\MR{1293974}},
}



\bib{Jac79}{article}{
   author={Jacquet, Herv\'{e}},
   title={Principal $L$-functions of the linear group},
   conference={
      title={Automorphic forms, representations and $L$-functions},
      address={Proc. Sympos. Pure Math., Oregon State Univ., Corvallis,
      Ore.},
      date={1977},
   },
   book={
      series={Proc. Sympos. Pure Math., XXXIII},
      publisher={Amer. Math. Soc., Providence, R.I.},
   },
   date={1979},
   pages={63--86},
   review={\MR{546609}},
}

\bib{Jan88}{article}{
   author={Jannsen, Uwe},
   title={Continuous \'{e}tale cohomology},
   journal={Math. Ann.},
   volume={280},
   date={1988},
   number={2},
   pages={207--245},
   issn={0025-5831},
   review={\MR{929536}},
   doi={10.1007/BF01456052},
}


\bib{Jan00}{article}{
   author={Jannsen, Uwe},
   title={Letter from Jannsen to Gross on higher Abel-Jacobi maps},
   conference={
      title={The arithmetic and geometry of algebraic cycles},
      address={Banff, AB},
      date={1998},
   },
   book={
      series={NATO Sci. Ser. C Math. Phys. Sci.},
      volume={548},
      publisher={Kluwer Acad. Publ., Dordrecht},
   },
   isbn={0-7923-6193-8},
   date={2000},
   pages={261--275},
   review={\MR{1744948}},
}


\bib{KMSW}{article}{
   author={Kaletha, Tasho},
   author={Minguez, Alberto},
   author={Shin, Sug Woo},
   author={White, Paul-James},
   title={Endoscopic Classification of Representations: Inner Forms of Unitary Groups},
   note={\href{https://arxiv.org/abs/1409.3731}{arXiv:1409.3731}},
}


\bib{Kar79}{article}{
   author={Karel, Martin L.},
   title={Functional equations of Whittaker functions on $p$-adic groups},
   journal={Amer. J. Math.},
   volume={101},
   date={1979},
   number={6},
   pages={1303--1325},
   issn={0002-9327},
   review={\MR{548883}},
   doi={10.2307/2374142},
}



\bib{Kat94}{article}{
   author={Kato, Kazuya},
   title={Semi-stable reduction and $p$-adic \'{e}tale cohomology},
   note={P\'{e}riodes $p$-adiques (Bures-sur-Yvette, 1988)},
   journal={Ast\'{e}risque},
   number={223},
   date={1994},
   pages={269--293},
   issn={0303-1179},
   review={\MR{1293975}},
}



\bib{KM74}{article}{
   author={Katz, Nicholas M.},
   author={Messing, William},
   title={Some consequences of the Riemann hypothesis for varieties over
   finite fields},
   journal={Invent. Math.},
   volume={23},
   date={1974},
   pages={73--77},
   issn={0020-9910},
   review={\MR{332791}},
   doi={10.1007/BF01405203},
}


\bib{KSZ}{article}{
    author={Kisin, Mark},
    author={Shin, Sug Woo},
    author={Zhu, Yihang},
    title={The stable trace formula for Shimura varieties of abelian type},
   note={\href{https://arxiv.org/abs/2110.05381}{arXiv:2110.05381}},
}


\bib{Kob13}{article}{
   author={Kobayashi, Shinichi},
   title={The $p$-adic Gross-Zagier formula for elliptic curves at
   supersingular primes},
   journal={Invent. Math.},
   volume={191},
   date={2013},
   number={3},
   pages={527--629},
   issn={0020-9910},
   review={\MR{3020170}},
   doi={10.1007/s00222-012-0400-9},
}



\bib{Kud97}{article}{
   author={Kudla, Stephen S.},
   title={Algebraic cycles on Shimura varieties of orthogonal type},
   journal={Duke Math. J.},
   volume={86},
   date={1997},
   number={1},
   pages={39--78},
   issn={0012-7094},
   review={\MR{1427845}},
   doi={10.1215/S0012-7094-97-08602-6},
}



\bib{Kud02}{article}{
   author={Kudla, Stephen S.},
   title={Derivatives of Eisenstein series and generating functions for
   arithmetic cycles},
   note={S\'{e}minaire Bourbaki, Vol. 1999/2000},
   journal={Ast\'{e}risque},
   number={276},
   date={2002},
   pages={341--368},
   issn={0303-1179},
   review={\MR{1886765}},
}


\bib{Kud03}{article}{
   author={Kudla, Stephen S.},
   title={Modular forms and arithmetic geometry},
   conference={
      title={Current developments in mathematics, 2002},
   },
   book={
      publisher={Int. Press, Somerville, MA},
   },
   date={2003},
   pages={135--179},
   review={\MR{2062318}},
}


\bib{Kud04}{article}{
   author={Kudla, Stephen S.},
   title={Special cycles and derivatives of Eisenstein series},
   conference={
      title={Heegner points and Rankin $L$-series},
   },
   book={
      series={Math. Sci. Res. Inst. Publ.},
      volume={49},
      publisher={Cambridge Univ. Press, Cambridge},
   },
   date={2004},
   pages={243--270},
   review={\MR{2083214}},
   doi={10.1017/CBO9780511756375.009},
}





\bib{KS97}{article}{
   author={Kudla, Stephen S.},
   author={Sweet, W. Jay, Jr.},
   title={Degenerate principal series representations for ${\rU}(n,n)$},
   journal={Israel J. Math.},
   volume={98},
   date={1997},
   pages={253--306},
   issn={0021-2172},
   review={\MR{1459856}},
   doi={10.1007/BF02937337},
}




\bib{Lan12}{article}{
   author={Lan, Kai-Wen},
   title={Comparison between analytic and algebraic constructions of
   toroidal compactifications of PEL-type Shimura varieties},
   journal={J. Reine Angew. Math.},
   volume={664},
   date={2012},
   pages={163--228},
   issn={0075-4102},
   review={\MR{2980135}},
   doi={10.1515/crelle.2011.099},
}


\bib{Lan99}{article}{
   author={Langer, Andreas},
   title={Local points of motives in semistable reduction},
   journal={Compositio Math.},
   volume={116},
   date={1999},
   number={2},
   pages={189--217},
   issn={0010-437X},
   review={\MR{1686773}},
   doi={10.1023/A:1000829923416},
}



\bib{LL}{article}{
   author={Li, Chao},
   author={Liu, Yifeng},
   title={Chow groups and $L$-derivatives of automorphic motives for unitary groups},
   journal={Ann. of Math. (2)},
   volume={194},
   date={2021},
   number={3},
   pages={817--901},
   issn={0003-486X},
   review={\MR{4334978}},
   doi={10.4007/annals.2021.194.3.6},
}


\bib{LL2}{article}{
   author={Li, Chao},
   author={Liu, Yifeng},
   title={Chow groups and $L$-derivatives of automorphic motives for unitary
   groups, II},
   journal={Forum Math. Pi},
   volume={10},
   date={2022},
   pages={Paper No. e5},
   review={\MR{4390300}},
   doi={10.1017/fmp.2022.2},
}


\bib{LZ}{article}{
   author={Li, Chao},
   author={Zhang, Wei},
   title={Kudla-Rapoport cycles and derivatives of local densities},
   journal={J. Amer. Math. Soc.},
   volume={35},
   date={2022},
   number={3},
   pages={705--797},
   issn={0894-0347},
   review={\MR{4433078}},
   doi={10.1090/jams/988},
}




\bib{Li92}{article}{
   author={Li, Jian-Shu},
   title={Nonvanishing theorems for the cohomology of certain arithmetic
   quotients},
   journal={J. Reine Angew. Math.},
   volume={428},
   date={1992},
   pages={177--217},
   issn={0075-4102},
   review={\MR{1166512}},
   doi={10.1515/crll.1992.428.177},
}

\bib{Liu11}{article}{
   author={Liu, Yifeng},
   title={Arithmetic theta lifting and $L$-derivatives for unitary groups,
   I},
   journal={Algebra Number Theory},
   volume={5},
   date={2011},
   number={7},
   pages={849--921},
   issn={1937-0652},
   review={\MR{2928563}},
   doi={10.2140/ant.2011.5.849},
}


\bib{Liu12}{article}{
   author={Liu, Yifeng},
   title={Arithmetic theta lifting and $L$-derivatives for unitary groups,
   II},
   journal={Algebra Number Theory},
   volume={5},
   date={2011},
   number={7},
   pages={923--1000},
   issn={1937-0652},
   review={\MR{2928564}},
   doi={10.2140/ant.2011.5.923},
}



\bib{Liu19}{article}{
   author={Liu, Yifeng},
   title={Bounding cubic-triple product Selmer groups of elliptic curves},
   journal={J. Eur. Math. Soc. (JEMS)},
   volume={21},
   date={2019},
   number={5},
   pages={1411--1508},
   issn={1435-9855},
   review={\MR{3941496}},
   doi={10.4171/JEMS/865},
}



\bib{Liu22}{article}{
   author={Liu, Yifeng},
   title={Theta correspondence for almost unramified representations of
   unitary groups},
   journal={J. Number Theory},
   volume={230},
   date={2022},
   pages={196--224},
   issn={0022-314X},
   review={\MR{4327954}},
   doi={10.1016/j.jnt.2021.03.027},
}



\bib{LTXZZ}{article}{
   author={Liu, Yifeng},
   author={Tian, Yichao},
   author={Xiao, Liang},
   author={Zhang, Wei},
   author={Zhu, Xinwen},
   title={On the Beilinson-Bloch-Kato conjecture for Rankin-Selberg motives},
   journal={Invent. Math.},
   volume={228},
   date={2022},
   number={1},
   pages={107--375},
   issn={0020-9910},
   review={\MR{4392458}},
   doi={10.1007/s00222-021-01088-4},
}


\bib{Mok15}{article}{
   author={Mok, Chung Pang},
   title={Endoscopic classification of representations of quasi-split
   unitary groups},
   journal={Mem. Amer. Math. Soc.},
   volume={235},
   date={2015},
   number={1108},
   pages={vi+248},
   issn={0065-9266},
   isbn={978-1-4704-1041-4},
   isbn={978-1-4704-2226-4},
   review={\MR{3338302}},
   doi={10.1090/memo/1108},
}


\bib{Mok93}{article}{
   author={Mokrane, A.},
   title={La suite spectrale des poids en cohomologie de Hyodo-Kato},
   language={French},
   journal={Duke Math. J.},
   volume={72},
   date={1993},
   number={2},
   pages={301--337},
   issn={0012-7094},
   review={\MR{1248675}},
   doi={10.1215/S0012-7094-93-07211-0},
}



\bib{Nek93}{article}{
   author={Nekov\'{a}\v{r}, Jan},
   title={On $p$-adic height pairings},
   conference={
      title={S\'{e}minaire de Th\'{e}orie des Nombres, Paris, 1990--91},
   },
   book={
      series={Progr. Math.},
      volume={108},
      publisher={Birkh\"{a}user Boston, Boston, MA},
   },
   date={1993},
   pages={127--202},
   review={\MR{1263527}},
   doi={10.1007/s10107-005-0696-y},
}


\bib{Nek95}{article}{
   author={Nekov\'{a}\v{r}, Jan},
   title={On the $p$-adic height of Heegner cycles},
   journal={Math. Ann.},
   volume={302},
   date={1995},
   number={4},
   pages={609--686},
   issn={0025-5831},
   review={\MR{1343644}},
   doi={10.1007/BF01444511},
}



\bib{Nek00}{article}{
   author={Nekov\'{a}\v{r}, Jan},
   title={$p$-adic Abel-Jacobi maps and $p$-adic heights},
   conference={
      title={The arithmetic and geometry of algebraic cycles},
      address={Banff, AB},
      date={1998},
   },
   book={
      series={CRM Proc. Lecture Notes},
      volume={24},
      publisher={Amer. Math. Soc., Providence, RI},
   },
   date={2000},
   pages={367--379},
   review={\MR{1738867}},
   doi={10.1090/crmp/024/18},
}




\bib{NT}{article}{
   author={Newton, James},
   author={Thorne, Jack A.},
   title={Symmetric power functoriality for Hilbert modular forms},
   note={\href{https://arxiv.org/abs/2212.03595}{arXiv:2212.03595}},
}



\bib{PR87}{article}{
   author={Perrin-Riou, Bernadette},
   title={Points de Heegner et d\'{e}riv\'{e}es de fonctions $L$ $p$-adiques},
   language={French},
   journal={Invent. Math.},
   volume={89},
   date={1987},
   number={3},
   pages={455--510},
   issn={0020-9910},
   review={\MR{903381}},
   doi={10.1007/BF01388982},
}


\bib{Ral82}{article}{
   author={Rallis, Stephen},
   title={Langlands' functoriality and the Weil representation},
   journal={Amer. J. Math.},
   volume={104},
   date={1982},
   number={3},
   pages={469--515},
   issn={0002-9327},
   review={\MR{658543}},
   doi={10.2307/2374151},
}




\bib{Ram}{article}{
   author={Ramakrishnan, D.},
   title={A theorem on $\GL(n)$ a la Tchebotarev},
   note={\href{https://arxiv.org/abs/1806.08429}{arXiv:1806.08429}},
}



\bib{RSZ20}{article}{
   author={Rapoport, M.},
   author={Smithling, B.},
   author={Zhang, W.},
   title={Arithmetic diagonal cycles on unitary Shimura varieties},
   journal={Compos. Math.},
   volume={156},
   date={2020},
   number={9},
   pages={1745--1824},
   issn={0010-437X},
   review={\MR{4167594}},
   doi={10.1112/s0010437x20007289},
}




\bib{Sat07}{article}{
   author={Sato, Kanetomo},
   title={$p$-adic \'{e}tale Tate twists and arithmetic duality},
   language={English, with English and French summaries},
   note={With an appendix by Kei Hagihara},
   journal={Ann. Sci. \'{E}cole Norm. Sup. (4)},
   volume={40},
   date={2007},
   number={4},
   pages={519--588},
   issn={0012-9593},
   review={\MR{2191526}},
   doi={10.1016/j.ansens.2007.04.002},
}


\bib{Sat13}{article}{
   author={Sato, Kanetomo},
   title={Cycle classes for $p$-adic \'{e}tale Tate twists and the image of
   $p$-adic regulators},
   journal={Doc. Math.},
   volume={18},
   date={2013},
   pages={177--247},
   issn={1431-0635},
   review={\MR{3064983}},
}



\bib{Sch94}{article}{
   author={Scholl, A. J.},
   title={Height pairings and special values of $L$-functions},
   conference={
      title={Motives},
      address={Seattle, WA},
      date={1991},
   },
   book={
      series={Proc. Sympos. Pure Math.},
      volume={55},
      publisher={Amer. Math. Soc., Providence, RI},
   },
   date={1994},
   pages={571--598},
   review={\MR{1265545}},
   doi={10.1090/pspum/055.1/1265545},
}


\bib{Sha}{article}{
   author={Shah, Syed Waqar Ali},
    title={Explicit Hecke descent for special cycles},
   note={\href{https://arxiv.org/abs/2310.01677}{arXiv:2310.01677}},
}


\bib{Shn16}{article}{
   author={Shnidman, Ariel},
   title={$p$-adic heights of generalized Heegner cycles},
   language={English, with English and French summaries},
   journal={Ann. Inst. Fourier (Grenoble)},
   volume={66},
   date={2016},
   number={3},
   pages={1117--1174},
   issn={0373-0956},
   review={\MR{3494168}},
}


\bib{Tan99}{article}{
   author={Tan, Victor},
   title={Poles of Siegel Eisenstein series on $\mathrm{U}(n,n)$},
   journal={Canad. J. Math.},
   volume={51},
   date={1999},
   number={1},
   pages={164--175},
   issn={0008-414X},
   review={\MR{1692899}},
   doi={10.4153/CJM-1999-010-4},
}


\bib{TY07}{article}{
   author={Taylor, Richard},
   author={Yoshida, Teruyoshi},
   title={Compatibility of local and global Langlands correspondences},
   journal={J. Amer. Math. Soc.},
   volume={20},
   date={2007},
   number={2},
   pages={467--493},
   issn={0894-0347},
   review={\MR{2276777}},
   doi={10.1090/S0894-0347-06-00542-X},
}


\bib{Tsu99}{article}{
   author={Tsuji, Takeshi},
   title={$p$-adic \'{e}tale cohomology and crystalline cohomology in the
   semi-stable reduction case},
   journal={Invent. Math.},
   volume={137},
   date={1999},
   number={2},
   pages={233--411},
   issn={0020-9910},
   review={\MR{1705837}},
   doi={10.1007/s002220050330},
}




\bib{Tsu00}{article}{
   author={Tsuji, Takeshi},
   title={On $p$-adic nearby cycles of log smooth families},
   language={English, with English and French summaries},
   journal={Bull. Soc. Math. France},
   volume={128},
   date={2000},
   number={4},
   pages={529--575},
   issn={0037-9484},
   review={\MR{1815397}},
}


\bib{Wal88}{article}{
   author={Wallach, Nolan R.},
   title={Lie algebra cohomology and holomorphic continuation of generalized
   Jacquet integrals},
   conference={
      title={Representations of Lie groups, Kyoto, Hiroshima, 1986},
   },
   book={
      series={Adv. Stud. Pure Math.},
      volume={14},
      publisher={Academic Press, Boston, MA},
   },
   date={1988},
   pages={123--151},
   review={\MR{1039836}},
   doi={10.2969/aspm/01410123},
}

\bib{Yam14}{article}{
   author={Yamana, Shunsuke},
   title={L-functions and theta correspondence for classical groups},
   journal={Invent. Math.},
   volume={196},
   date={2014},
   number={3},
   pages={651--732},
   issn={0020-9910},
   review={\MR{3211043}},
   doi={10.1007/s00222-013-0476-x},
}


\bib{Yam11}{article}{
   author={Yamashita, Go},
   title={$p$-adic Hodge theory for open varieties},
   language={English, with English and French summaries},
   journal={C. R. Math. Acad. Sci. Paris},
   volume={349},
   date={2011},
   number={21-22},
   pages={1127--1130},
   issn={1631-073X},
   review={\MR{2855488}},
   doi={10.1016/j.crma.2011.10.016},
}



\bib{YZZ}{book}{
   author={Yuan, Xinyi},
   author={Zhang, Shou-Wu},
   author={Zhang, Wei},
   title={The Gross-Zagier formula on Shimura curves},
   series={Annals of Mathematics Studies},
   volume={184},
   publisher={Princeton University Press, Princeton, NJ},
   date={2013},
   pages={x+256},
   isbn={978-0-691-15592-0},
   review={\MR{3237437}},
}


\end{biblist}
\end{bibdiv}
\end{document}